\documentclass[11pt,reqno]{amsart}
\usepackage[utf8]{inputenc}
\usepackage[T1]{fontenc}
\usepackage{lmodern}
\usepackage[english]{babel}
\usepackage{amsmath,a4wide}
\usepackage{extarrows}
\usepackage{dutchcal}
\usepackage{xfrac}
\usepackage{esint}
\usepackage{mathtools}

\mathtoolsset{showonlyrefs}
\usepackage{stmaryrd,mathrsfs,amsthm,mathtools,yfonts,amssymb,color,braket,booktabs,graphicx,graphics,amsfonts}
\usepackage{latexsym,microtype,indentfirst,hyperref}
\usepackage{xcolor}
\usepackage{courier}
\usepackage{todonotes}
\usepackage{tikz}
\usepackage[margin=1.8cm]{geometry}

\newtheorem{theorem}{Theorem}[section]
\newtheorem{corollary}[theorem]{Corollary}

\newtheorem{lemma}[theorem]{Lemma}
\newtheorem{proposition}[theorem]{Proposition}

\theoremstyle{definition}
\newtheorem{definition}[theorem]{Definition}
\newtheorem{remark}[theorem]{Remark}

\newtheorem{theoremletter}{Theorem}

\numberwithin{equation}{section}
\numberwithin{subsection}{section}

\newcommand{\R}{\mathbb{R}}

\newcommand{\Sf}{\mathbb{S}}

\newcommand{\spt}{\mathrm{spt}}

\newcommand{\dist}{\mathrm{dist}}

\newcommand{\Lip}{\mathrm{Lip}}

\newcommand{\Ha}{\mathcal{H}}

\newcommand{\Leb}{\mathcal{L}}

\newcommand{\pa}{\partial}

\newcommand{\loc}{\mathrm{loc}}

\newcommand{\mres}{\mathbin{\vrule height 1.6ex depth 0pt width
0.13ex\vrule height 0.13ex depth 0pt width 1.3ex}}

\newcommand{\weakstar}{\stackrel{*}{\rightharpoonup}}

\newcommand{\abs}[1]{\lvert#1\rvert}

\newcommand{\var}{\mathbf{Var}}

\renewcommand{\t}{~}
\newcommand{\beq}{\begin{equation}}
\newcommand{\eeq}{\end{equation}}

\newcommand{\bu}{\mathbf{b}}

\newcommand{\res}
{\mathop{\hbox{\vrule height 8pt width .5pt depth 0pt \vrule height .5pt width 4pt depth 0pt}}\nolimits}

\def\R{\mathbb R}
\def\N{\mathbb N}

\def\cal{\mathcal}

\def\H{{\cal H}}

\def\de{\delta}
\def\e{\varepsilon}

\def\vphi{\varphi}

\newcommand{\medint}{-\kern -,375cm\int}
\newcommand{\medintinrigo}{-\kern -,315cm\int}

\renewcommand{\S}{\mathbb{S}}

\def\pa{\partial}
\def\Id{{\rm Id}\,}
\def\loc{{\rm loc}}
\def\Div{{\rm div}}
\def\B{\mathcal{B}}

\newcommand{\udens}[3]{\Theta^*_{#1}(#2,#3)}
\newcommand{\ldens}[3]{\Theta_*^{#1}(#2,#3)}
\newcommand{\dens}[3]{\Theta^{#1}(#2,#3)}

\title{Elliptic regularization for anisotropic mean curvature flow}

\author[Antonini]{Carlo Alberto Antonini}
\address{{\href{mailto:carlo.antonini@unimi.it}{carlo.antonini@unimi.it}} \hfill \scriptsize{Department of Mathematics ``Federigo Enriques'', University of Milan}}
\author[De Rosa]{Antonio De Rosa}
\address{{\href{mailto:antonio.derosa@unibocconi.it}{antonio.derosa@unibocconi.it}} \hfill \scriptsize{Department of Decision Sciences and BIDSA, Bocconi University}}
\author[Morini]{Massimiliano Morini}\address{{\href{mailto:massimiliano.morini@unipr.it}{massimiliano.morini@unipr.it}} \hfill \scriptsize{Dipartimento di Scienze Matematiche Fisiche e Informatiche, Universit\`a di Parma}}
\author[Stuvard]{Salvatore Stuvard}
\address{{\href{mailto:salvatore.stuvard@unimi.it}{salvatore.stuvard@unimi.it}} \hfill \scriptsize{Department of Mathematics ``Federigo Enriques'', University of Milan}}

\date{\today}

\begin{document}

\begin{abstract}
We develop an elliptic regularization scheme for anisotropic mean curvature flow in arbitrary dimension, with general norms serving as surface tension and mobility. Starting from a
non-fattening set of finite perimeter, we prove that the vertically rescaled minimizers and translating solitons converge in the Kuratowski sense to the unique generalized
level-set motion, and locally in BV to the flat flow. Furthermore, by implementing a restarting and matching procedure prior to curvature identification, we establish for the limit
flat flow a coupled phase-transport and surface-energy dissipation balance valid for arbitrary, possibly crystalline, anisotropies. When the surface tension is continuously
differentiable (away from the origin) and strictly convex, we further prove that the flat flow is a BV flow and that its reduced boundary defines a unit-density
anisotropic Brakke flow. To our knowledge, this is the first construction of unit-density anisotropic Brakke flows valid for every strictly convex differentiable surface tension, in arbitrary dimension and for prescribed non-fattening initial data, without assuming perimeter convergence. 
\end{abstract}

\maketitle

\tableofcontents

\section{Introduction}

In this paper, we study the evolution of sets $E(t)\subset\R^{N}$ governed by the formal geometric motion law
\begin{equation}\label{intro:general-law}
 v_t=-\mathfrak m(\nu_{E(t)})H^\phi_{\partial E(t)}
 \qquad\text{on }\partial E(t)\,,
\end{equation}
where $v_t$ denotes the outward normal velocity of the evolving boundary $\partial E(t)$ along its outer unit normal $\nu_{E(t)}$, and $\phi,\mathfrak m$ are norms on $\R^N$
playing the roles of surface tension and mobility, respectively. When the interface $\partial E$ is of class $C^2$ and the norm $\phi$ is of class $C^2$ away from the origin, the scalar anisotropic mean curvature
$H^\phi_{\partial E}$ is given by
\[
 H^\phi_{\partial E}
 =\operatorname{div}_{\partial E}\bigl(\nabla\phi(\nu_E)\bigr),
\]
where $\operatorname{div}_{\partial E}$ denotes the tangential divergence on $\partial E$. Variationally, $H^\phi_{\partial E}$ represents the first variation of the
$\phi$-anisotropic perimeter
\begin{equation}\label{intro:surface-energy}
 P^\phi(E)=\int_{\partial^*E}\phi(\nu_E)\,d\mathcal H^{N-1},
\end{equation}
written here in terms of the reduced boundary $\partial^*E$ so as to extend the interfacial energy to sets of finite perimeter in $\R^N$.

Curvature-driven interface motion is a fundamental model for the capillarity-driven relaxation of microstructures. The role of orientation-dependent surface tension goes back to
Herring's work on crystalline interfaces, while Mullins' idealized grain-boundary model made the relation between curvature and interface migration explicit
\cite{Herring1951,Mullins1956}. In the vector thermodynamics of Hoffman and Cahn \cite{HoffmanCahn1972,CahnHoffman1974}, the Euclidean unit normal is naturally replaced in
anisotropic capillarity by the Cahn--Hoffman vector $\nabla\phi(\nu)$. In a local two-phase setting, the evolution law \eqref{intro:general-law} balances this capillary driving
force against the kinetic resistance to interfacial migration; the same law also arises as the sharp-interface limit of phase-field models for the relaxation of a nonconserved order
parameter, as in the Allen--Cahn theory of antiphase boundary motion \cite{AllenCahn1979}.

From a physical standpoint, the two norms in \eqref{intro:general-law} encode distinct mechanisms: $\phi$ measures the energetic cost of a given interface orientation, whereas
$\mathfrak m$ governs the kinetic response rate to a thermodynamic driving force. Because their angular dependences need not coincide — as confirmed, for instance, by atomistic
simulations of grain-boundary migration \cite{JanssensEtAl2006} — it is natural to prescribe $\phi$ and $\mathfrak m$ independently. Throughout this work, we assume both $\phi$ and
$\mathfrak m$ to be arbitrary norms on $\R^N$ (that is, even, convex, and positively one-homogeneous functions, positive away from the origin). This generality encompasses
crystalline surface tensions, whose lack of differentiability gives rise to faceted Wulff shapes.

The distinct roles of surface tension and mobility are clearly reflected in the formal energy-dissipation identity for smooth evolutions:
\begin{equation}\label{intro:smooth-dissipation}
 \frac{d}{dt}P^\phi(E(t))
 =-\int_{\partial E(t)}\mathfrak m(\nu_{E(t)})\bigl(H^\phi_{\partial E(t)}\bigr)^2\,d\mathcal H^{N-1}
 =-\int_{\partial E(t)}\frac{v_t^2}{\mathfrak m(\nu_{E(t)})}\,d\mathcal H^{N-1}\,.
\end{equation}
This identity exhibits \eqref{intro:general-law} as the formal gradient flow of the anisotropic perimeter $P^\phi$ with respect to the $L^2$ metric on normal velocities weighted by
$\mathfrak m(\nu_{E(t)})^{-1}$ relative to the surface measure $\mathcal H^{N-1}\mres\partial E(t)$.

 Continuing the evolution past singularities, or even giving meaning to \eqref{intro:general-law} when $\phi$ is crystalline and the pointwise expression for $H^\phi_{\partial E}$
 breaks down, requires suitable weak formulations.

\subsection{Weak formulations and variational constructions}

Different weak theories of mean curvature flow capture complementary features of the smooth evolution. In the level-set approach, the evolving sets are realized as sublevel sets of
a function $u:\R^N\times[0,\infty)\to\R$ solving, formally where $\nabla u\neq0$, the degenerate parabolic equation
\[
 \partial_tu=\mathfrak m(\nabla u)
             \operatorname{div}\bigl(\nabla\phi(\nabla u)\bigr).
\]
The viscosity solution theory of Evans--Spruck \cite{EvansSpruckI} and Chen--Giga--Goto \cite{CGG} provides a well-posed, comparison-based framework for mean curvature flow, both in the isotropic case and for anisotropic surface tensions of class $C^2$ away from the origin. This formulation allows for topological changes of the evolving interface. Uniqueness of the level-set function $u$, however, does not prevent the phenomenon of \emph{fattening}: a level set
$\{u(\cdot,t)=\lambda\}$ may develop non-empty interior in space-time, causing the corresponding open and closed sublevel evolutions $\{u<\lambda\}$ and $\{u\leq\lambda\}$ to
separate. However, the phenomenon is rare: the set of levels $\lambda$ such that the corresponding level set fattens is at most countable.

When the surface tension $\phi$ is crystalline, the Cahn--Hoffman vector must be replaced at nondifferentiable directions by a selection in the subdifferential $\partial\phi(\nu)$,
making the speed of a facet nonlocal in character. A viscosity theory based on calibrated facet test functions was developed by Giga and Po\v{z}\'ar \cite{GigaPozar18}. In
parallel, Chambolle, Morini, and Ponsiglione \cite{CMP17} introduced a distributional comparison framework based on anisotropic distance functions, which was subsequently extended
by Chambolle, Morini, Novaga, and Ponsiglione \cite{CMNP19,CMNP19bis} to establish existence, uniqueness, and stability of level-set flows for arbitrary pairs of norms
$(\phi,\mathfrak m)$, including stability under smooth approximations. In addition to covering the crystalline case, \cite{CMNP19} provides, to our knowledge, the only available theory that establishes well-posedness in arbitrary dimension for general differentiable (but not necessarily $C^2$) surface tensions and arbitrary mobility norms.

On the geometric measure theory side, Brakke's formulation \cite{Brakkephd} models the evolving interface as a family of varifolds satisfying an inequality version of
\eqref{intro:smooth-dissipation} tested against smooth localization functions. While Brakke flows naturally track interfacial area and multiplicity, the varifold inequality by
itself does not directly track the enclosed phase. Conversely, the $BV$ formulation introduced by Luckhaus and Sturzenhecker \cite{LS} tracks the evolving phase $E(t)$ via its
characteristic function $\chi_{E(t)}$ and encodes the normal velocity in a distributional transport equation. These two descriptions can be combined into the notion of generalized
$BV$ flow introduced by Stuvard and Tonegawa in \cite{StuvardTonegawa2024} (see also work by Tashiro in \cite{Tashiro,Tashiro2,Tashiro3}), where the evolving phase is coupled to a Brakke flow whose weight measure dominates the perimeter measure $\mu_{E(t)}$.
Showing that the varifold weight actually coincides with the reduced boundary measure of the phase — so that no hidden interfaces or higher multiplicities persist — is a nontrivial
problem that does not follow from convergence of the phases alone.

Several works connect the viscosity and energy-measure viewpoints. Evans and Spruck \cite{EvansSpruckIV} proved that almost every level set of the isotropic viscosity solution
yields a unit-density Brakke flow, a result extended to smooth anisotropic flows by Tonegawa \cite{Tonegawa2000}. More recently, Ullrich and Laux \cite{UllrichLaux2024} established
the $BV$ formulation and optimal energy dissipation for generic level sets in the isotropic case. Because these results rely on coarea and generic-level arguments across a family
of level sets, they apply to almost every level $\lambda\in\R$ rather than to a single prescribed non-fattening initial set $E_0$.

A closely related variational perspective builds directly on the gradient-flow structure of \eqref{intro:smooth-dissipation}. Laux and Otto \cite{LauxOtto2020} proved conditional
convergence of the thresholding scheme to Brakke flow, while Hensel and Laux \cite{HenselLaux2025} introduced De Giorgi varifold solutions—coupling phase transport to an
energy-dissipation inequality—and established their unconditional existence via Allen--Cahn approximation as well as weak-strong uniqueness. Laux and Poiatti \cite{LauxPoiatti2026}
recently obtained De Giorgi varifold solutions 
%in dimensions $N=2,3$
 via a modified minimizing-movements scheme employing a regularized distance proxy, and Laux, Stinson, and
Ullrich \cite{LauxStinsonUllrich2025} proved conditional convergence of the anisotropic Allen--Cahn equation and weak-strong uniqueness via relative entropy. We mention that, in the isotropic case, very recently Laux and Takasao \cite{LauxTakasao} have proved unconditional convergence of the Allen-Cahn equation to multi-phase Brakke flows.

In the classical time-discrete minimizing-movements scheme of Almgren--Taylor--Wang \cite{ATW} and Luckhaus--Sturzenhecker \cite{LS}, one constructs discrete evolutions by
iteratively minimizing the sum of the $\phi$-perimeter and a bulk dissipation term proportional to the signed $\mathfrak m^o$-distance from the previous step (here, $\mathfrak m^o$ denotes, as customary, the polar norm of $\mathfrak m$). Compactness in $BV$
yields $L^1_{\loc}$ limits as the time step vanishes, known as \emph{flat flows}. Through Chambolle's dual formulation \cite{Chambolle} and the distributional theory of
\cite{CMP17,CMNP19,CMNP19bis}, flat flows are known to coincide, up to spatially negligible sets for almost every time, with the unique level-set evolution whenever the initial
datum is non-fattening in measure in the sense of Definition~\ref{def:fat}.

The central analytical challenge in the flat-flow program lies in passing from the time-discrete variational problems to the limiting evolution equation. Because $L^1_{\loc}$
convergence of characteristic functions only guarantees lower semicontinuity of the perimeter, parts of the approximating boundaries may cancel in the limit, causing a drop in
surface energy. In \cite{LS} (and in its anisotropic and inhomogeneous counterpart by Chambolle, De Gennaro, and Morini \cite{CDGM}), the $BV$ evolution equation is derived under
the conditional assumption that the time-integrated perimeters of the discrete approximations converge to the perimeter of the limit — an assumption that rules out boundary
cancellation and provides the necessary control on normals and first variations.

In the isotropic setting, De Gennaro and Stuvard \cite{DeGennaroStuvard2026} have recently removed both the energy-convergence and the non-fattening assumptions in dimensions
$2\leq N\leq 7$ (that is, for hypersurfaces of dimension $1\leq N-1\leq 6$), proving that every flat flow starting from a bounded set of finite perimeter with $\mathcal
L^N(\partial E_0)=0$ admits an integral Brakke enhancement (with possibly higher multiplicity) and satisfies the corresponding generalized $BV$ identity. Their proof, however, relies heavily on the quantitative
regularity theory for isotropic perimeter quasi-minimizers and Bernstein-type dimension bounds at the scale of the time step. Extending such unconditional arguments to general
anisotropic integrands faces fundamental geometric obstacles, since singular anisotropic area minimizers already occur in $\R^4$ even for smooth elliptic integrands, as shown by
Morgan \cite{Morgan1991}. In the present paper, we instead work under the non-fattening assumption, which allows us to treat arbitrary dimensions $N\geq 2$ and general norms
without assuming perimeter convergence.

\subsection{Elliptic regularization and the compactness obstruction}

For isotropic mean curvature flow, Ilmanen's elliptic regularization \cite{Ilm94} provides a powerful bridge between the comparison and variational theories. By minimizing an
exponentially weighted area functional in the extended space $\R^{N+1}=\R^N\times\R$, one obtains translating solitons that are themselves Brakke flows in $\R^{N+1}$. Passing to
the limit via Brakke's compactness theorem produces a limiting varifold flow, which, under non-fattening, can be matched with the level-set evolution to yield a unit-density Brakke
flow; see also Tashiro \cite{Tashiro} for the associated generalized $BV$ flow structure.

Here we address the long-standing problem of developing Ilmanen's elliptic regularization program for anisotropic mean curvature flows. In carrying out this program  a
major obstacle is compactness even for a strictly convex \(C^1\) surface tension. In the isotropic theory, the compactness of Brakke flows rests on two pillars: time-compactness
of the Radon weight measures, and Allard's rectifiability and compactness theorem \cite{Allard1972} on the time slices. Allard's theorem requires uniform lower density bounds to
prevent rectifiable surface measures from diffusing into the bulk in the weak limit—bounds that in the Euclidean setting are guaranteed by the monotonicity formula. For general
anisotropic integrands, however, the monotonicity formula fails \cite{Allard1974}. While the anisotropic rectifiability theorem of De Philippis, De Rosa, and Ghiraldin \cite{DPDRG}
shows that a varifold with locally bounded $\phi$-first variation is rectifiable provided $\phi\in C^1(\R^N\setminus\{0\})$ is strictly convex (the codimension-one atomic
condition) \emph{and} the varifold weight has positive lower density, such density estimate does not hold in general. Under the atomic condition, De Rosa
\cite[Theorem~4.1]{DeRosa2018Varifolds} also proved that a varifold limit of integral varifolds with locally uniformly bounded anisotropic first variations remains integral,
again provided its weight has positive lower density.

This missing density control is equally critical when attempting to pass to the limit in Brakke's inequality. Proving the limiting dissipation inequality requires the upper
semicontinuity of the Brakke functional, which in turn relies on the rectifiability of the limit varifold (cf. \cite[4.28]{Brakkephd}) in order to identify a unique approximate
tangent plane almost everywhere and pass to the limit in the transport term $(\nu\cdot\mathbf H^\phi)(\nabla\phi(\nu)\cdot\nabla\varphi)$. Without a priori rectifiability of the
limit, one cannot first establish the Brakke flow property and only subsequently match the flow to the level-set evolution.

To overcome this difficulty, we reverse the architecture of Ilmanen's proof and implement the \emph{restarting and matching} procedure \emph{before} proving the Brakke inequality
or identifying the mean curvature. Specifically:
\begin{enumerate}
 \item Using anisotropic paraboloid barriers and the distributional comparison theory of \cite{CMP17,CMNP19}, we first show directly at the level of sets that the vertically rescaled minimizers and translating solitons converge in the Kuratowski sense (and in $BV_{\loc}$ under \emph{non-fattening in measure}) to the unique generalized $(\phi,\mathfrak m)$-motion.
 \item From the Euler--Lagrange equations of the elliptic minimizers, we extract a limiting surface-energy measure $\lambda_t\,dt$ dominating the anisotropic perimeter $\mu^\phi_{\widehat E(t)}\,dt$ of the flat flow, together with a phase-transport identity and a surface-energy balance with dissipation measure $\Lambda$.
 \item Whenever the total mass $\lambda_t(\R^N)$ exceeds $P^\phi(\widehat E(t))$ by a threshold $\delta>0$, we restart the elliptic regularization scheme from the flat-flow slice $\widehat E(t)$ and absorb the downward mass jumps into the dissipation measure. Letting $\delta\downarrow 0$, the matched surface-energy measure is forced to coincide identically with the flat-flow perimeter measure $\sigma = \mu^\phi_{\widehat E(t)}\,dt$.
\end{enumerate}
The time slices $\mu^\phi_{\widehat E(t)}$ of $\sigma=\mu^\phi_{\widehat E(t)}\,dt$ then have positive $(N-1)$-dimensional density at $\mu^\phi_{\widehat E(t)}$-almost every point, for almost every time.
At that stage — and only when $\phi\in C^1(\R^N\setminus\{0\})$ is strictly convex — we can invoke the anisotropic rectifiability theorem of \cite{DPDRG}; once we have rectifiability and, evidently, unit density, we can finally identify the transport vector field in terms of the generalized anisotropic mean curvature $\mathbf H^\phi_{\partial^*\widehat E(t)}$, and obtain the
full $(\phi,\mathfrak m)$-Brakke and canonical $BV$ flow equations.

Very recently, the general compactness problem for anisotropic varifolds has seen major breakthroughs: De Philippis and Pigati \cite{DePhilippisPigati2026} established the
anisotropic Michael--Simon Sobolev inequality for two-dimensional varifolds in $\mathbb R^3$ for norms $C^1$-close to the Euclidean one (and its equivalence with
rectifiable-varifold compactness under the atomic condition), Firester and Tsiamis \cite{FiresterTsiamis2026} extended the inequality to all positive, convex, even $C^1$
hypersurface anisotropies, and De Philippis, Gennaioli, Pigati and Rindler \cite{DePhilippisGennaioliPigatiRindler2026} proved it in arbitrary dimension and codimension for
positive $C^1$ anisotropic integrands satisfying the weak atomic condition \textup{(AC1)}. While these results now provide the lower-density input needed to run a classical
varifold-compactness argument in the $C^1$ strictly convex setting, our reversed matching strategy has a crucial further advantage: steps (1)--(3) require \emph{no
differentiability or strict convexity} of $\phi$ or $\mathfrak m$. As a consequence, our approach yields a complete local surface-energy transport and dissipation balance for
arbitrary crystalline norms, which lie outside the $C^1$ framework of the above anisotropic varifold and Michael--Simon results.

\subsection{Main results}

Throughout the paper, we adopt the non-fattening notions of Definition~\ref{def:fat}: an initial set $E_0\subset\R^N$ is \emph{$(\phi,\mathfrak m)$-non-fattening} if the open and
closed zero-level sets of the $(\phi,\mathfrak m)$-level-set flow $u$ starting from a representation $u_0$ of $E_0$ satisfy $\overline{\{u<0\}}=\{u\leq 0\}$ and
$\{u<0\}=\mathrm{int}(\{u\leq 0\})$ in $\R^N\times(0,\infty)$; it is \emph{non-fattening in measure} if, in addition, $\mathcal L^{N+1}(\{u=0\}\cap(\R^N\times(0,\infty)))=0$. It is important to recall that for any uniformly continuous function $u_0$, all but countably many sublevel sets $\{u_0\leq\lambda\}$ are non-fattening in measure\cite{CMNP19}. When $E_0$ has finite perimeter and is non-fattening in measure, all $(\phi,\mathfrak m)$-flat flows starting from $E_0$ coincide for every $t\geq 0$ up to $\mathcal L^N$-negligible sets, and we denote this unique flat flow by
$\{\widehat E(t)\}_{t\geq 0}$, with $\widehat E(0)$ understood as the $L^1_{\loc}$ right trace as $t\downarrow 0$.

Before stating our first main result, let us introduce the notation in use in the anisotropic elliptic regularization scheme. Consider first the choice $\mathfrak m = \phi$. This is often referred to  as the {\em natural mobility:}  indeed, for this choice the Wulff shape of $\phi$ shrinks homothetically. In this case, the evolution law reads
\begin{equation}\label{evol law}
 v_t=-\phi(\nu_{E(t)})H^\phi_{\partial E(t)}\,,
\end{equation}
and the elliptic regularization scheme is defined in $\R^{N+1}=\R^N_x\times\R_z$ using the single extended norm
\[
 \Phi(x,z)=\sqrt{\phi(x)^2+z^2}.
\]
For $\varepsilon>0$, we then minimize the weighted anisotropic area functional
\[
 \mathcal A^\Phi_\varepsilon(F)
 =\int_{\partial^*F\cap\{z>-1\}}
    e^{-z/\varepsilon}\Phi(\nu_F)\,d\mathcal H^{N}
\]
among sets of locally finite perimeter $F\subset\R^{N+1}$ satisfying the bulk constraint $F\cap\{z<0\}\stackrel{a.e.}{=}E_0\times(-\infty,0)$. From a minimizer $E_\varepsilon$, we
form the translating solitons $\mathcal S_\varepsilon(t):=\overline{E_\varepsilon}-(t/\varepsilon)e_{_{N+1}}$ (with space-time tube $\mathcal S_\varepsilon$) and the vertically
rescaled sets $\mathcal R_\varepsilon(E_\varepsilon):=\{(x,\varepsilon z)\in\R^N\times(0,\infty):(x,z)\in E_\varepsilon\}$. For a general mobility $\mathfrak m$, we instead minimize
$\mathcal A^{\Phi^\kappa}_\varepsilon$ for a family of compatible norms $(\Phi^\kappa,\mathfrak M^\kappa)$ on $\R^{N+1}$ whose horizontal restrictions $\phi^\kappa:=\Phi^\kappa(\cdot,0)$ and $\mathfrak m^\kappa:=\mathfrak M^\kappa(\cdot,0)$ converge to $\phi$ and $\mathfrak m$,  and
pass to the limit along a suitable diagonal $\kappa=\kappa(\varepsilon)\to 0^+$. We denote by   $E_\varepsilon^\kappa$, $\mathcal S_\varepsilon^\kappa(t)$ and
$\mathcal S_\varepsilon^\kappa$ the corresponding  minimizers, translating solitons and tubes, respectively. The technical details concerning the need, in the general mobility case, of a one-parameter family of extensions as opposed to one single extension will be fully addressed in subsection \ref{intro:seq-ext}. These additional difficulties are the reason for presenting first the details of the regularization scheme in the natural mobility case (Sections \ref{sec:kur:conv} and \ref{sec:ellip:varifold}). 

\begin{theoremletter}[Geometric and $BV$ identification of elliptic regularization]
\label{intro:kuratowski}
Let $\phi,\mathfrak m$ be arbitrary norms on $\R^{N}$, and let $E_0\subset\R^N$ be a set of finite perimeter that is non-fattening for $(\phi,\mathfrak m)$. Denote by
$(E(t))_{t>0}$ the unique generalized $(\phi,\mathfrak m)$-motion starting from $E_0$, by $E\subset\R^N\times(0,\infty)$ its space-time tube, and by $\mathcal
E\subset\R^{N+1}\times(0,\infty)$ the cylindrical tube with slices $\mathcal E(t)=E(t)\times\R$. Then there exists a diagonal regime $\kappa=\kappa(\varepsilon)\to 0^+$ as
$\varepsilon\to 0^+$ along which the following Kuratowski convergences hold as $\varepsilon\to 0^+$:
\begin{itemize}
 \item[(a)] $\mathcal S_\varepsilon^{\kappa(\varepsilon)} \stackrel{\mathcal K}{\longrightarrow} \mathcal E$ in $\R^{N+1}\times(0,\infty)$;
 \item[(b)] $\mathcal S_\varepsilon^{\kappa(\varepsilon)}(t) \stackrel{\mathcal K}{\longrightarrow} E(t)\times\R$ in $\R^{N+1}$ for all but countably many $t>0$;
 \item[(c)] $\mathcal R_\varepsilon\bigl(E_\varepsilon^{\kappa(\varepsilon)}\bigr) \stackrel{\mathcal K}{\longrightarrow} E$ in $\R^N\times(0,\infty)$,
\end{itemize}
together with the corresponding Kuratowski convergences for the complementary sets and boundaries. In the natural mobility $\mathfrak m=\phi$ case, one simply takes
the fixed extension $\Phi$ and may drop all the $\kappa(\varepsilon)$-superscripts throughout.

If, in addition, $E_0$ is non-fattening in measure and $\{\widehat E(t)\}_{t\geq 0}$ denotes the unique $(\phi,\mathfrak m)$-flat flow starting from $E_0$ (with space-time tube
$\widehat E$), then as $\varepsilon\to 0^+$:
\begin{itemize}
 \item[(d)] $\mathcal S_\varepsilon^{\kappa(\varepsilon)}(t) \longrightarrow \widehat E(t)\times\R$ weakly-* in $BV_{\loc}(\R^{N+1})$ for $\mathcal L^1$-a.e.\ $t>0$;
 \item[(e)] $\mathcal R_\varepsilon\bigl(E_\varepsilon^{\kappa(\varepsilon)}\bigr) \longrightarrow \widehat E$ weakly-* in $BV_{\loc}(\R^N\times(0,\infty))$.
\end{itemize}
\end{theoremletter}

For natural mobility, Theorem~\ref{intro:kuratowski} corresponds to Theorems~\ref{thm:main1} and~\ref{thm:altro:main}; the general mobility case is established in
Theorem~\ref{mob:thm:main1}, Corollary~\ref{mob:cor:diagonal}, and Theorem~\ref{mob:completion:diagonal}. The key step in establishing the result is to show, by a
delicate approximation procedure,  that the solitons produced by the scheme are distributional superflows  in the sense of \cite{CMP17, CMNP19}. Because the proof relies on
anisotropic  paraboloid barriers, Wulff-shape comparison, and the uniform perimeter bounds on $\{z>0\}$, it requires no smoothness on $\phi$ or $\mathfrak m$. To the best of our
knowledge, the full Kuratowski convergence of the rescaled sets $\mathcal R_\varepsilon(E_\varepsilon)$, the translating soliton tubes, and their complements is new even in the
classical isotropic setting.

We now turn to the connection between the set-theoretic and energy-measure descriptions of the flow. Our second main result provides a local energy and transport description of the flat flow $\{\widehat E(t)\}_{t\geq 0}$ for arbitrary (possibly crystalline) norms $\phi$ and
$\mathfrak m$. Set
\[
 \widehat E=\{(x,t)\in\R^N\times(0,\infty):x\in\widehat E(t)\},\qquad
 d\mu=d\mu_{\widehat E(t)}\,dt,\qquad
 d\sigma=d\mu^\phi_{\widehat E(t)}\,dt=\phi(\nu_{\widehat E(t)})\,d\mu,
\]
where $\mu_{\widehat E(t)}=\mathcal H^{N-1}\mres\partial^*\widehat E(t)$ and $\nu=\nu_{\widehat E(t)}$ is the measure-theoretic outer unit normal to $\widehat E(t)$.

\begin{theoremletter}[Phase transport and surface-energy dissipation]
\label{intro:transport}
Let $\phi,\mathfrak m$ be arbitrary norms on $\R^N$, and let $E_0\subset\R^N$ be a set of finite perimeter that is non-fattening in measure. Then there exist a normal velocity $v\in
L^2(\mu)$, a transport vector field $\mathbf b\in L^2(\sigma;\R^{N})$, a mobility weight $h\in L^\infty(\sigma)$, and a finite nonnegative Radon measure $\Lambda$ on
$\R^N\times(0,\infty)$ such that
\begin{equation}\label{intro:phase-energy}
 \begin{gathered}
 D_{x,t}\chi_{\widehat E}=(-\nu,v)\mu,
 \qquad \mathbf b\cdot\nu=v \quad\mu\text{-a.e.},\\
 \partial_t\sigma+\operatorname{div}_x(\mathbf b\sigma)=-\Lambda \quad\text{in }\mathcal D'\bigl(\R^N\times(0,\infty)\bigr),
 \end{gathered}
\end{equation}
and
\begin{equation}\label{intro:general-cost}
 \frac{\mathfrak m(\nu)}{\phi(\nu)}\leq h\leq
 \max_{\mathbb S^{N-1}}\frac{\mathfrak m}{\phi} \quad \sigma-\mbox{a.e.}, \qquad
 \frac{\bigl(\phi^o(\mathbf b)\bigr)^2}{h}\,\sigma\leq\Lambda,
 \qquad \Lambda\bigl(\R^{N}\times(0,\infty)\bigr)\leq P^\phi(E_0).
\end{equation}
More precisely, for every $0\leq t_1<t_2$ and every test function $\eta\in C^1(\R^N\times[0,\infty))$ with compact spatial support uniformly on $[t_1,t_2]$, one has the BV identity
\begin{equation} \label{intro:BV-id}
\int_{\widehat E(t_2)} \eta (x,t_2) \, dx - \int_{\widehat E(t_1)} \eta (x,t_1) \, dx = \int_{t_1}^{t_2} \int_{\widehat E(t)} \partial_t \eta(x,t) \, dx dt + \iint_{\R^N \times (t_1,t_2)} \eta v \, d\mu
\end{equation}
for $v=\mathbf b \cdot \nu \in L^2(\mu)$ and $\mathbf b$ satisfying the transport identity
\begin{equation} \label{intro:energy-balance}
\mu^\phi_{\widehat E(t_2)}(\eta (\cdot,t_2)) -\mu^\phi_{\widehat E(t_1)}(\eta (\cdot,t_1)) = \iint_{\R^N\times (t_1,t_2)} (\partial_t \eta + \mathbf b \cdot \nabla_x \eta) \, d\sigma -  \iint_{\R^N\times (t_1,t_2]} \eta \, d\Lambda\,.
\end{equation}
Finally, if $\mathcal L^{N}(\partial E_0)=0$, then $\widehat E(0)\stackrel{a.e.}{=}E_0$.
\end{theoremletter}

Theorem~\ref{intro:transport} summarizes Theorems~\ref{trasportobv} and~\ref{mob:trasportobv}, together with Lemma~\ref{lem:flat-initial-trace}. The first identity in
\eqref{intro:phase-energy} states that the indicator function $\chi_{\widehat E}$ is BV in space-time, and it therefore expresses the motion of the phase through its distributional
normal velocity $v$. Its integrated form \eqref{intro:BV-id} is the weak counterpart of the classical transport formula for integrals over a moving domain.

The second line of \eqref{intro:phase-energy} is a local balance between the transport of the surface energy $\sigma=\mu^\phi_{\widehat E(t)}\,dt$ along the velocity field $\mathbf
b$ and its dissipation, encoded by the nonnegative Radon measure $\Lambda$. The compatibility relation $\mathbf b\cdot\nu=v$ ensures that the normal component of the energy flux
coincides with the normal velocity of the phase boundary, even though $\mathbf b$ itself need not be parallel to $\nu$ (reflecting the fact that, in the anisotropic setting, the
Cahn--Hoffman vector generally tilts away from the Euclidean normal). Note also that $\Lambda$ need not be absolutely continuous with respect to $\sigma$; in particular, it records
possible instantaneous drops in surface energy at singular times.

Through the quadratic estimate in \eqref{intro:general-cost}, the dissipation measure $\Lambda$ also yields quantitative $L^2$ control on the normal velocity $v$. Here, the weight
$h$ represents the local ratio between the limiting mobility and surface-tension measures on the reduced boundary: for general crystalline norms, $h$ is sandwiched between the
bounds in \eqref{intro:general-cost}, without necessarily coinciding pointwise with $\mathfrak m(\nu)/\phi(\nu)$. Combining the compatibility condition $\mathbf b\cdot\nu=v$ with
the duality inequality $|v|=|\mathbf b\cdot\nu|\leq\phi(\nu)\phi^o(\mathbf b)$, we obtain the localized kinetic bound
\begin{equation} \label{eq:kinetic-estimate}
 \frac{v^2}{h\,\phi(\nu)}\,\mu
 =
 \frac{(\mathbf b\cdot\nu)^2}{h\,\phi(\nu)^2}\,\sigma
 \leq
 \frac{\bigl(\phi^o(\mathbf b)\bigr)^2}{h}\,\sigma
 \leq\Lambda.
\end{equation}
Notice also that testing the integrated energy balance \eqref{intro:energy-balance} on
$[t_1,t_2]$ with spatial cutoff functions invading $\R^N$ gives the identity 
\begin{equation} \label{eq:totLambda=diss}
\Lambda\bigl(\R^N\times(t_1,t_2]\bigr)=P^\phi(\widehat E(t_1))-P^\phi(\widehat E(t_2))\,.
\end{equation}
In the natural mobility case $\mathfrak m=\phi$, the bounds in \eqref{intro:general-cost} force $h\equiv1$, so that \eqref{eq:kinetic-estimate} and \eqref{eq:totLambda=diss} can be combined to give the sharp energy-dissipation relation
\begin{equation}\label{intro:natural-action}
 P^\phi(\widehat E(t_2))+
 \int_{t_1}^{t_2}\int\frac{v^2}{\phi(\nu)}
 \,d\mu_{\widehat E(t)}\,dt
 \leq P^\phi(\widehat E(t_1)),
 \qquad 0\leq t_1<t_2.
\end{equation}
In this way, Theorem~\ref{intro:transport} captures the local kinematic and energetic structure of the flat flow directly on its reduced boundary, even for crystalline anisotropies
where a varifold first variation is not available.

When the surface tension $\phi$ is $C^1$ and strictly convex, we can identify both the weight $h$ and the transport vector field $\mathbf b$ in terms of the generalized anisotropic
mean curvature of the slices, obtaining a unit-density Brakke flow and a canonical $BV$ flow.

\begin{theoremletter}[Brakke and canonical $BV$ realization]
\label{intro:brakke}
Under the hypotheses of Theorem~\ref{intro:transport}, assume in addition that $\phi\in C^1(\R^{N}\setminus\{0\})$ is strictly convex. Then the family of unit-density varifolds
\[
 V_t=\var\bigl(\partial^*\widehat E(t)\bigr),\qquad t\geq0,
\]
is a $(\phi,\mathfrak m)$-anisotropic Brakke flow in $\R^N$. For almost every $t>0$, $V_t$ admits a generalized $\phi$-anisotropic mean curvature $\mathbf
H^\phi_{\partial^*\widehat E(t)}\in L^2(\mu_{\widehat E(t)}\,dt)$ which is orthogonal to $\partial^*\widehat E(t)$, and the fields constructed in the proof of
Theorem~\ref{intro:transport} satisfy
\begin{equation}\label{intro:curvature-identification}
 h=\frac{\mathfrak m(\nu)}{\phi(\nu)},\qquad
 \mathbf b=\frac{\mathfrak m(\nu)}{\phi(\nu)}\bigl(\nu\cdot\mathbf H^\phi_{\partial^*\widehat E(t)}\bigr)\nabla\phi(\nu),\qquad
 v=\mathfrak m(\nu)\bigl(\nu\cdot\mathbf H^\phi_{\partial^*\widehat E(t)}\bigr)
 \quad\mu\text{-a.e.},
\end{equation}
together with the localized dissipation bound $\mathfrak m(\nu)\bigl|\mathbf H^\phi_{\partial^*\widehat E(t)}\bigr|^2\mu\leq\Lambda$.

In particular, for every $0\leq t_1<t_2$ and every $\eta\in C^1(\R^N\times[0,\infty);\R_+)$ with compact spatial support uniformly on $[t_1,t_2]$, the $(\phi,\mathfrak m)$-Brakke
inequality holds:
\begin{equation}\label{intro:brakke-ineq}
\begin{aligned}
 &\left.\int\eta(x,t)\phi(\nu)\,d\mu_{\widehat E(t)}\right|_{t=t_1}^{t_2}\\
 &\quad\leq\int_{t_1}^{t_2}\int
 \left[\phi(\nu)\partial_t\eta+
 \mathfrak m(\nu)\Bigl(\bigl(\nu\cdot\mathbf H^\phi_{\partial^*\widehat E(t)}\bigr)
 \nabla\phi(\nu)\cdot\nabla_x\eta-\eta\bigl|\mathbf H^\phi_{\partial^*\widehat E(t)}\bigr|^2\Bigr)\right]
 \,d\mu_{\widehat E(t)}\,dt,
\end{aligned}
\end{equation}
and, for arbitrary signed test functions $\eta$ in the same class, the evolving phase satisfies the canonical $BV$ identity
\begin{equation}\label{intro:bv-identity}
\begin{aligned}
 \left.\int_{\widehat E(t)}\eta(x,t)\,dx\right|_{t=t_1}^{t_2}
 &=\int_{t_1}^{t_2}\int_{\widehat E(t)}\partial_t\eta\,dx\,dt
 +\int_{t_1}^{t_2}\int
 \eta\,\mathfrak m(\nu)\bigl(\nu\cdot\mathbf H^\phi_{\partial^*\widehat E(t)}\bigr)
 \,d\mu_{\widehat E(t)}\,dt,
\end{aligned}
\end{equation}
as well as the optimal energy-dissipation inequality
\begin{equation}\label{intro:optimal-action}
 P^\phi(\widehat E(t_2))+
 \int_{t_1}^{t_2}\int\mathfrak m(\nu)\bigl|\mathbf H^\phi_{\partial^*\widehat E(t)}\bigr|^2
 \,d\mu_{\widehat E(t)}\,dt
 \leq P^\phi(\widehat E(t_1)).
\end{equation}
Finally, if $\mathcal L^{N}(\partial E_0)=0$, then $\chi_{\widehat E(t)}\to\chi_{E_0}$ in $L^1_{\loc}(\R^N)$ as $t\downarrow0$ and $V_0=\var(\partial^*E_0)$.
\end{theoremletter}

Theorem~\ref{intro:brakke} combines Theorems~\ref{brakkebv} and~\ref{mob:brakkebv} with Lemma~\ref{lem:flat-initial-trace}. The additional assumptions on $\phi$ allow us to
identify the normal velocity in Theorem~\ref{intro:transport} with the anisotropic curvature law \eqref{intro:general-law}. They also identify the surface-energy transport field.
Indeed, \eqref{intro:curvature-identification} gives
\[
 \mathbf b=\frac{v}{\phi(\nu)}\nabla\phi(\nu).
\]
Thus $\mathbf b$ is directed along the Cahn--Hoffman vector, which need not be parallel to $\nu$. Euler's identity $\nabla\phi(\nu)\cdot\nu=\phi(\nu)$ ensures that $\mathbf
b\cdot\nu=v$, as required by the phase transport identity.

The identification of $h$ recovers the kinetic cost dictated by the prescribed mobility. Using $\phi^o(\nabla\phi(\nu))=1$ and $d\sigma=\phi(\nu)\,d\mu$, we obtain
\[
\frac{\bigl(\phi^o(\mathbf b)\bigr)^2}{h}\,d\sigma
 =\frac{v^2}{h\,\phi(\nu)^2}\,d\sigma=\frac{v^2}{\mathfrak m(\nu)}\,d\mu
 =\mathfrak m(\nu)
   \bigl|\mathbf H^\phi_{\partial^*\widehat E(t)}\bigr|^2
   \,d\mu.
\]
The last equality uses the perpendicularity of $\mathbf H^\phi_{\partial^*\widehat E(t)}$ and the formula for $v$ in \eqref{intro:curvature-identification}. Consequently,
\eqref{intro:optimal-action} is the weak energy-dissipation counterpart of \eqref{intro:smooth-dissipation}, with the same quadratic cost on normal velocities.

The Brakke inequality \eqref{intro:brakke-ineq} localizes this dissipation law: the term containing $\nabla_x\eta$ accounts for the transport of surface energy through the region
selected by the test function, while the negative curvature term accounts for dissipation. The identity \eqref{intro:bv-identity} describes the motion of the enclosed phase using
the same curvature field. The two conclusions therefore couple the evolution of the interface to that of the phase it bounds. Following the terminology of
\cite{StuvardTonegawa2024}, this coupling between the Brakke flow and the evolution of the phases gives rise to a canonical Brakke flow. In the general setting of
\cite{StuvardTonegawa2024}, the Brakke flow may carry higher multiplicities and interfacial measure beyond the reduced boundaries of the phases. In the present setting, non-fattening in measure and the matching procedure yield the further identification of the varifold weight with $\mu_{\widehat E(t)}$, with multiplicity one. Hence, there is no additional
interfacial measure beyond the reduced boundary of the phase, and the BV identity \eqref{intro:bv-identity} is expressed through the generalized anisotropic mean curvature of
$\partial^*\widehat E(t)$ itself, rather than through the curvature of a Brakke enhancement whose weight merely dominates the perimeter measure of the flat flow.

We emphasize that no $C^2$ regularity or uniform convexity is required of $\phi$, that the mobility $\mathfrak m$ can be an arbitrary, possibly crystalline norm, and that the
result holds in every dimension $N\geq2$ without imposing an energy-convergence hypothesis on the elliptic approximations.  The perpendicularity of the generalized anisotropic mean
curvature is obtained directly from the construction. To our knowledge,  this is the first unconditional construction of unit density anisotropic Brakke flows holding
for general $C^1$ strictly convex surface tensions. Unit density is the key condition to establishing Allard- and Brakke-type regularity results, which are known in both the isotropic case \cite{Allard1972,Brakkephd,KasaiTonegawa2014,StuvardTonegawa2024b}, and in the anisotropic time-independent one \cite{allard-regularity,DRT,DeRosaFiresterTsiamis2026}. The validity of Brakke's regularity theorem for unit density anisotropic Brakke flows is the subject of work in preparation.

\subsection{Why general mobility requires a sequence of extensions} \label{intro:seq-ext}

In the natural mobility case $\mathfrak m=\phi$, the single Euclidean extension $\Phi(\xi,z)=\sqrt{\phi(\xi)^2+z^2}$ satisfies $\partial_z\Phi(\xi,z)=z/\Phi(\xi,z)$, so that the
surface tension and mobility weights coincide on $\R^{N+1}$. To treat an independent mobility $\mathfrak m$ via elliptic regularization, one would ideally seek a fixed pair of
norms $(\Phi,\mathfrak M)$ on $\R^{N+1}$ extending $(\phi,\mathfrak m)$ and satisfying the compatibility relations
\begin{equation}\label{intro:fixed-compatible}
 \Phi(\xi,0)=\phi(\xi),\qquad
 \mathfrak M(\xi,0)=\mathfrak m(\xi),\qquad
 \partial_z\Phi(\xi,z)=\frac{z}{\mathfrak M(\xi,z)}.
\end{equation}
Here, the first two conditions in \eqref{intro:fixed-compatible} prescribe the surface tension and the mobility along horizontal directions, ensuring that vertical cylinders
$E(t)\times\R$ evolve according to the original pair $(\phi,\mathfrak m)$ in the limit $\varepsilon\to 0^+$. The third identity is instead dictated by the translating-soliton
mechanism at the level of the $\varepsilon$-approximations. Indeed, by the Euler--Lagrange equation \eqref{euler:lagrange1} for the weighted functional $\mathcal
A^\Phi_\varepsilon$, the generalized $\Phi$-anisotropic mean curvature of a minimizer $E_\varepsilon$ in $\{z>0\}$ satisfies
\[
 \mathbf H^\Phi_{\partial^*E_\varepsilon}
 =
 -\frac{\partial_z\Phi(\nu_{E_\varepsilon})}{\varepsilon}
 \,\nu_{E_\varepsilon}.
\]
On the other hand, translating $E_\varepsilon$ vertically with velocity $-\varepsilon^{-1}e_{_{N+1}}$ imparts the normal speed $v_t = -\varepsilon^{-1}\nu_{E_\varepsilon}\cdot
e_{_{N+1}}$ to the moving boundary $\partial^*\mathcal S_\varepsilon(t)$. Requiring the translating solitons $\mathcal S_\varepsilon(t)$ to obey the $(\Phi,\mathfrak M)$-mean
curvature flow in $\R^{N+1}$ already before passing to the cylindrical limit — and thus also away from the horizontal equator, where $\nu_{E_\varepsilon}$ has a nonzero vertical
component — amounts to having
\[
 -\frac{\nu_{E_\varepsilon}\cdot e_{_{N+1}}}{\varepsilon}
 =
 \mathfrak M(\nu_{E_\varepsilon})
 \,\nu_{E_\varepsilon}\cdot
 \mathbf H^\Phi_{\partial^*E_\varepsilon}
 =
 -\frac{\mathfrak M(\nu_{E_\varepsilon})
 \,\partial_z\Phi(\nu_{E_\varepsilon})}{\varepsilon},
\]
namely imposing $\partial_z\Phi(\nu)=\mathfrak M^{-1}(\nu)\, \nu \cdot e_{N+1}$ for every unit vector $\nu \in \R^{N+1}$: by the $1$-homogeneity of $\mathfrak M$ and the $0$-homogeneity of $\partial_z\Phi$, this yields the third relation in \eqref{intro:fixed-compatible} for every
$(\xi,z)\in\R^{N+1}\setminus\{0\}$.

However, such a fixed pair of convex norms $(\Phi,\mathfrak M)$ does not exist in general. Consider, for instance, the pair $\phi(\xi)=\|\xi\|_\infty$ and $\mathfrak m(\xi)=|\xi|$
in $\R^N$ with $N\geq 2$. If norms $\Phi,\mathfrak M$ satisfying \eqref{intro:fixed-compatible} existed, integrating in $z$ would yield the expansion
\[
 \Phi(\xi,z)=\phi(\xi)+\frac{z^2}{2\mathfrak m(\xi)}+o(z^2)
 \qquad\text{as }z\to 0
\]
for every $\xi\neq0$. Choosing $\xi_\pm=(1,\pm a,0,\ldots,0)$ and $\xi_0=\frac12(\xi_++\xi_-)=(1,0,\ldots,0)$ with $a\in(0,1)$, we have $\phi(\xi_\pm)=\phi(\xi_0)=1$, whereas
$\mathfrak m(\xi_0)=1$ and $\mathfrak m(\xi_\pm)=\sqrt{1+a^2}$. Convexity of $\Phi(\cdot,z)$ would then require $\Phi(\xi_0,z)\leq\frac12\bigl(\Phi(\xi_+,z)+\Phi(\xi_-,z)\bigr)$,
which for small $z\neq 0$ reads
\[
 1+\frac{z^2}{2}+o(z^2)
 \leq 1+\frac{z^2}{2\sqrt{1+a^2}}+o(z^2),
\]
a contradiction.

To bypass this obstruction, we construct in Appendix~\ref{app:ext:norms} a one-parameter family of smooth, uniformly convex compatible extensions $(\Phi^\kappa,\mathfrak
M^\kappa)_{\kappa\in(0,\kappa_0)}$ on $\R^{N+1}$ satisfying
\begin{equation}\label{intro:compatible-extension}
 \partial_z\Phi^\kappa(\xi,z)=\frac{z}{\mathfrak M^\kappa(\xi,z)},\qquad
 \Phi^\kappa(\xi,z)\xrightarrow{\kappa\to 0^+}\phi(\xi),\qquad
 \mathfrak M^\kappa(\xi,0)\xrightarrow{\kappa\to 0^+}\mathfrak m(\xi).
\end{equation}
Because $\Phi^\kappa(\xi,z)\to\phi(\xi)$, the extended surface tension $\Phi^\kappa$ necessarily degenerates in the vertical direction $e_{_{N+1}}$ as $\kappa\downarrow 0$, while
the compatible extended mobility $\mathfrak M^\kappa$ blows up away from the horizontal equator $\{z=0\}$. Consequently, the $\Phi^\kappa$-perimeter and $\mathfrak
M^\kappa$-mobility measures on $\R^{N+1}$ are no longer uniformly comparable as $\kappa\downarrow 0$.

In Section~\ref{sec:gen:mob}, we overcome this degeneracy by comparing the full $(N+1)$-dimensional measures with their \emph{equatorial} counterparts, which evaluate
$\phi^\kappa(\nu')$ and $\mathfrak m^\kappa(\nu')$ on the horizontal projection $\nu'$ of the normal. Using quantitative quasi-cylindricity estimates derived from the soliton
equation, we show that along a suitably slow diagonal regime $\kappa=\kappa(\varepsilon)\to 0^+$, the full and equatorial measures share the same weak-* limits. 

\medskip
\noindent\textbf{Organization of the paper.} In Section~\ref{s:preliminaries} we fix the notation and review sets of finite perimeter, varifolds, distributional and level-set
solutions of anisotropic curvature flows, flat flows, and anisotropic Brakke flows. Section~\ref{sec:waf} is devoted to the weighted anisotropic area functional in $\R^{N+1}$, its
Euler--Lagrange equations, boundary contact properties at $\{z=0\}$, and anisotropic paraboloid barriers. In Section~\ref{sec:kur:conv} we study the elliptic regularization scheme
for the mobility $\mathfrak m=\phi$ and prove the Kuratowski convergence results of Theorem~\ref{intro:kuratowski}(a)--(c). In Section~\ref{sec:ellip:varifold} we establish the
uniform energy estimates, the $BV_{\loc}$ convergence of Theorem~\ref{intro:kuratowski}(d)--(e), and the restarting-and-matching construction leading to
Theorems~\ref{intro:transport} and~\ref{intro:brakke} for natural mobility. Section~\ref{sec:gen:mob} treats general mobilities using the compatible extensions
$(\Phi^\kappa,\mathfrak M^\kappa)$. Finally, the Appendices collect several technical properties of uniformly convex norms and the construction of the compatible extensions
$(\Phi^\kappa,\mathfrak M^\kappa)$.

\smallskip \noindent \textbf{AI Statement.} The vast majority of this work had been completed already prior to June 2026, with no input from AI, and in fact already presented by the authors in conferences and seminars. At the time, we had the complete result in the natural mobility regime, and the strategy for attacking the general mobility case. We had a complete understanding of the properties needed for the family of extended norms, but we were lacking an explicit construction. This was suggested to us by ChatGPT 5.6 Sol and carefully re-written and elaborated by us: the content of our interaction with AI is then confined to Appendix \ref{app:ext:norms}.

\smallskip \noindent \textbf{Acknowledgments.}
The work of CAA and SS is supported by the Italian Ministry of University and Research (MUR) through the FIS 2 StG \textbf{SiGmA}: “\textit{Singularities in Geometric Analysis:
Minimal Surfaces and Mean Curvature Flows}”, project code FIS-2023-02962 (CUP G53C25000120001). ADR was funded by the European Union: the European Research Council (ERC), through
StG \textbf{ANGEVA}, project number: 101076411. MM 
has been partially supported by PRIN 2022 Project “\textit{Geometric Evolution Problems and Shape Optimization (GEPSO)}”, PNRR Italia Domani, financed by European Union via the Program NextGenerationEU, CUP D53D23005820006. SS and MM are members of Gruppo Nazionale per l’Analisi Matematica, la Probabilità e le loro Applicazioni (GNAMPA), which is part of the Istituto Nazionale di
Alta Matematica (INdAM). Views and opinions expressed are however those of the authors only and do not necessarily reflect those of the European Union or the
European Research Council. Neither the European Union nor the granting authority can be held responsible for them.

\section{Notation and Preliminaries} \label{s:preliminaries}

The ambient space for our anisotropic flows is Euclidean space $\R^{N}$, with coordinates $x=(x_1,\ldots,x_N)$. We write $\R_+$ for the ``time interval'' $[0,\infty)$, and the
parameter $t \in \R_+$ will be referred to as ``time''. As in \cite{Ilm94}, the minimization problem for elliptic regularization will be set in the \textit{extended space}
$\R^{N+1}=\R^N_x\times\R_z$. For this reason, the first few sections of the paper will be set in an ``abstract'' Euclidean space $\R^d$, and then we will specialize to the cases
$d=N$ or $d=N+1$ upon introducing elliptic regularization.

We denote by $e_i=(0,\ldots,0,1,0,\ldots,0)$, with $1$ in the $i$-th position, the $i$-th vector of the standard basis of $\R^{N+1}$. In particular, $e_{_{N+1}}=(0,\dots,0,1)\in
\R^{N+1}$.

We denote by $B_r(x)$ the $d$-dimensional open ball of radius $r$ centered at $x$ in the abstract space $\R^d$, and the same notation will be in use later when $d=N$; instead, we
will use $\B_r(x,z)$ for the $(N+1)$-dimensional open ball of radius $r$ centered at $(x,z)$ in the extended space $\R^{N+1}$. In all cases, when we omit the center it is
understood that the ball is centered at the origin.

For $A\subset\mathbb R^{d}$ (with $d=N$ or $d=N+1$), ${\rm clos}\,A$ (or $\overline A$) is the topological closure of $A$ in $\R^d$ and ${\rm int}\,A$ is the set of interior points
of $A$. The standard Euclidean inner product between vectors in $\R^{d}$ is denoted $u \cdot v$, and $\abs{u} := \sqrt{u \cdot u}$. We also denote by $\mathbb{{S}}^{d-1}=\{x\in
\R^d:|x|=1\}$ the unit sphere in $\R^d$.

If $L,S \in \mathscr{L}(\R^{d};\R^{d})$ are linear operators in $\R^{d}$, their (Hilbert-Schmidt) inner product is $L \colon S := {\rm trace}(L^{^{\top}} \circ S)$, where
$L^{^{\top}}$ is the transpose (adjoint) of $L$ and $\circ$ denotes composition. The corresponding (Euclidean) norm in $\mathscr{L}(\R^{d};\R^{d})$ is then $\abs{L} := \sqrt{L
\colon L}$, whereas the operator norm in $\mathscr{L}(\R^{d};\R^{d})$ is $\|L\| := \sup\left\lbrace \abs{L(u)}\,\colon\, u \in \R^d\,, |u| \leq 1\right \rbrace$. If $u,v \in
\R^{d}$ then $u \otimes v \in \mathscr{L}(\R^{d};\R^{d})$ is defined by $(u \otimes v)(w) := (w \cdot v)\, u$, so that $\| u \otimes v \| = \abs{u}\,\abs{v}$. The Lebesgue measure
of a set $A \subset \R^{d}$ is denoted $\Leb^{d}(A)$ or $|A|$. We will set $\omega_d := \Leb^d(\{u \in \R^d\,\colon\, |u|\leq 1\})$. For $0 \leq k \leq d$, the symbol $\Ha^k$
denotes the $k$-dimensional Hausdorff measure in $\R^{d}$, normalized in such a way that $\Ha^{d}$ and $\Leb^{d}$ coincide as measures.

We will also use standard notation concerning Radon measures on $\R^d$: in particular, if $\mu$ is a Radon measure on $\R^d$ and $F \subset \R^d$ is a Borel subset then $\mu \mres
F$ denotes the restriction of $\mu$ to $F$. Given two sets $E, F \subset \R^d$, the writing $E \stackrel{\mu}{=} F$ signifies that $\mu(E \Delta F)=0$, where $E \Delta F := (E
\setminus F) \cup (F \setminus E)$ is the symmetric difference between $E$ and $F$. When $\mu$ is clear from the context, we will also use $E \stackrel{a.e.}{=} F$ with the same
significance. We will use standard notation to identify spaces of functions: in particular, $\mathrm{Lip}(X;Y)$ will denote the space of Lipschitz continuous functions between
metric spaces $X$ and $Y$, whereas $C^k$ will denote spaces of functions with continuous derivatives of all orders $0\leq \ell \leq k$; a subscript $_\mathrm{c}$ means that
functions are required to have compact support in their domain; finally, $L^p(X;d\mu)$ consists of all $\mu$-measurable functions $f$ on $X$ so that $|f|^p$ is $\mu$-summable.

We shall sometimes resort to the use of regularization by convolution: we will denote by $\rho$ the standard mollifier in $\R^d$, so that
\[
\rho(x)=
\begin{cases}
   c_d\,\exp\left(-\frac{1}{1-|x|^2} \right) & \mbox{if $|x| < 1$} \\
    0 & \mbox{if $|x| \geq 1$}\,,
\end{cases}
\]
with $c_d$ so that $\int_{\R^d} \rho \,dx=1$; we will also set $\rho_\e(x):=\e^{-d}\rho(x/\e)$ whenever $\e> 0$.

\subsection{Rectifiable sets and sets of finite perimeter}\label{sec:definiscorect}

In this paper, a set $M \subset \R^d$ is called (countably) $k$-rectifiable if it is $\Ha^k$-measurable, $\Ha^k(M\cap K)<\infty$ for all compact subsets $K\subset \R^d$, and $M$
can be covered, up to a $\H^k$-null set, by countably many Lipschitz images of $\R^k$. If $M$ is $k$-rectifiable, an approximate tangent space to $M$ exists for $\H^k$-a.e. $x \in
M$, and it is denoted $\mathrm{Tan}(M,x)$. A Lebesgue measurable set $E \subset \R^d$ is a set of locally finite perimeter if for every bounded open set $A \subset {\R^{d}}$ the
relative perimeter of $E$ in $A$ is finite, namely
\[
P(E;A):=\sup \left\{\int_E \text{div} X(x)\, d x\ : \ X \in C_{\mathrm{c}}^1({A}; \R^{d}),\, \sup |X| \leq 1\right\}<\infty\, .
\]
If $P(E;\R^d)<\infty$, we say that $E$ is a set of finite perimeter. If $E$ is of locally finite perimeter, then the indicator function $\chi_E \in
\mathrm{BV}_{\mathrm{loc}}(\R^d)$, and thus the distributional derivative $\nabla\chi_E$ is a vector-valued Radon measure such that
\[
\int_E \text{div} X\, d x = -\int_{\R^d} X \cdot d\nabla\chi_E \qquad \mbox{for every $X \in C^1_{\mathrm c}(\R^d;\R^d)$}\,.
\]
The relative perimeter of $E$ in $A$ is then
\begin{equation} \label{def:perimeter}
P(E;A) = |\nabla\chi_E|(A)\,,
\end{equation}
where $\mu_E:=|\nabla\chi_E|$ is the total variation measure of $\nabla\chi_E$. We will call $\mu_E$ the perimeter measure of $E$.

For $E$ of locally finite perimeter in $\R^d$, the {reduced boundary $\partial^* E$ of $E$} and the {(generalized) outer unit normal $\nu_E(x)$ to $E$ at $x$} are then defined by
\[
    \partial^* E := \left \{x\in\spt(\mu_E): \nu_E(x) := -\lim_{r\to 0^+}\frac{\nabla\chi_E(B_r(x))}{\mu_E(B_r(x))}\mbox{ exists and belongs to $\mathbb S^{d-1}$} \right\}.
\]
By De Giorgi's structure theorem, for a set $E$ of locally finite perimeter, $\partial^* E$ is $(d-1)$-rectifiable and
\begin{equation}
    -\nabla\chi_E = \nu_E\,\H^{d-1}\res\partial^*E,\quad \mu_E =\H^{d-1}\res\partial^*E;
\end{equation}
 see \cite[Remark 15.3 and Corollary 16.1]{MaggiBook}. For brevity, we will refer to $\nu_E$ as the \textit{normal} to $E$ or, equivalently, to $\partial^*E$, whereas
 $-\nabla\chi_E=\nu_E\,\Ha^{d-1} \mres \partial^*E$ will be called the \textit{Gauss-Green measure} associated with $E$. Notice that, for any Borel set $A\subset \R^d$, we have
 \begin{equation}\label{def:perclassico}
 P(E;A)= \mu_E(A)=\H^{d-1}(\partial^* E\cap A)\,.
 \end{equation}
For a set $E$ of locally finite perimeter, we recall the definition of the sets $E^{(\theta)}$ of points of density $\theta$, namely
\begin{equation} \label{def:density-pts}
   E^{(\theta)} = \left\lbrace x \in \R^d\,\colon\, \lim_{r \to 0^+} \frac{\Leb^d(E\cap B_r(x))}{\omega_d\, r^d}=\theta\right\rbrace\,,
\end{equation}
as well as the definition of the essential boundary
\begin{equation} \label{def:essential-bdry}
    \partial^e E:=\R^d \setminus (E^{(1)} \cup E^{(0)})\,,
\end{equation}
so that
\begin{equation} \label{bdry-inclusions}
    \partial^*E \subset E^{(1/2)} \subset \partial^e E \subset \partial E\,,
\end{equation}
and Federer's theorem holds:
\begin{equation} \label{eq:Federer}
    \H^{d-1}(\partial^e E \setminus \partial^* E) = 0\,.
\end{equation}

 Given a norm $\psi:\R^d\to [0,\infty)$ and a set $E$ of locally finite perimeter, we define the $\psi$-\textit{anisotropic perimeter measure $\mu_E^\psi$ of $E$} and the
 \textit{$\psi$-anisotropic perimeter $P^\psi(E;\cdot)$}  by the formula
\begin{equation}\label{def:peranis}
    d\mu^\psi_E:=\psi(\nu_E)\,d\mu_E\,,\quad  P^\psi(E;A):=\mu^\psi_E(A) =\int_{\partial^* E\cap A}\psi(\nu_E)\,d\H^{d-1}
\end{equation}
for any Borel set $A\subset \R^d$. Hence we have
\begin{equation}\label{perimetri:Equivalent}
\begin{split}
\Big(\inf_{\S^{d-1}} \psi\Big)\,\mu_E&\leq \mu^\psi_E\leq \Big(\sup_{\S^{d-1}} \psi\Big)\,\mu_E\quad\text{as (positive) Radon measures,} \\
&\text{and thus} \\
    \Big(\inf_{\S^{d-1}} \psi\Big)\,P(E;A)&\leq P^\psi(E;A)\leq \Big(\sup_{\S^{d-1}} \psi\Big)\,P(E;A)\quad\text{for any Borel set $A\subset \R^d$.}
    \end{split}
\end{equation}

When $A=\R^d$, we will simply write $P(E),\,P^\psi(E)$.

\subsection{Varifolds}\label{s:var:def}
In this paper, we will always work with codimension one varifolds only. Let $\Omega \subset \R^d$ be an open set. A $(d-1)$-varifold on $\Omega$ is a positive Radon measure $V$ on
the Grassmannian bundle $\Omega \times G(d,d-1)$, where $G(d,d-1)$ denotes the Grassmannian of (unoriented) hyperplanes in $\R^d$. As in \cite[Section 5]{DPDRG}, we shall always
identify $G(d,d-1)$ with the real projective space $\R\mathbb{P}^{d-1}$ via the map
\[
[\nu] \in \mathbb{RP}^{d-1} \mapsto \nu^{\perp}\in G(d,d-1)\,.
\]
For any measurable set $A\subset \Omega$, we will write
\[
G_{d-1}(A):=A\times G(d,d-1)\simeq A\times \mathbb{RP}^{d-1}.
\]
Given a unit vector $\nu\in \S^{d-1}$, we will denote by $[\nu]\in \mathbb{RP}^{d-1}$ its projection onto $\mathbb{RP}^{d-1}$. With some abuse of notation, given a function
$\varphi\in C^0_c(G_{d-1}(\Omega))$ we will write $\varphi=\varphi(x,\nu)$ in place of the more cumbersome $\varphi\big(x,[\nu]\big)=\varphi(x,\nu^{\perp})$, and thus when $V$ is a
$(d-1)$-varifold on $\Omega$ its action as a continuous linear functional will be denoted
\begin{equation}\label{this:is:the:notation}
    V(\varphi)=\int_{ G_{d-1}(\Omega)}\varphi(x,\nu)\,dV(x,\nu) \qquad \mbox{for $\varphi \in C^0_{\mathrm{c}}(G_{d-1}(\Omega))$}\,.
\end{equation}
When it is clear from the context, we shall also omit the domain of integration of the varifold.

The weight $\|V\|$ of $V$ is the Radon measure on $\Omega$ defined by
\begin{equation}\label{def:weightV}
    \|V\|(A) := V\big( G_{d-1}(A)\big), \quad \mbox{for every Borel set } A\subset \Omega.
\end{equation}
We will write $\mathbf{V}_{d-1}(\Omega)$ for the space of $(d-1)$-varifolds on $\Omega$. A $(d-1)$-varifold $V$ is \textit{rectifiable}, and we write $V \in
\mathbf{RV}_{d-1}(\Omega)$ if there exist a $(d-1)$-rectifiable set $M$ and a function $\theta \in L^1_{\mathrm{loc}}(\Ha^{d-1}\mres (M \cap \Omega);(0,\infty))$, which will be
referred to as \emph{multiplicity}, such that
\begin{equation}\label{if:ecos}
    V \equiv \mathbf{Var}(M,\theta)  := \theta\,\H^{d-1}\res (M\cap \Omega) \otimes \delta_{\mathrm{Tan}(M,x)}\,,
\end{equation}
where $\delta_a$ denotes the Dirac delta concentrated at $a \in G(d,d-1)$. Upon identifying $\mathrm{Tan}(M,x)$ with $[\nu_M(x)]\in \mathbb{RP}^{d-1}$, we call $\nu_M(x)$ the
generalized unit normal to $M$. We warn the reader that we are not choosing an orientation on $M$, as $\nu_M(x)$ is simply a short-hand notation for its projection onto
$\mathbb{RP}^{d-1}$: throughout the paper, $\nu_M$ will only appear in formulas which descend to the quotient space, that is, which are invariant with respect to flipping $\nu_M
\mapsto -\nu_M$. For brevity, when $V=\mathbf{Var}(M,\theta)$ we call $\nu_M$ the \textit{normal to $V$}.

In the case when $\theta\equiv 1$, we will simply write $\mathbf{Var}(M)$, i.e.
\[
\mathbf{Var}(M)=\H^{d-1}\res (M\cap \Omega) \otimes \delta_{\nu^\perp_M(x)}.
\]

Let us observe that, if $V=\var(M,\theta)\in \mathbf{RV}_{d-1}(\Omega)$, its weight measure is
\begin{equation}\label{weigh:Vrect}
\|V\|=\theta\,\H^{d-1}\mres (M\cap \Omega),
\end{equation}
and for any $\varphi\in C^0_c\big(G_{d-1}(\Omega)\big)$, we have
\begin{equation}\label{test:Vrect}
    V(\varphi)=\int_{M}\theta(x)\,\varphi\big(x,\nu_M(x)\big)\,d\mathcal{H}^{d-1}(x)=\int_\Omega\varphi\big(x,\nu_M(x)\big)\,d\|V\|(x).
\end{equation}
If $V$ is rectifiable and the multiplicity $\theta$ is integer-valued then we say that $V$ is integral and we write $V \in \mathbf{IV}_{d-1}(\Omega)$.

When $\Omega$ and $\Omega'$ are open subsets of $\R^d$, $f \colon \Omega \to \Omega'$ is a $C^1$ diffeomorphism, and $V \in \mathbf{V}_{d-1}(\Omega)$, the {push-forward $f^\sharp
V$ of $V$  with respect to $f$} is the $(d-1)$-varifold on $\Omega'$ defined by

\begin{equation}\label{E:pushforward-varifold}
    \int_{G_{d-1}(\Omega')}\varphi(x,\nu)\,d(f^\sharp V)(x,\nu) = \int_{G_{d-1}(\Omega)}\varphi\left(f(x),\tfrac{(Df_x)^{^-\top}(\nu)}{|(Df_x)^{^-\top}(\nu)|}\right)Jf(x,\nu) \, dV(x,\nu)\,,
\end{equation}
for all $\varphi\in C^0_c(G_{d-1}(\Omega'))$. Here $Df_x$ is the differential of $f$ at $x$, $(Df_x)^{^-\top}$ is the transpose of its inverse and $Jf(x,\nu)$ is the Jacobian
determinant of the differential $Df_x$ restricted to $\nu^\perp$, i.e.,
\[
Jf(x,\nu) :=\sqrt{\det\Big(\big(\left.Df_x\right|_{\nu^\perp}\big)^{^\top} \circ \left.Df_x\right|_{\nu^\perp}\Big)}\,.
\]
The upper and lower $(d-1)$-dimensional density of a Radon measure $\mu$ on $\R^{d}$ at $x$ are,  respectively,
\begin{equation*}
    \udens{d-1}{\mu}{x} := \limsup_{r\to 0^+} \frac{\mu(B_r(x))}{\omega_{d-1} r^{d-1}}\qquad \mbox{ and }\qquad \ldens{d-1}{\mu}{x} := \liminf_{r\to 0^+} \frac{\mu(B_r(x))}{\omega_{d-1} r^{d-1}}.
\end{equation*}
If $\udens{d-1}{\mu}{x} = \ldens{d-1}{\mu}{x},$ we call this number the {$(d-1)$-density of $\mu$ at $x$} and denote it by $\dens{d-1}{\mu}{x}$. { For a $(d-1)$-varifold, we let
$\udens{d-1}{V}{x}:=\udens{d-1}{\|V\|}{x}$ and for a $(d-1)$-rectifiable set $M$ we let $\udens{d-1}{M}{x}:=\udens{d-1}{\H^{d-1}\res M}{x}$ (and likewise for the lower density and
the density). We will omit the dimension when there is no chance of confusion. } \vspace{0.2cm}

Let now \(\psi:\R^d\to [0,\infty)\) be an (even) norm on \(\R^d\). Then \(\psi\) induces a well-defined function on the quotient \(\mathbb{RP}^{d-1}\), which we still denote by
\(\psi\), via
\[
\psi([\nu]) \equiv \psi (\nu) \,, \qquad [\nu] \in \mathbb{RP}^{d-1}\,.
\]

 For a given $V \in \mathbf V_{d-1}(\Omega)$, we define the $\psi$-anisotropic weight measure of $V$ by
\begin{equation} \label{e:an_en_var}
\|V\|^\psi(\varphi) := \int_{G_{d-1}(\Omega)} \varphi(x)\psi(\nu) \, d V(x,\nu),\quad \varphi \in C^0_c(\Omega).
\end{equation}
In particular, we have
\begin{equation}\label{weight:comparison}
    \Big(\inf_{\S^{d-1}}\psi\Big)\,\|V\|\leq \|V\|^\psi\leq\Big(\sup_{\S^{d-1}}\psi\Big)\, \|V\|\quad\text{as Radon measures on $\Omega$,}
\end{equation}
and for $V\in \mathbf{RV}_{d-1}(\Omega)$, we may write
\begin{equation}\label{two:weights}
    d\|V\|^\psi=\psi(\nu)\,d\|V\|,
\end{equation}
where $\nu=\nu_V$ denotes the normal to $V$.

Let us also observe that, given a set $E\subset \Omega$ of locally finite perimeter, and setting $V=\mathbf{Var}(\partial^* E,\theta)$, then according to the notation of Section
\ref{sec:definiscorect}, we have $\nu_V=\nu_E$ (as elements of $\mathbb{RP}^{d-1}$), and
\begin{equation}\label{varifold:perimeter}
    \|V\|=\theta\,\mu_E,\quad\text{and}\quad\|V\|^\psi=\theta\,\mu_E^\psi
\end{equation}
as Radon measures on $\Omega$.

Assume now, in addition, that $\psi \in C^1(\R^d\setminus\{0\})$.  Throughout the paper, $\nabla\psi$ will always denote the gradient of $\psi=\psi(\xi)$ w.r.t.  its original
variable $\xi$. For instance, whenever $v=v(x)$ is a non-vanishing vector field, we will always use the convention $\nabla\psi(v(x))\equiv \nabla\psi(\xi)|_{\xi=v(x)}$.

By letting $\mathbf{F}^\psi$ denote the $\psi$-anisotropic energy
\[
\mathbf{F}^\psi(V):=\|V\|^\psi(\Omega)=\int_{G_{d-1}(\Omega)}\psi(\nu)\,dV(x,\nu),
\]
the first variation $\delta_\psi V$ of $V$ with respect to the energy $\mathbf{F}^\psi$ is the order one distribution defined as follows. Given any vector field $X \in
C^1_c(\Omega;\R^d)$ with $\spt\,X\subset \widetilde\Omega \Subset \Omega$, we let $f_s(x):=x+s\,X(x)$ for small $|s|$, and set $$\delta_{\psi} V(X):=\left.\frac{d}{ds}\right|_{s=0}
\mathbf{F}^\psi\big( (f_s)^\sharp (V \mres G_{d-1}(\widetilde\Omega))\big)\,. $$ A classical computation {(cf. \cite[Section 5]{DPDRG})} shows that
$\delta_\psi V$ is independent of the choice of $\widetilde\Omega$ and
\begin{equation}\label{E:first-variation-Ani0}
\begin{aligned}
	\delta_{\psi} V(X) &= \int_{G_{d-1}(\Omega)} B_\psi(\nu): DX(x)\, d V(x,\nu) \,,
\end{aligned}
\end{equation}
where
\begin{equation}\label{def:Bpsi}
    B_\psi(\nu) := \psi(\nu)\, \text{Id}  -\nu \otimes \nabla\psi(\nu).
\end{equation}
We record here that, by definition and by the properties of $\psi$,
\begin{equation}\label{prop-Bpsi}
    B_\psi(-\nu)=B_\psi(\nu)\,, \quad \mbox{and} \quad B_\psi(\nu)\nu=0\,.
\end{equation}
In particular, $B_\psi$ is well defined on $\mathbb{RP}^{d-1}$. If $\delta_{\psi} V$ can be extended to a linear and continuous functional on $C_c(\Omega;\R^{d})$, we say that $V$
has {\it bounded first variation with respect to $\mathbf{F}^\psi$} in $\Omega$. In this case, $\delta_{\psi} V$ is naturally associated with a unique $\R^{d}$-valued measure on
$\Omega$ by means of the Riesz representation theorem. If such a measure is absolutely continuous with respect to the weight $\|V\|$, then there exists a $\|V\|$-measurable and
locally $\|V\|$-integrable vector field $\mathbf{H}^\psi_V$ such that
\begin{equation} \label{def:generalized mean curvature}
\delta_\psi V(X) = - \int_{\Omega} X(x) \cdot \mathbf{H}^\psi_V \, d\|V\|(x) \qquad \mbox{for every $X \in C^0_c(\Omega;\R^{d})$}
\end{equation}
by the Radon-Nikod\'ym theorem. The vector field $\mathbf{H}^\psi_V$ is called the {\it generalized $\psi$-anisotropic mean curvature vector} of $V$.

If $E\subset \Omega$ is a set of locally finite perimeter, we say that it admits (generalized) $\psi$-anisotropic mean curvature, denoted by $\mathbf{H}^\psi_{\partial^*E}$,
whenever $\mathbf{Var}(\partial^* E)$ does. In case $E\subset\Omega$ has boundary of class $C^2$ and also $\psi \in C^2(\R^d \setminus \{0\})$, we also define its \emph{scalar}
$\psi$-anisotropic mean curvature by
\begin{equation} \label{scalarH}
    H^\psi_{\partial E} :=\Div_{\partial E}\bigl(\nabla\psi(\nu_E)\bigr).
\end{equation}
It is then classical to see that
\[
    \delta_\psi\mathbf{Var}(\partial E)(X) = \int_{\partial E} H^\psi_{\partial E}\,(X\cdot\nu_E)\, d\mathcal H^{d-1},
\]
and hence, consistently with \eqref{def:generalized mean curvature},
\begin{equation}\label{vecscalH}
    \mathbf H^\psi_{\partial E} =-H^\psi_{\partial E}\nu_E\,.
\end{equation}
For a general set of finite perimeter, the generalized vector $\mathbf H^\psi_{\partial^*E}$ is not known a priori to be normal to $\partial^*E$. We will establish this property
for the varifolds arising in our construction: see \eqref{identifica}.

Whenever this orthogonality property holds, we can extend the definition of scalar $\psi$-anisotropic mean curvature and set
\begin{equation} \label{scalH-SOFP}
    H^\psi_{\partial^*E} :=-\mathbf H^\psi_{\partial^*E}\cdot\nu_E.
\end{equation}

Other than the first variation $\delta_{\psi} V$ discussed above, we shall also use a {\it weighted first variation} with respect to $\mathbf{F}^\psi$, defined as follows. For
given $\varphi\in C^1_c(\Omega)$, $V\in {\bf V}_{d-1}(\Omega)$, and $X \in C^1_c(\Omega;\R^{d})$, we introduce the $\varphi$-weighted first variation of $V$ with respect to
$\mathbf{F}^\psi$ in the direction of $X$, denoted $\delta_{\psi}(V,\varphi)(X)$, by setting
\begin{equation} \label{defFV_modified}
\delta_\psi(V,\varphi)(X) := \left.\frac{d}{ds}\right|_{s=0}  \mathbf{F} ^\psi\big((\mathrm{Id}+sX)^\sharp V,\varphi\big)\,,
\end{equation}
having set the $\varphi$-weighted anisotropic energy
\begin{equation}\label{def:Fpesata}
\mathbf{F}^\psi(V,\varphi) :=\|V\|^\psi(\varphi)= \int_{G_{d-1}(\Omega)} \varphi(x)\, \psi(\nu) \, dV(x,\nu)\,.
\end{equation}
Proceeding as in the derivation of \eqref{E:first-variation-Ani0}, one then obtains the expression (with integrals on $G_{d-1}(\Omega)$),
\begin{equation}\label{defFV2}
\delta_{\psi}(V,\varphi)(X)=\int \varphi(x)\, B_\psi(\nu)\,\colon\,DX(x)\,dV(x,\nu)+ \int X(x)\cdot\nabla\varphi(x)\,\psi(\nu)dV(x,\nu)\,.
\end{equation}
Using $\varphi D X= D(\varphi X)- X\otimes\nabla\varphi$, we obtain, using the explicit expression in \eqref{def:Bpsi},
\begin{equation}\label{defFV3}
\begin{split}
\delta_{\psi}(V,\varphi)(X)&=\delta_{\psi} V(\varphi X)+\int X(x)\cdot(\nabla\varphi(x)\psi(\nu)\,-B_\psi(\nu)( \nabla\varphi(x)))\,dV(x,\nu) \\
&=\delta_{\psi} V(\varphi X)+\int \big(\nu\otimes\nabla\psi(\nu)\big):\big(X(x)\otimes \nabla\varphi(x)\big)\,dV(x,\nu)\,.
\end{split}
\end{equation}
If $V$ has generalized $\psi$-mean curvature $\mathbf{H}^\psi_V$, then we may use \eqref{def:generalized mean curvature} in \eqref{defFV3} to obtain
\begin{equation}\label{defFV4}
\delta_{\psi}(V,\varphi)(X)=-\int \varphi(x)X(x)\cdot \mathbf{H}^\psi_V \, d\|V\|(x)+\int \big(\nabla\varphi(x)\cdot \nabla\psi(\nu)\big)\,(X(x)\cdot \nu)\,dV(x,\nu)\,.
\end{equation}
The right-hand side of \eqref{defFV4} defines an extension of $\delta_\psi(V,\varphi)(X)$ to every pair $\varphi \in C^1(\R^d)$ and $X \in C^0(\R^d;\R^d)$ such that $\spt\,(\varphi
X) \Subset \Omega$.

\begin{remark}\label{remark:minimality}
    Let $\varphi \in C^1_c(\Omega)$, $g\in C^0_c(\Omega)$, and suppose $E\subset \Omega$ is a set of locally finite perimeter which is a local minimizer of the functional
    \begin{equation*}
        \mathcal{F}^\psi_{g}(F;\varphi):=\int_{\partial^*F}\varphi(x)\,\psi(\nu_F)\,d\H^{d-1}+\int_{F} g\,dx,
    \end{equation*}
    that is
    \[
    \mathcal{F}^\psi_{g}(E;\varphi)\leq \mathcal{F}_{g}^\psi(F;\varphi),\quad\text{for all $F\Delta E\Subset \Omega$.}
    \]
Then by the above minimality property we may apply a standard first variation argument which together with \cite[Proposition 17.8]{MaggiBook}, \eqref{varifold:perimeter} and the
definition \eqref{defFV_modified}, yields that $$\delta_\psi \big(\mathbf{Var}(\partial^* E),\varphi\big)(X)=-\int_{\partial^* E}g\,(X\cdot\nu_{E})\,d\H^{d-1}\quad\text{for all
$X\in C^1_c(\Omega;\R^d)$},$$ where $\nu_E$ is the outward unit normal to $\partial^* E$. Hence from \eqref{defFV2} and \eqref{varifold:perimeter}, we deduce
\begin{equation}\label{stationarity:fromminimality}
-\int_{\partial^* E}g\,(X\cdot\nu_{E})\,d\H^{d-1}=\int_{\partial^* E} \varphi(x)\,B_\psi(\nu_E): DX(x)\,d\mathcal{H}^{d-1}+\int_{\partial^* E} X(x)\cdot\nabla\varphi(x)\,\psi(\nu_E)\,d\mathcal{H}^{d-1}.
\end{equation}
The same identity holds in fact for $\varphi\in C^1(\Omega)$ and $g\in C^0(\Omega)$
if $E$ is locally minimizing for these weights and $X\in C^1_c(\Omega;\R^d)$.
Indeed, the variation may be computed in any open set compactly contained
in $\Omega$ and containing $\spt X$.
\end{remark}

\subsection{Distributional anisotropic  curvature flows}
\label{subsec:cf}

In this section we review the distributional and level-set settings for anisotropic (possibly crystalline) mean curvature flows with general mobility, as developed
in~\cite{CMP17,CMNP19}. Throughout the section, $\Psi$ and $\mathfrak M$ denote two norms on $\R^d$, playing the role of surface tension and mobility, respectively. The
corresponding geometric law is formally given by
\[
v_t=-\mathfrak M(\nu_{E(t)})H^\Psi_{\pa E(t)}\qquad\text{on }\pa E(t),
\]
with $v_t$ and $H^\Psi_{\pa E(t)}$ denoting the normal outer velocity and the scalar $\Psi$-anisotropic curvature of the evolving boundary, respectively. We first recall the
distributional theory when the mobility is suitably compatible with the anisotropy, and then the stability result of~\cite{CMNP19}, which provides a unique well-posed level-set
flow for arbitrary pairs $(\Psi,\mathfrak M)$.

For a norm $\vartheta$ on $\R^d$, we denote by $\vartheta^o$ its polar norm,
\begin{equation}\label{polar:norm}
    \vartheta^o(y):=\sup\{y\cdot x:\vartheta(x)\leq1\},
\end{equation}
and by
\begin{equation}\label{wulff:shape}
    W_r^\vartheta(x_0)\coloneqq \{x\in\R^d:\vartheta^o(x-x_0)<r\}
\end{equation}
the corresponding Wulff shape. For a set $F\subseteq\R^d$ we set
\begin{equation}\label{distphipallino}
    \dist^{\vartheta^o}(x,F):=\inf\{\vartheta^o(x-y):y\in F\},
\end{equation}
with the convention that the distance from the empty set is $+\infty$. The distance between two sets is defined analogously. {For $F\neq\varnothing$,} we recall that $
\dist^{\vartheta^o}$ is $1$-Lipschitz with respect to $\vartheta^o$, so it is differentiable almost everywhere, and there holds the \textit{anisotropic Eikonal equation} (see
\cite[Section 3]{Shen:finsler})
\begin{equation}\label{Eikonal}
    \vartheta\big(\nabla d_F(x) \big)=1\quad\text{for every $x \in \R^d \setminus \mathrm{clos}\, F$ of differentiability of $d_F(\cdot)= \dist^{\vartheta^o}(\cdot,F)$.}
\end{equation}
For $p\in\R^d$, the subdifferential of $\Psi$ at $p$ is $\partial\Psi(p):=\{z\in\R^d:\Psi(q)\geq\Psi(p)+z\cdot(q-p)\text{ for every }q\in\R^d\}$. Equivalently,
\[
\partial\Psi(p)=\{z\in\R^d:\Psi^o(z)\leq1,\ z\cdot p=\Psi(p)\}.
\]
In particular, $\partial\Psi(0)=\{\Psi^o\leq1\}$, whereas for $p\neq0$ the condition $\Psi^o(z)\leq1$ upgrades to $\Psi^o(z)=1$ once $z\cdot p=\Psi(p)$ is imposed.

Let $U\subset\R^d$ be an open subset. We recall that a sequence of  closed sets $(E_n)_{n\geq1}$ of $\R^d$ converges to a closed set $E$ in the \emph{Kuratowski sense} in $U$ if:
(i) whenever, for $n_j \to \infty$,  $x_{n_j}\in E_{n_j}$, and $x_{n_j}\to x\in U$, then $x\in E$; (ii) for every $x\in E\cap U$, there exists $x_n\in E_n$ such that $x_n\to x$. We
write
 \[
 E_n\stackrel{\mathcal K}{\longrightarrow}E \qquad\text{in }U.
 \]

 When $U=\R^d$, this is equivalent to the locally uniform convergence of the distance functions associated with any fixed norm. The Ascoli--Arzel\`a theorem ensures that every
 sequence of closed sets admits a subsequence converging in this sense, possibly to the empty set when the distance functions diverge to $+\infty$. In what follows, when applied to
 time-dependent functions, the differential operators $\Div$, $\nabla$, $\nabla^2$ always refer to the spatial variables only.

\medskip
We now recall the compatibility condition introduced in~\cite{CMNP19}. The mobility $\mathfrak M$ is said to be \emph{$\Psi$-regular} if there exists $\varepsilon_0>0$ such that
for every $x\in\pa W_1^{\mathfrak M}$ there exists $y\in W_1^{\mathfrak M}$ with
\[
W_{\varepsilon_0}^{\Psi}(y)\subset W_1^{\mathfrak M}, \qquad x\in\pa W_{\varepsilon_0}^{\Psi}(y).
\]
This condition yields the uniform interior Wulff-shape property needed for the distributional comparison theory.

The special case $\mathfrak M=\Psi$ will be referred to as the case of \emph{natural mobility}. In this case the $\Psi$-Wulff shapes evolve homothetically: the flow starting from
$W_R^\Psi(x_0)$ is $W_{R(t)}^\Psi(x_0)$, with $R(t)=\sqrt{R^2-2(d-1)t}$, up to the extinction time. In the following definition, the anisotropic $\mathfrak M^o$–distance is used to
express the geometric motion in a metric form.

\begin{definition}[Definition of $(\Psi,\mathfrak M)$-flows, \cite{CMNP19}]\label{Defsol}
Assume that $\mathfrak M$ is $\Psi$-regular and let $F_0\subseteq\R^d$ be a closed initial set. A closed subset $E\subseteq\R^d\times[0,+\infty)$ is called a \emph{$(\Psi,\mathfrak
M)$-superflow} with initial datum $F_0$ if, setting $E(t):=\{x\in\R^d:(x,t)\in E\}$, the following conditions hold:
\begin{itemize}
  \item[(a)] $E(0)\subseteq F_0$;
  \item[(b)] $E(s)\stackrel{\mathcal K}{\longrightarrow}E(t)$ in $\R^d$ as $s\nearrow t$ for every $t>0$;
  \item[(c)] if $E(t)=\varnothing$ for some $t\geq0$, then $E(s)=\varnothing$ for all later times $s>t$;
  \item[(d)] defining
  \[
  T^*:=\inf\{t>0:E(s)=\varnothing\ \text{for all }s\geq t\}, \qquad d(x,t):=\dist^{\mathfrak M^o}(x,E(t)),
  \]
  there exists $\mathcal z\in L^\infty(\R^d\times(0,T^*);\R^d)$ such that $\mathcal z\in\partial\Psi(\nabla d)$ almost everywhere in $(\R^d\times(0,T^*))\setminus E$ and
  \begin{equation}\label{eq:supersol}
  \partial_t d\geq\Div\mathcal z
  \end{equation}
  distributionally there. In addition, $\Div\mathcal z$ is a Radon measure on $(\R^d\times(0,T^*))\setminus E$ and its positive part satisfies $(\Div\mathcal z)^+\in
  L^\infty(\{(x,t):0<t<T^*,\ d(x,t)\geq\delta\})$ for every $\delta\in(0,1)$.
\end{itemize}

An open set $A\subseteq\R^d\times[0,+\infty)$ is called a \emph{$(\Psi,\mathfrak M)$-subflow} with initial datum an open set $A_0$ if the complement $\R^d\times[0,+\infty)\setminus
A$ is a $(\Psi,\mathfrak M)$-superflow with initial datum $\R^d\setminus A_0$.

Finally, a closed set $E\subseteq\R^d\times[0,+\infty)$ is said to be a \emph{$(\Psi,\mathfrak M)$-flow} with initial datum $F_0$ if it is a $(\Psi,\mathfrak M)$-superflow with the
same initial datum and its interior, taken in $\R^d\times[0,+\infty)$, is a $(\Psi,\mathfrak M)$-subflow with initial datum $\mathrm{int}(F_0)$.
\end{definition}

\begin{remark}\label{phiC2}
Assume $\Psi,\mathfrak M^o\in C^2(\R^d\setminus\{0\})$. Then condition \textup{(d)} in the above definition can be replaced by
\begin{itemize}
\item[(d')] there holds
\begin{equation}\label{eq:supersol'}
  \partial_t d\geq\Div(\nabla\Psi(\nabla d))
\end{equation}
in the distributional sense on $(\R^d\times(0,T^*))\setminus E$, with $T^*$ and $d$ defined as before.
\end{itemize}
Indeed, since $d$ is an $\mathfrak M^o$-distance, outside of $E(t)$ $d(\cdot,t)$ satisfies $\mathfrak M(\nabla d)=1$ at every point of differentiability of $d$, and thus $\nabla d\neq0$ a.e. and $\partial\Psi(\nabla
d)=\{\nabla\Psi(\nabla d)\}$ a.e. Moreover, the $L^\infty$ bound from above on $\Div(\nabla\Psi(\nabla d))$ follows from the standard semiconcavity estimate for the $\mathfrak
M^o$-distance, which can be proven as follows. For any $t \in (0,T^*)$ and for any $x\in \{d(\cdot, t)> \delta\}$, let $y\in \pa E(t)$ be such that $d(x,t)=\mathfrak M^o(x-y)$.
Then,
\[
\begin{split}
d(x+h,t)-2d(x,t)+d(x-h,t)&\le {\mathfrak M^o}(x+h-y)-2{\mathfrak M^o}(x-y)+ {\mathfrak M^o}(x-h-y)\\
&= \int_{-1}^1(1-|s|)(\nabla^2{\mathfrak M^o}(x-y+sh)h) \cdot h \, ds \\
&\le C |h|^2 \int_{-1}^1 \frac{1-|s|}{{\mathfrak M^o}(x-y+sh)} \, ds \le \frac{C}{\delta}|h|^2,
\end{split}
\]
where in the second last inequality we used that the map $\xi \mapsto \mathfrak M^o(\xi)\nabla^2\mathfrak M^o(\xi) $ is $0$-homogeneous and thus bounded, while the last one follows
provided that $|h|$ is sufficiently small.

Thus, $d(\cdot,t)$ is locally semiconcave outside $E(t)$, and for every $\delta>0$,
\begin{equation}\label{dsemi}
D^2d(\cdot,t)\leq\frac C\delta\mathrm{Id}\,\Leb^d
\qquad\text{as matrix-valued measures in }
\{d(\cdot,t)>\delta\},
\end{equation}
where $C$ depends only on $\mathfrak M$. In particular, the singular part of $D^2d(\cdot,t)$ is negative semidefinite. By the chain rule for $BV$ maps and the convexity of $\Psi$,
$\Div(\nabla\Psi(\nabla d))$ is locally a signed Radon measure with nonpositive singular part, while the absolutely continuous part
\[
\big[\Div(\nabla\Psi(\nabla d))\big]_{\mathrm{ac}}
=
\operatorname{tr}\!\left(
\nabla^2\Psi(\nabla d)(D^2d)_{\mathrm{ac}}
\right)
\leq \frac{C}{\delta}\sup_{\mathfrak M(\xi)=1} \mathrm{tr}(\nabla^2 \Psi(\xi))\leq \frac{C'}{\delta}
\quad\text{a.e. in }\{d(\cdot,t)>\delta\}.
\]
Here $C'$ depends only on $\Psi$ and $\mathfrak M$, since $\mathfrak M(\nabla d)=1$ almost everywhere outside $E$. These estimates are uniform in $t$ and yield the Radon measure
property and the bound on the positive part required in \textup{(d)}.

\end{remark}

The next result is the fundamental comparison property established in~\cite{CMNP19}.

\begin{theorem}\label{th:compar}
Assume that $\mathfrak M$ is $\Psi$-regular. Let $E$ be a $(\Psi,\mathfrak M)$-superflow with initial datum $F_0$ and $A$ a $(\Psi,\mathfrak M)$-subflow with initial datum $A_0$.
If
\[
\dist^{\mathfrak M^o}(F_0,\R^d\setminus A_0)=:\Delta>0,
\]
then
\[
\dist^{\mathfrak M^o}(E(t),\R^d\setminus A(t))\geq\Delta \qquad\text{for every }t\geq0.
\]
\end{theorem}

\medskip
We now move to the level-set description. Let $\mathrm{UC}(\R^d)$ denote the space of uniformly continuous functions on $\R^d$.

\begin{definition}[Level-set sub- and supersolutions]\label{deflevelset1}
Assume that $\mathfrak M$ is $\Psi$-regular and let $u_0\in\mathrm{UC}(\R^d)$. A lower-semicontinuous function $u:\R^d\times[0,+\infty)\to\R$ is called a \emph{$(\Psi,\mathfrak
M)$-level-set superflow} with initial datum $u_0$ if $u(\cdot,0)\geq u_0$ and, for almost every $\lambda\in\R$, the sublevel set $\{u\leq\lambda\}$ is a $(\Psi,\mathfrak
M)$-superflow with initial set $\{u_0\leq\lambda\}$. An upper-semicontinuous function $u$ is a \emph{$(\Psi,\mathfrak M)$-level-set subflow} with initial datum $u_0$ if $-u$ is a
$(\Psi,\mathfrak M)$-level-set superflow with initial datum $-u_0$. A continuous function $u$ is a \emph{$(\Psi,\mathfrak M)$-level-set flow} if it is both a level-set subflow and
a level-set superflow.
\end{definition}

The following comparison principle is a direct consequence of Theorem~\ref{th:compar}; see~\cite[Theorem 2.8]{CMNP19}.

\begin{theorem}\label{th:lscomp}
Assume that $\mathfrak M$ is $\Psi$-regular. Let $u_0,v_0\in\mathrm{UC}(\R^d)$ and let $u$ and $v$ denote, respectively, a $(\Psi,\mathfrak M)$-level-set subflow and superflow
starting from $u_0$ and $v_0$. If $u_0\leq v_0$, then $u\leq v$ on $\R^d\times[0,+\infty)$.
\end{theorem}

\medskip
Existence and uniqueness for $\Psi$-regular mobilities follow from the minimizing-movements construction of~\cite{CMNP19}. We retain the formulation of the natural-mobility case,
with the corresponding notions of $(\Psi,\mathfrak M)$-superflow and subflow.

\begin{theorem}[{\cite[Theorem 4.8]{CMNP19}}]\label{th:phiregularlevelset}
Assume that $\mathfrak M$ is $\Psi$-regular and let $u_0\in\mathrm{UC}(\R^d)$. Then the following properties hold:
\begin{itemize}
  \item[(i)] There exists a unique $(\Psi,\mathfrak M)$-level-set flow $u$ in the sense of Definition~\ref{deflevelset1} starting from $u_0$.
  \item[(ii)] For every $\lambda\in\R$, the set $\{(x,t)\in\R^d\times[0,+\infty):u(x,t)\leq\lambda\}$ (resp. $\{(x,t)\in\R^d\times[0,+\infty):u(x,t)<\lambda\}$) coincides with the maximal $(\Psi,\mathfrak M)$-superflow (resp. the minimal $(\Psi,\mathfrak M)$-subflow) with initial data $\{u_0\leq\lambda\}$ (resp. $\{u_0<\lambda\}$).
  \item[(iii)] For all but countably many $\lambda\in\R$, no fattening occurs at level $\lambda$; namely,
  \begin{equation}\label{eq:nonfattening}
  \begin{split}
  &\overline{\{(x,t)\in\R^d\times[0,+\infty):u(x,t)<\lambda\}}=\{(x,t)\in\R^d\times[0,+\infty):u(x,t)\leq\lambda\},\\
  &\{(x,t)\in\R^d\times[0,+\infty):u(x,t)<\lambda\}=\mathrm{int}\big(\{(x,t)\in\R^d\times[0,+\infty):u(x,t)\leq\lambda\}\big),
  \end{split}
  \end{equation}
  where the closure and interior operations are taken in the space-time  topology of $\R^d\times[0,+\infty)$.
\end{itemize}
For all such values of $\lambda$, the set $\{(x,t)\in\R^d\times[0,+\infty):u(x,t)\leq\lambda\}$ represents the unique $(\Psi,\mathfrak M)$-flow starting from $\{u_0\leq\lambda\}$.
\end{theorem}

Note that by the maximality (minimality) properties of $\{u\leq \lambda\}$ (resp. $\{u< \lambda\}$), the evolution $\{u(\cdot, t)\leq \lambda\}$ (resp. $\{u(\cdot, t)<\lambda\}$)
depends only on \(\{u_0\le\lambda\}\) (resp. $\{u_0 < \lambda\}$) and not on the choice of $u_0$. This authorizes one to consider $\{u(\cdot, t)\leq \lambda\}$ as a uniquely
defined {\em generalized anisotropic mean curvature flow} starting from \(F_0:=\{u_0\le\lambda\}\).

We now recall the stability result of~\cite{CMNP19}, which removes the $\Psi$-regularity assumption on the mobility. Let $\Psi$ and $\mathfrak M$ be arbitrary norms and let
$\{\mathfrak M_n\}$ be any sequence of $\Psi$-regular mobilities converging locally uniformly to $\mathfrak M$; for instance, one may take $\mathfrak M_n:=\mathfrak M+n^{-1}\Psi$.
If $u_n$ denotes the $(\Psi,\mathfrak M_n)$-level-set flow starting from $u_0\in\mathrm{UC}(\R^d)$, then $u_n$ converges locally uniformly in $\R^d\times[0,+\infty)$ to a limit
which is independent of the chosen regular approximation.

\begin{definition}\label{def:general-levelset}
For arbitrary norms $\Psi$ and $\mathfrak M$, the limit described above is called the \emph{$(\Psi,\mathfrak M)$-level-set flow} starting from $u_0$.
\end{definition}

The resulting evolution is well posed, satisfies the comparison principle, and is stable under simultaneous approximation of anisotropy and mobility. More precisely, we shall use
the following consequence of~\cite[Proposition 5.12]{CMNP19}.

\begin{theorem}[Stability of level-set flows]\label{th:levelset-stability}
Let $\Psi_n,\mathfrak M_n$ be norms on $\R^d$ converging locally uniformly to the norms $\Psi,\mathfrak M$, respectively. Given $u_0\in\mathrm{UC}(\R^d)$, let $u_n$ and $u$ be the
corresponding $(\Psi_n,\mathfrak M_n)$- and $(\Psi,\mathfrak M)$-level-set flows starting from $u_0$. Then $u_n\to u$ locally uniformly in $\R^d\times[0,+\infty)$.
\end{theorem}

In particular, from now on the expression $(\Psi,\mathfrak M)$-level-set flow will be used without any compatibility assumption. Notice instead that the distributional notions of
$(\Psi,\mathfrak M)$-superflow, subflow and flow introduced in Definition~\ref{Defsol} refer to the $\Psi$-regular case.

\medskip
We now define the notions of non-fattening used throughout the paper. They refer to the unique $(\Psi,\mathfrak M)$-level-set flow just introduced and hence make sense for
arbitrary pairs of norms.

\begin{definition}[Non-fattening conditions]\label{def:fat}
Let $E_0\subseteq\R^d$ be an arbitrary set. We say that $u_0\in\mathrm{UC}(\R^d)$ is a \emph{level-set representation} of $E_0$ at level zero if
\begin{equation}\label{nonfat1}
    \{u_0<0\}\subseteq E_0\subseteq\{u_0\leq0\}.
\end{equation}
We say that $E_0$ is \emph{$(\Psi,\mathfrak M)$-non-fattening}, or simply \emph{non-fattening} when the pair is understood, if it admits a level-set representation $u_0$ such that,
if $u$ denotes the corresponding $(\Psi,\mathfrak M)$-level-set flow, then no fattening occurs at level zero, that is
\begin{equation}\label{nonfat2}
\begin{split}
&\overline{\{(x,t)\in\R^d\times(0,+\infty):u(x,t)<0\}}=\{(x,t)\in\R^d\times(0,+\infty):u(x,t)\leq0\},\\
&\{(x,t)\in\R^d\times(0,+\infty):u(x,t)<0\}=\mathrm{int}\big(\{(x,t)\in\R^d\times(0,+\infty):u(x,t)\leq0\}\big),
\end{split}
\end{equation}
where the closure and interior operations are understood with respect to the relative topology of $\R^d\times(0,+\infty)$.

We additionally say that $E_0$ is \emph{non-fattening in measure} if it is non-fattening and, furthermore,
\begin{equation}\label{strong:nonfat}
    \mathcal L^{d+1}\big(\{(x,t)\in\R^d\times(0,+\infty):u(x,t)=0\}\big)=0.
\end{equation}
\end{definition}

The positive-time open and closed zero-level evolutions associated with a non-fattening representation depend only on $E_0$ and not on the chosen representation; see~\cite[Theorem
5.9 (v)]{CMNP19}.

\begin{proposition}\label{prop:independence-representation}
Let $E_0\subseteq\R^d$, and let $u_0,v_0\in\mathrm{UC}(\R^d)$ be two non-fattening level-set representations of $E_0$. Let $u$ and $v$ be the corresponding $(\Psi,\mathfrak
M)$-level-set flows. Then
\[
\begin{split}
\{(x,t)\in\R^d\times(0,+\infty):u(x,t)\leq0\} &=\{(x,t)\in\R^d\times(0,+\infty):v(x,t)\leq0\},\\
\{(x,t)\in\R^d\times(0,+\infty):u(x,t)<0\} &=\{(x,t)\in\R^d\times(0,+\infty):v(x,t)<0\}.
\end{split}
\]
In particular, an evolution that is non-fattening (in particular, one that is non-fattening in measure) is independent of the chosen non-fattening level-set representation.
\end{proposition}

\begin{definition}\label{def:gen-motion}
Let $E_0$ be non-fattening and let $u_0$ be any non-fattening level-set representation of $E_0$, with corresponding $(\Psi,\mathfrak M)$-level-set flow $u$. The mapping
\[
t\in(0,+\infty)\longmapsto\{u(\cdot,t)\leq0\}
\]
is called the \emph{generalized $(\Psi,\mathfrak M)$-motion starting from $E_0$}. In view of Proposition~\ref{prop:independence-representation}, this definition is independent of
the chosen non-fattening representation $u_0$.
\end{definition}

\begin{remark}
For a fixed representation, the generalized motion extends to time $t=0$ as $\{u_0\leq0\}$. Notice, however, that this set need not coincide with $E_0$ and may depend on the chosen
representation.
\end{remark}

We conclude by recalling that, given $u_0\in\mathrm{UC}(\R^d)$, for all but countably many $\lambda\in\R$, every set $E_0$ satisfying $\{u_0<\lambda\}\subseteq
E_0\subseteq\{u_0\leq\lambda\}$ is non-fattening in measure.

\subsection{Flat flows}\label{subses:flat:flow}

We now recall the construction of flat-flow solutions for the $(\Psi,\mathfrak M)$-flow, originally introduced in the pioneering works~\cite{ATW,AlmTay95,LS}, and their
relationship with the level-set solution described above. We use the global variational formulation of the ATW scheme. For its connection with Chambolle's differential formulation,
we refer to \cite{Chambolle,CMP17,CMNP19}. In this scheme the anisotropy $\Psi$ determines the perimeter term, while the mobility $\mathfrak M$ enters through the signed $\mathfrak
M^o$-distance.

 Let $E_0\subseteq\R^d$ be a Lebesgue measurable set such that $P^\Psi(E_0)<+\infty$, and which we identify with its density-one representative $E_0=E_0^{(1)}$.

For a nonempty proper subset $E\subset\R^d$, we denote by $\mathrm{sd}_E$ its signed $\mathfrak M^o$-distance:
\begin{equation}\label{signeddist}
\mathrm{sd}_E(x)
:=
\dist^{\mathfrak M^o}(x,E)
-
\dist^{\mathfrak M^o}(x,\R^d\setminus E),
\qquad x\in\R^d.
\end{equation}

Fix $h>0$ and set $E_0^h:=E_0^{(1)}$. Starting from $E_k^h$, we choose $E_{k+1}^h$ as a global minimizer of
\begin{equation}\label{ATWfunctional}
F\longmapsto
P^\Psi(F)
+\frac1h\int_{F\Delta E_k^h}
|\mathrm{sd}_{E_k^h}(x)|\,dx
\end{equation}
among all Lebesgue measurable sets $F\subseteq\R^d$ of finite perimeter. If $E_k^h\in\{\varnothing,\R^d\}$, we instead set $E_{k+1}^h:=E_k^h$.

At every positive discrete time, we identify the minimizer with its closed representative $\overline{(E_{k+1}^h)^{(1)}}$. This changes the set only on a negligible set, by the
local density estimates for perimeter quasi-minimizers; see \cite[Lemma 1.3 and Remark 1.4]{LS} and the anisotropic argument in \cite[Lemma 3.8 and Remark 3.9]{CMNP19}.

We then define
\[
E^h(t):=E_{\lfloor t/h\rfloor}^h,
\qquad t\geq0
\]
 the piecewise constant interpolation of the discrete scheme, and let
 \[
 E^h:=\overline{\{(x,t)\in\R^d\times[0,+\infty):x\in E^h(t)\}}
 \]
 be the corresponding closed space-time tube.

 The standard Luckhaus--Sturzenhecker/Almgren--Taylor--Wang compactness argument~\cite{ATW,LS}, based on~\eqref{ATWfunctional},  yields a sequence $h_j\to0^+$ and  a family of sets
 of finite perimeter $(\widehat E(t))_{t\geq0}$, called a \emph{$(\Psi,\mathfrak M)$-flat flow}, such that
\begin{equation}\label{eq:flatflow}
    \chi_{E^{h_j}(t)}\longrightarrow\chi_{\widehat E(t)} \qquad\text{in }L^1_{\mathrm{loc}}(\R^d)
\end{equation}
for every $t> 0$,
\[
\widehat E(t) \stackrel{t \to 0^+}{\longrightarrow} \widehat E(0)\qquad\text{in }L^1_{\mathrm{loc}}(\R^d)\,
\]
the perimeter estimate

\begin{equation}\label{flat:tutto:per}
\sup_{t\geq 0}P^\Psi(\widehat E(t))\leq P^\Psi(E_0)
\end{equation}
holds, and we have
\begin{equation}\label{eq:flatflowholder}
\mathcal L^d\big(\widehat E(t)\Delta\widehat E(s)\big)
\leq C|t-s|^{1/2}
\qquad\text{for every }s,t\geq 0,
\end{equation}
where $C=C(d,\Psi,\mathfrak M,P^\Psi(E_0))$.

The relation with the level-set flow is as follows: let $u_0\in\mathrm{UC}(\R^d)$ be a level-set representation of $E_0$ in the sense of Definition~\ref{def:fat}, and let $u$ be
the corresponding $(\Psi,\mathfrak M)$-level-set flow.

If $E$ is any space-time Kuratowski cluster point of $E^h$, namely if
\begin{equation}\label{Kuratowskicluster}
    E^{h_j}\stackrel{\mathcal K}{\longrightarrow}E \qquad\text{in }\R^{d+1}
\end{equation}
along a subsequence $h_j\to0^+$, then
\[
\{(x,t)\in\R^d\times[0,+\infty):u(x,t)<0\} \subseteq E\subseteq \{(x,t)\in\R^d\times[0,+\infty):u(x,t)\leq0\}.
\]
Consequently, if $E_0$ is non-fattening and $u_0$ is a non-fattening level-set representation of $E_0$, then every Kuratowski cluster point $E$ coincides in positive time with the
generalized $(\Psi,\mathfrak M)$-motion $\{u\leq0\}$.

Denoting by $\widehat E(t)^{(1)}$ the density-one representative of $\widehat E(t)$ and by $E(t):=\{x:(x,t)\in E\}$, with $E$ the Kuratowski limit in~\eqref{Kuratowskicluster}, one
has
\begin{equation}\label{eq:flat-levelset-sandwich}
    \{u(\cdot,t)<0\}\subseteq\widehat E(t)^{(1)}\subseteq E(t)\subseteq\{u(\cdot,t)\leq0\}
\end{equation}
for every $t\geq0$.  In particular, if $E_0$ is non-fattening in measure and $u_0$ is a level-set representation of $E_0$ that is non-fattening in measure, then
\begin{equation}\label{flat-uguale}
    \widehat E(t)=E(t)=\{u(\cdot,t)\leq0\} \qquad\text{up to $\mathcal L^d$-negligible sets for almost every $t\geq0$.}
\end{equation}
Under the same assumption of non-fattening in measure, the $1/2$-H\"older continuity of flat flows in time \eqref{eq:flatflowholder} implies that any two $(\Psi,\mathfrak M)$-flat flows
obtained from the scheme agree up to negligible sets for every $t \geq 0$.

If the boundary of $E_0$ has zero Lebesgue measure, then the flat flow attains the prescribed initial datum in the $L^1_{\loc}$ sense.

\begin{lemma}
\label{lem:flat-initial-trace}
Let  $E_0$ be a nonempty, proper subset of $\R^d$ of finite perimeter with $\mathcal L^d(\partial E_0)=0$. Let $\{\widehat E(t)\}$ be the flat flow starting from $E_0$. Then
$\widehat E(0)\stackrel{a.e.}{=}E_0$ on $\R^d$.
\end{lemma}

\begin{proof}
Choose $u_0=\mathrm{sd}_{E_0}$, so that $\{u_0<0\}=\mathrm{int}(E_0)$ and $\{u_0\leq0\}=\overline{E_0}$. Let $u$ be the corresponding level-set flow. By the continuity of $u$ and
the $L^1_{\loc}$-continuity of $\widehat E$ in \eqref{eq:flatflowholder}, letting $t\downarrow0$ in \eqref{eq:flat-levelset-sandwich}, we obtain
\[
\mathrm{int}(E_0)\subseteq\widehat E(0)
\subseteq\overline{E_0}
\qquad\text{up to negligible sets}.
\]
Since $\mathcal L^d(\partial E_0)=0$, we have $\mathrm{int}(E_0)\stackrel{a.e.}{=}E_0 \stackrel{a.e.}{=}\overline{E_0}$, which proves the claim.
\end{proof}

\subsection{Anisotropic Brakke flow}

Let $\psi,\mathfrak{m}$ be norms on $\R^d$ such that $\psi\in C^1(\R^d \setminus \{0\})$.

The following definition of $(\psi,\mathfrak{m})$-anisotropic Brakke flow is motivated by \eqref{defFV4}: formally, we require that the rate of change in time of the anisotropic
energies of the form
\[
\mathbf{F}^\psi(V_t,\varphi(\cdot,t)) =\|V_t\|^\psi\big(\varphi(\cdot,t) \big)= \int\varphi(x,t) \,\psi(\nu)\,dV_t(x,\nu)
\]
is at most given by the weighted first variation $\delta_\psi(V_t,\varphi(\cdot,t))$ along the rough vector field $X(x)=\mathfrak{m}\big(\nu(x)\big)\,\mathbf{H}^\psi_{V_t}(x)$,
where here $\nu=\nu_{V_t}$ is the normal to $V_t$.

\begin{definition}[Anisotropic $(\psi,\mathfrak{m})$-Brakke flow]\label{def:brakke}
    Let $\Omega\subset \R^{d}$ be open, $\psi,\mathfrak{m}$ be norms as above, and let $\{V_t\}_{t\in [0,T)}$ be a (measurable) 1-parameter family of $(d-1)$-varifolds on $\Omega$.
    We say that $\{V_t\}_{t\in [0,T)}$ is an anisotropic Brakke flow  w.r.t. $\psi$ with mobility $\mathfrak{m}$ if it satisfies
    \begin{enumerate}
        \item for a.e. $t$, $V_t \in \mathbf{IV}_{d-1}(\Omega)$;
        \item for any compact $K\Subset \Omega$ and $t<T$, $\sup_{s\in [0,t]}\|V_s\|(K)<\infty$;
        \item for a.e. $t$, $V_t$ has locally bounded $\psi$-anisotropic first variation, $\delta_{\psi} V_t$ is absolutely continuous with respect to $\|V_t\|$ (i.e.  the anisotropic mean curvature $\mathbf{H}^\psi_{V_t}$ exists), and  the field $(x,t)\mapsto\mathbf{H}^\psi_{V_t}(x)\in L^2_{loc}(\|V_t\| dt)$;
        \item For all $0\leq t_1<t_2<T$ and $\varphi\in C^1(\Omega\times [0,T);\mathbb R^+)$ such that
        \[
                 {\spt(\varphi)} \quad \mbox{is compact in $\Omega \times [0,T)$}
        \]
  the following \textit{$(\psi,\mathfrak{m})$-anisotropic Brakke's inequality} holds
\begin{equation}\label{brakineq}
\begin{split}
\qquad &\left.\|V_t\|^\psi\big(\varphi(\cdot,t) \big)\right|_{t=t_1}^{t_2} :=\|V_{t_2}\|^\psi\big(\varphi(\cdot,t_2)\big)-\|V_{t_1}\|^\psi\big(\varphi(\cdot,t_1) \big) \\
&\leq \int_{t_1}^{t_2}\int\bigg\{\Big[-\varphi(x,t) \big|\mathbf{H}^\psi_{V_t}\big|^2 + \big(\nu\cdot \mathbf{H}^\psi_{V_t}\big)\,\big(\nabla\varphi(x,t)\cdot \nabla\psi(\nu)\big)\Big]\mathfrak{m}(\nu) + \partial_t\varphi(x,t)\,\psi(\nu)\bigg\} \,dV_t(x,\nu)\,dt
\,.
\end{split}
\end{equation}
    \end{enumerate}

In the case $\psi=\mathfrak{m}$, we simply write $\psi$-anisotropic Brakke flow in place of $(\psi,\mathfrak{m})$-anisotropic Brakke flow.

\end{definition}

We note explicitly that, when $\psi=\mathfrak m= |\cdot|$ then Definition \ref{def:brakke} reduces to the classical definition of Brakke flow as, for instance, in
\cite{TonegawaBook}. Additionally, we make the following important observation.

\begin{remark}\label{rem:brineq}
    The inner integral in \eqref{brakineq} is well defined; indeed, the second integrand may be written as $ [\nu\otimes \nabla\psi(\nu)]:[ \mathbf{H}^
    \psi_{V_t}\otimes\nabla\varphi(x,t)]$, and as the tensor function
    \[
   \S^{d-1}\ni \nu\mapsto\nu \otimes \nabla\psi(\nu) \] is invariant w.r.t. flipping $\nu \mapsto -\nu$, it descends to the quotient $\mathbb{RP}^{d-1}\simeq G(d,d-1)$.

  Moreover, as $V_t$ is rectifiable for a.e. $t\in (0,T)$, then for all such times it is well defined the ``Brakke term'' {\small
    \begin{equation}\label{def:Brakketerm}
\mathfrak{B}\big(\varphi(\cdot,t);t \big) :=\int_\Omega\left\lbrace -\varphi(x,t) \big|\mathbf{H}^\psi_{V_t}\big|^2 + \big(\nu(x)\cdot \mathbf{H}^\psi_{V_t}\big)\,\big(\nabla\varphi(x,t)\cdot \nabla\psi\big(\nu(x)\big)\big)\right\rbrace \,d\|V_t\|^{\mathfrak{m}}(x) + \|V_t\|^\psi\big(\partial_t\varphi(\cdot,t)\big)
    \end{equation}
    } with $\|V_t\|^{\mathfrak{m}}$  as in \eqref{two:weights} with $\mathfrak{m}$ in place of $\psi$, and with $\nu(x)=\nu_{V_t}(x)$ denoting the normal to $V_t$ at $x$. Then by
    means of \eqref{test:Vrect} and \eqref{two:weights}, inequality \eqref{brakineq} can be rewritten as
    \begin{equation}\label{neq:brakke}\left.\|V_t\|^\psi\big(\varphi(\cdot,t) \big)\right|_{t=t_1}^{t_2}
\leq \int_{t_1}^{t_2}\mathfrak{B}(\varphi(\cdot,t);t)\,dt.
\end{equation}
\end{remark}

\section{The weighted anisotropic area functional}\label{sec:waf}
Here we study the minimum problem which defines the anisotropic elliptic regularization scheme. As anticipated in Section \ref{s:preliminaries}, we will work in the extended space
$\R^{N+1}=\R^N_{x} \times \R_{z}$. All notation introduced for the abstract $d$-dimensional space specializes here with $d=N+1$. Let $\Psi$ be a norm on $\R^{N+1}$, and let $E_0
\subset \R^N_x$ be a set of finite perimeter.

 For every $\e>0$ we consider the following \textit{weighted anisotropic area functional} \beq\label{defAF}
 \mathcal{A}^\Psi_\e(F):=\int_{\pa^*F\cap\{z>-1\}}{e}^{-{z}/\e}\,\Psi(\nu)\, d\H^{N}(x,z), \eeq together  with the following localized version: given $A\subset \{z>-1\}$, we set
\begin{equation}\label{localized:Apsi}
 \mathcal{A}^\Psi_\e(F; A):=\int_{\pa^*F\cap A}{e}^{-{z}/\e}\,\Psi(\nu)\, d\H^{N}(x,z),
\end{equation}
for any set $F\subset \R^{N+1}$ of locally finite perimeter.

 Given $g\in \mathrm{Lip}_c(\R^{N+1})$ we are interested in the following minimization problem:
\begin{equation}\label{minpb}
\inf\Biggl\{\mathcal{A}^\Psi_\e(F)+\int_{F}g(x,z)\, dxdz: F\subset\R^{N+1}\text{ of loc. finite perimeter s.t. }F\cap H^{-} \stackrel{a.e.}{=} E_0\times (-\infty,0)\Biggr\}\,,
\end{equation}
where we have set
\begin{equation}\label{hyperplane}
    H^\pm:=\big\{(x,z): z \gtrless 0\big \}.
\end{equation}
The existence of a minimizer $E_\e$ follows by standard methods using the lower semicontinuity of the functional and the compactness properties of sets with locally bounded
perimeter.

Up to choosing an equivalent representative, we may, and shall assume that
\[
\partial E_\e=\operatorname{spt}\mu_{E_\e};
\]
see \cite[Proposition 12.19]{MaggiBook}.

\begin{proposition}\label{prop:Euler_lagrange}
    Let $\Psi\in C^1(\R^{N+1}\setminus\{0\})$ be a norm, and let $E_\e$ be a minimizer of \eqref{minpb}. Then it satisfies the Euler-Lagrange equations

    \begin{equation}\label{euler:lagrange00}
                     \int_{\partial^* E_\e} e^{-z/\e}B_\Psi(\nu_{E_\e}):DXd\H^N=\frac{1}{\e} \int_{\partial^* E_\e}e^{-z/\e}(X\cdot e_{_{N+1}})\,\Psi(\nu_{E_\e})d\H^N-\int_{\partial^* E_\e} g\, (X\cdot \nu_{E_\e})d\H^{N},
    \end{equation}
   for all $X\in C^1_c(\{z>0\};\R^{N+1})$, and
    \begin{equation}\label{euler:lagrange0}
        \int_{\partial^* E_\e} B_{\Psi}(\nu_{E_\e}):DY\,d\H^N=\int_{\partial^* E_\e} \left(\frac{\nabla \Psi(\nu_{E_\e})\cdot e_{_{N+1}}}{\e}-e^{z/\e}\,g\right)\,(Y\cdot \nu_{E_\e})\,d\H^N(x,z)
    \end{equation}
   for all $Y\in C^1_c(\{z>0\};\R^{N+1})$, where $\nu_{E_\e}$ is the outer normal to $E_\e$, and $B_\Psi$  is defined in \eqref{def:Bpsi}, with $\Psi$ in place of $\psi$. In
   particular,  $\partial^* E_\e$ admits $\Psi$-anisotropic mean curvature  which satisfies
\begin{equation}\label{euler:lagrange1}
        \mathbf{H}_{\partial^* E_\e}^\Psi=-\left(\frac{\nabla \Psi(\nu_{E_\e})\cdot e_{_{N+1}}}{\e}-e^{z/\e}\,g\right)\,\nu_{E_\e},\quad\text{on $\partial^* E_\e\cap \{z>0\}$.}
    \end{equation}
\end{proposition}

\begin{proof}
    We apply the argument of Remark \ref{remark:minimality} with $\Omega=\{z>0\}$ and $\varphi(x,z)=e^{-z/\e}$ so that, for all $X\in C^1_c(\{z>0\};\R^{N+1})$, we have
    \begin{equation*}
        \begin{split}
            -\int_{\partial^* E_\e} g\, (X\cdot \nu_{E_\e})d\H^{N}= \int_{\partial^* E_\e} e^{-z/\e}B_\Psi(\nu_{E_\e}):DXd\H^N-\frac{1}{\e} \int_{\partial^* E_\e}e^{-z/\e}(X\cdot e_{_{N+1}})\,\Psi(\nu_{E_\e})d\H^N,
        \end{split}
    \end{equation*}
    that is \eqref{euler:lagrange00}. Next Equation \eqref{euler:lagrange0} immediately follows by setting $Y=e^{-z/\e}X$, by using that
    $DY=e^{-z/\e}\,DX-\frac{1}{\e}e^{-z/\e}X\otimes e_{_{N+1}}=e^{-z/\e}\,DX-\frac{1}{\e}Y\otimes e_{_{N+1}}$, and that, from \eqref{def:Bpsi} we have
\begin{equation*}
    B_\Psi(\nu):(Y\otimes e_{_{N+1}})=\Psi(\nu)\,(Y\cdot e_{_{N+1}})-(\nabla\Psi(\nu)\cdot e_{_{N+1}})\,(Y\cdot \nu).
\end{equation*}
The proof is complete.
\end{proof}

\begin{remark}\label{rm:complementare}
If \(E_\e\) is a minimizer of \eqref{minpb}, then \(\R^{N+1}\setminus E_\e\) is a minimizer of
\[
\min\left\{ \mathcal{A}^\Psi_\e(F)-\int_F g(x,z)\,dxdz: F\cap H^-=(\R^N\setminus E_0)\times(-\infty,0) \right\}.
\]
\end{remark}

We recall the following important notion of $\Lambda$-minimizer for the functional $\mathcal{A}^\Psi_\e$.
\begin{definition}\label{lambda-min}
Let  $A\subset \{z>-1\}$ be an open set and $\Lambda\geq 0$ be a constant. We say that a set of locally finite perimeter $E\subset \R^{N+1}$ is a $\Lambda$-minimizer for
$\mathcal{A}^\Psi_\e$ in $A$ if for any ball $\B_r(x,z)\subset A$ and for any $F$ such that $E\Delta F\Subset \B_r(x,z)$ we have
\[
\mathcal{A}^\Psi_\e(E; \B_r(x,z))\leq \mathcal{A}^\Psi_\e(F; \B_r(x,z))+ \Lambda |E\Delta F|.
\]
\end{definition}

\begin{remark}\label{rem:Lmin}
Let \(E\) be a \(\Lambda\)-minimizer for \(\mathcal A^\Psi_\e\) in \(A\). By the standard characterization of the support of the perimeter measure, $\partial\bigl(E^{(1)}\bigr)\cap
A= \spt\,\mu_E\cap A\,.$ Moreover, standard local density estimates for \(\Lambda\)-minimizers (see e.g. \cite[Lemma 2.8]{DPhMa15}) imply that $\spt\,\mu_E\cap\bigl(E^{(1)}\cup
E^{(0)}\bigr)\cap A=\varnothing$. Thus, \(E^{(1)}\) and \(E^{(0)}\) are relatively open in \(A\), \(\mathcal L^{N+1}(\spt\,\mu_E\cap A)=0\), and $\partial\overline{E^{(1)}}\cap A =
\partial\bigl(E^{(1)}\bigr)\cap A = \spt\,\mu_E\cap A$. Thus both \(E^{(1)}\) and \(\overline{E^{(1)}}\) are representatives of \(E\) in \(A\). Whenever a closed representative is
needed in $A$, we shall identify \(E\) with \(\overline{E^{(1)}}\) in $A$.
\end{remark}

The following result is classical, but we record the proof for the sake of completeness.
\begin{lemma}[Regularity at contact points]\label{lm:contact}
Assume $\Psi, \Psi^o\in C^\infty(\R^{N+1}\setminus\{0\})$.  Let  $A\subset \{z > -1\}$ be an open set  and let  $E\subset \R^{N+1}$ be a (closed) $\Lambda$-minimizer for
$\mathcal{A}^\Psi_\e$ in $A$ according to the previous definition. Assume also that  there exists an open set $U\subset \R^{N+1}\setminus E$ of class $C^1$ and a point
$(x_0,z_0)\in \pa U\cap\pa E\cap A$. Then $\pa E$ is of class $C^{1,\alpha}$  in a neighborhood of $(x_0, z_0)$ for some $\alpha\in (0,1)$.
\end{lemma}
\begin{proof}
 For any $\eta_n\to 0^+$ set \(E_n:=(E-(x_0,z_0))/\eta_n\), $U_n:=(U-(x_0, z_0))/\eta_n$, and note that $E_n$ is an $(\eta_n\Lambda)$-minimizer for
\[
 \mathcal{A}^\Psi_{\e,n}(F):=\int_{\pa^*F\cap\{z>-(z_0+1)/\eta_n\}}{e}^{-(z_0+\eta_n z)/\e}\,\Psi(\nu)\, d\H^{N}(x,z)
\]
in the open set $A_n:=(A-(x_0,z_0))/\eta_n$, and that $U_n\subset \R^{N+1}\setminus E_n$, with $(0,0)\in \pa U_n\cap\pa E_n\cap A_n$. By the well-known compactness properties of
$\Lambda$-minimizers, see for instance \cite[Theorem~2.9]{DPhMa15}, up to a non-relabeled subsequence, $E_n\to E_\infty$ in $L^1_{loc}( \R^{N+1})$, with $E_\infty$ a local
minimizer for the anisotropic perimeter functional $P^\Psi(\cdot)$. in the whole $\R^{N+1}$. Moreover, $\pa E_n \stackrel{\mathcal K}{\longrightarrow}\pa E_\infty$ in the
Kuratowski sense  and $\mu_{E_n} \to \mu_{E_\infty}$ as (Radon) measures in $\R^{N+1}$. In particular, since $(0,0)\in \pa E_n$ for all $n$, we have $(0,0)\in \pa E_\infty$. Note
also that by our assumptions  $\R^{N+1}\setminus U_n\stackrel{\mathcal K}{\longrightarrow} H_\infty:=\{(x,z): (x,z)\cdot \nu_0\geq 0\}$, where $\nu_0$ denotes the outer normal to
$\pa U$ at $(x_0,z_0)$. Thus, $E_\infty\subset H_\infty$. We may now invoke the Strong Maximum Principle \cite[Lemma~2.13]{DPhMa15}, \cite[Theorem, p.~686]{SolomonWhite89} to
conclude that $E_\infty$ in fact coincides with the half-space $H_\infty$. It easily follows that $(x_0,z_0)\in \pa^* E$. The thesis of the lemma is now a consequence of the
classical partial regularity theory for $\Lambda$-minimizers of elliptic integrands, see  \cite{SSA77, Duzaar:Steffen}.
\end{proof}

We now show that under the same regularity assumption on $\Psi$ as in previous lemma and by taking $\e$ small enough, the trace of the boundary of any minimizer $E_\e$ of
\eqref{minpb} on $\{z=0\}$ must be included in $\partial E_0\times \{0\}$: in particular, no points of the boundary of $E_\e$ can touch the hyperplane in the interior of $E_0$.

\begin{proposition}[Boundary contact with the hyperplane $\{z=0\}$] \label{prop:traccia}
Assume $\Psi, \Psi^o\in C^\infty\big(\R^{N+1}\setminus\{0\}\big)$ and let $E_\e$ be a minimizer of \eqref{minpb}. Assume also that $\frac1\e \nabla\Psi(e_{_{N+1}})\cdot
e_{_{N+1}}>\|g\|_\infty$. Then, $\pa E_\e\cap (\{z=0\})\subset\pa E_0\times \{0\}$.
\end{proposition}
\begin{proof}
First observe that, since $E_\e$ is a $\|g\|_\infty$-minimizer of $\mathcal{A}^\Psi_\e$ in $\{z>0\}$, by Remark \ref{rem:Lmin} we may identify $E_\e$ with its closed representative
$\overline{E_\e^{(1)}}$ in $\{z>  0\}$.

Fix $x_0 \in \mathrm{int} (E_0)$, and assume that  $B_r(x_0)\subset E_0$. We claim that
\begin{equation}\label{eq:SW}
\mathcal{A}^\Psi_\e(E_\e)\leq \mathcal{A}^\Psi_\e(F)+\|g\|_\infty|E_\e\Delta F| \qquad\text{whenever } F\text{ is s.t. } F\Delta E_\e\Subset B_r(x_0)\times (-1,+\infty),
\end{equation}
i.e., $E_\e$ is a $\|g\|_\infty$-minimizer for $\mathcal{A}^\Psi_\e$ in $ B_r(x_0)\times (-1,+\infty)$. Indeed, let $F$ be a competitor as above  and set $A:=(E_0\times
(-1,0))\setminus F$. Note that the set $F^\sharp:=F \cup (E_0 \times (-1,0))$ is a competitor in \eqref{minpb}. Also note that the set $A$ has finite perimeter and is separated
from the lateral and bottom boundaries of $B_r(x_0) \times (-1,0)$. Then we can estimate
\[
\begin{split}
\mathcal{A}^\Psi_\e(F)-\mathcal{A}^\Psi_\e(F^\sharp)&\geq \int_{\pa^*A\cap \big(B_r(x_0)\times (-1,0)\big)}{e}^{-{z}/\e}\,\Psi(\nu)\, d\H^{N}(x,z)- \int_{\pa^*A\cap \{z=0\}}\Psi(\nu)\, d\H^{N}(x,z)\\
&\geq \int_{\pa^*A\cap \big(B_r(x_0)\times (-1,0)\big)}\Psi(\nu)\, d\H^{N}(x,z)- \int_{\pa^*A\cap \{z=0\}}\Psi(\nu)\, d\H^{N}(x,z)\geq 0,
\end{split}
\]
where the last inequality follows for instance from \cite[Remark~20.3]{MaggiBook}. Thus, by minimality
\[
\begin{split}
\mathcal{A}^\Psi_\e(E_\e)+\int_{E_\e}g\, dxdz&\leq \mathcal{A}^\Psi_\e(F^\sharp)+\int_{F^\sharp} g\, dxdz \\
&\leq \mathcal{A}^\Psi_\e(F)+\int_{F\cup \big(E_0\times (-1,0)\big)} g\, dxdz\,.
\end{split}
\]
Property  \eqref{eq:SW} now easily follows. Thus, both $E_\e$ and $\R^{N+1}\setminus E_\e$ are $\|g\|_\infty$-minimizers for $\mathcal{A}^\Psi_\e$ in $B_r(x_0) \times (-1,\infty)$,
hence we can apply Remark \ref{rem:Lmin} in this cylinder and use canonical representatives there.

Assume by contradiction that $(x_0,0)\in \partial E_\e$. Then, we can apply Lemma~\ref{lm:contact} to $\R^{N+1}\setminus E_\e$, with $U=H^-$ near $(x_0,0)$, in order to obtain that
$\pa E_\e$ is of class $C^{1,\alpha}$ (for some $\alpha\in (0,1)$) in a suitable neighborhood of $(x_0, 0)$.  More precisely, we may find $\eta\in (0, r)$ and a non-negative
function $f\in C^{1,\alpha}(\overline B_\eta (x_0))$ such that $\pa E_\e\cap \big(B_\eta(x_0)\times (-\eta, \eta)\big)$ and $E_\e\cap \big(B_\eta(x_0)\times (-\eta, \eta)\big)$
coincide, respectively, with the graph and subgraph of $f$,  with $f(x_0)=0$, $\|f\|_{\infty}\leq \eta$, and  outer normal $\nu(x_0,0)=e_{_{N+1}}$.

\begin{center}
    \includegraphics[width=6cm]{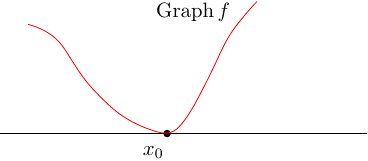}
\end{center}

We now claim that, upon setting
\begin{equation}\label{temp:psi0}
    \psi_0(\xi)\coloneqq \Psi(-\xi,1),\quad \xi\in \R^N,
\end{equation}
there exists a bounded function $h$, with $\|h\|_\infty\leq  e^{\eta/\e}\|g\|_\infty$, such that
\begin{equation}\label{claim:euler:f}
H^\Psi_{\pa E}(\cdot , f(\cdot))=-\Div\big(\nabla{\psi}_0(\nabla f)\big) =\frac1\e \nabla\Psi(\nu)\cdot e_{_{N+1}}+h \quad\text{distributionally in $B_\eta(x_0)$}.
\end{equation}
where $\nu=\nu(x,f(x))=\frac{(-\nabla f(x),1)}{|(-\nabla f(x),1)|}$. Indeed, by the area formula and the homogeneity of $\Psi$, for all $0<\rho<\eta$ we have
\begin{equation*}
    \mathcal{A}_\e^\Psi\big(\mathrm{Subgraph}(f);B_\rho(x_0)\times (-\eta,+\eta)\big)=\int_{B_\rho(x_0)}e^{-f(x)/\e}\,\Psi(-\nabla f,1)\,dx=\int_{B_\rho(x_0)}e^{-f(x)/\e}\,\psi_0(\nabla f)\,dx\,.
\end{equation*}
 Since for every $\vphi\in C^1_c(B_\rho(x_0))$, we have
\begin{equation*}
    \mathcal{L}^{N+1}\Big(\mathrm{Subgraph}(f)\Delta \,\mathrm{Subgraph}(f+t\varphi) \Big)=|t|\,\int_{B_\rho(x_0)}|\vphi(x)|\,dx,
\end{equation*}
then, owing to \eqref{eq:SW}, via a first variation argument we deduce
\begin{equation*}
    \bigg|\int_{B_\rho(x_0)}  \bigg[e^{-f/\e}\nabla \psi_0(\nabla f)\cdot \nabla \vphi-\frac{e^{-f/\e}}{\e}\psi_0(\nabla f)\,\vphi\bigg]\,dx\bigg|\leq \|g\|_\infty\,\int_{B_\rho(x_0)} |\vphi|\,dx.
\end{equation*}
Setting $\tilde\vphi=e^{-f/\e}\vphi$ in the above expression, we get
\begin{equation*}
    \begin{split}
    \|g\|_\infty\,\int_{B_\rho(x_0)}e^{f/\e}\,|\tilde\vphi|\,dx&\geq \bigg|\int_{B_\rho(x_0)} \nabla \psi_0(\nabla f)\cdot \nabla \tilde \vphi\,dx+\frac{1}{\e}\int_{B_\rho(x_0)}\Big[\nabla\psi_0(\nabla f)\cdot \nabla f-\psi_0(\nabla f)\Big]\,\tilde\vphi\,dx\bigg| \\
    &=\bigg|\int_{B_\rho(x_0)} \nabla \psi_0(\nabla f)\cdot \nabla \tilde \vphi-\frac{1}{\e}\int_{B_\rho(x_0)}\big[\nabla \Psi((-\nabla f,1))\cdot e_{_{N+1}}\big]\,\tilde\vphi\,dx\bigg|,
    \end{split}
\end{equation*}
 where in the last identity we used that $\nabla \psi_0(\nabla f)=-\nabla_x \Psi((-\nabla f,1))$, and  $\psi_0(\nabla f)=\Psi((-\nabla f,1))=\nabla \Psi((-\nabla f,1))\cdot
 (-\nabla f,1)$ by homogeneity. Recalling that $\|f\|_\infty\leq \eta$, by the arbitrariness of $\rho$ and $\tilde\vphi\in C^1_c(B_\rho(x_0))$, Claim \eqref{claim:euler:f}
 follows.

Now, recalling that $\nu(x_0,0)=e_{_{N+1}}$ and  by assumption $\frac1\e \nabla\Psi(e_{_{N+1}})\cdot e_{_{N+1}}>\|g\|_\infty$, by taking $\eta$ smaller if needed, by continuity we
deduce
\[
-\Div\big(\nabla{\psi}_0(\nabla f)\big)\geq c>0 \quad\text{distributionally in $B_\eta(x_0)$},
\]
which by strong minimum principle contradicts the fact that $f$ attains its minimum at $x_0$.

Finally, observe that if there exists $x_0\in \mathrm{int\, }(\R^N\setminus  E_0)$ such that $(x_0, 0)\in \pa E_\e$, we reach a contradiction by  the  same argument, taking into
account Remark~\ref{rm:complementare}.
\end{proof}

\begin{remark}\label{remark:elliptic:regul}
Let $\Psi$ be a norm on $\R^{N+1}$ such that $\Psi,\Psi^o\in C^2(\R^{N+1}\setminus\{0\})$, and let $\psi_0(\xi)$ be given by \eqref{temp:psi0}. Then, we have that $\psi_0\in
C^2(\R^N)$, and
\begin{equation}\label{lasua:unifell}
    \lambda_0\,\frac{|v|^2}{(1+|\xi|^2)^{3/2}}\leq \nabla^2\psi_0(\xi)\,v\cdot v \leq \frac{\Lambda_0}{\sqrt{1+|\xi|^2}}\,|v|^2\quad\text{for all $\xi,v\in \R^N$.}
\end{equation}
for some $0<\lambda_0\leq \Lambda_0$. To see this,  let $\tilde v\in \R^{N+1}$ be the orthogonal projection of $(v,0)\in \R^N\times \R$ onto the orthogonal complement of
$(-\xi,1)$, that is $\tilde v=(v,0)+\frac{(v\cdot \xi)}{1+|\xi|^2}(-\xi,1)$.  By the characterization of Lemma \ref{lem:varie:unifconv},  we have
\begin{equation*}
    \begin{split}
        \nabla^2\psi_0(\xi)\,v\cdot v&=\nabla^2 \Psi(-\xi,1)\,(v,0)\cdot (v,0)=\nabla^2\Psi(-\xi,1)\,\tilde v\cdot \tilde v \\
        \frac{\lambda_0}{\sqrt{1+|\xi|^2}}\,|\tilde v|^2&\leq \nabla^2\Psi(-\xi,1)\,\tilde v\cdot \tilde v\leq \frac{\Lambda_0}{\sqrt{1+|\xi|^2}}\,|\tilde v|^2  ,
    \end{split}
\end{equation*}
and since   $ \frac{| v|^2}{1+|\xi|^2}\leq|\tilde v|^2\leq |v|^2$, \eqref{lasua:unifell} follows.

Therefore, if $f\in C^{1,\alpha}(\overline B_R)$ is a weak solution to  the anisotropic mean curvature equation
\[
-\mathrm{div}(\nabla \psi_0(\nabla f))=h\quad\text{in $B_R$,}
\]
 then this is a uniformly elliptic equation as, by the local boundedness of $\nabla f$, the condition \eqref{lasua:unifell} only needs to be evaluated  for $\xi$ in the range of
 $\nabla f$, so that
 \[
 \frac{\lambda_0}{(1+\|\nabla f\|_\infty^2)^{3/2}}\,\mathrm{Id} \leq\nabla^2\psi_0(\xi)\leq \Lambda_0\,\mathrm{Id},\quad\text{for all $\xi\in \R^N$ s.t. $|\xi|\leq \|\nabla f\|_\infty$.}
 \]
  Furthermore, if $h\in L^\infty(B_R)$, then via a difference quotient argument (see \cite[Theorem 8.1]{Ben:Fre})  one gets  $f\in W^{2,2}_{\rm{loc}}(B_R)$, and thus by the chain
  rule $f$ satisfies the equation of trace form
 \[
 -\mathrm{tr}\big(\nabla^2\psi_0(\nabla f\big)\,\nabla^2 f)=h\quad\text{a.e. in $B_R$.}
 \]
Since the coefficients $\nabla^2\psi_0(\nabla f)$ are locally continuous and uniformly elliptic, by standard regularity theory we have $f\in W^{2,p}_{\mathrm{loc}}(B_R)$ for every
$1<p<\infty$, and thus $f\in C^{1,\alpha}_{\mathrm{loc}}(B_R)$ for every $\alpha\in(0,1)$.

\end{remark}

We now state  two useful geometric lemmas.

\begin{lemma}[Anisotropic curvature estimate]\label{lem:curvature}
    Let $\Psi\in C^2(\R^{N+1}\setminus\{0\})$ be a uniformly convex norm and set
    \[
        G:=\frac12\Psi(\cdot,-1)^2.
    \]
    Given $\lambda>0$ and $a\in\R$, let
    \[
        g(x):=\frac{a-G^*(-\lambda x)}{\lambda}.
    \]
    Then, with
    \[
        p:=(\nabla g,-1),
    \]
    one has
    \begin{equation}\label{eq:curvature-estimate}
        \Psi(p)\nabla^2_{\xi\xi}\Psi(p):\nabla^2g \geq -\lambda N.
    \end{equation}
\end{lemma}

\begin{proof}
By Lemma~\ref{lem:varie:unifconv}, there exists $c>0$ such that $\nabla^2G\geq cI_N$ on $\R^N$, where $I_N$ is the identity matrix on $\R^N$. Thus $G$ is uniformly convex and
superlinear, and Lemma~\ref{lem:legendre} yields $G^*\in C^2(\R^N)$.

    From the definition of $g$,
    \[
        \nabla g(x)=\nabla G^*(-\lambda x), \qquad \nabla^2g(x)=-\lambda\nabla^2G^*(-\lambda x).
    \]
By the identities for the Legendre transform in Lemma~\ref{lem:legendre},
    \[
        \nabla^2G(\nabla g)\nabla^2g=-\lambda I_N\quad\text{and hence}\quad  \nabla^2G(\nabla g):\nabla^2g=-\lambda N.
    \]
On the other hand,
    \[
        \nabla^2G(\xi) = \Psi(\xi,-1)\nabla^2_{\xi\xi}\Psi(\xi,-1) + \nabla_\xi\Psi(\xi,-1)\otimes\nabla_\xi\Psi(\xi,-1),
    \]
    and since $\nabla^2g\leq0$, we get
    \[
    \begin{aligned}
        \Psi(p)\nabla^2_{\xi\xi}\Psi(p):\nabla^2g &=-\lambda N -\bigl(\nabla_\xi\Psi(p)\otimes\nabla_\xi\Psi(p)\bigr):\nabla^2g\geq-\lambda N,
    \end{aligned}
    \]
    which proves the claim.
\end{proof}

In what follows, given  a norm $\vartheta$ on $\R^N$, we denote by $W_r^{\vartheta}(x_0)$ the open $\vartheta$-Wulff shape of radius $r$ centered at $x_0$ defined in
\eqref{wulff:shape}, and we also write $W_r^\vartheta:=W_r^\vartheta(0)$.

We will also  make use of the following well-known properties of Legendre transforms of convex functions. Given \(F:\R^d\to\R\), we denote with $F^*$ its Legendre transform
	\[
	F^*(p):=\sup_{y\in\R^d} \bigl\{\langle p,y\rangle-F(y)\bigr\}.
	\]
\begin{lemma}[Properties of the Legendre transform]\label{lem:legendre}
The set of functions \(F:\R^d\to\R\) that are  strictly convex, superlinear and $C^1(\R^d)$ is self dual with respect to the Legendre transform. That is, $F$ is strictly convex,
superlinear and $C^1(\R^d)$ if and only if $F^*$ has the same properties.

    Moreover, under these assumptions on $F$, then $\nabla F^*$ is the inverse map of $\nabla F$; that is, the identity
	\begin{equation}\label{id:grad}
		\nabla F^*(\nabla F(y))=y
	\end{equation}
	holds for every $y$. Finally, if in addition $F\in C^k(U)$, $k\geq 2$, for some open set $U\subset\R^d$ and
	\begin{equation}\label{ip:hess}
	\nabla^2F(y)>0 \qquad\text{for every }y\in U,
\end{equation}
then    $F^*\in C^k(V)$, where $V:=\nabla F(U)$, and
	\begin{equation}\label{id:hess}
	\nabla^2F^*(\nabla F(y)) = \bigl(\nabla^2F(y)\bigr)^{-1}.
	\end{equation}
	for all  $y\in U$.
\end{lemma}
\begin{proof} The above properties are classical:  the fact that $\nabla F$ is invertible and \eqref{id:grad} follow for instance from \cite[Theorems~26.5 and~26.6]{RockafellarConvexAnalysis}. The regularity properties of $F^*$  are  an easy consequence of the Inverse Function Theorem  while \eqref{id:hess} can be derived by differentiating  \eqref{id:grad}.
\end{proof}

Next, for  $\Psi,\Psi^o\in C^\infty(\R^{N+1}\setminus\{0\})$, with $\Psi$ uniformly convex and symmetric with respect to the last variable, we set
\[
    \psi(\xi):=\Psi(\xi,0), \qquad G:=\frac12\Psi(\cdot,-1)^2.
\]
Then $G$ and $G^*$ are even and $G^*$ attains its minimum at the origin. Given $r,\lambda>0$, set
\begin{equation}\label{abs:defa}
    a^\Psi_{\lambda,r} := \min_{\psi^o(y)=r} \bigl\{G^*(\lambda y)-G^*(0)\bigr\}>0.
\end{equation}
The minimum exists since $\{\psi^o=r\}$ is compact, and it is positive by the strict convexity of $G^*$. For $x_0\in\R^N$ and $r,\eta,\e>0$, define the intrinsic paraboloid
\begin{equation}\label{abs:paraboloidi}
    \mathcal C^\Psi_{r,\eta,\e}(x_0) := \left\{(x,z)\in\R^{N+1}:\,0<z\leq g^\Psi_{r,\eta,\e}(x-x_0)\right\},
\end{equation}
where
\begin{equation}\label{abs:defg}
    g^\Psi_{r,\eta,\e}(x) := \frac{G^*(0)+a^\Psi_{\lambda,r}-G^*(-\lambda x)}{\lambda}, \qquad \lambda:=\frac{2\eta}{\e}.
\end{equation}
By convexity of $G^*$ and the fact that it attains its minimum at the origin,
\begin{equation}\label{abs:footprint}
    \{g^\Psi_{r,\eta,\e}>0\}\subset W^\psi_r.
\end{equation}
Indeed, if $\psi^o(x)\geq r$, writing $x=ty$ with $\psi^o(y)=r$ and $t\geq1$, the function $s\mapsto G^*(sy)$ is non-decreasing for $s\geq0$, and hence
\[
    G^*(\lambda x)-G^*(0) =G^*(\lambda ty)-G^*(0) \geq a^\Psi_{\lambda,r}.
\]

Moreover, there exist constants $c,C>0$, depending only on $N,\Psi$, such that, for every $r,\eta,\e>0$,
\begin{equation}\label{abs:height-comparison}
c\,\frac{\eta r^2}{\e}
\leq g^\Psi_{r,\eta,\e}(0)
\leq C\,\frac{\eta r^2}{\e}.
\end{equation}
Indeed, by Lemma~\ref{lem:varie:unifconv}, using the smoothness and uniform convexity of $\Psi$, together with the zero-homogeneity of $\nabla^2(\Psi^2)$, we obtain
\[
c\,\mathrm{Id}_N
\leq\nabla^2G(\xi)
\leq C\,\mathrm{Id}_N
\qquad\text{for every }\xi\in\R^N.
\]
By \eqref{id:hess}, the Hessian of $G^*$ satisfies analogous bounds. Since $G^*$ is even, $\nabla G^*(0)=0$, and Taylor's formula therefore yields
\[
c|\xi|^2
\leq
G^*(\xi)-G^*(0)
=
\int_0^1(1-t)\,
\nabla^2G^*(t\xi)\xi\cdot\xi\,dt
\leq C|\xi|^2,
\]
after possibly changing the constants. Evaluating this estimate at $\xi=\lambda y$, taking the minimum over $\{\psi^o(y)=r\}$, and using equivalence of norms, we obtain
\[
c\lambda^2r^2
\leq a^\Psi_{\lambda,r}
\leq C\lambda^2r^2.
\]
Finally, by \eqref{abs:defg} we obtain
\[
g^\Psi_{r,\eta,\e}(0)
=\frac{a^\Psi_{\lambda,r}}{\lambda},
\qquad
\lambda=\frac{2\eta}{\e},
\]
which proves \eqref{abs:height-comparison}. In particular, since $G^*$ has a minimum at the origin,
\begin{equation}\label{abs:smallheight}
    g^\Psi_{r,\eta,\e}(x)\leq \frac{a^\Psi_{\lambda,r}}{\lambda}\approx \frac{\eta\, r^2}{\e} \longrightarrow0 \qquad\text{as $\eta\to0^+$,}
\end{equation}

\medskip
\begin{lemma}[Anisotropic paraboloid barriers]\label{lm:paraboloidi}
Let $g\in\mathrm{Lip}_c(\R^{N+1})$ and let $E_\e$ be a minimizer of \eqref{minpb}. Assume that
\begin{equation}\label{abs:trace-hyp}
    \frac1\e\Psi(e_{_{N+1}})>\|g\|_\infty.
\end{equation}
Let $r>0$, $x_0\in\R^N$, and assume that
\begin{equation}\label{abs:wulff-in}
    \overline{W^\psi_r(x_0)}\subset\mathrm{int}(E_0).
\end{equation}
If there exists $\eta_0>0$ such that, for every $\eta\in(0,\eta_0]$ and every $y$ with $g^\Psi_{r,\eta,\e}(y)>0$, upon setting
\[
    p:=\bigl(\nabla g^\Psi_{r,\eta,\e}(y),-1\bigr),
\]
one has
\begin{equation}\label{abs:vertical}
\begin{split}
    \frac1\e\Psi(p)\partial_\zeta\Psi(p) &+ e^{g^\Psi_{r,\eta,\e}(y)/\e} g\bigl(x_0+y,g^\Psi_{r,\eta,\e}(y)\bigr)\Psi(p) <-\frac{2N\eta_0}{\e},
\end{split}
\end{equation}
then
\[
    \mathcal C^\Psi_{r,\eta,\e}(x_0)\subset E_\e \qquad\text{for all $\eta\in(0,\eta_0]$.}
\]
Analogously, if
\[
\overline{W^\psi_r(x_0)} \subset\R^N\setminus\overline{E_0},
\]
the corresponding inclusion holds for $\R^{N+1}\setminus E_\e$ provided that \eqref{abs:vertical} holds with $g$ replaced by $-g$.
\end{lemma}

\begin{proof}
By \eqref{abs:footprint}, \eqref{abs:wulff-in}, Proposition~\ref{prop:traccia}, the definition \eqref{abs:paraboloidi} and  \eqref{abs:smallheight}, we have
\[
    \mathcal C^\Psi_{r,\eta,\e}(x_0)\subset E_\e
\]
for $\eta>0$ sufficiently small. Assume by contradiction that the conclusion fails and set
\[
    \bar\eta := \sup\left\{s\in(0,\eta_0]:\, \mathcal C^\Psi_{r,\eta,\e}(x_0)\subset E_\e \text{ for every }\eta\in(0,s]\right\}.
\]
Therefore there exists a touching point
\[
    (\bar x,\bar z) \in \pa\mathcal C^\Psi_{r,\bar\eta,\e}(x_0)\cap\pa E_\e,
\]
with $\bar z>0$. Applying Lemma~\ref{lm:contact} to the complementary minimizer, and then the standard elliptic regularity for the Euler--Lagrange equation (see Remark
\ref{remark:elliptic:regul}), we infer that in a neighborhood of $(\bar x,\bar z)$ the boundary $\pa E_\e$ coincides with the graph of a function $f=f(x)$ of class $C^{2,\alpha}$
solving, with $p_f=(\nabla f,-1)$,
\begin{equation}\label{abs:EL}
    \Psi(p_f)\nabla^2_{\xi\xi}\Psi(p_f):\nabla^2f = \frac1\e\Psi(p_f)\partial_\zeta\Psi(p_f) + e^{f/\e}g(x,f(x))\Psi(p_f).
\end{equation}
At the touching point, setting
\[
    p:= \bigl( \nabla g^\Psi_{r,\bar\eta,\e}(\bar x-x_0),-1 \bigr),
\]
we have \(\bar z=g^\Psi_{r,\bar\eta,\e}(\bar x-x_0), \) and
\[
    \nabla g^\Psi_{r,\bar\eta,\e}(\bar x-x_0) =\nabla f(\bar x), \qquad \nabla^2g^\Psi_{r,\bar\eta,\e}(\bar x-x_0) \leq\nabla^2f(\bar x).
\]

\begin{center}
\includegraphics[width=9cm]{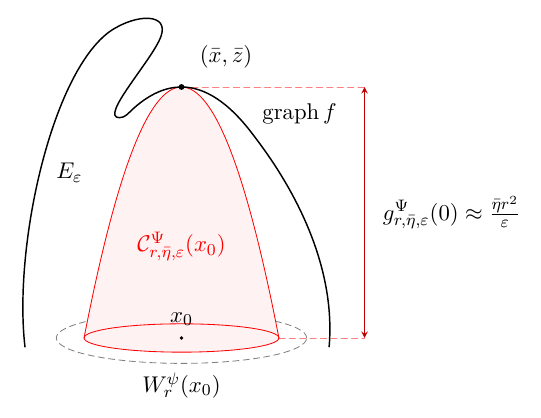}
\end{center}

Thus, by Lemma~\ref{lem:curvature}, \eqref{abs:EL}, and \eqref{abs:vertical},
\[
\begin{aligned}
    -\frac{2\bar\eta N}{\e} &\leq \Psi(p)\nabla^2_{\xi\xi}\Psi(p): \nabla^2g^\Psi_{r,\bar\eta,\e}(\bar x-x_0) \\
    &\leq \Psi(p)\nabla^2_{\xi\xi}\Psi(p):\nabla^2f(\bar x) \\
    &= \frac1\e\Psi(p)\partial_\zeta\Psi(p) + e^{\bar z/\e}g(\bar x,\bar z)\Psi(p) <-\frac{2N\eta_0}{\e}.
\end{aligned}
\]
This implies $\bar\eta>\eta_0$, a contradiction. The complementary statement follows from Remark~\ref{rm:complementare}.
\end{proof}

\section{Elliptic regularization with natural mobility: Kuratowski convergence}\label{sec:kur:conv}

Throughout this section and the next one, we restrict to the important special case of \emph{natural mobility}, namely $\mathfrak m=\phi$. Thus, the geometric evolution law
considered throughout these two sections is
\begin{equation}\label{eq:natural-mobility-flow}
    V=-\phi(\nu)\,H^\phi.
\end{equation}
In accordance with the terminology introduced in Section~\ref{subsec:cf}, the corresponding superflows and subflows would be $(\phi,\phi)$-superflows and $(\phi,\phi)$-subflows.
Since no other mobility will be considered in these two sections, we shall simply call them \emph{superflows} and \emph{subflows} for~\eqref{eq:natural-mobility-flow}.

Let $\phi$ be an \emph{arbitrary norm} on $\R^N$.  We consider  its Euclidean extension  to $\R^{N+1}=\R^N\times \R$, defined by \beq\label{Phi}
\Phi(\xi,\zeta):=\sqrt{\phi^2(\xi)+\zeta^2}, \eeq for all $\xi\in \R^N$ and $\zeta\in \R$. We assume
\begin{equation}\label{mis:deE0}
   E_0\subset\R^N\,\text{  is a set of finite perimeter.}
\end{equation}
\emph{Here and in what follows, we always identify $E_0$  with its precise representative, namely the set of its points of density one.}

The anisotropic elliptic regularisation scheme for the flow \eqref{evol law}  starting from $E_0$ will consist in  approximating the flow with \beq\label{minpbell} E_\e\in
\operatorname*{argmin}\Bigl\{\mathcal{A}^\Phi_\e(F): F\subset\R^{N+1}\text{ s.t. }F\cap H^{-}=E_0\times (-\infty,0)\Bigr\}, \eeq where $\mathcal{A}^\Phi_\e$ is defined as in
\eqref{defAF}, with $\Phi$ in place of $\Psi$. Note that by the results of Section~\ref{sec:waf}, we may assume without loss of generality that $E_\e\cap \{z>0\}$ is a (relatively)
closed set, and
\[
\partial E_\e=\operatorname{spt}\mu_{E_\e}.
\]

Following \cite{Ilm94}, we introduce the {\em translative solitons} \beq\label{solitoni} \mathcal S_\e(t):=\overline{E}_\e-\frac{t}\e e_{_{N+1}} \eeq for all $t\geq 0$ and  the
associated space-time tube \beq\label{tubosolitoni} \mathcal S_\e:=\big\{(x,z, t)\in \R^{N+1}\times[0,+\infty):\, (x,z)\in \mathcal S_\e(t)\big\}. \eeq Moreover, for $t\geq 0$, set
\begin{equation}\label{def:Schius}
\mathcal S_\e^\complement(t)
:=\overline{\R^{N+1}\setminus\mathcal S_\e(t)}
=\R^{N+1}\setminus
\operatorname{int}\bigl(\mathcal S_\e(t)\bigr),
\end{equation}
 and
\begin{equation}\label{def:Schiu1}
\mathcal{S}_{\e}^\complement:=\left\{(x,z,t)\in \R^{N+1}\times[0,\infty):\, (x,z)\in \mathcal{S}_{\e}^\complement(t)\right\}.
\end{equation}
Given $F\subset \R^{N+1}$ we denote by $\mathcal R_\e(F)$ the vertically rescaled and truncated set defined as
\[
\mathcal R_\e(F):=\{(x,\e z):\, (x,z)\in F\text{ and }z>0\}.
\]
Our first goal is to prove the following  theorem. In what follows, we will interpret the rescaled variable $\e z$ as a time variable denoted by $t$.
\begin{theorem}\label{thm:main1}
Let $\phi, \Phi, E_0$, $E_\e$, $\mathcal S_\e$ and $\mathcal S_\e^\complement$ be as above. Assume also that $E_0$ is non-fattening for \eqref{eq:natural-mobility-flow} in the
sense of Definition~\ref{def:fat}  and let $(E(t))_{t> 0}$ be the (unique) generalised motion by \eqref{eq:natural-mobility-flow} starting from $E_0$  as defined in
Definition~\ref{def:gen-motion}.

Denote by $E$ the associated space-time tube
\[
E:=\big\{(x,t)\in \R^{N}\times(0,+\infty):\, x\in E(t)\big\}
\]
 and by $\mathcal E$ the space-time tube built upon  the cylindrical evolution  $\mathcal E(t):=E(t)\times \R$; i.e.,
 \[
 \mathcal E:=\big\{(x,z, t)\in \R^{N+1}\times(0,+\infty):\, x\in E(t)\big\}.
 \]

Then,  the following Kuratowski limits:
\begin{itemize}
\item[(a)]
\[
\mathcal S_\e\stackrel{\mathcal K}{\longrightarrow} \mathcal E \quad\text{in } \R^{N+1}\times(0,+\infty);
\]
\item[(b)]  for all but countably many $t>0$
\[
\mathcal S_{\e}(t)\stackrel{\mathcal K}{\longrightarrow}  E(t)\times \R\quad\text{in }\R^{N+1};
\]
\item[(c)]
\[
\mathcal R_\e(E_\e) \stackrel{\mathcal K}{\longrightarrow}E    \quad\text{in }\R^{N}\times (0,+\infty)
\]
\end{itemize}
hold true as $\e\to 0^+$.
 \end{theorem}

 \begin{remark}\label{remark:complementary}
Under the assumptions of Theorem~\ref{thm:main1}, let $u$ be the level-set flow associated with a non-fattening level-set representation of $E_0$, so that
$E(t)=\{u(\cdot,t)\leq0\}$, where $u$ is as in Definition~\ref{def:fat}. We observe that $-u$ is the level-set flow associated with a non-fattening representation of $\R^N\setminus
E_0$.

Therefore, by applying Theorem~\ref{thm:main1} to the complementary initial datum $\R^N\setminus E_0$ in place of $E_0$, we clearly deduce the following Kuratowski limits:
   \[
    \begin{split}
        &\mathcal S_\e^\complement\stackrel{\mathcal K}{\longrightarrow} \big(\R^{N+1}\times (0,\infty)\big)\setminus \mathrm{int}\big(\mathcal E \big)\quad\text{in }\R^{N+1} \times (0,+\infty),\\
        &\mathcal S_{\e}^\complement(t)\stackrel{\mathcal K} {\longrightarrow}\Big(\R^{N}\setminus \mathrm{int}\big({E}(t)\big)\Big)\times \R\quad\text{in $\R^{N+1}$ for all but countably many $t>0$},\\
        &   \big(\R^{N}\times (0,\infty)\big)\setminus \mathcal R_\e\big(\mathrm{int} (E_\e)\big) \stackrel{\mathcal K}{\longrightarrow}\big(\R^{N}\times (0,\infty)\big)\setminus \mathrm{int}(E)  \quad\text{in }\R^{N}\times (0,+\infty)
    \end{split}
    \]
    as $\e\to 0^+$. Here we used that, by the non-fattening condition,
    \[
\operatorname{int}E(t)=\{u(\cdot,t)<0\},
\qquad
\partial E(t)=\{u(\cdot,t)=0\}\quad\text{for all but countably many $t>0$.}
\]

    By combining the convergences above with the results of Theorem~\ref{thm:main1} and setting
\[
\Gamma_\e:=\mathcal R_\e(\partial E_\e),
\qquad
\Gamma:=\{(x,t)\in\R^N\times(0,+\infty):u(x,t)=0\},
\]
 we deduce also the following Kuratowski limits:
    \[
    \begin{array}{rcll}
         \pa \mathcal S_{\e}(t) &\stackrel{\mathcal K} {\longrightarrow}  &\pa E(t)\times \R &\quad\text{in $\R^{N+1}$ for all but countably many $t>0$},\\
        \Gamma_\e &\stackrel{\mathcal K}{\longrightarrow} & \Gamma  &\quad\text{in }\R^{N}\times (0,+\infty)
    \end{array}
          \]
    as $\e\to 0^+$.
  \end{remark}

The proof of Theorem~\ref{thm:main1} will be a consequence of several intermediate results that will be established  throughout the following subsections.

 {The cases $E_0=\varnothing$ and $E_0=\R^N$ are immediate; in the proof below we therefore assume that $E_0$ is nonempty and $E_0\neq \R^N$.}

\subsection{The approximate translative solitons}\label{sec:asol}
In this section we fix a sequence $\{\Phi_n\}$ of  norms on $\R^{N+1}$ such that
\begin{equation}\label{eq:Phin}
    \begin{split}
    & \{\Phi_n\}\subset C^\infty(\R^{N+1}\setminus\{0\}),\\
    &\Phi_n\to \Phi \quad \text{locally uniformly in  $\R^{N+1}$}, \vspace{1cm}\\
    &c_n\, \mathrm{Id} \leq \frac12\nabla^2 \Phi_n^2\leq C_n\, \mathrm{Id}\quad\text{in $\R^{N+1}\setminus\{0\}$  for some $C_n, c_n>0$.}
\end{split}
\end{equation}
Such a sequence $\Phi_n$ can be constructed, for instance, by taking a sequence $\tilde{\Phi}_n$ of $C^\infty(\R^{N+1}\setminus\{0\})$-smooth norms which converge locally uniformly
to $\Phi$ (see, e.g., Section \ref{sec:normeunifconvesse}), and by taking $\Phi_n(\xi,\zeta)\coloneqq \sqrt{\tilde{\Phi}_n^2(\xi,\zeta)+\tfrac{1}{n}|(\xi,\zeta)|^2}$ which
satisfies the third condition of \eqref{eq:Phin} by the discussion of \cite[Example 2.1]{ACF23}.

Since $\Phi$ is symmetric with respect to the last variable, we may also assume without loss of generality that
\begin{equation}\label{eq:Phin-sym}
    \Phi_n(\xi,\zeta)=\Phi_n(\xi,-\zeta) \qquad\text{for all $(\xi,\zeta)\in\R^{N+1}$.}
\end{equation}
Indeed, if necessary, it is enough to replace $\Phi_n$ by $\sqrt{\frac12\bigl(\Phi_n^2(\xi,\zeta)+\Phi_n^2(\xi,-\zeta)\bigr)}$ which still satisfies \eqref{eq:Phin}, with possibly
different constants $c_n,C_n$.

\begin{remark}\label{rm:Phin}
	Set $F_n:=\frac12(\Phi_n)^2$. Then one can easily check that the Legendre transform $F_n^*$ is given by         $F_n^*=\frac12(\Phi^o_n)^2$, with $\Phi_n^o$ denoting the dual norm
	of $\Phi_n$. Thus, by Lemma~\ref{lem:legendre} (applied with $U=\R^{N+1}\setminus\{0\}$) we deduce that $\frac12(\Phi^o_n)^2$ and, in turn,  $\Phi^o_n$ belong to
	$C^{\infty}(\R^{N+1}\setminus\{0\})$ and that
	\begin{equation}\label{eq:Phino}
	\frac1{C_n}\, \mathrm{Id} \leq \frac12\nabla^2 (\Phi^o_n)^2\leq \frac1{c_n}\, \mathrm{Id}\quad\text{in $\R^{N+1}\setminus\{0\}$}.
	\end{equation}
	Note that the last condition in \eqref{eq:Phin} and \eqref{eq:Phino} imply that  the $\Phi_n$-ball $\mathcal{B}^{\Phi_n}$ and the $\Phi_n$-Wulff shape $\mathcal{W}^{\Phi_n}$
	(=$\mathcal{B}^{\Phi^o_n}$) are uniformly convex and have positive principal curvatures at every boundary point.  Moreover, from the second condition in \eqref{eq:Phin} and the
	definition of dual norms one easily infers that
	\begin{equation}\label{convpolari}
	\Phi^o_n\to \Phi^o \quad \text{locally uniformly in  $\R^{N+1}$.}
	\end{equation}
	This is equivalent to  the existence of $\delta_n\searrow 0$ such that
	\[
	 (1-\delta_n)\Phi_n^o\leq \Phi^o\leq (1+\delta_n)\Phi_n^o,
	\]
	which in turn implies
	\begin{equation}\label{convdist}
	(1-\delta_n)\dist^{\Phi_n^o}(\cdot, C)\leq \dist^{\Phi^o}(\cdot, C)\leq 	(1+\delta_n)\dist^{\Phi_n^o}(\cdot, C)
	\end{equation}
	 for any closed set $C\subset\R^{N+1}$. Finally, from \eqref{id:grad} we deduce that
     \[
	\Phi_n ^o(\xi,\zeta)\nabla\Phi_n(\nabla\Phi_n^o(\xi, \zeta))=\big(\nabla(\tfrac12\Phi_n^2)\circ \nabla(\tfrac12(\Phi^o_n)^2)\big)(\xi, \zeta)=(\xi,\zeta)
	\]
	for all $(\xi, \zeta)\in (\R^N\times \R)\setminus\{(0,0)\}$.  By differentiating the above identity we get
	\begin{multline}\label{usefulid}
\nabla^2\Phi_n(\nabla\Phi_n^o): \nabla^2\Phi_n^o=\Div_{\xi, \zeta}( \nabla\Phi_n(\nabla\Phi_n^o))=	\Div_{\xi, \zeta}\Big(\frac{(\xi,\zeta)}{\Phi_n^o}\Big)
\\
=\frac{N+1}{\Phi_n^o}-\frac{\nabla\Phi_n^o(\xi,\zeta)\cdot(\xi,\zeta)}{(\Phi_n^o)^2}=\frac{N}{\Phi_n^o}\quad \text{in $(\R^N\times \R)\setminus\{(0,0)\}$,}
	\end{multline}
	where in the last equality we used the fact that $\nabla\Phi_n^o(\xi,\zeta)\cdot(\xi,\zeta)=\Phi_n^o(\xi,\zeta)$ by one-homogeneity.
\end{remark}

We will also need the following property, which shows that $\partial_{\zeta}\Phi_n(\xi,\zeta)\to \partial_\zeta \Phi(\xi,\zeta)=\frac{\zeta}{\Phi(\xi,\zeta)}$.

\begin{lemma}\label{lem:vertical} Let $\Phi_n=\Phi_n(\xi,\zeta)$ be as above. Then we have that
	\begin{equation}\label{eq:vertical-limit}
		\nabla \Phi_n(\xi,\zeta)\cdot e_{_{N+1}}=\partial_{\zeta}\Phi_n(\xi, \zeta)\to \frac{\zeta}{\Phi(\xi, \zeta)} \quad\text{locally uniformly in  $(\R^{N}\times\R)\setminus\{(0,0)\}$}.
	\end{equation}
\end{lemma}

\begin{proof}
	It is enough to prove that  if $(\xi_n, \zeta_n)\to (\xi, \zeta)\neq (0,0)$ then
	\begin{equation}\label{enough}
\partial_{\zeta}\Phi_n(\xi_n, \zeta_n)\to \frac{\zeta}{\Phi(\xi, \zeta)}.
	\end{equation}
	To this aim, set $v_n:=\nabla\Phi_n (\xi_n, \zeta_n)$. We may assume (up to extracting a subsequence if needed) that $v_n\to v$. Since $\Phi_n \to \Phi$  and $(\xi_n, \zeta_n)\to
	(\xi, \zeta)$ we get         $v\in \partial \Phi(\xi, \zeta)$.  Notice that by \eqref{Phi}, $\Phi$ is differentiable with respect to the variable $\zeta$, with
	$\partial_{\zeta}\Phi(\xi, \zeta)=\frac{\zeta}{\Phi(\xi, \zeta)}=w\cdot e_{_{N+1}}$ for any $w\in \partial \Phi(\xi, \zeta)$. Thus,
	\[
	\partial_{\zeta}\Phi_n(\xi_n, \zeta_n)=v_n\cdot e_{_{N+1}}\to v\cdot e_{_{N+1}}=\frac{\zeta}{\Phi(\xi, \zeta)},
	\]
	which establishes \eqref{enough}.
\end{proof}

We  fix a radius $R>0$ arbitrarily large; let $\zeta\in C^\infty_c(\B_2)$ be a cut-off function such that $0\leq \zeta\leq 1$,  $\zeta\equiv 1$ on $\B_1$, and we set
$\zeta_R(x,z)=\zeta\big((x,z)/R\big)$ for  $(x,z)\in \R^{N+1}$. For $\e\in (0,1)$ fixed, let
\begin{equation}\label{def:gR}
    g_R:=\mathrm{dist}(\cdot,\partial E_\e)\,\zeta_R,\quad \mathrm{s}g_R=\mathrm{sd}_{E_\e}(\cdot)\,\zeta_R
\end{equation}
where $\mathrm{sd}_{E_\e}$ is the signed distance function defined by \eqref{signeddist} with $\mathfrak M^o=|\cdot|$. We take a small parameter  such that
\begin{equation}\label{smallde}
    \delta\in (0,1)\quad \text{s.t.}\quad \delta\,\|\mathrm{s}g_R\|_{L^\infty(\R^{N+1})}<1,
\end{equation}
 and we consider a minimizer
\begin{equation}\label{min:Een}
    E_{n,\de,R}\in \operatorname*{argmin}\Bigl\{\mathcal{A}^{\Phi_n}_\e(F)+\de\int_{F\Delta E_\e}g_R\,dxdz: F\subset\R^{N+1}\text{ s.t. }F\cap H^{-}=E_0\times (-\infty,0)\Bigr\},
\end{equation}
where $\mathcal{A}^{\Phi_n}_\e$ is defined as in \eqref{defAF}, with $\Phi_n$ in place of $\Psi$.  Note that to ease the notation we dropped the dependence on $\e$ of the minimizer
$E_{n,\de,R}$.
\begin{remark}\label{rm:complementarebis}
Note that  $\R^{N+1}\setminus E_{n,\de,R}$ is a minimizer for the problem
\begin{equation*}
\min\Bigl\{\mathcal{A}^{\Phi_n}_\e(F)+\de\int_{F\Delta (\R^{N+1}\setminus E_\e)}g_R\,dxdz: F\subset\R^{N+1} \text{ s.t. }F\cap H^{-}=(\R^N\setminus  E_0)\times (-\infty,0)\Bigr\}.
\end{equation*}
\end{remark}

One can easily check that
\begin{equation*}
    \int_{F\Delta E_\e}g_R\,dxdz=\int_F \mathrm{s}g_R\,dxdz-\int_{E_\e}\mathrm{s}g_R\,dxdz,
\end{equation*}
so the minimization problem \eqref{min:Een} is equivalent to
\begin{equation}\label{min:Een1}
     E_{n,\de,R}\in\operatorname*{argmin}\Bigl\{\mathcal{A}^{\Phi_n}_\e(F)+\de\int_{F}\mathrm{s}g_R\,dxdz\,: F\subset\R^{N+1}\text{ s.t. }F\cap H^{-}=E_0\times (-\infty,0)\Bigr\}.
\end{equation}
Observe that, as $\nabla\Phi_n(e_{_{N+1}})\cdot e_{_{N+1}}\xrightarrow{n\to\infty} 1$ by Lemma \ref{lem:vertical} and $\Phi(e_{_{N+1}})=1$, from \eqref{smallde}  we have
$\frac{1}{\e}\nabla\Phi_n(e_{_{N+1}})\cdot e_{_{N+1}}>\delta\,\|\mathrm{s}g_R\|_{L^\infty(\R^{N+1})}$ for $n>n_0$ large enough, so  by \eqref{min:Een1} and
Proposition~\ref{prop:traccia} we have that
\begin{equation}\label{traccia:zero}
    \partial E_{n,\de,R}\cap (\{z=0\})\subset \partial E_0\times \{0\}.
\end{equation}
By the classical partial regularity for minimizers of elliptic integrands \cite{SSA77,Duzaar:Steffen}, we have that $\partial^* E_{n,\de,R}\in C^{1,\alpha}$ on $\{z>0\}$ for some
$\alpha\in (0,1)$, and by Proposition \ref{prop:Euler_lagrange}, we have $\mathbf{H}^{\Phi_n}_{\partial^* E_{n,\de,R}}=-H^{\Phi_n}_{\partial^* E_{n,\de,R}}\,\nu_{E_{n,\de,R}}$,
with
\begin{equation}\label{ELn}
    H^{\Phi_n}_{\partial E_{n,\de,R}}=\frac{1}{\e} \partial_\zeta\Phi_n(\nu)- \delta\, e^{z/\e}\mathrm{s}g_R,\quad\text{on $\partial^* E_{n,\de,R}\cap \{z>0\}$.}
\end{equation}
In turn, the Lipschitz continuity of $\mathrm{s}g_R$,  Remark  \ref{remark:elliptic:regul} and  standard elliptic regularity theory \cite[Theorem 9.19]{Gilbarg:Trudinger} entail
that
\begin{equation}\label{regolarita}
    \partial^* E_{n,\de,R}\,\text{ is of class }\, C^{2,\beta}\,\text{ on $\{z>0\}$, for all $\beta\in (0,1)$.}
\end{equation}
Hence by \eqref{ELn} and Lemma~\ref{lem:vertical}, we infer
\begin{equation}\label{ineq:AMC}
     H^{\Phi_n}_{\partial E_{n,\de,R}}\leq\frac{1}{\e} \frac{\nu\cdot e_{_{N+1}}}{\Phi_n(\nu)}+\delta+ \delta\, e^{2R/\e}\|\mathrm{s}g_R\|_\infty\quad\text{on $\partial^*E_{n,\delta,R}\cap \{z>0\}$}
\end{equation}
for $n>n_0$ sufficiently large. Next, by taking $F=E_0\times \R$ as a competitor in \eqref{min:Een}, and recalling \eqref{eq:Phin}, one can easily check that
$\sup_{n}\mathcal{A}_\e^{\Phi_n}(E_{n,\de,R})<\infty$, and therefore
\begin{equation}\label{per:nde}
    \sup_{n>n_0} P(E_{n,\de,R};\B_M)<C_M\quad \text{for every $M>0$}
\end{equation}
Moreover, the sets $E_{n,\de,R}$ are evidently $\Lambda$-minimizers of the functional $\mathcal{A}^{\Phi_n}_\e$ in $\R^N\times \{z>0\}$, with
\begin{equation}\label{lambda1}
\Lambda=\delta\,\|\mathrm{s}g_R\|_\infty<1
\end{equation}
by our choice of $\delta$. Thus we may assume that $E_{n,\de,R}$ coincides with its closed precise representative in $\{z>0\}$--see Remark \ref{rem:Lmin}.

Additionally, since $\Phi_n\to\Phi$ locally uniformly by \eqref{eq:Phin}, there exist constants $c,C>0$, independent of $n$, such that $c|v|\leq\Phi_n(v)\leq C|v|$ for all
$v\in\R^{N+1}$ and all $n\in\N$. Hence one can check that the following (uniform in $n\in \N$) volume density estimates hold true: for every $M>0$, there exists $r_0=r_0(M)$ and
$c_0=c_0(M)$ such that,
\begin{equation}\label{density:nd}
    \mathcal{L}^{N+1}\big(E_{n,\de,R}\cap \B_{r}(x,z)\big)\geq c_0\,r^{N+1}
\end{equation}
for all $0<r<r_0$, and all $(x,z)\in E_{n,\de,R}$ such that $1/M<z<M$-- see, e.g., the argument of \cite[Lemma~2.8]{DPhMa15}. Then, again standard arguments stemming from
\eqref{per:nde} and \eqref{density:nd} yield that (up to a non-relabelled subsequence)
\begin{equation}\label{conv:ndelta}
    E_{n,\de,R}\xrightarrow{n\to \infty} E_{\delta,R}\quad \text{in the Kuratowski sense on $\{z>0\}$   and in $L^1_{loc}(\{z> 0\})$ },
\end{equation}
for some set $E_{\delta,R}$ of locally finite perimeter, relatively closed in $\{z>0\}$. We extend it to $\{z<0\}$ by setting
\[
E_{\delta,R}\cap H^-=E_0\times(-\infty,0).
\]
Since the approximating sets $E_{n,\delta,R}\equiv E_0\times (-\infty,0]$ in $\{z\leq 0\}$  by \eqref{min:Een}, the convergence also holds in $L^1_{\mathrm{loc}}(\R^{N+1})$.

We now claim that
\begin{equation}\label{ugualeBR}
    E_{\de,R}\cap \B_R=E_\e\cap \B_R.
\end{equation}
Indeed, by the minimality \eqref{minpbell}, and since  $E_{\de,R}\cap H^-=E_0\times(-\infty,0)$, we have
\begin{equation*}
    \mathcal{A}^{\Phi}_\e(E_\e)\leq \mathcal{A}^{\Phi}_\e(E_{\delta,R})+\de\int_{E_{\delta,R}\Delta E_\e}g_R\,dxdz,
\end{equation*}
with strict inequality if $\de\int_{E_{\delta,R}\Delta E_\e}g_R\,dxdz\neq 0$, while by the minimality \eqref{min:Een}, \eqref{conv:ndelta}, \eqref{eq:Phin} and the lower
semicontinuity of  the functional $\mathcal{A}^{\Phi_n}_\e$, we have
\begin{equation*}
    \begin{split}
     \mathcal{A}^{\Phi}_\e(E_{\delta,R})+\de\int_{E_{\delta,R}\Delta E_\e}g_R\,dxdz&\leq   \liminf_{n\to\infty}\bigg( \mathcal{A}^{\Phi_n}_\e(E_{n,\de,R})+\de\int_{E_{n,\de,R}\Delta E_\e}g_R\,dxdz\bigg) \\
     &\leq \lim_{n\to\infty} \mathcal{A}^{\Phi_n}_\e(E_\e)=\mathcal{A}^{\Phi}_\e(E_\e).
    \end{split}
\end{equation*}
Combining the two inequalities above, we infer that  $\int_{E_{\delta,R}\Delta E_\e}g_R\,dxdz=0$, from which claim \eqref{ugualeBR} follows.

We now introduce the approximate translative solitons
\begin{equation}\label{approxtubosolitoni}
    \mathcal{S}_{n,\de,R}(t)=E_{n,\de,R}-\frac{t}{\e}e_{_{N+1}}\,,
\end{equation}
and, correspondingly, the limiting translating solitons
\begin{equation}\label{approxsol2}
     \mathcal{S}_{\de,R}(t)=E_{\de,R}-\frac{t}{\e}e_{_{N+1}}\,.
\end{equation}
Clearly, from \eqref{conv:ndelta}, for every $t\geq 0$, we have
\begin{equation}\label{conv:soldelta}
    \mathcal{S}_{n,\de,R}(t)\xrightarrow{n\to\infty} \mathcal{S}_{\de,R}(t)\quad\text{in the Kuratowski sense in  $\{z>-t/\e\}$, and in $L^1_{loc}(\R^{N+1})$.}
\end{equation}

\subsection{Translative solitons are super-flows}\label{subse:transl}
Let us introduce the polar distance function
\begin{equation}\label{de:se}
    d^\e\big(x,z,t \big):=\mathrm{dist}^{\Phi^o}((x,z),\mathcal{S}_\e(t)),
\end{equation}
and for any $s\in \R$ we denote by $\mathcal{H}_s$ the ``horizontal'' hyperplane \(\{(x,z)\in \R^{N+1}:\, z=s\}\). Define
\begin{equation}\label{ue}
    \mathcal{U}_\e=\Big\{(x,z,t)\in \R^N\times \R\times (0,+\infty):\, z>-\frac{t}{\e},  0<d^\e(x,z,t)<\mathrm{dist}^{\Phi^o}\left((x,z),\mathcal{H}_{-t/\e}\right)\Big\}.
\end{equation}

\begin{remark}\label{rm:ue}
(a) By \eqref{Phi}, the polar norm is given by
\[
\Phi^o(\xi,\zeta)
=\sqrt{(\phi^o(\xi))^2+\zeta^2}.
\]
In particular,
\[
\Phi^o(e_{_{N+1}})=1,
\qquad
\mathrm{dist}^{\Phi^o}((x,z),\mathcal H_s)=|z-s|.
\]
Moreover, by the translation property \eqref{solitoni} we have
\[
d^\e(X,t)
=\mathrm{dist}^{\Phi^o}
\left(X+\frac{t}{\e}e_{_{N+1}},\overline E_\e\right).
\]
Hence, by the Lipschitz continuity of the distance function,
\begin{equation}\label{lipintime}
\big|d^\e\big((x,z),t+\tau\big)-d^\e\big((x,z),t\big)\big|
\leq\frac{|\tau|}{\e}
\end{equation}
whenever $t,t+\tau\geq0$.

(b) Given $\eta>0$, set
\[
\e_0(\eta):=\frac{\eta^2}{16}
\]
and define
\begin{equation}\label{ueeta}
\mathcal U_{\e,\eta}
:=
\left\{(x,z,t)\in
\R^N\times\left(-\frac{\eta}{2\e},+\infty\right)
\times(\eta,+\infty):
0<d^\e(x,z,t)<\frac1{\sqrt\e}\right\}.
\end{equation}
If $0<\e\leq\e_0(\eta)$, $(x,z,t)\in\mathcal U_{\e,\eta}$ and $\tau\geq0$, then \eqref{lipintime} yields
\[
\begin{split}
d^\e(x,z,t+\tau)
&\leq d^\e(x,z,t)+\frac{\tau}{\e}<\frac1{\sqrt\e}+\frac{\tau}{\e}
\leq\frac{\eta}{4\e}+\frac{\tau}{\e}\\
&<z+\frac{t+\tau}{\e}
=\mathrm{dist}^{\Phi^o}
\left((x,z),\mathcal H_{-(t+\tau)/\e}\right).
\end{split}
\]
Therefore, writing $e_t=(0,0,1)\in\R^{N+2}$, we obtain
\[
(\mathcal U_{\e,\eta}+\tau e_t)\cap\{d^\e>0\}
\subset\mathcal U_\e
\qquad\text{for every }\tau\geq0.
\]
In particular, $\mathcal U_{\e,\eta}\subset\mathcal U_\e$.
\end{remark}

We now prove that the solitons \eqref{solitoni} are superflows of \eqref{eq:natural-mobility-flow} with anisotropy $\Phi$ in place of $\phi$ in the open sets $\mathcal{U}_\e$. This
is made precise by the next proposition.

\begin{proposition}[Translative solitons are super-flows]\label{prop:superflow}
    For any $0<\e<1$   the function $d^\e$  in \eqref{de:se} satisfies
    \begin{equation}\label{step5}
        \partial_t d^\e\geq \Div \mathcal{z}_\e \quad\text{in $\mathcal{U}_\e$}
    \end{equation}
    in the sense of distributions,  for a suitable vectorfield $\mathcal{z}_\e\in L^\infty(\mathcal{U}_\e; \R^{N+1})$ such that
    \begin{equation}\label{step5bis}
    \mathcal{z}_\e\in \pa \Phi(\nabla d^\e)\quad\text{a.e. in $\mathcal{U}_\e $}.
    \end{equation}
 Moreover,     for every $\sigma>0$
    \begin{equation}\label{step5cor2}
         \Div \mathcal{z}_\e \leq \frac{N}{\sigma}\quad\text{in the sense of distributions in  $\mathcal{U}_\e\cap\{d^\e>\sigma\}$.}
    \end{equation}
\end{proposition}

\begin{proof}
     We fix $0<\e<1$ and  we  define
\begin{equation}\label{dndr}
   d_{n,\delta,R}(x,z,t)=\mathrm{dist}^{\Phi^o_n}((x,z),\mathcal{S}_{n,\delta,R}(t)),
\end{equation}
where $\Phi_n, \mathcal{S}_{n,\delta,R}(t)$ are defined in Section \ref{sec:asol}. We also  let $\mathrm{sd}_{n,\delta,R}$ be the corresponding signed distance function.

Let us introduce the sets
\begin{multline}\label{setsUnd}
   \mathcal{U}_{n,\delta,R}=\big\{(x,z,t)\in \R^N\times \R\times (0,+\infty):\, z>-\frac{t}{\e}, 0<d_{n,\delta,R}\big(x,z,t\big)< \mathrm{dist}^{\Phi^o_n}\left((x,z),\mathcal{H}_{-t/\e}\right)\big\}.
\end{multline}

\noindent \textit{Step 1.} We claim that
\begin{equation}\label{step1}
    \partial_t d_{n,\delta,R}\geq \mathrm{div}\big(\nabla \Phi_n(\nabla d_{n,\delta,R}) \big)-\delta-\delta\,e^{2R/\e}\|sg_R\|_\infty\quad\text{in the viscosity sense in $\mathcal{U}_{n,\delta,R}$,}
\end{equation}
for $n$ large enough.

To this end, let $(x_0,z_0,t_0)\in \mathcal{U}_{n,\delta,R}$, and let $\varphi\in C^2$ be a viscosity test function such that $d_{n,\delta,R}-\varphi$ has a local minimum at
$(x_0,z_0,t_0)$.

Consider $(x'_0,z'_0)\in \mathcal{S}_{n,\delta,R}(t_0)$ the point of $\Phi^o_n$-minimal distance from $(x_0,z_0)$ and note that by the very definition of $\mathcal{U}_{n,\delta,R}$
we have
\begin{equation}\label{lauso1}
    z'_0>-\frac{t_0}{\e}.
\end{equation}Let also
\begin{equation*}
    \bar \varphi(x,z,t)=\varphi(x+x_0-x_0',z+z_0-z_0',t),
\end{equation*}
and we note that $\mathrm{s}d_{n,\delta,R}-\bar \varphi$ has a local minimum at $(x'_0,z'_0,t_0)$. Indeed for $(x,z,t)$ sufficiently close to $(x'_0,z'_0,t_0)$, we have
\begin{equation}\label{locmin}
    \begin{split}
        (\mathrm{s}d_{n,\delta,R}-\bar \varphi)\big((x,z),t\big)\geq &\,(d_{n,\delta,R}- \varphi)\big(x+x_0-x_0',z+z_0-z_0',t\big) -\Phi^o_n(x_0-x'_0,z_0-z'_0) \\
       \geq & (d_{n,\delta,R}- \varphi)\big(x_0,z_0,t_0\big)-\Phi^o_n(x_0-x'_0,z_0-z'_0) \\
       =& -\varphi(x_0,z_0,t_0)=(\mathrm{s}d_{n,\delta,R}- \bar\varphi)(x'_0,z'_0,t_0),
    \end{split}
\end{equation}
where in the first inequality we used the definition of $\bar \varphi$, the 1-Lipschitz continuity of $\mathrm{s}d_{n,\delta,R}$ with respect to $\Phi_n^o$, and that
$\mathrm{s}d_{n,\delta,R}=d_{n,\delta,R}$ at the point $(x+x_0-x_0',z+z_0-z_0',t)$; in the second inequality we used the local minimality of $d_{n,\delta,R}-\varphi$ at
$(x_0,z_0,t_0)$, in the third one that $d_{n,\delta,R}(x_0,z_0,t_0)=\Phi_n^o(x_0-x_0',z_0-z_0')$.

Note now that by \eqref{lauso1} we may apply Lemma \ref{lm:contact} to the point $(x'_0,z'_0)\in \partial^*\mathcal{S}_{n,\delta,R}(t_0)=\partial^*E_{n,\delta,R}-\frac{t_0}{\e}
e_{_{N+1}}$; in turn, the regularity \eqref{regolarita} yields that in a neighbourhood of $(x_0',z_0')$, the set $\partial\mathcal{S}_{n,\delta,R}(t_0)$ is of class $C^{2,\beta}$.

Therefore, by \cite[Theorem~4.16]{CrastaMalusa07}, we may find a ball $\B_r(x_0',z_0')$ and $\gamma\in (0,1)$ such that the signed distance function $\mathrm{s}d_{n,\delta,R}\in
C^2\big(\B_r(x_0',z_0')\times (t_0-\gamma,t_0+\gamma)\big)$.

\begin{center}
    \includegraphics[width=9cm]{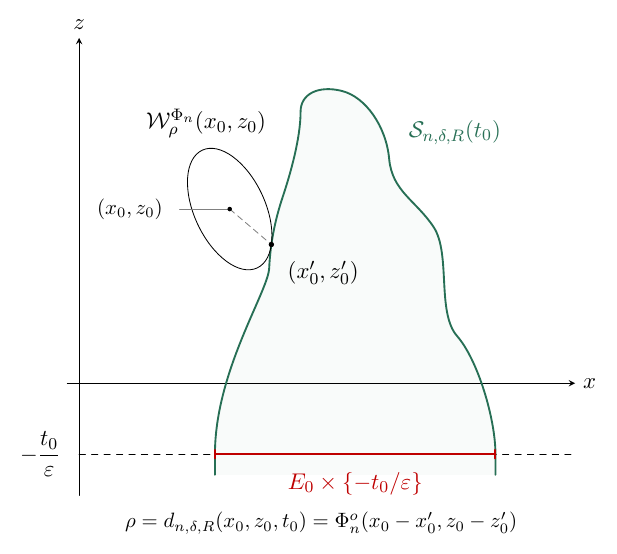}
\end{center}

Taking into account \eqref{approxtubosolitoni}, \eqref{dndr}, \eqref{eq:Phin}, inequality \eqref{ineq:AMC},  and the local minimality property \eqref{locmin}, we have
\begin{equation*}
    \begin{split}
        \partial_t\varphi(x_0, z_0, t_0)&=\partial_t\bar \varphi(x'_0, z'_0, t_0)=\partial_t\mathrm{s}d_{n,\delta,R}(x'_0, z'_0, t_0)\\
        &=\frac{1}{\e} \frac{\nu\cdot e_{_{N+1}}}{\Phi_{n}(\nu)}(x'_0,z'_0) \geq H^{\Phi_n}_{\partial^*\mathcal{S}_{n,\delta,R}(t_0)}-\delta-\delta\,e^{2R/\e}\|\mathrm{s}g_R\|_\infty\\
        &=\mathrm{div}\big(\nabla \Phi_n(\nabla \mathrm{s}d_{n,\delta,R}) \big)(x'_0, z'_0, t_0)-\delta-\delta\,e^{2R/\e}\|sg_R\|_\infty\\
        &\geq\mathrm{div}\big(\nabla \Phi_n(\nabla \bar \vphi) \big)(x'_0, z'_0, t_0)-\delta-\delta\,e^{2R/\e}\|sg_R\|_\infty\\
        &=\mathrm{div}\big(\nabla \Phi_n(\nabla \vphi) \big)(x_0, z_0, t_0)-\delta-\delta\,e^{2R/\e}\|\mathrm{s}g_R\|_\infty.
    \end{split}
\end{equation*}
Observe that, in the third identity we used that, at smooth contact points on $\partial \mathcal{S}_{n,\delta,R}(t_0)$, we have the transport identity
\[
\partial_t\mathrm{sd}_{n,\delta,R}=(\nu\cdot e_{_{N+1}})/(\e\Phi_n(\nu)).
\]
This follows from the translative property $\mathrm{sd}_{n,\delta,R}\big((x,z),t\big)=\mathrm{sd}_{n,\delta,R}\big((x,z)+t e_{_{N+1}}/\e,0\big)$, and that, at the smooth contact
point, $\nabla\mathrm{sd}_{n,\delta,R}=\nu/\Phi_n(\nu)$.

This establishes claim \eqref{step1}.
\\

\noindent \textit{Step 2.} We claim that the function
\begin{equation*}
    d_{n,\delta,R}\,\text{ is locally semiconcave in $\{d_{n,\delta,R}>\sigma\}$}
\end{equation*}
in the space-time variables, for every $\sigma>0$.

Fix $(x,z,t)\in\{d_{n,\delta,R}>\sigma\}$ and take $h\in\R^{N+1}$ and $\tau\in\R$ sufficiently small that both $\big((x,z)+h,t+\tau\big)$ and $\big((x,z)-h,t-\tau\big)$ belong to $\{d_{n,\delta,R}>\sigma\}$.

Let $(\hat{x},\hat{z})\in \partial \mathcal{S}_{n,\delta,R}(t)$ be a projection point of $(x,z)$, i.e., $\Phi_n^o\big((x-\hat x,z-\hat z) \big)=d_{n,\delta,R}(x,z,t)$. By
\eqref{approxtubosolitoni}, we have that $(\hat{x},\hat{z})\pm\frac{\tau}{\e}\,e_{_{N+1}}\in \partial\mathcal{S}_{n,\delta,R}(t\mp\tau)$. Therefore, letting $\bar h=h+\frac\tau \e
e_{_{N+1}}$, we have
\begin{equation}
    \begin{split}
       & d_{n,\delta,R}((x,z)+h,t+\tau)-2d_{n,\delta,R}(x,z,t)+d_{n,\delta,R}((x,z)-h,t-\tau) \\
       &\qquad \leq \Phi^o_n\big((x-\hat x,z-\hat z)+h+\tfrac{\tau}{\e}e_{_{N+1}}\big)-2\,\Phi^o_n\big((x-\hat x,z-\hat z)\big)+\Phi^o_n\big((x-\hat x,z-\hat z)-h-\tfrac{\tau}{\e}e_{_{N+1}}\big) \\
       &\qquad=\int_{-1}^1(1-|s|)(\nabla^2\Phi_n^o((x-\hat x,z-\hat z)+s\bar h)\bar h) \cdot \bar h \, ds \\
&\qquad\le c_n|\bar h|^2 \int_{-1}^1 \frac{1-|s|}{\Phi_n^o((x-\hat x,z-\hat z)+s\bar h)} \, ds \le \frac{C_n}{\sigma}\big(|h|^2+\e^{-2}|\tau|^2\big),
    \end{split}
\end{equation}
with $C_n$ depending only on $\|\Phi^o_n\nabla^2\Phi^o_n\|_\infty$. This proves the claim of Step 2. \\

\noindent \textit{Step 3.} We claim that
\begin{equation}\label{step3}
    \partial_t d_{n,\delta,R}\geq \mathrm{div}\big(\nabla \Phi_n(\nabla d_{n,\delta,R}) \big)-\delta-\delta\,e^{2R/\e}\|\mathrm{s}g_{R}\|_\infty\quad\text{distributionally  in $\mathcal{U}_{n,\delta,R}$,}
\end{equation}
for $n$ large enough. Using Aleksandrov's theorem, for a.e. $(x,z,t)\in \mathcal{U}_{n,\delta,R}$  and $h\in \R^{N+1}$, $\tau\in \R$ small enough, we have a second-order parabolic
jet:
\begin{align}
d_{n,\delta,R}((x,z)&+h,t+\tau)\nonumber \\
=\ & d_{n,\delta,R}(x,z,t)+\partial_t d_{n,\delta,R}(x,z,t)\tau +  \nabla d_{n,\delta,R}(x,z,t)\cdot h + \frac{1}{2} (\nabla^2 d_{n,\delta,R}(x,z,t)h)\cdot h \nonumber
\\ &+ \frac{1}{2}\partial^2_{tt} d_{n,\delta,R}(x,z,t)\tau^2 + \tau \partial_t \nabla d_{n,\delta,R}(x,z,t)\cdot h + o(\tau^2+|h|^2)
\label{eq:Aleksandrov}
\\\ =\ & d_{n,\delta,R}(x,z,t)+\partial_t d_{n,\delta,R}(x,z,t)\tau +  \nabla d_{n,\delta,R}(x,z,t)\cdot h + \frac{1}{2} (\nabla^2 d_{n,\delta,R}(x,z,t)h)\cdot h \nonumber \\ & + o(|\tau|+|h|^2). \nonumber
\end{align}

Being $d_{n,\delta,R}$ a viscosity supersolution of \eqref{step1}, it follows  that
\begin{equation}\label{eq:viscodeps}
    \partial_t d_{n,\delta,R} \geq  \nabla^2 \Phi_n(\nabla d_{n,\delta,R}):\nabla^2 d_{n,\delta,R}-\delta-\delta\,e^{2R/\e}\|\mathrm{s}g_R\|_\infty\quad\text{a.e. in $\mathcal{U}_{n,\delta,R}$.}
\end{equation}

We now show, via a regularization procedure, that \eqref{eq:viscodeps} also holds distributionally. Choose a smooth, convex, nondecreasing function $\chi:[0,+\infty)\to\R$ such
that $\chi$ is constant on $[0,1/4]$ and $\chi(s)=s$ for $s\geq1/2$, and set
\[
\widetilde\Phi_n:=\chi\circ\Phi_n.
\]
Then $\widetilde\Phi_n$ is smooth and convex on the whole $\R^{N+1}$ and coincides with $\Phi_n$ on $\{\Phi_n>1/2\}$.

Let $\rho_\gamma$ be a space-time mollifier and set
\[
d_\gamma:=\rho_\gamma*d_{n,\delta,R}.
\]
All the following arguments are local in $\mathcal U_{n,\delta,R}$, with $\gamma$ sufficiently small.

By semiconcavity, the spatial distributional Hessian $D^2d_{n,\delta,R}$ is a locally finite matrix-valued Radon measure whose singular part is negative semidefinite. Moreover, its
absolutely continuous part is given  by the Alexandrov Hessian $\nabla^2 d_{n,\de,R}$. Thus,
\[
D^2d_{n,\delta,R}
=
\nabla^2d_{n,\delta,R}\,dxdzdt
+(D^2d_{n,\delta,R})^s
\leq
\nabla^2d_{n,\delta,R}\,dxdzdt.
\]
Using the properties of convolution and the convexity of $\widetilde\Phi_n$, we infer that
\[
\begin{split}
\Div\bigl(\nabla\widetilde\Phi_n(\nabla d_\gamma)\bigr)
&=
\nabla^2\widetilde\Phi_n(\nabla d_\gamma):
\nabla^2d_\gamma=
\nabla^2\widetilde\Phi_n(\nabla d_\gamma):
(\rho_\gamma*D^2d_{n,\delta,R})\\
&\leq
\nabla^2\widetilde\Phi_n(\nabla d_\gamma):
(\rho_\gamma*\nabla^2d_{n,\delta,R}).
\end{split}
\]

As $\gamma\to0$, we have $\nabla d_\gamma\to\nabla d_{n,\delta,R}$ almost everywhere, with locally uniformly bounded gradients, and
\[
\rho_\gamma*\nabla^2d_{n,\delta,R}
\to\nabla^2d_{n,\delta,R}
\quad\text{in }L^1_{\mathrm{loc}}.
\]
The smoothness of $\widetilde\Phi_n$ therefore implies that the right-hand side of the preceding inequality converges in $L^1_{\mathrm{loc}}$, while
\[
\nabla\widetilde\Phi_n(\nabla d_\gamma)
\to
\nabla\widetilde\Phi_n(\nabla d_{n,\delta,R})
\quad\text{in }L^p_{\mathrm{loc}}
\quad\text{for every }p<\infty.
\]
Passing to the limit yields
\[
\Div\bigl(\nabla\widetilde\Phi_n
(\nabla d_{n,\delta,R})\bigr)
\leq
\nabla^2\widetilde\Phi_n(\nabla d_{n,\delta,R}):
\nabla^2d_{n,\delta,R}
\]
in the sense of distributions in $\mathcal U_{n,\delta,R}$.

Finally, since $\Phi_n(\nabla d_{n,\delta,R})=1$ almost everywhere, $\widetilde\Phi_n$ and $\Phi_n$, together with their first and second derivatives, coincide at $\nabla
d_{n,\delta,R}$. Hence,
\begin{equation}\label{stimadiv}
\Div\bigl(\nabla\Phi_n(\nabla d_{n,\delta,R})\bigr)
\leq
\nabla^2\Phi_n(\nabla d_{n,\delta,R}):
\nabla^2d_{n,\delta,R}
\quad\text{distributionally in }\mathcal U_{n,\delta,R}.
\end{equation}
Combining this inequality with \eqref{eq:viscodeps} establishes \eqref{step3}.
 \\

\noindent \textit{Step 4.} We claim that for every $\sigma>0$
\begin{equation}\label{eq:nuovostep4}
\Div \big(\nabla \Phi_n(\nabla d_{n,\delta,R})\big) \le \frac{N}{\sigma}\quad\text{in the sense of distributions in $\mathcal{U}_{n,\delta,R}\cap\{d_{n,\delta,R}>\sigma\}$.}
\end{equation}
To this aim, let $(x,z,t)\in \mathcal{U}_{n,\delta,R}\cap\{d_{n,\delta,R}>\sigma\}$ be such that the spatial Hessian $\nabla^2 d_{n,\delta,R}(x,z,t)$ exists in the sense of
Alexandrov. Denote by $(x',z')\in \mathcal{S}_{n,\delta,R}(t)$ the point of $\Phi^o_n$-minimal distance  from $(x,z,t)$. Then, taking into account \eqref{eq:Aleksandrov} we have
for  $h\in \R^{N+1}$ sufficiently small
\begin{align*}
\Phi_n^o((x-x', z-z')+h)&\geq	d_{n,\delta,R}((x,z)+h,t)\nonumber \\
	=\ & d_{n,\delta,R}(x,z,t) +  \nabla d_{n,\delta,R}(x,z,t)\cdot h + \frac{1}{2} \nabla^2 d_{n,\delta,R}(x,z,t)h\cdot h+ o(|h|^2),
\end{align*}
with equality holding at $h=0$. Thus,
\[
\nabla \Phi_n^o(x-x', z-z')=\nabla d_{n,\delta,R}(x,z,t)\quad\text{and}\quad \nabla^2 \Phi_n^o(x-x', z-z')\geq \nabla^2 d_{n,\delta,R}(x,z,t).
\]
Thus,
\[
\begin{split}
\nabla^2\Phi_n(\nabla &d_{n,\delta,R}(x,z,t)):\nabla^2d_{n,\delta,R}(x,z,t)= \nabla^2\Phi_n(\nabla \Phi_n^o(x-x', z-z')):\nabla^2d_{n,\delta,R}(x,z,t)\\
&\leq \nabla^2\Phi_n(\nabla \Phi_n^o(x-x', z-z')):\nabla^2 \Phi_n^o(x-x', z-z')=\frac{N}{\Phi_n^o(x-x', z-z')}\leq \frac{N}{\sigma},
\end{split}
\]
where in the last equality we used \eqref{usefulid}. Since the above chain of inequalities holds for almost every  $(x,z,t)\in
\mathcal{U}_{n,\delta,R}\cap\{d_{n,\delta,R}>\sigma\}$, claim \eqref{eq:nuovostep4} follows recalling \eqref{stimadiv}.
\\

\noindent \textit{Step 5.}  Let us set
\[
\mathcal{z}_{n,\delta,R}:=
\begin{cases}
	\nabla \Phi_n(\nabla d_{n,\delta,R}) & \text{if $\exists\, \nabla d_{n,\delta,R}\neq 0$,}\\
	0 & \text{otherwise.}
\end{cases}
\]
 Since $\Phi_n^o(\mathcal{z}_{n,\delta,R})\leq1$  a.e. by \eqref{Eikonal} and recalling \eqref{convpolari}, we may assume that up to a non-relabeled subsequence
\begin{equation}\label{convznd}
\mathcal{z}_{n,\delta,R}\to \mathcal{z}_{\delta,R} \quad\text{weakly* in $L^\infty(\R^{N+1}\times (0,+\infty))$, with $\Phi^o(\mathcal{z}_{\delta, R})\leq 1$.}
\end{equation}
We now claim that
 \begin{equation}\label{step4.0}
 	\mathcal{z}_{\delta,R}\in \partial \Phi(\nabla d_{\delta,R})\quad\text{a.e in  $\mathcal{U}_{\delta,R}$,}
\end{equation}
 \begin{equation}\label{step4.01}
	\Div \mathcal{z}_{\delta,R}\leq \frac{N}\sigma \quad\text{distributionally  in $\mathcal{U}_{\delta,R}\cap\{d_{\delta,R}>\sigma\}$ for any $\sigma>0$,}
\end{equation}
and
 \begin{equation}\label{step4}
    \partial_t d_{\delta,R}\geq	\Div \mathcal{z}_{\delta,R}-\delta-\delta\,e^{2R/\e}\|sg_R\|_\infty\quad\text{distributionally  in $\mathcal{U}_{\delta,R}$,}
\end{equation}
where we defined the distance function
\begin{equation*}
    d_{\delta,R}(x,z,t)\coloneqq \mathrm{dist}^{\Phi^o}((x,z),\mathcal{S}_{\delta,R}(t)),
\end{equation*}
and the sets $\mathcal{U}_{\delta,R}$ as in \eqref{setsUnd}, with $d_{\delta,R}$ and $\Phi^o$ in place of $d_{n,\delta,R}$ and $\Phi_n^o$, respectively.

First, thanks to \eqref{conv:soldelta},  \eqref{convdist}, and the fact that $\mathcal{S}_{n,\delta,R}(t)$ and $\mathcal{S}_{\delta,R}(t)$ are translating solitons, we have
\begin{equation}\label{dist:convunif}
d_{n,\delta,R}\longrightarrow d_{\delta,R}
\quad\text{locally uniformly in $\mathcal U_{\delta,R}$.}
\end{equation}

Now let $\varphi\in C^\infty_c(\mathcal{U}_{\delta,R})$ be a test function such that $\varphi\geq 0$. By \eqref{dist:convunif} and the definition of the sets \eqref{setsUnd}, for
$n$ large enough, it holds $\varphi\in C^\infty_c(\mathcal{U}_{n,\delta,R})$; and thus by \eqref{step3}, we get
\begin{equation}\label{testn1}
    -\iiint\partial_t\varphi\,d_{n,\delta,R}\,dxdzdt\geq -\iiint \nabla \Phi_n(\nabla d_{n,\delta,R})\cdot \nabla \varphi +\big(\delta+\delta\,e^{2R/\e}\|\mathrm{s}g_R\|_\infty\big)\,\varphi\, dxdzdt
\end{equation}
Owing to \eqref{convznd} and \eqref{dist:convunif}, we can pass to the limit in the above inequality and establish \eqref{step4}. Claim \eqref{step4.01} follows similarly taking
into account \eqref{eq:nuovostep4}.

In order to prove \eqref{step4.0}, we exploit the bound \eqref{eq:nuovostep4} as done in the final part of the proof of \cite[Theorem~4.5]{CMP17}. We reproduce the argument for
completeness. Consider as before  a nonnegative test function $\varphi\in C^\infty_c(\mathcal{U}_{\delta,R})$ Then,  by \eqref{eq:Phin}, \eqref{dist:convunif} and by lower
semicontinuity
\begin{multline*}
\iiint \Phi(\nabla d_{\delta, R})\varphi\, dxdzdt\le \liminf_{n}\iiint \Phi_n(\nabla d_{n, \delta, R})\varphi\, dxdzdt = \liminf_n \iiint (\mathcal{z}_{n,\delta, R}\cdot\nabla d_{n, \delta, R})\varphi\, dxdzdt.
\end{multline*}
Writing
\begin{multline*}
 \iiint (\mathcal{z}_{n,\delta, R}\cdot\nabla d_{n, \delta, R})\varphi\, dxdzdt = \iiint (\mathcal{z}_{n,\delta, R}\cdot\nabla d_{\delta, R})\varphi\, dxdzdt+ \iiint \mathcal{z}_{n,\delta, R}\cdot\nabla (d_{n, \delta, R}-d_{\delta, R})\varphi\, dxdzdt,
\end{multline*}
where
\[
 \iiint (\mathcal{z}_{n,\delta, R}\cdot\nabla d_{\delta, R})\varphi\, dxdzdt \stackrel{n\to\infty}{\longrightarrow} \iiint (\mathcal{z}_{\delta, R}\cdot\nabla d_{\delta, R})\varphi\, dxdzdt,
\]
we obtain
\begin{equation}\label{eq:verygood}
\iiint\Phi(\nabla d_{\delta, R})\varphi\, dxdzdt\le\iiint(\mathcal{z}_{\delta, R}\cdot\nabla d_{\delta, R})\varphi\, dxdzdt,
\end{equation}
if we show that
\begin{equation}\label{eq:good}
	 \iiint \mathcal{z}_{n,\delta, R}\cdot\nabla (d_{n, \delta, R}-d_{\delta, R})\varphi\, dxdzdt\to 0 \quad\text{as $n\to\infty$.}
\end{equation}
For each $t$, set
\[
\begin{split}
m_n(t)&:=\min_{(x,z)\in\spt \varphi(\cdot,t)} \bigl(d_{n,\delta, R}(x,z,t)-d_{\delta, R}(x,z,t)\bigr),\\
M_n(t)&:=\max_{(x,z)\in\spt \varphi(\cdot,t)} \bigl(d_{n,\delta, R}(x,z,t)-d_{\delta, R}(x,z,t)\bigr).
\end{split}
\]
These quantities are uniformly bounded and converge to $0$ at all $t$ (in fact uniformly) thanks to \eqref{dist:convunif}. Then, by integration by parts,  we have for $n$ large
enough
\begin{equation}\label{yeswecan}
	\begin{split}
 \iiint \mathcal{z}_{n,\delta, R}&\cdot\nabla (d_{n, \delta, R}-d_{\delta, R})\varphi\, dxdzdt =\iiint \mathcal{z}_{n,\delta, R}\cdot\nabla (d_{n, \delta, R}-d_{\delta, R}-m_n)\varphi\, dxdzdt \\
	&= -\iiint (d_{n, \delta, R}-d_{\delta, R}-m_n)(\mathcal{z}_{n,\delta, R}\cdot \nabla \varphi)\, dxdzdt
    \\
    &\qquad-\iiint (d_{n, \delta, R}-d_{\delta, R}-m_n)\varphi\,d(\mathrm{div} \, \mathcal{z}_{n,\delta, R})
    \\
	&\geq -\iiint (d_{n, \delta, R}-d_{\delta, R}-m_n)\left(\mathcal{z}_{n,\delta, R}\cdot \nabla \varphi + \frac{N}\sigma\varphi\right)\, dxdzdt\stackrel{n\to\infty}{\longrightarrow}0,
\end{split}
\end{equation}
where  $\mathrm{div} \mathcal{z}_{n,\delta, R}$ is the  distributional divergence of $\mathcal{z}_{n,\delta, R}$, in the last inequality we have used that, for $n$ large enough,
the support of $\vphi$ is contained in $\{d_{n,\delta, R}>\sigma\}$ for some $\sigma>0$, and the bound \eqref{eq:nuovostep4}. Thus we have shown that
\[
\liminf_n \iiint \mathcal{z}_{n,\delta, R}\cdot\nabla (d_{n, \delta, R}-d_{\delta, R})\varphi\, dxdzdt \ge 0.
\]
 Writing now
\[
\iiint \mathcal{z}_{n,\delta, R}\cdot\nabla (d_{n, \delta, R}-d_{\delta, R})\varphi\, dxdzdt =\iiint \mathcal{z}_{n,\delta, R}\cdot\nabla (d_{n, \delta, R}-d_{\delta, R}-M_n)\varphi\, dxdzdt
\]
and arguing similarly, one can show that
\[
\limsup_n \iiint \mathcal{z}_{n,\delta, R}\cdot\nabla (d_{n, \delta, R}-d_{\delta, R})\varphi\, dxdzdt \le 0
\]
so that \eqref{eq:good} and in turn, \eqref{eq:verygood} follows for any non-negative $C^\infty_c(\mathcal{U}_{\delta,R})$. Thus,
\[
\Phi(\nabla d_{\delta, R})\leq \mathcal{z}_{\delta, R}\cdot \nabla d_{\delta, R}\leq \Phi^o(\mathcal{z}_{\delta, R})\Phi(\nabla d_{\delta, R})\leq \Phi(\nabla d_{\delta, R}) \quad\text{a.e. in $\mathcal{U}_{\delta,R}$,}
\]
where we used that $\Phi^o(\mathcal{z}_{\delta, R})\leq 1$ a.e.. Hence,  necessarily $\Phi(\nabla d_{\delta, R})= \mathcal{z}_{\delta, R}\cdot \nabla d_{\delta, R}$ a.e. in
$\mathcal{U}_{\delta,R}$, which establishes also \eqref{step4.0}.
\\

\noindent\textit{Step 6.} (Conclusion) Thanks to \eqref{ugualeBR},  for every $M>0$ fixed, we have
\begin{equation}
    d_{\delta,R}\equiv d^\e\quad\text{in $\mathcal{U}_\e\cap \big(\mathcal{B}_{M}\times (0,M)\big)=\mathcal{U}_{\delta,R}\cap \big(\mathcal{B}_{M}\times (0,M)\big)$,}
\end{equation}
 provided that $R>R(M)$ is large enough and $\delta$ satisfies \eqref{lambda1}. Therefore, from \eqref{step4.0}--\eqref{step4}, $d^\e$ satisfies
\begin{equation}\label{quasiquasi}
	\begin{split}
    &\partial_t d^\e\geq \Div \mathcal{z}_{\de, R}-\delta-\delta\,e^{2R/\e}\|sg_R\|_\infty\quad\text{distributionally  in $\mathcal{U}_\e\cap \big(\mathcal{B}_{M}\times (0,M)\big)$,}\\
    & \mathcal{z}_{\delta,R}\in \partial \Phi(\nabla d^\e)\quad\text{a.e in $\mathcal{U}_\e\cap \big(\mathcal{B}_{M}\times (0,M)\big)$,}\\
    &\Div \mathcal{z}_{\delta,R}\leq \frac{N}\sigma \quad\text{distributionally  in $\mathcal{U}_\e\cap \big(\mathcal{B}_{M}\times (0,M)\big)\cap\{d^\e>\sigma\}$ for any $\sigma>0$,}
    \end{split}
\end{equation}
 For any $M>0$ and thus for any $R>R(M)$, we may find a sequence $\delta_n\searrow 0$ (possibly depending on $R$) such that
 \[
 \mathcal{z}_{\delta_n,R}\to \mathcal{z}_{R} \quad\text{weakly* in $L^\infty(\R^{N+1}\times (0,+\infty))$}
 \]
 and \eqref{quasiquasi} holds with $\delta=0$ and $\mathcal{z}_{\delta,R}$ replaced by $\mathcal{z}_{R}$. Finally, we  find $R_n\to \infty$ such that for a suitable vector-field
 $\mathcal{z}_\e$ we have
  \[
  \mathcal{z}_{R_n}\to \mathcal{z}_{\e} \quad\text{weakly* in $L^\infty(\R^{N+1}\times (0,+\infty))$.}
  \]
   Since for each $n$ \eqref{quasiquasi} holds with $\delta=0$, $\mathcal{z}_{R_n}$ in place of $\mathcal{z}_{\delta, R}$ and $M=M_n\to\infty$, we easily deduce that
   $\mathcal{z}_\e$ has all the required properties \eqref{step5}--\eqref{step5cor2}.
\end{proof}

In what follows, we aim to show that the  minimizers $E_\e$ of \eqref{minpbell} and their complement are almost cylindrical. To this aim,  given $x_0\in \R^N$, $\eta, r>0$ and
$\e\in (0,1)$, as in \eqref{abs:paraboloidi} we define  the following intrinsic paraboloids
\begin{equation}\label{paraboloidi}
    \mathcal{C}_{r,\eta, \e}(x_0):=\left\{(x,z)\in \R^{N+1}:\, 0<z\leq g_{r,\eta, \e}(x-x_0) \right\},\qquad  g_{r,\eta, \e}(x)\coloneqq\frac{\eta}{\e}\big(r^2-\phi^o(x)^2\big).
\end{equation}

Correspondingly, set
\[
    \phi_n(\xi):=\Phi_n(\xi,0), \qquad G_n:=\frac12\Phi_n(\cdot,-1)^2,
\]
and, with $\lambda:=2\eta/\e$,
\begin{equation}\label{defan}
    a^n_{\lambda,r} := \min_{\phi_n^o(y)=r} \bigl\{G_n^*(\lambda y)-G_n^*(0)\bigr\}.
\end{equation}
We consider the approximating paraboloids
\begin{equation}\label{approx-paraboloidi}
    \mathcal{C}^n_{r,\eta, \e}(x_0):=\left\{(x,z)\in \R^{N+1}:\, 0<z\leq g^n_{r,\eta, \e}(x-x_0) \right\},
\end{equation}
where
\begin{equation}\label{defgn}
    g^n_{r,\eta, \e}(x):=\frac{G_n^*(0)+a^n_{\lambda,r}-G_n^*(-\lambda x)}{\lambda}.
\end{equation}
By \eqref{abs:footprint},
\begin{equation}\label{foot:gn}
    \{g^n_{r,\eta,\e}>0\}\subset W^{\phi_n}_r.
\end{equation}
Since $G_n\to G:=\frac12\Phi(\cdot,-1)^2 =\frac12\phi^2+\frac12$ locally uniformly in $\R^N$, it follows (see for instance \cite[Theorems~7.11, 7.17 and 11.34]{RockafellarWets})
that
\[
G_n^*\to G^*= \frac12(\phi^o)^2-\frac12\quad\text{ locally uniformly in $\R^N$.}
\]
Moreover, $\phi_n\to\phi$ locally uniformly, and therefore $\phi_n^o\to\phi^o$ locally uniformly. Hence, for every fixed $\lambda>0$,
\begin{equation}\label{antolimit}
    a^n_{\lambda,r}\longrightarrow\frac{\lambda^2r^2}{2}.
\end{equation}
It follows that
\begin{equation}\label{gntog}
    g^n_{r,\eta, \e}\to g_{r,\eta, \e}\quad\text{locally uniformly in $\R^N$ as }n\to\infty.
\end{equation}
In turn, using also the fact that $\{g_{r,\eta, \e}\ge 0\}=\overline{\{g_{r,\eta, \e}> 0\}}$,
\begin{equation}\label{gntogzero}
        \{g^n_{r,\eta, \e}\ge 0\}\stackrel{\mathcal K}{\longrightarrow} \{g_{r,\eta, \e}\ge 0\}=\overline{W^\phi_r}\quad \text{in $\R^N$ as }n\to\infty.
\end{equation}

We now prove the following

\begin{proposition}[Intrinsic paraboloid inclusion]\label{prop:paraboloidi}
   Let $E_\e$ be as in \eqref{minpbell}, with $\e\in (0,\tfrac12)$. If $\overline{W^\phi_r(x_0)}\subset \mathrm{int}(E_0)$, then
    \[
    \mathcal{C}_{r,\eta, \e}(x_0)\subset E_\e
    \]
    for all $\eta\in (0, \frac1{4N}]$. Symmetrically, if $\overline{W^\phi_r(x_0)}\subset \R^N\setminus \overline{E_0}$, then
    \[
    \mathcal{C}_{r,\eta, \e}(x_0)\subset \overline{\R^{N+1}\setminus E_\e}
    \]
     for all $\eta\in (0, \frac1{4N}]$.
\end{proposition}

\begin{proof}
    By Remark~\ref{rm:complementare}, it is enough to consider the case $\overline{W^\phi_r(x_0)}\subset\mathrm{int}(E_0)$. Set $\eta_0:=\frac1{4N}$  and fix $R>0$. Since
    $\phi_n^o\to\phi^o$ locally uniformly, for $n$ sufficiently large
    \begin{equation}\label{Wn:inside}
        \overline{W^{\phi_n}_r(x_0)}\subset\mathrm{int}(E_0).
    \end{equation}

    We now claim that
    \begin{equation}\label{unif:gradpar}
        \sup_{0<\eta\leq\eta_0} \|\nabla g^n_{r,\eta,\e}\|_{L^\infty(W^{\phi_n}_r)} \leq C_{r,\e} \qquad\text{and}\qquad \sup_{0<\eta\leq\eta_0} \|\Phi_n(\nabla g^n_{r,\eta,\e},-1)\|_{L^\infty(W^{\phi_n}_r)} <M_{r,\e},
    \end{equation}
   hold for some $C_{r,\e}, M_{r,\e}>0$ independent of $n$. To prove them, we recall that, by the local uniform convergence of $\phi_n^{o}\to \phi^o$,  the sets $W_r^{\phi_n}$ are
   uniformly bounded in $n$; hence we may choose $L>0$ such that
\[
-\frac{2\eta}{\e}x\in B_L
\qquad\text{for all }n,\quad
0<\eta\leq\eta_0,\quad x\in W_r^{\phi_n}.
\]
The local uniform convergence of $G_n^*$ and its convexity \cite[Theorem 10.4]{RockafellarConvexAnalysis} yield
\[
\sup_n\|\nabla G_n^*\|_{L^\infty(B_L)}<+\infty,
\]
and since $\nabla g^n_{r,\eta,\e}(x) =\nabla G_n^*\left(-\frac{2\eta}{\e}x\right)$, we obtain the first bound in \eqref{unif:gradpar}. The second follows from the uniform estimate
$\Phi_n(v)\leq C|v|$, which is a consequence of \eqref{eq:Phin} and the norm equivalence.

   We now choose $\de>0$ so small that, in addition to \eqref{lambda1},
    \begin{equation}\label{supersmallde}
        \delta\,e^{2R/\e} \|\mathrm{s}g_R\|_{L^\infty(\R^{N+1})}\,M_{r,\e} <\frac12.
    \end{equation}
    Notice also that, for $n$ sufficiently large, $\frac1\e\Phi_n(e_{_{N+1}}) >\delta\|\mathrm{s}g_R\|_\infty$, so that the trace assumption \eqref{abs:trace-hyp} is satisfied for
    $E_{n,\de,R}$.

    By Lemma~\ref{lem:vertical} and \eqref{unif:gradpar}, for $n$ sufficiently large we have, uniformly for $0<\eta\leq\eta_0$ and $x\in W^{\phi_n}_r$,
    \[
        \Phi_n(p)\partial_\zeta\Phi_n(p)<-\frac34, \qquad p=(\nabla g^n_{r,\eta,\e}(x),-1).
    \]
    Moreover, since $\mathrm{s}g_R$ is supported in $\mathcal B_{2R}$,
    \[
        e^{g^n_{r,\eta,\e}(x)/\e} \left| \mathrm{s}g_R \bigl(x_0+x,g^n_{r,\eta,\e}(x)\bigr) \right| \leq e^{2R/\e}\|\mathrm{s}g_R\|_\infty.
    \]
    Hence, by \eqref{supersmallde} and the fact that $\e<\frac12$,
    \[
    \begin{aligned}
        \frac1\e\Phi_n(p)\partial_\zeta\Phi_n(p) &+ \delta e^{g^n_{r,\eta,\e}(x)/\e} \mathrm{s}g_R \bigl(x_0+x,g^n_{r,\eta,\e}(x)\bigr)\Phi_n(p) \\
        &<-\frac{3}{4\e}+\frac12 <-\frac1{2\e} =-\frac{2N\eta_0}{\e}.
    \end{aligned}
    \]
    Thus, applying Lemma~\ref{lm:paraboloidi} with $\Psi=\Phi_n$ and  $g=\delta\,\mathrm{s}g_R$, we obtain $\mathcal C^n_{r,\eta,\e}(x_0)\subset E_{n,\de,R}$ for every
    $\eta\in(0,\eta_0]$ and all $n$ sufficiently large. Passing to the limit as $n\to\infty$ and using \eqref{conv:ndelta} and \eqref{gntog}, we infer
    \[
        \mathcal C_{r,\eta,\e}(x_0)\subset E_{\de,R}.
    \]
    Finally, the thesis follows from \eqref{ugualeBR} by choosing $R>0$ (and $\de$ accordingly) so large that $\mathcal C_{r,\eta,\e}(x_0)\subset\mathcal B_R$.
\end{proof}

\subsection{Limiting analysis and proof of Theorem~\ref{thm:main1}}

We now analyze the behavior of $E_\e$ as $\e\to 0^+$. It is well known that the evolution of a {\em $\Phi$-Wulff shape}
\[
\mathcal{W}_R^{\Phi}(x_0, z_0):=\big\{(x,z):\, \Phi^o(x-x_0, z-z_0)<R\big\}
\]
under the flow $V=-\Phi\, H^\Phi $ is self similar; that is,   given by $\mathcal{W}_{R(t)}^{\Phi}(x_0, z_0)$, with
\[
R(t):=\sqrt{R^2-2Nt},
\]
until the extinction time $\bar t(R):=\frac{R^2}{2N}$, and the empty set afterwards (see for instance \cite{CaCha}). Denote by $\mathcal{W}_R$ the associated closed space-time
tube:
\[
\mathcal{W}_R:=\left\{(x,z, t)\in \R^{N+1}\times\big[0,\bar t(R)\big):\, (x,z)\in \overline{\mathcal{W}_{R(t)}^{\Phi}(x_0, z_0)} \right\}\cup (x_0, z_0, \bar t(R)).
\]
We now show the following comparison with the evolution of $\Phi$-Wulff shapes.
\begin{proposition}\label{prop:wulff}
    Let $\eta>0$ be fixed and let $\e_0(\eta)$ be the constant introduced in Remark~\ref{rm:ue}-(b). Assume
    \[
    \overline{\mathcal{W}_R^{\Phi}(x_0, z_0)}\cap \mathcal{S}_\e(t_0)=\varnothing\quad\text{and}\quad  \overline{\mathcal{W}_R^{\Phi}(x_0, z_0)}\subset\{(x,z):\, z>-\frac{\eta}{2\e}\},
    \]
    for some $\e \in (0, \e_0(\eta)\land \tfrac{1}{16R^2})$ and  $t_0>\eta>0$. Then,
    \[
    \mathcal{S}_\e(t_0+s)\cap \overline{\mathcal{W}_{R(s)}^{\Phi}(x_0, z_0)}=\varnothing\quad\text{for all }s\in \left(0, \tfrac{R^2}{2N}\right).
    \]
\end{proposition}
\begin{proof} Set $\Delta(s):=\dist^{\Phi^o}\Big(\mathcal{S}_\e(t_0+s), \overline{\mathcal{W}_{R(s)}^{\Phi}(x_0, z_0)})\Big)$ and note that $\Delta(\cdot)$ is continuous in $[0, \tfrac{R^2}{2N})$.  Set $\Delta:=\Delta(0)\land \frac1{2\sqrt{\e}}$.

We claim that  $\Delta(s)\geq \Delta$ for all $s\in \left(0, \tfrac{R^2}{2N}\right)$. To this aim, let $\bar \sigma:=\sup\{\sigma\geq 0: \Delta(s)\geq \Delta\text{ for all }0\leq s
\leq\sigma\}$ and assume by contradiction that $\bar \sigma<\tfrac{R^2}{2N}$. Then, by continuity $\Delta(\bar \sigma)=\Delta\leq \frac1{2\sqrt{\e}}$. The latter condition,
together with the assumptions $R<\frac{1}{4\sqrt{\e}}$, $t_0>\eta$ and  $\overline{\mathcal{W}_R^{\Phi}(x_0, z_0)}\subset\{(x,z):\, z>-\frac{\eta}{2\e}\}$, implies
$\overline{\mathcal{W}_{R(\bar \sigma)+r}^{\Phi}(x_0, z_0)}\times\{t_0+\bar \sigma\} \subset \mathcal{U}_{\e, \eta}$, with $\mathcal{U}_{\e, \eta}$ defined as in \eqref{ueeta},
provided that $r>0$ is small enough. By Remark~\ref{rm:ue}-(b), this compact set is contained in $\mathcal U_\e$; since $\mathcal U_\e$ is open, there exists $\tau_*>0$ such that
\begin{equation}\label{piusigma}
\overline{\mathcal W_{R(\bar\sigma)+r}^{\Phi}(x_0,z_0)}
\times[t_0+\bar\sigma,t_0+\bar\sigma+\tau_*]
\subset\mathcal U_\e.
\end{equation}
We may now argue as in \cite[Theorem~3.3]{CMP17}, with $\mathcal{W}_{R(\bar \sigma)}$ in place of $E$, with $F$ given by the space-time tube associate with the flow
$F(s):=\R^{N+1}\setminus \mathcal{S}_\e(t_0+\bar \sigma+s)$, and with the stripe $S$ in the first display of \cite[page 1095]{CMP17} replaced by
\[
S:=\{(x,z):\, \eta_1+R(\bar \sigma)<\Phi^o(x-x_0, z-z_0)<\eta_2+R(\bar \sigma)\},
\]
 with $0<\eta_1<\eta_2<\Delta\land r$. Notice that by \eqref{piusigma} and Proposition~\ref{prop:superflow}, $F$ is a subflow in $S\times(0,\tau_*)$. Thus, by the very same
 argument of \cite[Theorem~3.3]{CMP17}, we deduce that  there exists $s^*>0$ such that $\Delta(\bar\sigma+s)\geq\Delta$ for $s\in (0, s^*)$, contradicting the definition of
 $\bar\sigma$.
\end{proof}
Fix $\eta>0$ and assume that  $(x_0,z_0)\in \R^{N+1}\setminus \mathcal{S}_\e(t_0)$ for some $t_0>\eta$ and that
\begin{equation}\label{eq*}
    \mathcal{W}_{d^\e(x_0,z_0, t_0)}^{\Phi}(x_0, z_0)\subset\{(x,z):\, z>-\tfrac{\eta}{2\e}\}\quad\text{for some }\e \in \big(0, \e_0(\eta)\land \tfrac{1}{16(d^\e(x_0,z_0, t_0))^2}\big)
\end{equation}
Then, for any $0<R<d^\e(x_0,z_0, t_0)$, by the previous proposition we deduce that
\[
d^\e(x_0,z_0, t_0+s)\geq\sqrt{R^2-2Ns}\quad\text{for all }s\in \left(0, \tfrac{R^2}{2N}\right).
\]
Sending $R\nearrow d^\e(x_0,z_0, t_0)$, we get
\begin{equation}\label{eq:helly}
    d^\e(x_0,z_0, t_0+s)\geq\sqrt{(d^\e(x_0,z_0, t_0))^2-2Ns}\quad\text{for all }s\in \left(0, \tfrac{(d^\e(x_0,z_0, t_0))^2}{2N}\right),
\end{equation}
provided that \eqref{eq*} is satisfied, with $t_0>\eta$. We are  now ready to state  the following result, showing that any (subsequential) Kuratowski limit of the space-time
soliton tubes $(\mathcal{S_{\e}})$ is a superflow for $V=-\Phi\, H^\Phi$.
\begin{proposition}\label{prop:helly}
    Assume that
    \[
    \mathcal{S_{\e_j}} \stackrel{\mathcal K}{\longrightarrow}\mathcal{S}\quad\text{in the Kuratowski sense in }\R^{N+2}
    \]
    along a vanishing subsequence $\{\e_j\}$. Then, setting for $t\geq0 $
    \[
    \mathcal{S}(t):=\{(x,z):\, (x,z,t)\in \mathcal{S}\},
    \]
    the following properties hold true:
    \begin{itemize}
        \item [(a)] $\mathrm{int }(E_0)\times \R\subset \mathcal{S}(0)\subset \overline{E_0}\times \R$;
        \item [(b)] $\mathcal{S}(\cdot)$ is Kuratowski left-continuous;
        \item [(c)] $\mathcal{S}(t)=\varnothing$ implies $\mathcal{S}(s)=\varnothing$ for all $s>t$ so that it is well defined the extinction time
        \[
        T^*:=\inf\{t>0:\, \mathcal{S}(t)=\varnothing\};
        \]
        \item [(d)]  $\{d^{\e_j}\}$ is locally uniformly bounded in $\R^{N+1}\times (0,T^*)$ and there exists a countable set $\mathcal{N}\subset (0,T^*)$ such that for all $t\in (0,T^*)\setminus \mathcal{N}$
        \[
        d^{\e_j}(\cdot, \cdot, t)\to \dist^{\Phi^o}(\cdot, \cdot, \mathcal{S}(t) )\quad\text{locally uniformly},
        \]
        that is, \(\mathcal{S_{\e_j}}(t) \stackrel{\mathcal K}{\longrightarrow}\mathcal{S}(t)\) in the Kuratowski sense in $\R^{N+1}$;
        \item[(e)] for all $t>T^*$
        \[
        \mathcal{S_{\e_j}}(t) \stackrel{\mathcal K}{\longrightarrow}\mathcal{S}(t)=\varnothing\quad\text{in the Kuratowski sense in }\R^{N+1};
        \]
        \item[(f)]
            $\mathcal{S}$ is a superflow for $V=-\Phi\, H^\Phi$ in the sense of Definition~\ref{Defsol}.
                      \end{itemize}
\end{proposition}
\begin{proof}
We first prove (c). Assume that $\mathcal S(t_0)=\varnothing$. Fix $\eta\in(0,t_0)$. By the space-time Kuratowski convergence, for every fixed $R>0$ and $j$ sufficiently large we
have $\overline{\mathcal W_R^\Phi}\cap\mathcal S_{\e_j}(t_0)=\varnothing$, and the hypotheses of Proposition~\ref{prop:wulff} are satisfied. Hence
\[
\mathcal S_{\e_j}(t_0+s)\cap\overline{\mathcal W^\Phi_{\sqrt{R^2-2Ns}}}=\varnothing \qquad\text{for }0<s<\frac{R^2}{2N}.
\]
Passing to the limit, and then using the arbitrariness of $R$, we obtain $\mathcal S(t_0+s)=\varnothing$ for every $s>0$. This proves (c).

We next establish, via Wulff shape comparison, the time regularity of the limiting slices. For $t\in(0,T^*)$, set
\[
d(x,z,t):=\dist^{\Phi^o}((x,z),\mathcal S(t)).
\]
Fix $(x,z)\in\R^{N+1}$ and $0<t<s<T^*$. For every $0<R<d(x,z,t)$, the compact set
\[
\overline{\mathcal W_R^\Phi(x,z)}\times\{t\}
\]
is disjoint from $\mathcal S$. Hence, by the space-time Kuratowski convergence,
\[
\overline{\mathcal W_R^\Phi(x,z)}
\cap\mathcal S_{\e_j}(t)=\varnothing
\]
for all sufficiently large $j$. Since $R$ is fixed, the remaining assumptions of Proposition~\ref{prop:wulff} also hold for large $j$. Passing to the limit in the resulting
Wulff-shape avoidance, we obtain
\[
d(x,z,s)\geq\sqrt{R^2-2N(s-t)}
\qquad\text{if }2N(s-t)<R^2.
\]
Letting $R\nearrow d(x,z,t)$, and observing that the inequality below is immediate when its right-hand side is nonpositive, we obtain
\begin{equation}\label{eq:limit-monotonicity}
d(x,z,s)^2\geq d(x,z,t)^2-2N(s-t).
\end{equation}
Thus, for every fixed $(x,z)$, the map $t\longmapsto d(x,z,t)^2+2Nt$ is nondecreasing on $(0,T^*)$. This is similar to the local one-sided quasi-monotonicity used in the proof of
\cite[Proposition~4.4]{CMP17}. Applying this observation at the points of a countable dense subset of $\R^{N+1}$, and using the uniform $1$-Lipschitz continuity of
$d(\cdot,\cdot,t)$ with respect to $\Phi^o$, we find a countable set $\mathcal N\subset(0,T^*)$ such that
\begin{equation}\label{eq:limit-distance-cont}
d(\cdot,\cdot,t)\longrightarrow d(\cdot,\cdot,t_0)
\qquad\text{locally uniformly in }\R^{N+1}
\end{equation}
as $t\to t_0$, for every $t_0\in(0,T^*)\setminus\mathcal N$. Moreover, the closedness of the space-time tube $\mathcal S$ implies $d(x,z,t_0)\leq\liminf_{t\to t_0}d(x,z,t)$. On the
other hand, by \eqref{eq:limit-monotonicity} we have $\limsup_{t\nearrow t_0}d(x,z,t)\leq d(x,z,t_0)$. Consequently,
\[
d(\cdot,\cdot,t)\longrightarrow d(\cdot,\cdot,t_0)
\qquad\text{locally uniformly in }\R^{N+1}
\quad\text{as }t\nearrow t_0
\]
for every $t_0\in(0,T^*)$, proving the Kuratowski left-continuity before extinction.

If $T^*<+\infty$ and $\mathcal S(T^*)\neq\varnothing$, the same comparison argument applies with $s=T^*$ and establishes left-continuity at $T^*$ as well. If instead $\mathcal
S(T^*)=\varnothing$, the closedness of $\mathcal S$ directly yields $\mathcal S(t)\stackrel{\mathcal K}{\longrightarrow}\varnothing$ as $t\nearrow T^*$. Together with (c), this
proves (b).

We now prove (d). We first establish the local uniform boundedness of $\{d^{\e_j}\}$. Fix a compact set $K\subset\R^{N+1}$ and $0<a<b<T^*$. Choose $s\in(b,T^*)$ and
$(x_*,z_*)\in\mathcal S(s)$. By the space-time Kuratowski convergence, there exist
\[
(x_j,z_j,s_j)\in\mathcal S_{\e_j},
\qquad
(x_j,z_j,s_j)\longrightarrow(x_*,z_*,s).
\]
Choose $R>0$ sufficiently large that
\[
R^2>
\max_{(x,z)\in K}
\Phi^o(x-x_*,z-z_*)^2+2N(s-a).
\]
Then, for all sufficiently large $j$,
\[
\Phi^o(x-x_j,z-z_j)^2+2N(s_j-t)<R^2
\qquad
\text{for every }(x,z,t)\in K\times[a,b].
\]

We claim that $d^{\e_j}\leq R$ on $K\times[a,b]$ for all sufficiently large $j$. Otherwise, at some $(x,z,t)\in K\times[a,b]$ we would have
\[
\overline{\mathcal W_R^\Phi(x,z)}
\cap\mathcal S_{\e_j}(t)=\varnothing.
\]
Fixing $\eta\in(0,a)$, the remaining assumptions of Proposition~\ref{prop:wulff} hold uniformly for $(x,z,t)\in K\times[a,b]$ when $j$ is sufficiently large. By that proposition we
would therefore have
\[
\mathcal S_{\e_j}(s_j)\cap
\overline{\mathcal W^\Phi_{\sqrt{R^2-2N(s_j-t)}}(x,z)}
=\varnothing.
\]
This contradicts $(x_j,z_j)\in\mathcal S_{\e_j}(s_j)$ and the choice of $R$. Thus $\{d^{\e_j}\}$ is locally uniformly bounded in $\R^{N+1}\times(0,T^*)$.

Next, we prove the Kuratowski convergence of the slices. So fix $t\in(0,T^*)\setminus\mathcal N$. By the space-time Kuratowski convergence we immediately have
$\limsup_{j\to+\infty}\mathcal S_{\e_j}(t)\subset\mathcal S(t)$. We claim that the reverse inclusion holds as well.

Suppose otherwise. Then, up to passing to a subsequence, there exist $(x,z)\in\mathcal S(t)$ and $R>0$ such that $ \mathcal W_R^\Phi(x,z)\cap\mathcal S_{\e_j}(t)=\varnothing$ for
every $j$. For $j$ sufficiently large the hypotheses of Proposition~\ref{prop:wulff} are satisfied. Hence, if $s>t$ is sufficiently close to $t$, by comparison with the shrinking
$\Phi$-Wulff shape we have $\mathcal W^\Phi_{R/2}(x,z)\cap\mathcal S_{\e_j}(s)=\varnothing$ for all sufficiently large $j$. More precisely, by Proposition~\ref{prop:wulff} the
solitons avoid the whole corresponding Wulff space-time tube; therefore, after possibly taking $s-t$ smaller, a space-time neighborhood of $\overline{\mathcal
W^\Phi_{R/2}(x,z)}\times\{s\}$ is disjoint from $\mathcal S_{\e_j}$ for large $j$. Passing to the space-time Kuratowski limit yields $ \mathcal W^\Phi_{R/2}(x,z)\cap\mathcal
S(s)=\varnothing$ for every $s>t$ sufficiently close to $t$. On the other hand, since $t\notin\mathcal N$, the slices $\mathcal S(s)$ converge to $\mathcal S(t)$ in the Kuratowski
sense as $s\searrow t$. It follows that $ \mathcal W^\Phi_{R/2}(x,z)\cap\mathcal S(t)=\varnothing, $ contradicting $(x,z)\in\mathcal S(t)$. Thus
\[
\mathcal S_{\e_j}(t)\stackrel{\mathcal K}{\longrightarrow}\mathcal S(t) \qquad\text{in }\R^{N+1}
\]
for every $t\in(0,T^*)\setminus\mathcal N$, which proves (d).

Finally, let $t>T^*$. Choose $t_0\in(T^*,t)$. By (c), $\mathcal S(t_0)=\varnothing$. Starting at time $t_0$ with fixed Wulff shapes of arbitrarily large radius, the same avoidance
argument shows that $\mathcal S_{\e_j}(t)$ eventually avoids every compact subset of $\R^{N+1}$. Hence $\mathcal S_{\e_j}(t)\stackrel{\mathcal
K}{\longrightarrow}\varnothing=\mathcal S(t)$, proving (e).

We now  prove (a) and (f). For the former, let $\overline{W^\phi_r(x_0)}\subset \mathrm{int}(E_0)$. By Proposition~\ref{prop:paraboloidi} for $j$ large enough, we have that
$\mathcal{C}_{r,\frac1{4N}, \e_j}(x_0)\subset E_{\e_j}$. Passing to the limit we obtain
  \[
  \overline{W^\phi_r(x_0)}\times[0,+\infty)\subset \mathcal S(0)\quad\text{for any }\overline{W^\phi_r(x_0)}\subset \mathrm{int}(E_0)
  \]
  and thus, recalling that $E_{\e_j}\cap H^-=E_0\times(-\infty, 0)$, we deduce $\mathrm{int }(E_0)\times \R\subset \mathcal S(0)$.

  For the opposite inequality assume that  $\overline{W^\phi_r(x_0)}\subset \R^N\setminus \overline{E_0}$. Then, by the last part of Proposition~\ref{prop:paraboloidi} we have
  $\mathcal{C}_{r,\frac1{4N}, \e_j}(x_0)\subset  \overline{\R^{N+1}\setminus E_{\e_j}}$ for $j$ large. In particular, setting
\[
\widetilde{\mathcal{C}}_{r/2,\eta, \e_j}(x_0):=\mathcal{C}_{r/2,\eta, \e_j}(x_0)\cup \left(\overline{W^\phi_{r/2}(x_0)}\times(-\infty,0]\right)
\]
for any fixed $\eta\in (0, \frac1{4N})$, we have that  $\widetilde{\mathcal{C}}_{r/2,\eta, \e_j}(x_0)-\frac{t}{\e_j}e_{_{N+1}}$ has positive distance from
$E_{\e_j}-\frac{t}{\e_j}e_{_{N+1}}=\mathcal{S_{\e_j}}(t)$ bounded away from zero uniformly with respect to $j$ and $t>0$. Observe now that for any  sequence of times $t_j\searrow
0$  we have
  \[
\widetilde{\mathcal{C}}_{r/2,\eta, \e_j}(x_0)-\frac{t_j}{\e_j}e_{_{N+1}} \stackrel{\mathcal K}{\longrightarrow} \overline{W^\phi_{r/2}(x_0)}\times \R\quad\text{in the Kuratowski sense in }\R^{N+1}
\]
  and thus $\big(\overline{W^\phi_{r/2}(x_0)}\times \R\big)\cap \mathcal{S}(0)=\varnothing$. This concludes the proof of (a).

  It remains to prove  (f). Let $\mathcal{z}_\e$, $\mathcal{U}_\e$ be as in Proposition~\ref{prop:superflow} and note any compact subset $K$ of
  $\bigl(\R^{N+1}\times(0,T^*)\bigr)\setminus\mathcal S$ is eventually contained in $\mathcal{U}_\e$ for all $\e$ sufficiently small. Hence, up to extracting a non-relabelled
  subsequence, we may assume that
  \[
  \mathcal{z}_\e\to \mathcal{z}\quad\text{locally weakly* in }L^\infty\big(\R^{N+1}\times(0,\infty)\big).
  \]
  The fact that
  \[
  \partial_t d\geq \Div \mathcal{z}\quad\text{distributionally in }\R^{N+1}\times(0,T^*)\setminus \mathcal{S}
  \]
   and that
   \[
  \mathcal{z}\in \partial\Phi(\nabla d) \quad\text{a.e. in }\R^{N+1}\times(0,T^*)\setminus \mathcal{S}
   \]
   can be proven passing to the limit  properties \eqref{step5} and \eqref{step5bis} and arguing as in Step~5 of the proof of Proposition~\ref{prop:superflow}. Likewise, from
   \eqref{step5cor2} we deduce
   \[
     \Div \mathcal{z} \leq \frac{N}{\sigma}\quad\text{in the sense of distributions in  $\R^{N+1}\times(0,T^*)\cap \{d>\sigma\}$.}
   \]
   This concludes the proof of property (f) and thus of the proposition.
\end{proof}

We are now ready to prove the main theorem of the section.
\begin{proof}[Proof of Theorem~\ref{thm:main1}]  Let $\mathcal{S}$ be any subsequential Kuratowski limit of $(\mathcal{S_\e})$ as in Proposition~\ref{prop:helly} along a vanishing subsequence $\{\e_j\}$. Thus, in order to establish (a) and (b) of Theorem~\ref{thm:main1}, we need to show that $\mathcal{S}\cap\big( \R^{N+1}\times(0,+\infty)\big)=\mathcal{E}$.

To this aim, let $\mathcal{S}_{\e}^\complement(t), \mathcal{S}_{\e_j}^\complement$ be given by \eqref{def:Schius},\eqref{def:Schiu1}, respectively. Recalling
Remark~\ref{rm:complementare} and arguing symmetrically with $E_0$ replaced by $\R^N\setminus E_0$, we infer from Proposition~\ref{prop:helly} that up to passing to a further
(non-relabeled) subsequence
\[
\mathcal{S}_{\e_j}^\complement \stackrel{\mathcal K}{\longrightarrow}\widehat{\mathcal{S}}\quad\text{in the Kuratowski sense in }\R^{N+2},
\]
with $\widehat{\mathcal{S}}$ a superflow for $V=-\Phi\, H^\Phi$ such that
\[
(\R^N\setminus \overline{E_0})\times \R\subset \widehat{\mathcal{S}}(0)\subset (\R^N\setminus \mathrm{int }(E_0))\times \R.
\]
 Thus, setting $\mathcal A:= \left(\R^{N+1}\times[0,+\infty)\setminus \widehat{\mathcal{S}}\right)$, we have
\begin{equation}\label{aeesse}
\mathcal A\subset \mathcal{S}\text{ with $\mathcal A$ and $\mathcal{S}$ a subflow and a superflow, respectively, for $V=-\Phi\, H^\Phi$},
 \end{equation}
 with $\mathrm{int\, }(E_0)\times\R\subset \mathcal{A}(0)\subset  \mathcal{S}(0)\subset \overline{E_0}\times\R.$ We split the rest of the proof into two steps.

\medskip

\noindent{\it Step 1.} ($\mathcal{S}=\mathcal{E}$) Let $ u= u(x,t)$ be the level set flow for $V=-\phi\, H^\phi$, given  by the non-fattening assumption on $E_0$, see
Definition~\ref{def:fat}, and recall that $u$ satisfies  \eqref{nonfat1} and \eqref{nonfat2}. Let now $  \hat u_0=\hat u_0(x,z):= u(x,0)$, so that by \eqref{nonfat1} and the
previous step
\begin{multline}\label{aesse2}
\big\{(x,z)\in \R^{N+1}:  \hat u_0(x,z)<0\big\}\subset \mathrm{int\, }(E_0)\times \R \subset \mathcal A(0)\\
\subset \mathcal S(0) \subset \overline{E_0}\times \R\subset \big\{(x,z)\in \R^{N+1}:  \hat u_0(x,z)\leq 0\big\}.
\end{multline}
Since $\Phi(x,0)=\phi(x)$, it readily follows by uniqueness that $\hat u=\hat u(x,z,t):= u(x,t)$ is the (unique) level set flow for  $V=-\Phi\, H^\Phi$ in the sense of
Definition~\ref{deflevelset1}, with initial datum $\hat u_0$.

Thus, recalling that by Theorem~\ref{th:phiregularlevelset}-(ii) $\big\{(x,z,t)\in \R^{N+1}\times[0,\infty):  \hat u(x,z, t)<0\big\}$ and $\big\{(x,z,t)\in
\R^{N+1}\times[0,\infty):  \hat u(x,z, t)\le 0\big\}$ are the minimal subflow and the maximal superflow, respectively, with respect to their own initial data,  and taking into
account \eqref{aeesse} and \eqref{aesse2}, we have
\begin{multline}\label{eq:quasi}
\big\{(x,z,t)\in \R^{N+1}\times[0,\infty):  \hat u(x,z, t)<0\big\}\subset \mathcal{A}\\
\subset \mathcal{S}\subset \big\{(x,z,t)\in \R^{N+1}\times[0,\infty):  \hat u(x,z, t)\le 0\big\}.
\end{multline}
Note now that \eqref{nonfat2} is equivalent to
\[
\begin{split}
    &\overline{\big\{(x,z,t)\in \R^{N+1}\times(0,\infty): \vphantom{\tfrac12}\hat u(x,z,t)<0\big\}}=\big\{(x,z,t)\in \R^{N+1}\times(0,\infty): \hat u(x,z,t)\leq 0\big\},\\
    & \big\{(x,z, t)\in \R^{N+1}\times(0,\infty): \vphantom{\tfrac12}\hat u(x,z,t)<0\big\}=\mathrm{int\, }\big(\big\{(x,z,t)\in \R^{N+1}\times(0,\infty): \hat u(x,z,t)\leq 0\big\}\big),
    \end{split}
\]
where the closure is taken with respect to the relative topology of $\R^{N+1}\times(0,\infty)$. The latter display, combined with \eqref{eq:quasi} and the fact that $\mathcal{S}$
is closed, finally yields
\[
\mathcal{S}\cap\big( \R^{N+1}\times(0,+\infty)\big)= \big\{(x,z,t)\in \R^{N+1}\times(0,\infty):  \hat u(x,z, t)\le 0\big\}=\mathcal E.
\]
 Recalling also Proposition~\ref{prop:helly}, this establishes parts (a) and (b) of Theorem~\ref{thm:main1}.
\medskip

\noindent{\it Step 2. (Convergence of $\mathcal{R_\e}(E_\e)$).} Let $\e_j\to 0^+$ and let $\{(x_j, t_j)\}\subset \R^{N}\times(0,+\infty)$ be such that $(x_j, t_j)\in
\mathcal{R}_{\e_j}(E_{\e_j})$  and $(x_j, t_j)\to (x,t)$, with $t>0$. We claim that $(x,t)\in E$.

To this aim, observe that by definition of $\mathcal{R}_{\e_j}(E_{\e_j})$,  $t_j=\e_jz_j$  for  a suitable $z_j> 0$ such that $(x_j, z_j)\in E_{\e_j}$. But then, $(x_j, 0)=(x_j,
z_j)-\frac{t_j}{\e_j}e_{_{N+1}}\in \mathcal{S}_{\e_j}(t_j)$ or, equivalently, $(x_j, 0, t_j)\in \mathcal{S}_{\e_j}$. Since $\mathcal{S}_{\e_j}\stackrel{\mathcal
K}{\longrightarrow}\mathcal{E}$ in the Kuratowski sense in $\R^{N+1}\times(0,+\infty)$, it follows that $(x,0,t)\in \mathcal{E}$. Thus, $x\in E(t)$, that is,  $(x,t)\in E$.

Assume now that $(x,t)\in E$ and thus $(x,0,t)\in \mathcal{E}$, $t>0$. Then, again by  Kuratowski convergence of  $\mathcal{S}_{\e_j}$ to  $\mathcal{E}$, there exists $(x_j, z_j,
t_j)\to (x,0,t)$, with $(x_j, z_j, t_j)\in \mathcal{S}_{\e_j}$ for all $j$. In particular, $(x_j, z_j)\in E_{\e_j}-\frac{t_j}{\e_j}e_{_{N+1}}$, which is equivalent to $(x_j, \e_j
z_j+t_j)\in \mathcal{R}_{\e_j}(E_{\e_j})$. Since $(x_j, \e_j z_j+t_j)\to (x,t)$, this allows us to conclude that $\mathcal{R}_{\e_j}(E_{\e_j})\stackrel{\mathcal
K}{\longrightarrow}E$ in the Kuratowski sense in $\R^{N}\times(0,+\infty)$. The theorem is proven.
\end{proof}

\section{Elliptic regularization with natural mobility: BV and transport identities,  Brakke-flow realization}\label{sec:ellip:varifold}

\subsection{Energy estimates and geometric properties of the limiting measures}\label{s:limit}

Let $\phi$ be a norm on $\R^N$, and let $\Phi$ be its extension defined in \eqref{Phi}. We first recall the first variation identities in the case $\phi\in
C^1(\R^N\setminus\{0\})$. Under this assumption $\Phi\in C^1(\R^{N+1}\setminus\{0\})$ and satisfies
\begin{equation}\label{prop:extnorm}
	\nabla \Phi(v)\cdot e_{_{N+1}}=\frac{v\cdot e_{_{N+1}}}{\Phi(v)},\quad\text{for all $v\in \R^{N+1}\setminus\{0\}$.}
\end{equation}

Let $E_0\subset \R^N$ be a set of finite perimeter, and for every $\e>0$ we consider a minimizer $E_\e$ for $\mathcal{A}^\Phi_\e$ as in \eqref{minpbell} with initial datum $E_0$;
we denote with
\begin{equation}\label{def:abbrevio}
    \mu_\e:=\mu_{E_\e}=\H^N\res \partial^* E_\e.
\end{equation}
Under the same differentiability assumption, by Proposition~\ref{prop:Euler_lagrange} we have, for every vector field $Y\in C_c^1(\R^N\times(0,\infty);\R^{N+1})$,
\begin{equation}\label{E:first-variationEe}
	 \frac1\e \int e^{-z/\e}\, \Phi(\nu)\,Y\cdot e_{_{N+1}}\,d\mu_\e(x,z) =\int e^{-z/\e} B_\Phi(\nu) : DY d\mu_\e(x,z),
\end{equation}
with $B_\Phi$ defined in \eqref{def:Bpsi};  by \eqref{euler:lagrange1} and \eqref{prop:extnorm}, $\partial^* E_\e$ admits $\Phi$-anisotropic mean curvature satisfying the soliton
equation
\begin{equation}\label{eq:soliton}
    \mathbf{H}_{\partial^*E_\e}^\Phi +\frac{1}{\e}\frac{(\nu \cdot e_{_{N+1}})}{\Phi(\nu)} \nu=0 \quad\text{on $\partial^* E_\e\cap \{z>0\}$,}
\end{equation}
where $\nu=\nu(x,z)$ is the outer unit normal to $\partial^* E_\e$.

\subsection{Energy estimates}\label{subsec:encvmeas}
Throughout this subsection $\phi$ is an arbitrary norm and  when $\phi$ is of class $C^1$ away from the origin, we write $\mathbf H^\Phi=\mathbf
H^\Phi_{\partial^*E_\e}$.

Let $H^+=\{(x,z)\in \R^{N+1}:\,z>0\}$, and let $\mathcal{A}_\e^\Phi$ be as in \eqref{defAF}-\eqref{localized:Apsi} with $\Phi$ in place of $\Psi$. We start with the following
energy estimates.

\begin{lemma}\label{lem:stimEeHp}
    Let $E_\e$ be a minimizer of \eqref{minpbell}. Then, for every $\e>0$ we have
\begin{equation}\label{ti:rendi:conto}
    \int_{{\{z>-1\}}} e^{-z/\e}\,\Phi(\nu)\,d\mu_\e(x,z)=\mathcal{A}_\e^\Phi(E_\e)\leq \mathcal{A}_\e^\Phi(E_0\times \R),
\end{equation}
   where $\nu=\nu_{E_\e}$ is the normal to $E_\e$, and
\begin{equation} \label{en-uniform-bound}
    \frac{1}{\e}\mathcal{A}^\Phi_\e(E_\e;H^+)\leq P^\phi(E_0).
\end{equation}
\end{lemma}

\begin{proof}
    Since  $E_0\times \R$ is an admissible competitor for the minimization  problem \eqref{minpbell}, we have
\begin{equation*}
   \frac{1}{\e}\mathcal{A}^\Phi_{\e}(E_0\times \R)\geq \frac{1}{\e}\mathcal{A}^\Phi_{\e}(E_\e),
\end{equation*}
that is \eqref{ti:rendi:conto}; then observe that, as the generalized normal  $\nu_{\partial^* E_0\times \R}=(\nu_{\partial^* E_0},0)$ and thus
$\Phi(\nu)=\phi(\nu_{\partial^*{E_0}})$ by \eqref{Phi}, we have
\begin{equation*}
\begin{split}
     \frac{1}{\e}\mathcal{A}^\Phi_{\e}(E_0\times \R)&= \frac{1}{\e} \Big(\int_{-1}^{0}e^{-z/\e}dz\Big)\,\int_{\partial^*E_0} \phi(\nu)d\mathcal{H}^{N-1}+\frac{1}{\e}\Big(\int_0^\infty e^{-z/\e}dz\Big)\int_{\partial^*E_0} \phi(\nu)d\mathcal{H}^{N-1} \\
     & =\frac{1}{\e} \Big(\int_{-1}^{0}e^{-z/\e}dz\Big)\,P^\phi(E_0)+P^\phi(E_0),
     \end{split}
\end{equation*}
and, since $E_\e=E_0\times (-\infty,0)$ in $\{z<0\}$, we also have

\begin{equation*}
    \frac{1}{\e} \mathcal{A}^\Phi_{\e}(E_\e)\geq \frac{1}{\e}\Big(\int_{-1}^0e^{-z/\e}dz\Big)\,P^\phi(E_0)+\frac{1}{\e}\mathcal{A}^\Phi_{\e}(E_\e;H^+).
\end{equation*}
Combining the three expressions above yields \eqref{en-uniform-bound}, that is our thesis.
\end{proof}

In what follows we write $\nu=(\nu',\nu\cdot e_{_{N+1}})$.
\begin{lemma}\label{lem:possible:test}\label{lem:1}\label{lem:2}
Let $E_\e$ be a minimizer of \eqref{minpbell}. For every Borel set $A\subset(0,\infty)$ of finite measure, we have
\begin{equation}\label{eq:prima e seconda}
\int_{\R^N\times A}\frac{\phi(\nu')^2}{\Phi(\nu)}\,d\mu_\e\leq |A|P^\phi(E_0),\qquad
\frac1\e\int_{H^+}\frac{(\nu\cdot e_{_{N+1}})^2}{\Phi(\nu)}\,d\mu_\e\leq P^\phi(E_0).
\end{equation}
If $\phi\in C^1(\R^N\setminus\{0\})$, the second estimate reads, by \eqref{eq:soliton},
\begin{equation}\label{stima:eneH2}
    \int_{H^+} \e |\mathbf{H}^{\Phi}|^2  \,\Phi(\nu) d\mu_\e \leq P^\phi(E_0).
\end{equation}
\end{lemma}
\begin{proof}
Let $\xi\in C_c^1((0,\infty))$ and note that  for $|s|$ sufficiently small, the map $F_s(x,z)=(x,z+s\xi(z))$ is a diffeomorphism coinciding  with  the identity near $\{z\leq0\}$.
Note also that by  the area formula and \eqref{Phi},
\[
\mathcal A_\e^\Phi(F_s(E_\e);H^+)
=\int_{H^+}e^{-(z+s\xi(z))/\e}\sqrt{(1+s\xi'(z))^2\phi(\nu')^2+(\nu\cdot e_{_{N+1}})^2}\,d\mu_\e .
\]
Differentiating with respect to $s$ at $s=0$,  by minimality we obtain
\begin{equation}\label{altre:eulero:lagrange}
\int_{H^+}e^{-z/\e}\frac{\phi(\nu')^2}{\Phi(\nu)}\,\xi'(z)\,d\mu_\e
=\frac1\e\int_{H^+}e^{-z/\e}\Phi(\nu)\,\xi(z)\,d\mu_\e .
\end{equation}
The same identity holds for $\xi\in C^1([0,\infty))$ with $\xi(0)=0$ and both $\xi$ and $\xi'$ bounded. Indeed, take $0<2\delta<R$, $R\geq1$ and consider the test function
$\chi_{\delta,R}\xi$, where $0\leq\chi_{\delta,R}\leq1$ is a smooth cutoff function equal to one on $[2\delta,R]$ and vanishing outside $(\delta,2R)$, with  $|\chi'_{\delta,R}|\leq
C/\delta$ on $(\delta,2\delta)$ and $|\chi'_{\delta,R}|\leq C/R$ on $(R,2R)$. In $(\chi_{\delta,R}\xi)'=\chi_{\delta,R}\xi'+\chi'_{\delta,R}\xi$, only the second term requires an
estimate. Since $\xi(0)=0$, we have $|\xi(z)|\leq z\|\xi'\|_\infty$, and hence
\[
|\chi'_{\delta,R}(z)\xi(z)|\leq
\begin{cases}
(C/\delta)\,z\|\xi'\|_\infty\leq 2C\|\xi'\|_\infty, & \text{for }\delta<z<2\delta,\\
(C/R)\,\|\xi\|_\infty, &  \text{for } R<z<2R;
\end{cases}
\]
that is,  $\chi_{\delta,R}\xi$ and their derivatives are uniformly bounded. This allows one to test  \eqref{altre:eulero:lagrange} with $\chi_{\delta,R}\xi$ in place of $\xi$ and
to recover  \eqref{altre:eulero:lagrange} by letting $\delta\to 0+$ and $R\to\infty$ (and using dominated convergence and \eqref{en-uniform-bound}).

Given $\eta\in C_c^1([0,\infty))$ nonnegative, we may therefore take $\xi(z)=e^{z/\e}\eta(z)-\eta(0)$ in \eqref{altre:eulero:lagrange} and  using the identity
$\Phi(\nu)^2-\phi(\nu')^2=(\nu\cdot e_{_{N+1}})^2$, we find
\begin{equation}\label{altre:eulero:lagrange1}
\frac1\e\int_{H^+}\eta(z)\frac{(\nu\cdot e_{_{N+1}})^2}{\Phi(\nu)}\,d\mu_\e
-\int_{H^+}\eta'(z)\frac{\phi(\nu')^2}{\Phi(\nu)}\,d\mu_\e
=\frac{\eta(0)}{\e}\mathcal A_\e^\Phi(E_\e;H^+)\le \eta(0)P^\phi(E_0),
\end{equation}
where the last inequality follows from \eqref{en-uniform-bound}. Taking $\eta(z)=\int_z^\infty\psi(r)\,dr$ for any nonnegative $\psi\in C_c((0,\infty))$ and discarding the first
term in \eqref{altre:eulero:lagrange1} we obtain
\[
\int_{H^+}\psi(z)\frac{\phi(\nu')^2}{\Phi(\nu)}\,d\mu_\e\leq P^\phi(E_0)\int_0^\infty\psi(r)\,dr,
\]
which implies the first estimate in \eqref{eq:prima e seconda}. We may now take $\eta$ nonnegative and nonincreasing, equal to one on $[0,R]$ in \eqref{altre:eulero:lagrange1}.
Discarding  the second term and letting  $R\to\infty$ we also obtain the second estimate in \eqref{eq:prima e seconda}, thus concluding the proof.
\end{proof}

\begin{proposition}\label{prop:energyestimate}
    For every measurable set $A\subset (0,\infty)$ with $|A|<\infty$, it holds
    \begin{equation} \label{eq:energyestimate}
    P^\Phi(E_\e;\R^N\times A) \leq  (|A|+\e) \,P^\phi(E_0),
    \end{equation}
    and
    \begin{equation}\label{eq:rescaled:energy}
        P\big(\mathcal{R}_\e(E_\e);\R^N\times A\big)\leq C_\Phi\,\Big((\e^2+|A|)+(\e^2+|A|)^{\frac{1}{2}} \Big)\,P^\phi(E_0).
    \end{equation}
for all $\e>0$, where $C_\Phi$ is a constant depending only on $N,\Phi$.
\end{proposition}

\begin{proof}
Combining the identity $\Phi(\nu)=\phi(\nu')^2/\Phi(\nu)+(\nu\cdot e_{_{N+1}})^2/\Phi(\nu)$ with the two estimates of Lemma~\ref{lem:2}, we obtain \eqref{eq:energyestimate}. For the
rescaling $(x,z)\mapsto(x,\e z)$, the area formula yields
\[
P\big(\mathcal R_\e(E_\e);\R^N\times A\big)
=\int_{\R^N\times(A/\e)}\sqrt{\e^2|\nu'|^2+(\nu\cdot e_{_{N+1}})^2}\,d\mu_\e
\leq\int_{\R^N\times(A/\e)}\big(\e|\nu'|+|\nu\cdot e_{_{N+1}}|\big)\,d\mu_\e .
\]
Since $|\nu'|\leq C_\Phi\Phi(\nu)$, the contribution of the first term is bounded by $C_\Phi(|A|+\e^2)P^\phi(E_0)$ by \eqref{eq:energyestimate}. For the second, we have
\[
\int_{\R^N\times(A/\e)}|\nu\cdot e_{_{N+1}}|\,d\mu_\e
\leq\left(\int_{H^+}\frac{(\nu\cdot e_{_{N+1}})^2}{\Phi(\nu)}\,d\mu_\e\right)^{1/2}
\left(P^\Phi(E_\e;\R^N\times(A/\e))\right)^{1/2}
\leq \sqrt{|A|+\e^2}\,P^\phi(E_0),
\]
where  we used the second estimate in \eqref{eq:prima e seconda} and \eqref{eq:energyestimate} again. Thus we have shown \eqref{eq:rescaled:energy}.
\end{proof}

\subsection{Convergence of rescalings, translating solitons and space-time tracks}
\label{subsec:limit:L1}

From now on, we shall assume that
\begin{equation*}
    \text{$E_0$ is a set of finite perimeter that is \textit{non-fattening in measure} in the sense of Definition \ref{def:fat}}
\end{equation*}
with respect to the anisotropy $\Psi=\phi$ on $\R^N$, and natural mobility $\mathfrak{M}=\phi$; correspondingly, we consider the (unique) flat flow $\{\widehat{E}(t)\}_{t\geq 0}$
as in Section \ref{subses:flat:flow}. We also denote by
\begin{equation}\label{altre:def}
\begin{split}
    &\widehat{E}\coloneqq  \left\{ (x,t)\in \R^{N}\times(0,\infty): x\in \widehat{E}(t)\right\}.
    \end{split}
\end{equation}

By using the newfound perimeter estimates, we are able to prove the $L^1$-counterpart of Theorem \ref{thm:main1}.

\begin{theorem}\label{thm:altro:main}
  Under the notations and assumptions above, the following limits hold:
  \begin{itemize}
\item[(a)]  there exists a set $\mathcal{N}\subset (0,+\infty)$  of null measure such that
   \begin{equation}\label{BV}
        \mathcal{S}_{\e}(t)\xrightarrow{\e\to 0^+} \widehat{E}(t)\times \R\quad\text{weakly-* in $BV_{\mathrm{loc}}(\R^{N+1})$ for every $t\in(0,\infty)\setminus\mathcal N$.}
    \end{equation}
\item[(b)] the rescaled sets
\begin{equation}\label{conv:riscalati}
    \mathcal{R}_\e\big(E_\e \big)\xrightarrow{\e\to0^+} \widehat{E}\quad\text{weakly-* in $BV_{\rm{loc}}(\R^{N}\times (0,\infty))$,}
\end{equation}
and we have the estimate
\begin{equation}\label{Eabs:cont}
    P\big(\widehat{E};\,\R^N\times A \big)\leq C_\Phi\big(|A|+|A|^{1/2} \big)\,P^\phi(E_0).
\end{equation}
for all Borel sets $A\subset (0,\infty)$ with $|A|<\infty$;
\end{itemize}
\end{theorem}

Before proving Theorem \ref{thm:altro:main}, let us prove a basic fact concerning the relation between Kuratowski and $L^1_{\loc}$ convergence of sets.

\begin{lemma}\label{lemma:L1Kur}
    Let $\{E_j\}_{j\in \N}$ be a sequence of measurable sets in $\R^d$. Suppose that
    \begin{equation*}
            E_j\xrightarrow{L^1_{loc}} E\quad\text{and}\quad \overline E_j\xrightarrow{\mathcal{K}} C,
    \end{equation*}
    for some measurable set $E$ and a closed set $C$. Then $ E^{(1)}\subset C$, where $E^{(1)}$ denotes the sets of points of $E$ with density one as in
    \eqref{def:density-pts}. In particular, $|E\setminus C|=0$.
\end{lemma}

\begin{proof}
Let $x\notin C$. Since $C$ is closed, there exists $r>0$ such that $\overline{B_r(x)}\cap C=\varnothing$. We claim that
\[
\overline{E_j}\cap\overline{B_r(x)}=\varnothing\quad\text{for all sufficiently large $j$.}
\]
 Otherwise, after passing to a subsequence, we could find $x_j\in\overline{E_j}\cap\overline{B_r(x)}.$ for all $j$.

By compactness, up to a further subsequence, $x_j\to y$ for some $y\in\overline{B_r(x)}$, and by the characterization of the Kuratowski convergence we would then have  $y\in C$,
contradicting $\overline{B_r(x)}\cap C=\varnothing$. Hence $E_j\cap B_r(x)=\varnothing$ for all sufficiently large $j$, and since $E_j\to E$ in $L^1_{\mathrm{loc}}(\R^d)$, it
follows that $|E\cap B_r(x)|=0$. This implies that $x\notin E^{(1)}$, that is our thesis.
\end{proof}

\begin{proof}[Proof of Theorem \ref{thm:altro:main}]
  (a)  First let us recall that, by density estimates (see Remark \ref{rem:Lmin}), we have that $\mathrm{int}(E_\e)=\overline{E_\e}$ for $\mathcal{L}^{N+1}$-a.e. $(x,z)\in \{z\geq
  0\}$ (in short $\mathrm{int}(E_\e)\stackrel{a.e.}{=}\overline{E_\e}$ on $\{z\geq 0\}$), and thus
    \begin{equation}\label{int:clos}
        \mathrm{int}\,(\mathcal{S}_\e(t))\stackrel{a.e.}{=}\overline{\mathcal{S}_\e(t)}\quad\text{ on $\{z\geq -t/\e\}$.}
    \end{equation}
 Then, owing to \eqref{eq:energyestimate} and \eqref{perimetri:Equivalent}, for every $t \ge 0$, the sequence $\{\mathcal{S}_\varepsilon(t)\}_{\varepsilon>0}$ satisfies local
 perimeter estimates  on $\R^N\times (-t/\e,+\infty)$ which are uniform in $\varepsilon\in (0,1)$ and $t\geq 0$. Therefore, we may extract a subsequence $\varepsilon_j(t)$
 (possibly depending on $t$) such that
\begin{equation}\label{BVS}
    \mathcal{S}_{\varepsilon_j(t)}(t) \xrightarrow{j \to \infty} \mathcal{S}(t) \quad \text{weakly-* in } BV_{\mathrm{loc}}(\mathbb{R}^{N+1}).
\end{equation}
for all $t>0$, for some limit set $\mathcal{S}(t)$. Analogously, up to extracting a further subsequence, for all  $t>0$ we have
\begin{equation}\label{BVcomplement}
    \mathcal{S}_{\e_j(t)}^{\complement}(t)\stackrel{a.e.}{=}\R^{N+1}\setminus \mathcal{S}_{\e_j(t)}(t)\xrightarrow{j\to\infty} \R^{N+1}\setminus \mathcal{S}(t)\quad\text{weakly-* in $BV_{\rm{loc}}(\R^{N+1})$,}
\end{equation}
where $\mathcal{S}_{\e_j(t)}^{\complement}(t)$ are given by \eqref{def:Schius}, and the first equality stems from \eqref{int:clos}.

On the other hand,  by Theorem \ref{thm:main1}, we have $\mathcal{S}_\e(t)\xrightarrow{\e\to 0^+}  \mathcal{E}(t)$ in the Kuratowski sense, for all $t\not\in \mathcal{N}_0$, for
some countable set $\mathcal{N}_0\subset(0,\infty)$. Lemma \ref{lemma:L1Kur} thus yields
\begin{equation}\label{first:inclusion}
    \mathcal{S}(t)^{(1)}\subset \mathcal{E}(t)\quad \text{for all  $t\notin \mathcal{N}_0$.}
\end{equation}

Owing to Remark \ref{remark:complementary}, we also have
\[
\mathcal{S}^\complement_{\e}(t)\xrightarrow{\e\to 0^+} \R^{N+1}\setminus \mathrm{int}\big(\mathcal{E}(t) \big)\quad\text{ in the Kuratowski sense for all  $t\notin \mathcal{N}_0$.}
\]
Thus, by \eqref{BVcomplement} and Lemma~\ref{lemma:L1Kur},  we deduce
\begin{equation}\label{second:inclusion}
    \R^{N+1}\setminus \mathcal{S}(t)^{(1)}\stackrel{a.e.}{=}(\R^{N+1}\setminus \mathcal{S}(t))^{(1)}\subset  \R^{N+1}\setminus \mathrm{int}\big(\mathcal{E}(t) \big)\quad \text{for all  $t\notin \mathcal{N}_0$ .}
\end{equation}
By non-fattening in measure and \eqref{flat-uguale}, there exists a set $\mathcal{N}_1\subset (0,\infty)$  of null measure  such that
$\mathrm{int}\big(\mathcal{E}(t)\big)\stackrel{a.e}{=}\widehat E(t)\times \R\stackrel{a.e.}{=}\mathcal{E}(t)$ for all $t\not\in \mathcal{N}_1$. This piece of information coupled
with \eqref{first:inclusion} and \eqref{second:inclusion} yields $\mathcal{S}(t)\stackrel{a.e}{=}\widehat{E}(t)\times \R$    for all $t\notin (\mathcal{N}_0\cup \mathcal
{N}_1)\eqqcolon \mathcal{N}$. Combining this fact with \eqref{BVS}, we obtain the convergence in \eqref{BV} along the subsequence $\e_j(t)$. Since the same argument applies to any subsequence $\{\e_k\}$ and yields the same limit, the convergence in \eqref{BV} thus holds for the entire sequence. \\
 \noindent (b) The proof of \eqref{conv:riscalati} is very similar to that of \eqref{BV}. By the energy estimate \eqref{eq:rescaled:energy}, up to subsequences
    \begin{equation}
        \mathcal{R}_\e(E_\e)\to \mathcal{R}\quad\text{weakly-* in $BV_{\rm{loc}}(\R^N\times (0,\infty))$},
    \end{equation}
    for some set $\mathcal{R}$ of locally finite perimeter in $\R^N\times (0,\infty)$. On the other hand, by the Kuratowski convergence of Remark \ref{remark:complementary} and
    Lemma \ref{lemma:L1Kur}, we have that $\mathcal{R}^{(1)}\subset E$. By the same argument applied to the complementary, we get $(\R^{N}\times (0,\infty)\setminus
    \mathcal{R})^{(1)}\subset \R^{N}\times (0,\infty)\setminus\mathrm{int}(E)$, hence
    \[
    \mathrm{int}(E)\stackrel{a.e}{\subset}\mathcal{R}^{(1)}\subset E\stackrel{a.e}{=}\mathrm{int}(E),
    \]
    the last equality due to the non-fattening in measure. Hence $\mathcal{R}\stackrel{a.e}{=}E\stackrel{a.e}{=}\widehat{E}$ (the last equality still due to the non-fattening
    condition), from which \eqref{conv:riscalati} follows. Finally, the energy estimate \eqref{Eabs:cont} for $A$ open follows, via lower-semicontinuity, by passing to the limit in
    \eqref{eq:rescaled:energy}. For general Borel set $A\subset (0,\infty)$, it follows by approximation.
\end{proof}
\subsection{General transport and BV identities}

\begin{theorem}[Transport and BV identities]\label{trasportobv}
Let $\phi$ be a norm on $\R^N$ (possibly crystalline) and consider the natural mobility $\mathfrak m=\phi$. Let $E_0\subset\R^N$ have finite perimeter and be non-fattening in measure. Denote by $\{\widehat E(t)\}_{t\ge 0}$ the corresponding flat flow.

Set
\[
\widehat E = \{(x,t)\in\R^N\times(0,\infty):x\in\widehat E(t)\}, \qquad d\mu=d\mu_{\widehat E(t)}\,dt, \qquad d\sigma=d\mu^\phi_{\widehat E(t)}\,dt.
\]

Then, there exist a scalar field $v\in L^2(\mu)$, a vector field $\mathbf b\in L^2(\sigma;\R^N)$, and a finite nonnegative Radon measure $\Lambda$ on $\R^N\times(0,\infty)$ with
the following properties.
\begin{itemize}
\item[i)] The tube $\widehat E$ has locally finite perimeter, and
\begin{equation}\label{tubofase}
D_{x,t}\chi_{\widehat E} = (-\nu_{\widehat E(t)},v)\,\mu,
\qquad \mathbf b\cdot\nu_{\widehat E(t)}=v \quad\mu\text{-a.e.}
\end{equation}
In particular, for every $0\leq t_1<t_2$ and every $C^1$ test function $\eta$ with compact spatial support uniformly on $[t_1,t_2]$ one has
\begin{equation}\label{trasportofase}
\begin{split}
\int_{\widehat E(t_2)}\eta(x,t_2)\,dx -\int_{\widehat E(t_1)}\eta(x,t_1)\,dx &= \int_{t_1}^{t_2}\int_{\widehat E(t)}\partial_t\eta\,dx\,dt + \int_{t_1}^{t_2}\int\eta v\,d\mu_{\widehat E(t)}\,dt .
\end{split}
\end{equation}
\item[ii)] The  energy transport equation
\begin{equation}\label{bilancioenergia}
\partial_t\sigma+\operatorname{div}_x(\mathbf b\sigma)=-\Lambda
\end{equation}
holds in the sense of distributions, and
\begin{equation}\label{dissipazione}
\bigl(\phi^o(\mathbf b)\bigr)^2\,\sigma \leq\Lambda, \qquad \Lambda(\R^N\times(0,\infty))
\leq P^\phi(E_0).
\end{equation}
More precisely, for every $0\leq t_1<t_2$ and every $C^1$ admissible test function $\eta$ as before one has
\begin{equation}\label{bilanciointervalli}
\begin{aligned}
&\mu^\phi_{\widehat E(t_2)}(\eta(\cdot,t_2)) -\mu^\phi_{\widehat E(t_1)}(\eta(\cdot,t_1)) \\
&\qquad= \int_{t_1}^{t_2}\int (\partial_t\eta+\mathbf b\cdot\nabla\eta)\, d\mu^\phi_{\widehat E(t)}\,dt - \int_{\R^N\times(t_1,t_2]}\eta\,d\Lambda .
\end{aligned}
\end{equation}

\end{itemize}
\end{theorem}
\begin{remark}[Transport and dissipation of surface energy]\label{trasportoenergia}
Theorem~\ref{trasportobv} provides a local energetic description of the flat flow, \emph{valid also for crystalline anisotropies}. The measure $\sigma=\mu^\phi_{\widehat E(t)}\,dt$
represents the surface energy, and $\bu\sigma$ is its transport flux. The identity $\bu\cdot\nu=v$, where $\nu=\nu_{\widehat E(t)}$, ensures that the normal component of this
transport agrees with the normal velocity of the evolving phase. The nonnegative measure $\Lambda$ records the local dissipation of surface energy after transport has been taken
into account, while the inequality $(\phi^o(\bu))^2\sigma\leq\Lambda$ states that this loss controls the quadratic cost of transport.

In particular, the theorem yields the energy dissipation inequality
\begin{equation}\label{dissopt}
P^\phi(\widehat E(t_2))+\int_{\R^N\times(t_1,t_2)}\left(\frac{v}{\phi(\nu)}\right)^2\,d\sigma \leq P^\phi(\widehat E(t_1)), \qquad 0\leq t_1<t_2.
\end{equation}
To see this, we may apply \eqref{bilanciointervalli} to spatial cutoff functions $\eta$ with support in $B_{2R}$, tending to one as $R\to\infty$, and with gradients bounded by
$C/R$. The flux contribution involving $\mathbf b\cdot\nabla\eta$ vanishes as $R\to\infty$ and thus
\[
P^\phi(\widehat E(t_1))-P^\phi(\widehat E(t_2)) =\Lambda\bigl(\R^N\times(t_1,t_2]\bigr).
\]
On the other hand,  by the compatibility condition in~\eqref{tubofase} we have
\[
\frac{\abs{v}}{\phi(\nu)} =\frac{\abs{\bu\cdot\nu}}{\phi(\nu)} \leq\phi^o(\bu).
\]
By squaring and integrating, and then using~\eqref{dissipazione}, we obtain the claimed inequality. Thus, the decrease of the anisotropic perimeter bounds the time-integrated
squared normal velocity, weighted by the inverse mobility.

These conclusions establish a local balance between transport, dissipation, and phase motion without identifying the transport field through curvature. Such an identification is
the missing information  required to recover the full curvature-driven evolution law, see the following remark.
\end{remark}
\begin{remark}\label{rm:once}
Once the above theorem is established, in order to get the Brakke/BV formulation (in the $C^1$ strictly convex case) it is enough to identify the transport field $\mathbf b$,
namely, to show that
\[
\mathbf b = (\nu_{\widehat E(t)}\cdot\mathbf H^\phi_{\pa^*\widehat E(t)}) \,\nabla\phi(\nu_{\widehat E(t)})
\]
Indeed, plugging the above expression in \eqref{trasportofase} and using the second equation in \eqref{tubofase}, we get the BV flow property, while from
\eqref{bilanciointervalli}, together with first inequality in \eqref{dissipazione}, we get the Brakke flow property. Note also that plugging the above identification of $\bu$  in
\eqref{dissopt} yields the optimal energy/dissipation inequality
\begin{equation}\label{daje}
P^\phi(\widehat E(t_2))-P^\phi(\widehat E(t_1))\leq -\int_{\R^N\times(t_1,t_2)}| \mathbf H^\phi_{\pa^*\widehat E(t)}|^2\,d\sigma, \qquad 0\leq t_1<t_2.
\end{equation}
\end{remark}

\begin{proof}[Proof of Theorem~\ref{trasportobv}] We split the proof into several steps.

\medskip

\noindent\textit{Step 1: Approximation of norms and convergence of minimizers.}

Fix a starting time $t_0\geq0$, and set $G_0=\widehat E(t_0)$. Choose smooth norms $\phi_n$ such that, after taking a subsequence,
\[
(1-n^{-1})\phi\leq\phi_n\leq(1+n^{-1})\phi,
\]
and define $\Phi_n(\xi,\zeta)=\sqrt{\phi_n(\xi)^2+\zeta^2}$, as in \eqref{Phi}.  In particular, $\Phi_n\geq c|\cdot|$ with $c>0$ independent of $n$, and $P^{\phi_n}(G_0)\to
P^\phi(G_0)$. If $\phi$ is already $C^1$ away from the origin, we may instead take $\phi_n=\phi$ throughout and omit $n$: the Euler--Lagrange equation and the energy estimates used
below require only this regularity. Let $G_{\e,n}$ minimize \eqref{minpbell} with anisotropy $\Phi_n$ and datum $G_0$. Comparison with $G_0\times\R$ (see the proof of Lemma
\ref{lem:stimEeHp}) yields
\[
\mathcal A_\e^{\Phi_n}(G_{\e,n}) \leq \e\, e^{1/\e}P^{\phi_n}(G_0).
\]
For fixed $\e>0$, by coercivity and positivity of the weight on compact subsets of $\{z>-1\}$ we obtain local BV compactness. Thus, along a subsequence, $G_{\e,n}\to G_\e$ in
$L^1_{\loc}$, and the limit satisfies the same constraint.

For every admissible competitor $F$ of finite weighted energy,
\[
\begin{aligned}
\mathcal A_\e^\Phi(G_\e) &\leq \liminf_n\mathcal A_\e^{\Phi_n}(G_{\e,n}) \\
&\leq \limsup_n\mathcal A_\e^{\Phi_n}(G_{\e,n}) \leq \lim_n\mathcal A_\e^{\Phi_n}(F) = \mathcal A_\e^\Phi(F).
\end{aligned}
\]
Note that the first inequality follows from lower semicontinuity for the fixed anisotropy $\Phi$ and the fact that $\Phi_n\to \Phi$ locally uniformly. Hence $G_\e$ minimizes the
problem with anisotropy $\Phi$.

\medskip
\noindent\textit{Step 2: Uniform estimates and convergence of the rescaled sets.}

In computations on $G_{\e,n}$, we write $\nu=(\nu',\nu\cdot e_{_{N+1}})$ for its outer normal and abbreviate $H^{\Phi_n}=H^{\Phi_n}_{\partial^*G_{\e,n}}$. By \eqref{eq:soliton} and
\eqref{scalH-SOFP},
\begin{equation}\label{solitonbis}
H^{\Phi_n} = \frac{\nu\cdot e_{_{N+1}}}{\e\,\Phi_n(\nu)}.
\end{equation}
For any Borel set $A\subset (0,+\infty)$, by  \eqref{stima:eneH2} and \eqref{eq:energyestimate} we have
\begin{equation}\label{1}
\begin{split}
&\e\int_{H^+}|H^{\Phi_n}|^2\, d\mu_{G_{\e,n}}^{\Phi_n} \leq P^{\phi_n}(G_0), \\
&\mu_{G_{\e,n}}^{\Phi_n}(\R^N\times A) \leq (|A|+\e)P^{\phi_n}(G_0),
\end{split}
\end{equation}
while from \eqref{eq:rescaled:energy} we obtain
\begin{equation}\label{2}
\begin{split}
P(\mathcal R_\e(G_{\e,n});\R^N\times A) &\leq \left( C_{\Phi_n}(|A|+\e^2) +(|A|+\e^2)^{\frac12} \right)P^{\phi_n}(G_0)\\
&\leq \left( C(|A|+\e^2) +(|A|+\e^2)^{\frac12} \right)P^{\phi_n}(G_0)
\end{split}
\end{equation}
after noticing that, since $\Phi_n\to \Phi$, we have that $C_{\Phi_n}$ can be bounded uniformly with respect to $n$.

At fixed $\e$, sending $n\to\infty$ in \eqref{2} by lower semicontinuity of the perimeter we obtain the corresponding   bound for $\mathcal R_\e(G_\e)$ on sets of the form
$\R^N\times A$, with $A\subset(0,+\infty)$ open. This allows us to reproduce the argument of Theorem~\ref{thm:altro:main} to get
\[
\mathcal R_\e(G_\e)+(0,t_0) \longrightarrow\widehat E \quad\text{in }L^1_{\loc}(\R^N\times(t_0,\infty)).
\]
Choose $\e_j\searrow0$. By the above considerations together with a  diagonal argument, we may find  $n_j\to \infty$ so that
\begin{equation}\label{3}
\mathcal R_{\e_j}(G_{\e_j,n_j})+(0,t_0) \longrightarrow\widehat E \quad\text{in }L^1_{\loc}(\R^N\times(t_0,\infty)).
\end{equation}
Passing to the limit in estimate \eqref{2} also shows that the limiting tube has locally finite perimeter.

\medskip
\noindent\textit{Step 3: Volume and energy identities before passage to the limit.}

Let $\eta\in C_c^1(\R^N\times(t_0,\infty))$.  Changing variables and applying the Gauss-Green formula for sets of finite perimeter yields
\begin{equation}\label{4}
\begin{split}
-\partial_t\chi_{\mathcal R_\e(G_{\e,n})+(0,t_0)}(\eta) &=\int_{\mathcal R_\e(G_{\e,n})+(0,t_0)} \partial_t\eta\,dx\,dt = \int_{G_{\e,n}\cap H^+} \partial_z[\eta(x,t_0+\e z)]\,dx\,dz\\
&= \int_{\partial^*G_{\e,n}\cap H^+} \eta(x,t_0+\e z) (\nu_{G_{\e,n}}\cdot e_{_{N+1}})\,d\Ha^N= \e\int_{H^+} \eta(x,t_0+\e z)H^{\Phi_n}\, d\mu_{G_{\e,n}}^{\Phi_n},
\end{split}
\end{equation}
where the last equality follows from ~\eqref{solitonbis}. The above identity easily yields the following estimate on the total variation measure $\Big| \partial_t\chi_{\mathcal
R_\e(G_{\e,n})+(0,t_0)} \Big|$: for every nonnegative continuous compactly supported $\eta$ we have
\begin{equation}\label{5}
\begin{aligned}
&\Big| \partial_t\chi_{\mathcal R_\e(G_{\e,n})+(0,t_0)} \Big|(\eta) \\
&\quad\leq \left( \e\int_{H^+}\eta(x,t_0+\e z)\, d\mu_{G_{\e,n}}^{\Phi_n} \right)^{1/2} \left( \e\int_{H^+}\eta(x,t_0+\e z) |H^{\Phi_n}|^2\,d\mu_{G_{\e,n}}^{\Phi_n} \right)^{1/2}.
\end{aligned}
\end{equation}
Next, the first variation identity reads
\[
\int B_{\Phi_n}(\nu):DY\,d\mu_{G_{\e,n}}(x,z)=\int H^{\Phi_n}\, \nu\cdot Y d\mu_{G_{\e,n}}(x,z),\quad Y\in C^1_c(H^+;\R^{N+1});
\]
choosing $Y(x,z)=\eta(x,t_0+\e z)\,e_{_{N+1}}$, and using \eqref{solitonbis},  \eqref{def:Bpsi} and \eqref{prop:extnorm}, we obtain

\[
\begin{aligned}
\int_{H^+} \left[ \e\left( 1-\frac{(\nu\cdot e_{_{N+1}})^2}{\Phi_n(\nu)^2} \right)\partial_t\eta - \frac{(\nu\cdot e_{_{N+1}})}{\Phi_n(\nu)} \nabla_\xi\Phi_n(\nu)\cdot\nabla_x\eta \right]d\mu^{\Phi_n}_{G_{\e,n}}= \int_{H^+} \frac{(\nu\cdot e_{_{N+1}})^2} {\e\,\Phi_n(\nu)^2} \,\eta\,d\mu^{\Phi_n}_{G_{\e,n}}.
\end{aligned}
\]
Here $\eta$ and its derivatives are evaluated at $(x,t_0+\e z)$. Using again \eqref{solitonbis}, which yields $\frac{(\nu\cdot e_{_{N+1}})^2}{\Phi_n(\nu)^2} =\e^2|H^{\Phi_n}|^2$ and $
\frac{(\nu\cdot e_{_{N+1}})^2} {\e\Phi_n(\nu)^2} =\e|H^{\Phi_n}|^2$, and substituting in the above identity, we therefore obtain
\begin{equation}\label{6}
\begin{aligned}
0=\e\int_{H^+}\Big( &(1-\e^2|H^{\Phi_n}|^2)\partial_t\eta -H^{\Phi_n}\nabla_\xi\Phi_n(\nu)\cdot\nabla\eta -|H^{\Phi_n}|^2\eta \Big)\, d\mu_{G_{\e,n}}^{\Phi_n},
\end{aligned}
\end{equation}
here the derivatives of $\eta$ are evaluated at $(x,t_0+\e z)$. Finally, by the extension formula we have, when $\nu'\neq0$, $\nabla_\xi\Phi_n(\nu)
=\frac{\phi_n(\nu')}{\Phi_n(\nu)}\,\nabla\phi_n(\nu')$. Since $\phi_n^o(\nabla\phi_n(\nu'))=1$ and the horizontal gradient is zero at $\nu'=0$, we get
\begin{equation}\label{7}
\phi_n^o(\nabla_\xi\Phi_n(\nu))\leq1.
\end{equation}
In particular, by uniform equivalence of the norms we also have $|\nabla_\xi\Phi_n(\nu)|\leq C_\phi$, independently of $n$.

\medskip
\noindent\textit{Step 4: Limiting measures and transport identity before restarting.} We now formulate the compactness argument directly in terms of space--time measures. Set
$\mathcal T_j(x,z):=(x,t_0+\e_jz)$ and define the nonnegative Radon measures on $\R^N\times(t_0,\infty)$ by
\[
\lambda_j := \e_j(\mathcal T_j)_\# \left( \mu^{\Phi_{n_j}}_{G_{\e_j,n_j}}\mres H^+ \right), \qquad \Lambda_j := \e_j(\mathcal T_j)_\# \left( |H^{\Phi_{n_j}}|^2 \mu^{\Phi_{n_j}}_{G_{\e_j,n_j}}\mres H^+ \right).
\]
Equivalently, for every $\eta\in C_c(\R^N\times(t_0,\infty))$,
\[
\int\eta\,d\lambda_j = \e_j\int_{H^+} \eta(x,t_0+\e_jz)\, d\mu^{\Phi_{n_j}}_{G_{\e_j,n_j}}, \quad \int\eta\,d\Lambda_j = \e_j\int_{H^+} \eta(x,t_0+\e_jz) |H^{\Phi_{n_j}}|^2\, d\mu^{\Phi_{n_j}}_{G_{\e_j,n_j}}.
\]
  The estimates in~\eqref{1} read
\begin{equation}\label{questaserve}
\begin{split}
&\Lambda_j\bigl(\R^N\times(t_0,\infty)\bigr) = \e_j\int_{H^+} |H^{\Phi_{n_j}}|^2\, d\mu^{\Phi_{n_j}}_{G_{\e_j,n_j}} \leq P^{\phi_{n_j}}(G_0),\\
&\lambda_j(\R^N\times I) = \e_j \mu^{\Phi_{n_j}}_{G_{\e_j,n_j}} \left( \R^N\times\frac{I-t_0}{\e_j} \right) \leq \bigl(|I|+\e_j^2\bigr)P^{\phi_{n_j}}(G_0).
\end{split}
\end{equation}
for any bounded open interval $I\Subset (t_0, +\infty)$. Since $P^{\phi_{n_j}}(G_0)\to P^\phi(G_0)$, by weak compactness of Radon measures, we may extract  a non relabelled
subsequence and find nonnegative Radon measures $\lambda$ and $\Lambda$ on $\R^N\times(t_0,\infty)$ such that
\[
\lambda_j\weakstar\lambda, \qquad \Lambda_j\weakstar\Lambda \quad\text{ in }\R^N\times(t_0,\infty),
\]
with
\begin{equation}\label{forseserve}
\Lambda\bigl(\R^N\times(t_0,\infty)\bigr) \leq P^\phi(G_0), \qquad \lambda(\R^N\times I) \leq |I|P^\phi(G_0)
\end{equation}
for every bounded open interval $I\Subset(t_0,\infty)$.  It follows in particular that the time-marginal of $\lambda$, defined by $A\mapsto\lambda(\R^N\times A)$,  is absolutely
continuous with respect to Lebesgue measure and is bounded above by $P^\phi(G_0)\mathcal L^1\mres(t_0,\infty)$. The disintegration theorem therefore yields a measurable family of
nonnegative Radon measures $\{\lambda_t\}_{t>t_0}$ on $\R^N$ such that
\begin{equation}\label{9}
d\lambda=d\lambda_t(x)\,dt, \qquad \lambda_t(\R^N)\leq P^\phi(G_0) \quad\text{for a.e. }t>t_0.
\end{equation}
Thus,
\begin{equation}\label{8}
\begin{aligned}
&\lim_j\e_j\int_{H^+} \eta(x,t_0+\e_jz)\, d\mu_{G_{\e_j,n_j}}^{\Phi_{n_j}} = \int_{t_0}^\infty\int\eta\,d\lambda_t\,dt, \\
&\lim_j\e_j\int_{H^+} \eta(x,t_0+\e_jz)|H^{\Phi_{n_j}}|^2\, d\mu_{G_{\e_j,n_j}}^{\Phi_{n_j}} = \int\eta\,d\Lambda.
\end{aligned}
\end{equation}
Using the maps $\mathcal T_j$ introduced above, we now introduce the ``transport''  vector-valued space-time measures defined as
\[
J_j := -\e_j(\mathcal T_j)_\# \left( H^{\Phi_{n_j}} \nabla_\xi\Phi_{n_j}(\nu) \, \mu^{\Phi_{n_j}}_{G_{\e_j,n_j}}\mres H^+ \right),
\]
that is, for any $X\in C_c(\R^N\times(t_0,\infty);\R^N)$
\begin{equation}\label{defJj}
\int X\cdot dJ_j = -\e_j\int_{H^+} H^{\Phi_{n_j}} \nabla_\xi\Phi_{n_j}(\nu) \cdot X(x,t_0+\e_jz)\, d\mu^{\Phi_{n_j}}_{G_{\e_j,n_j}}.
\end{equation}
We first establish local compactness of these vector measures. By \eqref{7},
\[
\left| \nabla_\xi\Phi_{n_j}(\nu) \cdot X(x,t_0+\e_jz) \right| \leq \phi_{n_j}\bigl(X(x,t_0+\e_jz)\bigr).
\]
It follows from Cauchy--Schwarz that
\begin{equation}\label{CS}
\begin{aligned}
\left|\int X\cdot dJ_j\right| &\leq \left( \e_j\int_{H^+} \phi_{n_j}\bigl(X(x,t_0+\e_jz)\bigr)^2 \,d\mu^{\Phi_{n_j}}_{G_{\e_j,n_j}} \right)^{1/2} \\
&\qquad{}\times \left( \e_j\int_{H^+} |H^{\Phi_{n_j}}|^2 \mathbf \chi_{\spt X}(x,t_0+\e_jz) \,d\mu^{\Phi_{n_j}}_{G_{\e_j,n_j}} \right)^{1/2} \\
&= \left( \int\phi_{n_j}(X)^2\,d\lambda_j \right)^{1/2} \left( \Lambda_j(\spt X) \right)^{1/2}\leq C_\phi \left( \int|X|^2\,d\lambda_j \right)^{1/2} P^{\phi_{n_j}}(G_0)^{1/2},
\end{aligned}
\end{equation}
for $j$ large enough. Note that in the last inequality,  we have used that the norms $\phi_{n_j}$ are converging to $\phi$ and \eqref{questaserve}. Recalling also local mass bounds
for $\lambda_j$ we infer that the total variations $|J_j|$ are uniformly bounded on compact subsets of $\R^N\times(t_0,\infty)$. Thus,  up to a further subsequence, there is
therefore an $\R^N$-valued Radon measure $J$ such that $J_j\weakstar J$ in $\R^N\times(t_0,\infty)$; that is,
\[
\begin{aligned}
\int X\cdot dJ ={}& -\lim_{j\to\infty} \e_j\int_{H^+} H^{\Phi_{n_j}} \nabla_\xi\Phi_{n_j}(\nu)\cdot X(x,t_0+\e_jz)\, d\mu^{\Phi_{n_j}}_{G_{\e_j,n_j}}
\end{aligned}
\]
for every compactly supported continuous vector field $X$.

Since  $\lambda_j\weakstar\lambda$, passing to the limit in \eqref{CS} and using the first estimate in~\eqref{1}, we obtain
\[
\left|\int X\cdot dJ\right| \leq C_\phi P^\phi(G_0)^{1/2} \left( \int|X|^2\,d\lambda \right)^{1/2}
\]
for all  $X\in C_c(\R^N\times(t_0,\infty);\R^N)$. Thus the functional $X\mapsto\int X\cdot dJ$, $X\in C_c(\R^N\times(t_0,\infty);\R^N)$, is continuous with respect to the
$L^2(\lambda;\R^N)$ norm. Since compactly supported continuous vector fields are dense in this space, the functional extends uniquely to all of $L^2(\lambda;\R^N)$, and thus  by
Riesz Theorem there exist $\mathbf b\in L^2(\lambda;\R^N)$ such that
\begin{equation}\label{limitingJ}
J=\mathbf b\,\lambda=\mathbf b\,\lambda_t\,dt,
\end{equation}
as vector-valued Radon measures.

Observe now that by~\eqref{7} we have $\phi_{n_j}^o(H^{\Phi_{n_j}} \nabla_\xi\Phi_{n_j}(\nu)) \leq |H^{\Phi_{n_j}}|$. In turn,  Young inequality yields for every $x\in\R^N$
\[
- 2 H^{\Phi_{n_j}} \nabla_\xi\Phi_{n_j}(\nu)\cdot x-\phi_{n_j}(x)^2 \leq |H^{\Phi_{n_j}}|^2.
\]
Thus, for any $X\in C_c(\R^N\times(t_0,\infty);\R^N)$ and any nonnegative $\eta\in C_c(\R^N\times(t_0,\infty))$, integrating with respect to the rescaled energy measure, and using
the definitions of $J_j$, $\lambda_j$, and $\Lambda_j$, we obtain
\[
2\int\eta X\cdot dJ_j - \int\eta\,\phi_{n_j}(X)^2\,d\lambda_j \leq \int\eta\,d\Lambda_j.
\]
Passing to the limit, we obtain
\begin{equation}\label{mii}
\int \eta\bigl(2\mathbf b\cdot X-\phi(X)^2\bigr)\,d\lambda \leq \int\eta\,d\Lambda.
\end{equation}
Recall that $(\phi^o)^2/2$ is the Legendre transform of $\phi^2/2$. Equivalently,
\[
\sup_{x\in\R^N} \left\{ 2\mathbf b\cdot x-\phi(x)^2 \right\} = \bigl(\phi^o(\mathbf b)\bigr)^2.
\]
Thus, by \eqref{mii} we obtain
\[
\sup_{X\in C_c(\R^N\times(t_0,\infty);\R^N)} \int \eta\bigl(2\mathbf b\cdot X-\phi(X)^2\bigr)\,d\lambda= \int\eta\bigl(\phi^o(\mathbf b)\bigr)^2\,d\lambda \leq \int\eta\,d\Lambda.
\]
Note that the fact that the  supremum can be moved under the integral follows from standard measure theory arguments. By the arbitrariness of $\eta$ and using  the
disintegration~\eqref{9}, we have  established
\begin{equation}\label{10}
\bigl(\phi^o(\mathbf b)\bigr)^2\lambda_t\,dt\leq\Lambda.
\end{equation}

We now rewrite~\eqref{6} in terms of the rescaled measures. By the definitions of $\lambda_j$, $J_j$, and $\Lambda_j$,  \eqref{6} is equivalent to
\[
\int\partial_t\eta\,d\lambda_j +\int\nabla\eta\cdot dJ_j = \int\eta\,d\Lambda_j +\e_j^2\int\partial_t\eta\,d\Lambda_j.
\]
The last term tends to zero, since $\left| \e_j^2\int\partial_t\eta\,d\Lambda_j \right| \leq \e_j^2P^{\phi_{n_j}}(G_0)\|\partial_t\eta\|_\infty$. Thus, we may pass to the limit and
use the disintegration~\eqref{9} and \eqref{limitingJ} to obtain
\begin{equation}\label{11}
\int_{t_0}^\infty\int \bigl(\partial_t\eta+\mathbf b\cdot\nabla\eta\bigr) \,d\lambda_t\,dt = \int\eta\,d\Lambda
\end{equation}
for any $\eta\in C^1_c(\R^N\times (t_0,+\infty))$. Likewise, by \eqref{3}, lower semicontinuity of the total (time) variation, and \eqref{5} we obtain
\begin{equation}\label{12}
|\partial_t\chi_{\widehat E}|(\eta) \leq \left(\int_{t_0}^\infty\int\eta\,d\lambda_t\,dt\right)^{1/2} \left(\int\eta\,d\Lambda\right)^{1/2}
\end{equation}
for any nonnegative $\eta\in C_c(\R^N\times (t_0,+\infty))$.

We now show that the limiting energy dominates the anisotropic perimeter of the slices. Let $\eta\in C_c(\R^N\times(t_0,\infty))$ be nonnegative, and set $F_j=\mathcal
R_{\e_j}(G_{\e_j,n_j})+(0,t_0)$. Primes on normals denote their spatial components. Then
\begin{equation}\label{as:in}
\begin{split}
\int_{t_0}^\infty\int \eta\,d\mu^\phi_{\widehat E(t)}\,dt &= \int_{\partial^*\widehat E} \eta\,\phi(\nu'_{\widehat E})\,d\Ha^N\leq \liminf_{j\to\infty} \int_{\partial^*F_j} \eta\,\phi_{n_j}(\nu')\,d\Ha^N \\
&= \liminf_{j\to\infty}\e_j \int_{\partial^*G_{\e_j,n_j}\cap H^+} \eta(x,t_0+\e_jz)\, \phi_{n_j}(\nu')\,d\Ha^N \\
&\leq \lim_{j\to\infty}\e_j \int_{H^+} \eta(x,t_0+\e_jz)\, d\mu^{\Phi_{n_j}}_{G_{\e_j,n_j}}= \lim_{j\to\infty}\int\eta\,d\lambda_j = \int_{t_0}^\infty\int\eta\,d\lambda_t\,dt.
\end{split}
\end{equation}
Note that the  first equality follows from Coarea Formula,  the first inequality is a consequence of  the $L^1_{\loc}$ convergence of $F_j$ to $\widehat E$, together with
Reshetnyak's lower semicontinuity theorem  and the fact that $\phi_{n_j}\to \phi$, while the second equality follows from  change-of-variables and the one-homogeneity of
$\phi_{n_j}$. The arbitrariness of $\eta$ then yields
\begin{equation}\label{13}
\lambda_t\geq\mu^\phi_{\widehat E(t)} \quad\text{for a.e. }t>t_0.
\end{equation}

From now on, including at time zero, the starting datum is $G_0=\widehat E(t_0)$.

\medskip
\noindent\textit{Step 5: Choice of the precise representative, and monotonicity of the total mass.}

The family $\lambda_t$ obtained by disintegration is defined only for almost every $t>t_0$. We now construct a precise representative for all $t\geq t_0$, by considering the
essential right limit  at $t_0$ and the essential left limit  at every interior time.

To this aim, fix $\varphi\in C_c^1(\R^N)$  and let $\pi(x,t)=t$. Testing~\eqref{11} with $\eta(x,t)=\zeta(t)\varphi(x)$, $\zeta\in C^1_c(t_0,+\infty)$, we obtain
\begin{equation}\label{timeder}
D_t\bigl(\lambda_t(\varphi)\bigr) = \left( \int_{\R^N} \mathbf b(x,t)\cdot\nabla\varphi(x)\,d\lambda_t(x) \right)dt - \pi_\#(\varphi\Lambda).
\end{equation}
The right-hand side has finite total variation on any bounded  interval $(t_0,T)$. Indeed, \eqref{9}, \eqref{10}, \eqref{forseserve}, and equivalence of norms imply
\begin{equation}\label{stimetta}
\int_{t_0}^T\int|\mathbf b|\,d\lambda_t\,dt \leq \left( \int_{t_0}^T\lambda_t(\R^N)\,dt \right)^{1/2} \left( \int_{t_0}^T\int|\mathbf b|^2\,d\lambda_t\,dt \right)^{1/2} \leq C_\phi P^\phi(G_0)\sqrt{T-t_0}.
\end{equation}
In particular,
\[
\left|D_t\bigl(\lambda_t(\varphi)\bigr)\right|((t_0,T)) \leq C_\phi P^\phi(G_0)\|\nabla\varphi\|_\infty \sqrt{T-t_0}+\|\varphi\|_\infty \Lambda(\R^N\times(t_0,T)).
\]
Together with the bound $|\lambda_t(\varphi)|\leq P^\phi(G_0)\|\varphi\|_\infty$, this proves that $t\mapsto\lambda_t(\varphi)$ is locally BV and thus it admits essential left
limit at every $\tau\in(t_0,T)$ and an essential right limit at $t_0$.

Choose a countable uniformly dense family $\{\varphi_k\}_{k\in\N}\subset C_c^1(\R^N)$ in $C_0(\R^N)$. For each $k$, the scalar function $\lambda_t(\varphi_k)$ agrees outside a null
set with its left-continuous BV representative. Taking the union of these exceptional sets, together with the exceptional set in~\eqref{9}, we obtain a single null set $\mathcal
N$, independent of the test and of the time. Outside $\mathcal N$, each $\lambda_t(\varphi_k)$ has a left limit at every $\tau>t_0$ and a right limit at $t_0$. The uniform mass
bound $|\lambda_t(\varphi)|\leq P^\phi(G_0)\|\varphi\|_\infty$ (which follows from \eqref{9}) allows to extend these limits to every $\varphi\in C_0(\R^N)$:   More explicitly, if
$\varphi_k\to\varphi$ uniformly, the corresponding limiting values form a Cauchy sequence, and its limit is independent of the approximating sequence. The limits define positive
linear functionals bounded by $P^\phi(G_0)\|\varphi\|_\infty$. By Riesz representation, there are therefore finite nonnegative Radon measures $\lambda_{\tau-}$ and $\lambda_{t_0+}$
such that
\[
\begin{split}
&\lambda_t \weakstar\lambda_{\tau-} \quad \text{as }t\nearrow\tau,\ t\notin\mathcal N, \quad\text{for every }\tau>t_0,\\
&\lambda_t\weakstar\lambda_{t_0+}\quad \text{as }t\searrow t_0,\ t\notin\mathcal N,
\end{split}
\]
with $\lambda_{t_0+}(\R^N)$, $\lambda_{\tau-}(\R^N)\leq P^\phi(G_0)$. Moreover, by construction of $\mathcal N$, we have $\lambda_{\tau-}=\lambda_\tau$ for all $\tau\not\in
\mathcal N$.

We next identify $\lambda_{t_0+}$ and extend the domination \eqref{13} to each $\lambda_{\tau-}$. Fix a nonnegative $\varphi\in C_c(\R^N)$ and $\tau>t_0$. Approach $\tau$ from the
left through times $t$ for which~\eqref{13} holds and $\lambda_t(\varphi)$  converge to $\lambda_{\tau-}(\varphi)$. By the $L^1_{\loc}$ continuity of the flat flow and lower
semicontinuity of the (localized) anisotropic perimeter we obtain
\[
\mu^\phi_{\widehat E(\tau)}(\varphi) \leq \liminf_{s\to \tau-} \mu^\phi_{\widehat E(s)}(\varphi) \leq \lambda_{\tau-}(\varphi),
\]
where the liminf is taken along the chosen times. The same argument at $t_0$, approaching from the right and using $G_0=\widehat E(t_0)$, yields
$\mu^\phi_{G_0}(\varphi)\leq\lambda_{t_0+}(\varphi)$. Thus, we obtain
\begin{equation}\label{tracce}
\lambda_{\tau-}\geq\mu^\phi_{\widehat E(\tau)} \quad\text{for all }\tau>t_0, \qquad \lambda_{t_0+}\geq\mu^\phi_{G_0}.
\end{equation}
But $\lambda_{t_0+}(\R^N)\leq P^\phi(G_0)=\mu^\phi_{G_0}(\R^N)$ and therefore $\lambda_{t_0+}=\mu^\phi_{G_0}$.

We now fix the representative by setting
\[
\lambda_{t_0}:=\mu^\phi_{G_0}, \qquad \lambda_t:=\lambda_{t-} \quad\text{for }t>t_0.
\]
The chosen representative  is weakly* left-continuous on $(t_0,\infty)$ and converges weakly* to $\mu^\phi_{G_0}$ as $t \to t_0+$. Moreover, by \eqref{timeder} for every $t>t_0$
and $\varphi\in C_c^1(\R^N)$ it satisfies
\begin{equation}\label{thfund}
\lambda_t(\varphi) =\mu^\phi_{G_0}(\varphi)+ \int_{t_0}^t\int\mathbf b\cdot\nabla\varphi\,d\lambda_s\,ds -\int_{\R^N\times(t_0,t)}\varphi(x)\,d\Lambda(x,s).
\end{equation}
We can now compute the total mass directly from this identity. Choose $\varphi_R\in C_c^1(\R^N)$ such that $0\leq\varphi_R\leq1$, $\varphi_R=1$ on $B_R$,
$\operatorname{spt}\varphi_R\subset B_{2R}$, $|\nabla\varphi_R|\leq C/R$, and $\varphi_R\nearrow 1$ as $R\to+\infty$. For each fixed $t>t_0$, the flux term satisfies
\[
\left| \int_{t_0}^t\int \mathbf b\cdot\nabla\varphi_R\,d\lambda_s\,ds \right| \leq \frac{C}{R} \int_{t_0}^t\int|\mathbf b|\,d\lambda_s\,ds \longrightarrow0
\]
as $R\to+\infty$. Monotone convergence in the other terms, together with the domination already proved, yields
\begin{equation}\label{14}
\lambda_t\geq\mu^\phi_{\widehat E(t)}, \quad \lambda_t(\R^N) = P^\phi(G_0)-\Lambda(\R^N\times(t_0,t)) \qquad\text{for all }t\geq t_0.
\end{equation}
Since $\Lambda$ is nonnegative, the formula shows that $t\mapsto\lambda_t(\R^N)$ is nonincreasing,  left-continuous on $(t_0,\infty)$ and tends to $P^\phi(G_0)$ as $t\to t_0+$.

Note now that more generally, for every $T>t_0$ and $\eta\in C_c^1(\R^N\times\R)$, we have
\begin{equation}\label{thfund-general}
\int_{t_0}^{T}\int \bigl(\partial_t\eta+\mathbf b\cdot\nabla\eta\bigr) \,d\lambda_t\,dt = \int_{\R^N\times(t_0,T)}\eta\,d\Lambda + \lambda_T\bigl(\eta(\cdot,T)\bigr) -\mu^\phi_{G_0}\bigl(\eta(\cdot,t_0)\bigr).
\end{equation}
Indeed, testing \eqref{11} with $\zeta(t)\eta(x,t)$, where $\zeta\in C_c^1((t_0,T))$, and expanding the time derivative, we find
\begin{equation}\label{Dt:lt}
D_t\!\left[\lambda_t(\eta(\cdot,t))\right] = \left( \int_{\R^N} (\partial_t\eta+\mathbf b\cdot\nabla\eta)\,d\lambda_t \right)dt - \pi_\#(\eta\Lambda).
\end{equation}
We now record the following properties of the actual perimeter measures. For a nonnegative $\varphi\in C_c^1(\R^N)$, \eqref{thfund}, the positivity of $\Lambda$, and
\eqref{stimetta} imply
\[
\mu^\phi_{\widehat E(t)}(\varphi) \leq\lambda_t(\varphi)\leq \mu^\phi_{\widehat E(t_0)}(\varphi) + C_\phi P^\phi(\widehat E(t_0)) \|\nabla\varphi\|_\infty\sqrt{t-t_0}.
\]
On the other hand, by $L^1_{\loc}$ continuity of the flat flow and lower semicontinuity of the localized perimeter we obtain $\mu^\phi_{\widehat E(t_0)}(\varphi) \leq \liminf_{t\to
t_0+} \mu^\phi_{\widehat E(t)}(\varphi)$. Thus $\mu^\phi_{\widehat E(t)}\weakstar \mu^\phi_{\widehat E(t_0)}$ as $t\to t_0+$. In other words the actual perimeter measures are
weakly right-continuous at $t_0$.

Finally, by \eqref{14} we have $P^\phi(\widehat E(t))\leq P^\phi(\widehat E(t_0))$ for every $t\geq t_0$. Since the  construction is available with every starting time $t_0$, the
actual perimeter measures of the flat flow starting from a set that is non-fattening in measure are weakly right-continuous at every time, and the function $t\mapsto P^\phi(\widehat E(t))$ is nonincreasing. Its
lower semicontinuity then implies its right continuity as well.

\medskip
\noindent\textit{Step 6: Restarting and matching.}

Fix $\delta>0$, set $t_0^\delta=0$ and proceed by induction as follows. Assume that $t_i^\delta$ is defined and consider the elliptic regularization scheme with starting time
$t_0=t_i^\delta$ and initial datum $\widehat E(t_i^\delta)$ as in the previous steps. Denote  by $\lambda_t^{(i)}$, $t>t_i^\delta$ the corresponding disintegrated measures, its
transport vector by $\mathbf b^{(i)}$, and its dissipation measure by $\Lambda^{(i)}$. Define the restarting time $t_{i+1}^\delta$ as
\begin{equation}\label{15}
t_{i+1}^\delta = \inf\left\{ t>t_i^\delta: \lambda_t^{(i)}(\R^N) \geq P^\phi(\widehat E(t))+\delta \right\},
\end{equation}
with value $+\infty$ if the set is empty. Note that any finite restarting time is attained. Indeed, by approaching the infimum along times satisfying the inequality in \eqref{15}
and  using the lower semicontinuity of the perimeter (in fact right-continuity) and monotonicity of $\lambda^{(i)}_t(\R^N)$ we infer that the inequality holds at the infimum as
well. Moreover, since $\lambda^{(i)}_{t_i^\delta}(\R^N)=P^\phi(\widehat E(t_i^\delta))$, $\lambda^{(i)}_{t}(\R^N)$ is nonincreasing and $P^\phi(\widehat E(t))$ is right continuous,
we easily deduce that $t_{i+1}^\delta>t_{i}^\delta$. Observe also that
\[
P^\phi(\widehat E(t_{i+1}^\delta))+\delta \leq \lambda^{(i)}_{t_{i+1}^\delta}(\R^N)\leq P^\phi(\widehat E(t_i^\delta)).
\]
Thus the procedure has at most $N_\delta\leq P^\phi(\widehat E(0))/\delta$ restarts. Set $t_{N_\delta+1}^\delta=\infty$ and $I_i^\delta=[t_i^\delta,t_{i+1}^\delta)$.

On $I_i^\delta$, set $\lambda_t^\delta=\lambda_t^{(i)}$ and $\mathbf b^\delta=\mathbf b^{(i)}$. Define the total dissipation, including the losses at restarting times, by
\begin{equation}\label{16}
\int\eta\,d\Lambda^\delta \coloneqq \sum_{i=0}^{N_\delta} \int_{\R^N\times(t_i^\delta,t_{i+1}^\delta)} \eta\,d\Lambda^{(i)} + \sum_{i=1}^{N_\delta} \left( \lambda_{t_i^\delta}^{(i-1)}
-\mu^\phi_{\widehat E(t_i^\delta)} \right) \bigl(\eta(\cdot,t_i^\delta)\bigr),
\end{equation}
where each measure in the second sum is nonnegative by \eqref{14}.

Then, recalling \eqref{thfund-general},  we get
\begin{equation}\label{17}
\int_0^\infty\int (\partial_t\eta+\mathbf b^\delta\cdot\nabla\eta) \,d\lambda_t^\delta\,dt = \int\eta\,d\Lambda^\delta.
\end{equation}

For each finite block, by \eqref{14} we have
\[
\begin{aligned}
&\Lambda^{(i)}(\R^N\times(t_i^\delta,t_{i+1}^\delta)) +\lambda_{t_{i+1}^\delta}^{(i)}(\R^N) -P^\phi(\widehat E(t_{i+1}^\delta))= P^\phi(\widehat E(t_i^\delta)) -P^\phi(\widehat E(t_{i+1}^\delta)).
\end{aligned}
\]
Summing over $i$ and cancelling the intermediate terms, we obtain the first inequality below, while the second follows from \eqref{10}:
\begin{equation}\label{18}
\Lambda^\delta(\R^N\times(0,\infty)) \leq P^\phi(\widehat E(0)), \qquad \bigl(\phi^o(\mathbf b^\delta)\bigr)^2 \,\lambda_t^\delta\,dt \leq\Lambda^\delta.
\end{equation}
Applying \eqref{12} on each  block and then Cauchy-Schwarz yields
\begin{equation}\label{19}
|\partial_t\chi_{\widehat E}|(\eta) \leq \left(\int_0^\infty\int\eta\,d\lambda_t^\delta\,dt\right)^{1/2} \left(\int\eta\,d\Lambda^\delta\right)^{1/2}.
\end{equation}

By construction
\[
\mu^\phi_{\widehat E(t)}\leq\lambda_t^\delta, \qquad \lambda_t^\delta(\R^N) \leq P^\phi(\widehat E(t))+\delta.
\]
Thus, recalling the definition of $d\sigma=d\mu^{\phi}_{\hat{E}(t)}dt$, we have
\begin{equation}\label{20}
\left|\lambda_t^\delta\,dt-\sigma\right|\big(\R^N\times(0,T)\big) \leq T\delta \qquad\text{for any }T<\infty.
\end{equation}

\medskip
\noindent\textit{Step 7: The matched limit.}

Let $\delta\downarrow0$. By \eqref{18}, equivalence of norms, and Cauchy--Schwarz, for every $X\in C_c(\R^N\times(0,\infty);\R^N)$,
\[
\left| \int_0^\infty\int X\cdot\mathbf b^\delta\,d\lambda_t^\delta\,dt \right| \leq C_\phi\sqrt{P^\phi(\widehat E(0))} \left( \int_0^\infty\int|X|^2\,d\lambda_t^\delta\,dt \right)^{1/2}.
\]
Combining this estimate with the mass bound in \eqref{20}, we obtain locally uniform bounds on the total variations of the vector measures $\mathbf b^\delta\lambda_t^\delta\,dt$.
Thus, after extraction, these vector measures and the nonnegative measures $\Lambda^\delta$ admit weak-star limits.

By \eqref{20}, the integral on the right-hand side converges to $\int|X|^2\,d\sigma$. The limiting vector measure therefore defines a continuous linear functional on
$L^2(\sigma;\R^N)$. The density and Riesz representation argument of Step 4 identifies this measure with $\mathbf b\sigma$, for some $\mathbf b\in L^2(\sigma;\R^N)$. Hence, along
the chosen subsequence,
\[
\mathbf b^\delta\lambda_t^\delta\,dt \weakstar\mathbf b\sigma, \qquad \Lambda^\delta\weakstar\Lambda.
\]
The limit $\Lambda$ is nonnegative, and the total-mass bound in \eqref{18} passes to the limit by lower semicontinuity:
\[
\Lambda(\R^N\times(0,\infty)) \leq P^\phi(\widehat E(0)).
\]
To inherit the localized quadratic estimate, we repeat the Young inequality argument used in Step 4. For every nonnegative $\eta\in C_c(\R^N\times(0,\infty))$ and every compactly
supported continuous vector field $X$, by \eqref{18} we have
\[
\int_0^\infty\int \eta\bigl(2\mathbf b^\delta\cdot X-\phi(X)^2\bigr) \,d\lambda_t^\delta\,dt \leq \int\eta\,d\Lambda^\delta.
\]
The convergences just established and \eqref{20} allow us to pass to the limit in each term, obtaining
\[
\int\eta\bigl(2\mathbf b\cdot X-\phi(X)^2\bigr)\,d\sigma \leq \int\eta\,d\Lambda.
\]
Taking the supremum over $X$, exactly as in Step 4 yields $\bigl(\phi^o(\mathbf b)\bigr)^2\sigma\leq\Lambda$. We can also pass to the limit in \eqref{17} to get for every $\eta\in
C_c^1(\R^N\times(0,\infty))$
\begin{equation}\label{restart:l1}
\int \bigl(\partial_t\eta+\mathbf b\cdot\nabla\eta\bigr)\,d\sigma = \int\eta\,d\Lambda,
\end{equation}
that is, \eqref{bilancioenergia}. Finally, fix a nonnegative $\eta\in C_c(\R^N\times(0,\infty))$. Passing to the limit in  \eqref{19}, we get
\[
|\partial_t\chi_{\widehat E}|(\eta) \leq \left(\int\eta\,d\sigma\right)^{1/2} \left(\int\eta\,d\Lambda\right)^{1/2}.
\]
By an approximation argument we deduce $|\partial_t\chi_{\widehat E}|(A) \leq\big(\sigma(A)\Lambda(A)\big)^{1/2}$ for every bounded Borel set $A\subset \R^N\times (0,+\infty)$. In
particular, $\partial_t\chi_{\widehat E}\ll\sigma$. Since $\sigma=\phi(\nu_{\widehat E(t)})\mu$ and $\phi$ is bounded above and below on the unit sphere, there exists $v\in
L^1_{\loc}(\mu)$ such that
\[
\partial_t\chi_{\widehat E}=v\mu.
\]
By slicing, the Gauss--Green theorem for sets of finite perimeter (and De Giorgi's Structure Theorem) we obtain $\nabla\chi_{\widehat E}=-\nu_{\widehat E(t)}\mu$. Thus,
\[
D_{x,t}\chi_{\widehat E} = (-\nu_{\widehat E(t)},v)\mu, \qquad |D_{x,t}\chi_{\widehat E}| = \sqrt{1+v^2}\,\mu.
\]

\medskip
\noindent\textit{Step 8: conclusion.}

By the identities established in Step 7 and the representation of the perimeter measure of the tube,
\[
\sqrt{1+v^2}\,\mu = |D_{x,t}\chi_{\widehat E}| = \Ha^N\mres\partial^*\widehat E.
\]
Since $\sigma=\phi(\nu_{\widehat E(t)})\mu$, we obtain
\[
\sigma = \frac{\phi(\nu_{\widehat E(t)})}{\sqrt{1+v^2}}\, \Ha^N\mres\partial^*\widehat E.
\]
Thus $\sigma$ is an $N$-rectifiable measure with positive finite density. At $\sigma$-almost every point $y_0=(x_0,t_0)$, the measure $\sigma$ has tangent $\theta_0\Ha^N\mres P_0$,
where $P_0=\operatorname{Tan}(\partial^*\widehat E,y_0)$ and $\theta_0>0$ and also
\[
\lim_{r\downarrow0} \frac{1}{\sigma(B_r(y_0))} \int_{B_r(y_0)} |\mathbf b(y)-\mathbf b(y_0)|\,d\sigma(y)=0.
\]
In particular, $\sigma(B_r(y_0))/r^N\to\omega_N\theta_0$, where $\omega_N$ is the volume of the unit ball in $\R^N$. Moreover, writing $\Lambda=f\sigma+\Lambda^s$ for its
decomposition with respect to $\sigma$, the same theorem yields $\Lambda(B_r(y_0))/\sigma(B_r(y_0))\to f(y_0)<\infty$ at $\sigma$-almost every point. Consequently,
\[
\frac{\Lambda(B_r(y_0))}{r^N} = \frac{\Lambda(B_r(y_0))}{\sigma(B_r(y_0))} \frac{\sigma(B_r(y_0))}{r^N} \longrightarrow f(y_0)\omega_N\theta_0.
\]
We choose $y_0$ in the common set of full $\sigma$-measure where all these properties hold. For $\eta\in C_c^1(\R^{N+1})$ and $r$ small enough, test  \eqref{bilancioenergia}  with
$r^{1-N}\eta((\,\cdot-y_0)/r)$, obtaining
\[
r^{-N}\int (\mathbf b(y),1)\cdot \nabla_{x,t}\eta\!\left(\frac{y-y_0}{r}\right)\,d\sigma(y) = r^{1-N}\int \eta\!\left(\frac{y-y_0}{r}\right)\,d\Lambda(y).
\]
If $\spt\eta\subset B_R(0)$, the rescaled test is supported in $B_{Rr}(y_0)$. The right-hand side therefore satisfies
\[
\begin{aligned}
\left| r^{1-N}\int \eta\!\left(\frac{y-y_0}{r}\right)\,d\Lambda(y) \right| \leq r^{1-N}\|\eta\|_\infty\,\Lambda(B_{Rr}(y_0))\leq CR^N\|\eta\|_\infty\,r \longrightarrow0.
\end{aligned}
\]
On the left-hand side, if we replace $\mathbf b(y)$ by $\mathbf b(y_0)$ we  introduce an error bounded by
\[
\|\nabla_{x,t}\eta\|_\infty\,r^{-N} \int_{B_{Rr}(y_0)} |\mathbf b(y)-\mathbf b(y_0)|\,d\sigma(y),
\]
which tends to zero by the Lebesgue-point property and the bound $\sigma(B_{Rr}(y_0))=O(r^N)$. By tangent-measure convergence we then obtain
\[
\theta_0\int_{P_0} (\mathbf b(y_0),1)\cdot\nabla_{x,t}\eta\,d\Ha^N=0.
\]
This identity implies $(\mathbf b(y_0),1)\in P_0$. Indeed, if $\mathbf n$ is normal to $P_0$, choose $\eta(y)=(y\cdot \mathbf n)\psi(y)$, where $\psi$ is smooth, compactly
supported, nonnegative, and has positive integral on $P_0$. Since $\nabla_{x,t}\eta=\mathbf n\,\psi$ on $P_0$, we deduce from this identity that $(\mathbf b(y_0),1)\cdot \mathbf n=0$.
Thus $(\mathbf b,1)$ is tangent to the tube a.e. In turn, since the outer unit normal to the tube is $(\nu_{\widehat E(t)},-v)/\sqrt{1+v^2}$, we deduce
\[
\mathbf b\cdot\nu_{\widehat E(t)}=v \qquad\text{for }\mu\text{-a.e.}.
\]
In particular, the $L^2(\sigma)$ bound for $\mathbf b$ established in Step 7 implies that $v\in L^2(\mu)$ as well.

Finally, for a time-dependent test $\eta$, starting from \eqref{restart:l1}, and arguing as in the proof of \eqref{thfund-general}-\eqref{Dt:lt}, we find that
$t\mapsto\mu^\phi_{\widehat E(t)}(\eta(\cdot,t))$ has derivative equal to the time projection of
\[
(\partial_t\eta+\mathbf b\cdot\nabla\eta)\sigma-\eta\Lambda.
\]
In turn, we obtain \eqref{bilanciointervalli}. Analogously, testing $\partial_t\chi_{\widehat E}=v\mu$ with $\zeta(t)\eta(x,t)$, $\zeta\in C_c^1(0,\infty)$, we obtain
\[
D_t\!\left[\int_{\widehat E(t)}\eta(x,t)\,dx\right] = \left[ \int_{\widehat E(t)}\partial_t\eta\,dx + \int\eta v\,d\mu_{\widehat E(t)} \right]dt.
\]
The density on the right is integrable on finite time intervals by the estimates already established. By \eqref{eq:flatflowholder}, the integral is continuous with respect to time.
Thus the fundamental theorem for BV functions yields \eqref{trasportofase}.
\end{proof}

\begin{theorem}[Brakke and canonical BV realization]
\label{brakkebv} Under the assumptions of Theorem~\ref{trasportobv}, suppose additionally that $\phi\in C^1(\R^N\setminus\{0\})$ is strictly convex. Let $\bu$ and $\Lambda$ be the transport field and the dissipation measure, respectively,  provided by that theorem. Consider   the unit-density varifolds
\[
V_t=\var(\partial^*\widehat E(t)), \qquad t\geq0.
\]
Then,  for a.e.  $t>0$, $V_t$ has generalized anisotropic mean curvature, which we will denote by  $\mathbf H^\phi_{\pa^*\widehat E(t)}$, and satisfies
\begin{equation} \label{identifica}
\begin{split}
 \mathbf H^\phi_{\pa^*\widehat E(t)}\perp \pa^*\widehat E(t), \qquad \mathbf b = (\nu_{\widehat E(t)}\cdot\mathbf H^\phi_{\pa^*\widehat E(t)}) \,\nabla\phi(\nu_{\widehat E(t)}),
 \\
 v = \phi(\nu_{\widehat E(t)})
(\nu_{\widehat E(t)}\cdot\mathbf H^\phi_{\pa^*\widehat E(t)}) \quad \mathcal H^{N-1}\text{-a.e. on } \partial^*\widehat E(t).
\end{split}
\end{equation}
In particular,
\begin{equation} \label{curvatura}
\big|\mathbf H^\phi_{\pa^*\widehat E(t)}\big|^2\,\sigma\leq\Lambda.
\end{equation}
As a consequence, the family $\{V_t\}_{t\geq0}$ is an integral (unit-density) $\phi$-anisotropic Brakke flow in the sense of Definition~\ref{def:brakke}; that is,    for every
$0\leq t_1<t_2$  and for every nonnegative  $C^1$ test function $\eta$ with compact spatial support uniformly on $[t_1,t_2]$ one has
\begin{equation} \label{brakke}
\begin{aligned}
&\mu^\phi_{\widehat E(t_2)}(\eta(\cdot,t_2)) -\mu^\phi_{\widehat E(t_1)}(\eta(\cdot,t_1)) \\
&\quad\leq \int_{t_1}^{t_2}\int \left[ \partial_t\eta + (\nu_{\widehat E(t)}\cdot\mathbf H^\phi_{\pa^*\widehat E(t)}) \nabla\phi(\nu_{\widehat E(t)})\cdot\nabla\eta - \eta\big|\mathbf H^\phi_{\pa^*\widehat E(t)}\big|^2 \right]d\mu^\phi_{\widehat E(t)}\,dt .
\end{aligned}\end{equation}
In particular, the optimal energy/dissipation inequality \eqref{daje} holds. Finally, the flow is a canonical BV flow; that is,
\begin{equation} \label{fase}
\begin{split}
&\int_{\widehat E(t_2)}\eta(x,t_2)\,dx -\int_{\widehat E(t_1)}\eta(x,t_1)\,dx \\
&\quad= \int_{t_1}^{t_2}\int_{\widehat E(t)}\partial_t\eta\,dx\,dt + \int_{t_1}^{t_2}\int \eta\, (\nu_{\widehat E(t)}\cdot\mathbf H^\phi_{\pa^*\widehat E(t)}) \,d\mu^\phi_{\widehat E(t)}\,dt
\end{split}
\end{equation}
for any $0\leq t_1<t_2$ and any $\eta$ of class $C^1$ with uniformly bounded spatial support on $[t_1, t_2]$.
\end{theorem}

\begin{proof}
We start by observing that it is enough to identify the field $\bu$ obtained in Theorem~\ref{trasportobv}. More precisely, we will show that, for almost every $t>0$, the varifold
$\var(\partial^*\widehat E(t))$ admits generalized anisotropic mean curvature, that this curvature is normal, and that the formula for $\bu$ in \eqref{identifica} holds. Indeed,
combining the identity $v=\bu\cdot\nu_{\widehat E(t)}$ in \eqref{tubofase} with Euler's identity for $\phi$, we obtain the formula for $v$ in \eqref{identifica}. Moreover, the
identification of $\mathbf{b}$ together with  $\phi^o(\nabla\phi(\nu))=1$ yields
\[
\bigl(\phi^o(\bu)\bigr)^2 = \bigl|\nu_{\widehat E(t)}\cdot\mathbf H^\phi_{V_t}\bigr|^2 = |\mathbf H^\phi_{V_t}|^2,
\]
and thus \eqref{dissipazione} implies  \eqref{curvatura}. In turn, using the latter and  plugging the formula for $\bu$ into \eqref{bilanciointervalli} we obtain  \eqref{brakke}.
Finally, inserting the formula for $v$ into \eqref{trasportofase}, and recalling $d\mu^\phi_{\widehat E(t)}=\phi(\nu_{\widehat E(t)})\,d\mu_{\widehat E(t)}$,  we infer
\eqref{fase}.

In order to obtain \eqref{identifica}, we follow the construction of Theorem~\ref{trasportobv}, taking $\phi_n=\phi$ throughout, as described in Step~1 of its proof. Indeed, since
$\phi$ is already $C^1$ and strictly convex, no approximation of the anisotropy is required. We therefore work directly with the minimizers $G_\e$ associated with the fixed
extension $\Phi$. Let $\e_j\to 0$ be the sequence selected in the proof of Theorem~\ref{trasportobv}. For notational simplicity, in what follows we write $\e$ in place of $\e_j$;
all limits are taken along this sequence, or along further subsequences, which are denoted again by $\e$.

\medskip
\noindent\textit{Step 1: Varifolds associated with the slices of the rescaled minimizers.}

Fix a starting time $t_0\geq0$, set $G_0=\widehat E(t_0)$, and let $G_\e$ be the minimizers used in the proof of Theorem~\ref{trasportobv}. Throughout Steps 1--3 of the current
proof, limits as $\e\to0$ are taken along the subsequence already chosen in Step 4 of the proof of Theorem~\ref{trasportobv}, or along a further subsequence. In particular, the
limits in \eqref{8} and \eqref{limitingJ} are fixed.

For almost every $t>t_0$, the spatial slice of $\mathcal R_\e(G_\e)+(0,t_0)$ has finite perimeter. We denote its unit-density varifold by
\[
V_t^\e := \var\!\left(\partial^*\{x\in\R^N:(x,(t-t_0)/\e)\in G_\e\}\right).
\]
We denote as usual by $\nu'_{G_\e}$ the horizontal component of the outer normal $\nu_{G_\e}$ and set $H^\Phi=H^\Phi_{\partial^*G_\e}$. By the slicing theorem for sets of finite
perimeter, for almost every $z>0$ the reduced boundary of the slice agrees, up to an $\Ha^{N-1}$-null set, with $\{x:(x,z)\in\partial^*G_\e,\ |\nu'_{G_\e}|>0\}$, and its outer
normal is $\nu'_{G_\e}/|\nu'_{G_\e}|$; see \cite[Theorem~18.11]{MaggiBook}. The tangential Jacobian of $(x,z)\mapsto z$ on $\partial^*G_\e$ is $\sqrt{1-(\nu_{G_\e}\cdot
e_{_{N+1}})^2}=|\nu'_{G_\e}|$. Thus by the coarea formula and the change of variable $t=t_0+\e z$ we have
\begin{equation}\label{coarea}
\int_{t_0}^\infty\int f(x,t,\nu)\,dV_t^\e(x,\nu)\,dt = \e\int_{H^+} f\!\left(x,t_0+\e z,\frac{\nu'_{G_\e}}{|\nu'_{G_\e}|}\right)|\nu'_{G_\e}|\,d\mu_{G_\e}
\end{equation}
for every continuous compactly supported $f$ on $\R^N\times(t_0,\infty)\times G(N,N-1)$. Here and below the integrand on the right is understood as equal to zero whenever
$\nu'_{G_\e}=0$, and  as in \eqref{this:is:the:notation}, $\nu$ denotes  the projective normal. Thus, taking $f(x,t,\nu)=\eta(x,t)\phi(\nu)$, using one-homogeneity, and recalling
\eqref{e:an_en_var}, we obtain
\begin{equation}\label{gu}
\int_{t_0}^\infty\int\eta\,d\|V_t^\e\|^\phi\,dt = \e\int_{H^+}\eta(x,t_0+\e z)\phi(\nu'_{G_\e})\,d\mu_{G_\e}.
\end{equation}
We now claim that this measure and the first rescaled energy measure in \eqref{8} have the same limit. To this aim, observe that by \eqref{Phi} and \eqref{solitonbis},
\[
0\leq\Phi(\nu_{G_\e})-\phi(\nu'_{G_\e}) = \frac{(\nu_{G_\e}\cdot e_{_{N+1}})^2}{\Phi(\nu_{G_\e})+\phi(\nu'_{G_\e})} \leq \frac{(\nu_{G_\e}\cdot e_{_{N+1}})^2}{\Phi(\nu_{G_\e})} = \e^2|H^\Phi|^2\Phi(\nu_{G_\e}).
\]
Multiplying by $\e\eta(x,t_0+\e z)$ and integrating, we infer from \eqref{1} that
\begin{equation}\label{scarto}
0\leq \e\int_{H^+}\eta(x,t_0+\e z)\,d\mu_{G_\e}^\Phi - \int_{t_0}^\infty\int\eta\,d\|V_t^\e\|^\phi\,dt \leq \e^2\|\eta\|_\infty P^\phi(G_0)
\end{equation}
for every nonnegative $\eta\in C_c(\R^N\times(t_0,\infty))$. By \eqref{weight:comparison}, $\phi(\nu'_{G_\e})\leq\Phi(\nu_{G_\e})$, and \eqref{questaserve}, for every bounded open
interval $I\Subset(t_0,\infty)$ we also have
\[
\int_I\|V_t^\e\|(\R^N)\,dt \leq C_\phi\bigl(|I|+\e^2\bigr)P^\phi(G_0).
\]
Thus, up to a further subsequence, $V_t^\e\,dt\weakstar\mathcal V$ on $\R^N\times(t_0,\infty)\times G(N,N-1)$. Passing to the limit in the last estimate yields
\[
\mathcal V\bigl(\R^N\times I\times G(N,N-1)\bigr) \leq C_\phi|I|P^\phi(G_0).
\]
Thus the time-marginal of $\mathcal V$ is absolutely continuous with respect to Lebesgue measure, with bounded density. We may therefore disintegrate it into  measurable family of
varifolds $V_t$ such that $d\mathcal V=dV_t(x,\nu)\,dt$.

Since $\phi$ is continuous on the Grassmannian, \eqref{scarto} and \eqref{8} lead to
\[
\int_{t_0}^\infty\int\eta\,d\|V_t\|^\phi\,dt = \int_{t_0}^\infty\int\eta\,d\lambda_t\,dt.
\]
Hence,  by the arbitrariness of $\eta$ we conclude that
\begin{equation}\label{pesoblocco}
\|V_t\|^\phi=\lambda_t \qquad\text{for a.e. }t>t_0.
\end{equation}

\medskip
\noindent\textit{Step 2: The approximate slice curvature vector measure.}

We now introduce the $\R^N$-valued Radon measure $\mathcal J_\e$ on $\R^N\times(t_0,\infty)\times G(N,N-1)$ which may be interpreted as encoding information on the ``approximate''
mean curvature vectors of the slices of  $\mathcal R_\e(G_\e)+(0,t_0)$. We define it as
\begin{equation}\label{mathcalJe}
\int Y(x,t,\nu)\cdot d\mathcal J_\e(x,t,\nu) := -\e\int_{H^+}H^\Phi\frac{\nu'_{G_\e}}{|\nu'_{G_\e}|} \cdot Y\!\left(x,t_0+\e z,\frac{\nu'_{G_\e}}{|\nu'_{G_\e}|}\right)\phi(\nu'_{G_\e})\,d\mu_{G_\e}
\end{equation}
for continuous compactly supported vector fields $Y$ on this product space (above, in the evaluation of $Y$, the vector $\nu'_{G_\e}/|\nu'_{G_\e}|$ is identified with the
corresponding element of $G(N,N-1)$ as described in Section \ref{s:var:def}).

By \eqref{coarea}, its density with respect to $\phi(\nu)V_t^\e\,dt$ is $-H^\Phi\nu'_{G_\e}/|\nu'_{G_\e}|$, where $H^\Phi$ is evaluated at $(x,(t-t_0)/\e)$. By Cauchy-Schwarz and
\eqref{1} we have
\[
\left|\int Y\cdot d\mathcal J_\e\right| \leq P^\phi(G_0)^{1/2}\left(\int_{t_0}^\infty\int|Y|^2\phi(\nu)\,dV_t^\e\,dt\right)^{1/2}.
\]
Thus, recalling the mass bound in Step 1, we may assume,  along a further subsequence, $\mathcal J_\e\weakstar\mathcal J$ for an $\R^N$-valued Radon measure $\mathcal J$ on
$\R^N\times(t_0,\infty)\times G(N,N-1)$ . Passing to the limit in the last inequality, using the convergence of $V_t^\e\,dt$ and disintegration, we obtain
\[
\left|\int Y\cdot d\mathcal J\right| \leq P^\phi(G_0)^{1/2}\left(\int_{t_0}^\infty\int|Y|^2\phi(\nu)\,dV_t\,dt\right)^{1/2}.
\]
In particular, as in Step 4 of the proof of Theorem~\ref{trasportobv}, by the Riesz representation theorem, we have
\begin{equation}\label{limitecurva}
\mathcal J = \mathbf h\,\phi(\nu)V_t\,dt, \qquad \text{with }\mathbf h\in L^2\big(\phi(\nu)V_t\,dt;\R^N\big).
\end{equation}

We now  establish a quadratic estimate on $\mathbf h$. To this aim, integrating the elementary inequality $-2H^\Phi(\nu'_{G_\e}/|\nu'_{G_\e}|)\cdot Y-|Y|^2\leq|H^\Phi|^2$ against
$\e\eta(x,t_0+\e z)\phi(\nu'_{G_\e})\,d\mu_{G_\e}$, and recalling \eqref{gu}, we get
\[
2\int\eta Y\cdot d\mathcal J_\e - \int_{t_0}^\infty\int\eta|Y|^2\phi(\nu)\,dV_t^\e\,dt \leq \e\int_{H^+}\eta(x,t_0+\e z)|H^\Phi|^2\,d\mu_{G_\e}^\Phi.
\]
Passing to the limit and recalling \eqref{8}, we obtain
\[
\int_{t_0}^\infty\int\eta\bigl(2\mathbf h\cdot Y-|Y|^2\bigr)\phi(\nu)\,dV_t\,dt \leq \int\eta\,d\Lambda
\]
and taking the supremum over  the admissible fields $Y$ yields
\begin{equation}\label{disscurva}
\int_{t_0}^\infty\int\eta(x,t)|\mathbf h(x,t,\nu)|^2\phi(\nu)\,dV_t(x,\nu)\,dt \leq \int\eta\,d\Lambda
\end{equation}
for any nonnegative $\eta\in C(\R^N\times [t_0, \infty))$. Observe now that the vectors defining $\mathcal J_\e$ are parallel to the normal of the slices. Therefore, denoting by
$\mathrm{Id}_N$ the identity matrix on $\R^N$, by testing \eqref{mathcalJe} with $(\mathrm{Id}_N-\nu\otimes\nu)\,Y$ in place of $Y$, we obtain zero, and consequently, passing to
the limit, we get
\[
\int_{t_0}^\infty\int (\mathrm{Id}_N-\nu\otimes\nu)\mathbf h\cdot Y\,\phi(\nu)\,dV_t\,dt=0
\]
for any continuous compactly supported vector field $Y$ on $\R^N\times(0,\infty)\times G(N,N-1)$. This shows that
\begin{equation}\label{normalita}
(\mathrm{Id}_N-\nu\otimes\nu)\,\mathbf h(x,t,\nu)=0 \qquad V_t\,dt\text{-a.e.. }
\end{equation}
It remains to relate $\mathbf h$ to the transport field $\bu$ constructed in Theorem~\ref{trasportobv}. Differentiating \eqref{Phi} in the horizontal variables and using
zero-homogeneity of $\nabla\phi$, we obtain
\begin{equation}\label{Phiphi}
\Phi(\nu_{G_\e})\nabla_\xi\Phi(\nu_{G_\e}) = \phi(\nu'_{G_\e})\nabla\phi\!\left(\frac{\nu'_{G_\e}}{|\nu'_{G_\e}|}\right) \qquad\text{on }\{\nu'_{G_\e}\neq0\}.
\end{equation}
Both sides vanish on $\{\nu'_{G_\e}=0\}$. Thus, for any $X\in C_c(\R^N\times(t_0,\infty);\R^N)$, the vector measure $J_\e$ defined as in \eqref{defJj} (with $\e$, $\Phi$, and
$G_\e$ in place of $\e_j$, $\Phi_{n_j}$, and $G_{\e_j, n_j}$, respectively) satisfies
\[
\int X\cdot dJ_\e = \int \bigl(\nabla\phi(\nu)\cdot X(x,t)\bigr)\,\nu\cdot d\mathcal J_{\e}(x,t,\nu),
\]
so that, passing to the limit and recalling that $J_\e\weakstar \bu \lambda_t dt$ (along the extracted subsequence),
\begin{equation}\label{flussocurva}
\int_{t_0}^\infty\int\bu\cdot X\,d\lambda_t\,dt = \int_{t_0}^\infty\int (\mathbf h\cdot\nu)\nabla\phi(\nu)\cdot X\,\phi(\nu)\,dV_t\,dt.
\end{equation}

\medskip
\noindent\textit{Step 3: Horizontal first variation.}

We now test the first variation on horizontal fields. Let $X\in C_c^1(\R^N\times(t_0,\infty);\R^N)$. We apply \eqref{euler:lagrange0}, with $g=0$, to $Y(x,z)=(X(x,t_0+\e z),0)$. By
\eqref{def:Bpsi} and \eqref{prop:extnorm},
\[
B_\Phi(\nu_{G_\e}):DY = \bigl[\Phi(\nu_{G_\e})\mathrm{Id}_{N}-\nu'_{G_\e}\otimes\nabla_\xi\Phi(\nu_{G_\e})\bigr]:D_xX - \e\frac{\nu_{G_\e}\cdot e_{_{N+1}}}{\Phi(\nu_{G_\e})}\nu'_{G_\e}\cdot\partial_tX,
\]
where $X$ and its derivatives are evaluated at $(x,t_0+\e z)$. Since $(\nu_{G_\e}\cdot e_{_{N+1}})/\Phi(\nu_{G_\e})=\e H^\Phi$ by \eqref{solitonbis}, multiplying the
Euler--Lagrange equation by $\e$ yields
\begin{equation}\label{varorizz}
\begin{split}
&\e\int_{H^+}\bigl[\Phi(\nu_{G_\e})\,\mathrm{Id}_{N}-\nu'_{G_\e}\otimes\nabla_\xi\Phi(\nu_{G_\e})\bigr]:D_xX\,d\mu_{G_\e} - \e^3\int_{H^+}H^\Phi\nu'_{G_\e}\cdot\partial_tX\,d\mu_{G_\e}
= \e\int_{H^+}H^\Phi\nu'_{G_\e}\cdot X\,d\mu_{G_\e}.
\end{split}
\end{equation}

We start by comparing the matrix in the first integral with the one appearing in the first variation of the slice varifolds. By one-homogeneity and \eqref{def:Bpsi}, on
$\{\nu'_{G_\e}\neq0\}$ we have
\[
|\nu'_{G_\e}|B_\phi\!\left(\frac{\nu'_{G_\e}}{|\nu'_{G_\e}|}\right) = \phi(\nu'_{G_\e})\,\mathrm{Id}_{N}-\nu'_{G_\e}\otimes\nabla\phi\!\left(\frac{\nu'_{G_\e}}{|\nu'_{G_\e}|}\right)
\]
and recalling \eqref{Phiphi}, we can write
\[
\Phi(\nu_{G_\e})\,\mathrm{Id}_{N}-\nu'_{G_\e}\otimes\nabla_\xi\Phi(\nu_{G_\e}) -|\nu'_{G_\e}|B_\phi\!\left(\frac{\nu'_{G_\e}}{|\nu'_{G_\e}|}\right) = \bigl(\Phi(\nu_{G_\e})-\phi(\nu'_{G_\e})\bigr)
\left[\mathrm{Id}_{N}+\frac{\nu'_{G_\e}}{\Phi(\nu_{G_\e})}\otimes\nabla\phi\!\left(\frac{\nu'_{G_\e}}{|\nu'_{G_\e}|}\right)\right].
\]
Observe that, since $|\nu'|/\Phi(\nu)\leq |\nu'|/\phi(\nu')\leq C_\phi$ by \eqref{Phi} and the equivalence of norms, the matrix inside the square brackets is bounded by a constant
depending only on $\phi$. Thus, by \eqref{coarea} and the estimate leading to \eqref{scarto}, if we replace the first integral in \eqref{varorizz} by the slice first variation
\begin{equation}\label{slice1stvar}
\int_{t_0}^\infty\int B_\phi(\nu):D_xX\,dV_t^\e\,dt
\end{equation}
 we introduce an error bounded by $C_\phi\,\e^2\,P^\phi(G_0)\|D_xX\|_\infty$, which clearly tends to zero. On the other hand, also the second integral in \eqref{varorizz} also
 tends to zero. Indeed, if  $\spt X\subset\R^N\times(t_0,T)$, by  Cauchy-Schwarz and \eqref{1}, we have
\[
\begin{split}
\left|\e^3\int_{H^+}H^\Phi\nu'_{G_\e}\cdot\partial_tX\,d\mu_{G_\e}\right| &\leq C_\phi\e^2\|\partial_tX\|_\infty \left(\e\int_{H^+}|H^\Phi|^2\,d\mu_{G_\e}^\Phi\right)^{1/2}
\left(\e\mu_{G_\e}^\Phi\!\left(\R^N\times\left(0,\frac{T-t_0}{\e}\right)\right)\right)^{1/2}\\
&\leq C_\phi\,\e^2\|\partial_tX\|_\infty P^\phi(G_0)\sqrt{T-t_0+\e^2}.
\end{split}
\]
Finally, by the definition of $\mathcal J_\e$ in \eqref{mathcalJe} and $\phi(\nu'_{G_\e})=|\nu'_{G_\e}|\phi(\nu'_{G_\e}/|\nu'_{G_\e}|)$, the right-hand side of \eqref{varorizz} can
be written as
\[
-\int \frac{X(x,t)}{\phi(\nu)}\cdot d\mathcal J_\e\longrightarrow -\int\frac{X(x,t)}{\phi(\nu)}\cdot d\mathcal J=-\int_{t_0}^\infty\int\mathbf h\cdot X\,dV_t\,dt
\]
as $\e\to 0$, where in the last equality we have used  \eqref{limitecurva}. Since $B_\phi$ is continuous on the Grassmannian, we can also pass to the limit  also the slice first
variation \eqref{slice1stvar} and by all the above observations,  we arrive at
\begin{equation}\label{limvar}
\int_{t_0}^\infty\int B_\phi(\nu):D_xX\,dV_t\,dt = -\int_{t_0}^\infty\int\mathbf h\cdot X\,dV_t\,dt
\end{equation}
for every $X\in C_c^1(\R^N\times(t_0,\infty);\R^N)$.

\medskip
\noindent\textit{Step 4: Restarting and matching.}

In this step we carry the varifolds and the fields constructed above through the restarting procedure of Theorem~\ref{trasportobv}.

Fix $\delta>0$, and recall the times $0=t_0^\delta<\cdots<t_{N_\delta}^\delta<t_{N_\delta+1}^\delta=\infty$ and the intervals $I_i^\delta=[t_i^\delta,t_{i+1}^\delta)$ constructed
in Step 6 of the proof of Theorem~\ref{trasportobv}. For each $i$, the construction starting from $\widehat E(t_i^\delta)$ provides the measures $\lambda_t^{(i)}\,dt$ and
$\Lambda^{(i)}$ and the transport vector $\bu^{(i)}\lambda_t^{(i)}\,dt$, by passing to the limit along a sequence of rescaled minimizers with that initial datum. We may also apply
Steps 1--3 of the current proof to this same sequence, with $t_0=t_i^\delta$ and $G_0=\widehat E(t_i^\delta)$ and, possibly  passing to further subsequences, obtain the additional
limits $V_t^{(i)}$ and $\mathbf h^{(i)}$.

We can argue in this way on every interval, since their number is finite. Note that the subsequence of minimizers used on one interval may depend on $i$ and  the limit as $\e\to0$
is taken separately for each re-starting datum. We then join the resulting limiting families. More precisely, for $t_i^\delta<t<t_{i+1}^\delta$ set
\[
V_t^\delta:=V_t^{(i)}, \qquad \mathbf h^\delta(x,t,\nu):=\mathbf h^{(i)}(x,t,\nu).
\]
Values at the finitely many endpoints are immaterial for the measures obtained by integration with respect to $dt$. By \eqref{pesoblocco}, and recalling the definition of
$\lambda_t^\delta$ given in the proof of Theorem~\ref{trasportobv}, we have
\begin{equation}\label{pesobloccodelta}
\|V_t^\delta\|^\phi=\lambda_t^\delta \quad\text{for almost every $t>0$.}
\end{equation}
We now show that \eqref{flussocurva} and \eqref{limvar} remain valid after joining the blocks. Fix $X\in C_c^1(\R^N\times(0,\infty);\R^N)$. For each $i$, we may choose $\zeta_k\in
C_c^\infty((t_i^\delta,t_{i+1}^\delta))$, with $0\leq\zeta_k\leq1$ and $\zeta_k\to1$  in that interval, and apply \eqref{limvar} for the $i$-th block to $\zeta_kX$. Since
$D_x(\zeta_kX)=\zeta_kD_xX$, we can let $\zeta_k\to1$ and use the  mass bound on $dV_t^{(i)}\,dt$,  the quadratic estimate \eqref{disscurva},  and  dominated convergence to obtain
\[
\int_{t_i^\delta}^{t_{i+1}^\delta}\int B_\phi(\nu):D_xX\,dV_t^{(i)}\,dt = -\int_{t_i^\delta}^{t_{i+1}^\delta}\int\mathbf h^{(i)}\cdot X\,dV_t^{(i)}\,dt.
\]
This shows that the above identity now holds even when $X$ does not vanish near the endpoints of the interval. Summing over $i$ and recalling the definitions of $V_t^\delta$ and
$\mathbf h^\delta$, we obtain
\begin{equation}\label{1stvardelta}
\int_0^\infty\int B_\phi(\nu):D_xX\,dV_t^\delta\,dt = -\int_0^\infty\int\mathbf h^\delta\cdot X\,dV_t^\delta\,dt.
\end{equation}
 Applying a similar argument  to \eqref{flussocurva} (and  using also \eqref{10}), we also get
\begin{equation}\label{id-b-delta}
\int_0^\infty\int\bu^\delta\cdot X\,d\lambda_t^\delta\,dt = \int_0^\infty\int(\mathbf h^\delta\cdot\nu)\nabla\phi(\nu)\cdot X\,\phi(\nu)\,dV_t^\delta\,dt.
\end{equation}
 Finally, for a nonnegative $\eta\in C_c(\R^N\times(0,\infty))$, applying \eqref{disscurva} on the open intervals, after removing time cutoff functions as above, we can write
\begin{equation}\label{raccordivar}
\int_0^\infty\int\eta|\mathbf h^\delta|^2\phi(\nu)\,dV_t^\delta\,dt \leq \sum_{i=0}^{N_\delta} \int_{\R^N\times(t_i^\delta,t_{i+1}^\delta)}\eta\,d\Lambda^{(i)} \leq \int\eta\,d\Lambda^\delta,
\end{equation}
where the last inequality follows from \eqref{16}, since the additional measures at the restarting times are nonnegative.

Now let $\delta\to 0$ and recall that Step 7 of the proof of Theorem~\ref{trasportobv} provides a sequence along which
\begin{equation}\label{limitozzi}
\lambda_t^\delta\,dt\longrightarrow\sigma \quad\text{in total variation on }\R^N\times(0,T)\text{ for every }T<\infty, \qquad \bu^\delta\lambda_t^\delta\,dt\weakstar\bu\sigma, \qquad \Lambda^\delta\weakstar\Lambda.
\end{equation}
The first convergence is precisely \eqref{20}. In what follows we take further subsequences of this sequence.

Recall that on the $i$-th interval, the total mass of $\lambda_t^\delta$ is bounded by $P^\phi(\widehat E(t_i^\delta))\leq P^\phi(\widehat E(0))$. Hence, by
\eqref{weight:comparison} and \eqref{pesobloccodelta} we have  for every bounded time interval $I\subset(0,\infty)$
\[
\int_I\|V_t^\delta\|(\R^N)\,dt \leq C_\phi\int_I\lambda_t^\delta(\R^N)\,dt \leq C_\phi |I|P^\phi(\widehat E(0)).
\]
Thus, we may assume that along a further subsequence the measures $V_t^\delta\,dt\weakstar \mathcal V$ on $\R^N\times(0,\infty)\times G(N,N-1)$, with $\mathcal V(\R^N\times I\times
G(N,N-1))\leq C_\phi |I|P^\phi(\widehat E(0))$. Thus, exactly as in Step 1, the time marginal of $\mathcal V$ is absolutely continuous with respect to Lebesgue measure and we may
disintegrate it  as  $d\mathcal V=dV_t(x,\nu)\,dt$.

We now consider the vector measures $\mathbf h^\delta\phi(\nu)V_t^\delta\,dt$. By \eqref{raccordivar} and \eqref{18}, for every continuous compactly supported vector field
$Y(x,t,\nu)$ we have
\[
\left|\int_0^\infty\int Y\cdot\mathbf h^\delta\phi(\nu)\,dV_t^\delta\,dt\right| \leq P^\phi(\widehat E(0))^{1/2} \left(\int_0^\infty\int |Y|^2\phi(\nu)\,dV_t^\delta\,dt\right)^{1/2}.
\]
Thus,  after a further extraction, we can also assume $\mathbf h^\delta\phi(\nu)V_t^\delta\,dt\weakstar \mathcal J$, with
\[
\left|\int Y\cdot d\mathcal J\right| \leq P^\phi(\widehat E(0))^{1/2} \left(\int_0^\infty\int |Y|^2\phi(\nu)\,dV_t\,dt\right)^{1/2},
\]
and by the same Riesz representation argument of Step 2 we obtain the existence of $\mathbf h\in L^2(V_t\,dt;\R^N)$ such that $\mathcal J=\mathbf h\,\phi(\nu)V_t\,dt$. Summarizing,
we have  obtained
\[
V_t^\delta\,dt\weakstar V_t\,dt, \qquad \mathbf h^\delta\phi(\nu)V_t^\delta\,dt\weakstar\mathbf h\phi(\nu)V_t\,dt.
\]
We next identify the anisotropic mass of $V_t$. For every $\eta\in C_c(\R^N\times(0,\infty))$, \eqref{pesobloccodelta} and the first limit in \eqref{limitozzi} yield
\[
\int_0^\infty\int\eta(x,t)\phi(\nu)\,dV_t\,dt = \lim_{\delta\downarrow0}\int_0^\infty\int\eta(x,t)\phi(\nu)\,dV_t^\delta\,dt = \lim_{\delta\downarrow0}\int_0^\infty\int\eta\,d\lambda_t^\delta\,dt = \int\eta\,d\sigma.
\]
Thus $\|V_t\|^\phi\,dt=\sigma=\mu^\phi_{\widehat E(t)}\,dt$, and  in turn
\begin{equation}\label{massalimite}
\|V_t\|^\phi=\mu^\phi_{\widehat E(t)}\quad\text{ for almost every $t>0$.}
\end{equation}
Moreover, for  $\eta\in C_c(\R^N\times(0,\infty))$  nonnegative and for $Y\in C_c\big(\R^N\times(0,\infty)\times G(N,N-1) {;\R^N}\big)$, we obtain
\[
\int_0^\infty\int\eta \bigl(2\mathbf h\cdot Y-|Y|^2\bigr)\phi(\nu)\,dV_t\,dt\leq\int\eta\,d\Lambda
\]
by passing to the limit in the corresponding inequality for $\delta>0$, which in turn is a consequence of \eqref{raccordivar} and Young's inequality. Taking the supremum over all
such fields $Y$ as in Step 2, we have
\begin{equation}\label{pesolimite}
\int_0^\infty\int\eta|\mathbf h|^2\phi(\nu)\,dV_t\,dt \leq\int\eta\,d\Lambda \quad\text{for all $\eta\in C_c(\R^N\times(0,\infty);\R_+)$}.
\end{equation}
It remains to pass to the limit in the identities obtained before. By \eqref{normalita}, for any continuous compactly supported vector field $Y$ on $\R^N\times(0,\infty)\times
G(N,N-1)$ we have $\int_0^\infty\int\mathbf h^\delta\cdot(\mathrm{Id}_N-\nu\otimes\nu)Y\,\phi(\nu)\,dV_t^\delta\,dt=0$. Since the field $(\mathrm{Id}_N-\nu\otimes\nu)Y$ is
continuous and compactly supported on the same product space, by weak* convergence of {the vector measures $\mathbf h^\delta\phi(\nu)V_t^\delta\,dt$ to $\mathbf h\phi(\nu)V_t\,dt$}
we get {$\int_0^\infty\int\mathbf h\cdot(\mathrm{Id}_N-\nu\otimes\nu)Y\,\phi(\nu)\,dV_t\,dt=0$}, and the arbitrariness of $Y$ yields
\begin{equation}\label{matchingnormalita}
(\mathrm{Id}_N-\nu\otimes\nu)\mathbf h(x,t,\nu)=0 \quad \text{ for $V_t\,dt$-a.e. $(x,t,\nu)$.}
\end{equation}
Finally, using the convergences in \eqref{limitozzi}, we may pass to limit as $\delta\to 0$ in \eqref{1stvardelta} and \eqref{id-b-delta}, to obtain that for any $X\in
C_c^1(\R^N\times(0,\infty);\R^N)$
\begin{equation}\label{primavar}
\int_0^\infty\int B_\phi(\nu):D_xX\,dV_t\,dt = -\int_0^\infty\int\mathbf h\cdot X\,dV_t\,dt
\end{equation}
and
\begin{equation}\label{id-b}
\int\bu\cdot X\,d\sigma = \int_0^\infty\int(\mathbf h\cdot\nu)\nabla\phi(\nu)\cdot X\,\phi(\nu)\,dV_t\,dt.
\end{equation}

\medskip
\noindent\textit{Step 5: Rectifiability and identification of the curvature.}

We start by observing that \eqref{pesolimite}, applied to nonnegative cutoff functions increasing to one, implies
\[
\int_0^\infty\int|\mathbf h(x,t,\nu)|^2\phi(\nu)\,dV_t(x,\nu)\,dt \leq\Lambda(\R^N\times(0,\infty))\leq P^\phi(\widehat E(0)).
\]
Thus,  $\int|\mathbf h(x,t,\nu)|^2\,dV_t(x,\nu)<\infty$ for a.e.  $t>0$. Moreover, by \eqref{massalimite}  and the properties of the flat flow $\|V_t\|(\R^N)<\infty$ for a.e.
$t>0$.

Testing \eqref{primavar} with  $X(x,t)=\zeta(t)X_0(x)$, where $\zeta\in C_c^1(0,\infty)$ and $X_0\in C_c^1(\R^N;\R^N)$ and using the arbitrariness of $\zeta$ (and recalling
\eqref{E:first-variation-Ani0}), we obtain
\begin{equation}\label{varfette}
\delta_\phi V_t(X_0)=-\int\mathbf h(x,t,\nu)\cdot X_0(x)\,dV_t(x,\nu) \qquad\text{for a.e.  $t>0$.}
\end{equation}
Let $\mathcal N_0\subset(0,\infty)$ be a negligible set such that, for every $t\notin\mathcal N_0$, \eqref{massalimite} holds,
\begin{equation}\label{tempibelli}
\|V_t\|(\R^N)<\infty, \quad \int|\mathbf h(x,t,\nu)|^2\,dV_t(x,\nu)<\infty,\quad\text{and}\quad (\mathrm{Id}_N-\nu\otimes\nu)\mathbf h(x,t,\nu)=0 \text{ for $V_t$-a.e. $(x,\nu)$.}
\end{equation}

We now enlarge $\mathcal N_0$ so that \eqref{varfette} holds for every test field at every  time outside  $\mathcal N_0$. To this aim, choose a countable family
$\{X_k\}_{k\in\mathbb N}\subset C_c^1(\R^N;\R^N)$ with the  property that every $X_0\in C_c^1(\R^N;\R^N)$ is the $C^1$ limit of a sequence of members of $\{X_k\}$ whose supports
lie in a common compact set. For each $k$, let $\mathcal N_k$ be a negligible set outside which \eqref{varfette} holds with $X_0=X_k$, and set $\mathcal N:=\mathcal
N_0\cup\bigcup_{k=1}^{\infty}\mathcal N_k$. Then $\mathcal N$ is a negligible set with the property that  for every $t\notin\mathcal N$, identity \eqref{varfette} holds for all
$X_0\in C_c^1(\R^N;\R^N)$, together with~\eqref{tempibelli}.

Fix a time outside this exceptional set. For $\spt\, X_0\subset K$, with $K\subset\R^N$ compact, by \eqref{varfette} we have
\[
|\delta_\phi V_t(X_0)| \leq\|X_0\|_\infty\|V_t\|(K)^{1/2} \left(\int_{K\times G(N,N-1)}|\mathbf h(x,t,\nu)|^2\,dV_t(x,\nu)\right)^{1/2}<+\infty.
\]
Moreover, \eqref{massalimite} and \eqref{weight:comparison} yield
\[
c_\phi\mu_{\widehat E(t)}\leq\|V_t\|\leq C_\phi\mu_{\widehat E(t)}.
\]
Thus, since $\mu_{\widehat E(t)}$ has $(N-1)$-dimensional density one at $\mu_{\widehat E(t)}$-almost every point, we infer  that $\|V_t\|$ has positive lower $(N-1)$-dimensional
density at $\|V_t\|$-almost every point. We may therefore apply the rectifiability theorem proved in \cite{DPDRG}, see also \cite{rabasa:etal,DRK}. Indeed, since $\phi$ is $C^1$
away from the origin and strictly convex, the associated integrand on $G(N,N-1)$ satisfies the atomic condition by \cite[Theorem~1.3]{DPDRG}. The locally bounded anisotropic first
variation and the positive lower density established above allow us to apply \cite[Theorem~1.2]{DPDRG}, showing that $V_t$ is rectifiable. In turn, appealing again to
\eqref{massalimite}, we get
\begin{equation}\label{densitauno}
V_t=\mu_{\widehat E(t)}\otimes\delta_{\operatorname{Tan}(\partial^*\widehat E(t),x)}\quad\text{for all $t\not\in \mathcal N$.}
\end{equation}
Thus, for all such times, defining $\mathbf H^\phi_{\pa^*\widehat E(t)}(x):=\mathbf h\big(x,t,[\nu_{\widehat E(t)}(x)]\big)$, identity \eqref{varfette} becomes
\[
\delta_\phi V_t(X_0)=-\int X_0\cdot\mathbf H^\phi_{\pa^*\widehat E(t)}\,d\mu_{\widehat E(t)} \qquad\text{for every }X_0\in C_c^1(\R^N;\R^N);
\]
that is,  $\mathbf H^\phi_{\pa^*\widehat E(t)}$ is the generalized anisotropic mean curvature of $V_t$. Moreover, by the third equality in \eqref{tempibelli} we obtain
$(\mathrm{Id}_N-\nu_{\widehat E(t)}\otimes\nu_{\widehat E(t)})\mathbf H^\phi_{\pa^*\widehat E(t)}=0$ at $\mu_{\widehat E(t)}$-a.e. point, as stated in \eqref{identifica}.

Finally, by using \eqref{densitauno} in \eqref{id-b}, we obtain
\[
\int\bu\cdot X\,d\sigma =\int(\mathbf H^\phi_{\pa^*\widehat E(t)}\cdot\nu_{\widehat E(t)}) \nabla\phi(\nu_{\widehat E(t)})\cdot X\,d\sigma
\]
for every $X\in C_c^1(\R^N\times(0,\infty);\R^N)$. By the arbitrariness of $X$, this proves the formula for $\bu$ in \eqref{identifica} $\sigma$-a.e.  As explained at the beginning
of the proof (and by  re-defining $V_t=\var(\partial^*\widehat E(t))$ also for  $t\in \mathcal N$ and $t=0$), this concludes the proof of the theorem.
\end{proof}

\section{Elliptic regularization with general mobility} \label{sec:gen:mob}

In this section we are going to adapt the previous proofs to the elliptic regularization scheme with arbitrary mobility. We will only expand the parts that significantly differ
from the case of natural mobility, while we refer to the previous sections for all the other parts.

\subsection{Kuratowski convergence}
Throughout this section, we fix two (arbitrary) norms $\phi$ and $\mathfrak m$ on $\R^N$ and we consider the corresponding families of extended norms
$\{\Phi^\kappa\}_{\kappa\in(0,\kappa_0)}$ and $\{\mathfrak M^\kappa\}_{\kappa\in(0,\kappa_0)}$ provided by Proposition~\ref{prop:extended-norms}. In particular, for every
$\kappa\in(0,\kappa_0)$, the norms $\Phi^\kappa,\mathfrak M^\kappa\in C^\infty(\R^{N+1}\setminus\{0\})$ are uniformly convex, even with respect to the last variable, and satisfy
\begin{equation}\label{ijr}
\nabla\Phi^\kappa(\xi,z)\cdot e_{_{N+1}}= \partial_z\Phi^\kappa(\xi,z) = \frac{z}{\mathfrak M^\kappa(\xi,z)} \qquad \text{for every }(\xi,z)\in\R^{N+1}\setminus\{0\}.
\end{equation}
We also recall that, by the standard duality properties of smooth uniformly convex norms, $(\Phi^\kappa)^o$ and $(\mathfrak M^\kappa)^o$ belong to
$C^\infty(\R^{N+1}\setminus\{0\})$ and are uniformly convex as well-- see Lemma \ref{lem:varie:unifconv}.

We first establish the analogue of Proposition~\ref{prop:superflow} for the elliptic regularization associated with the pair $(\Phi^\kappa,\mathfrak M^\kappa)$.

Fix $\kappa\in(0,\kappa_0)$ and $\e>0$, and let
\begin{equation}\label{mob:minpb}
E^\kappa_\e\in\operatorname*{argmin}\{\mathcal A^{\Phi^\kappa}_\e(F):F\subset\R^{N+1},\;F\cap H^-=E_0\times(-\infty,0)\}.
\end{equation}
As in the natural-mobility case, we choose the representative of $E^\kappa_\e$ to be relatively closed in $\{z>0\}$--see Remark \ref{rem:Lmin}.

We define the corresponding translative solitons by
\begin{equation}\label{Ssolliton}
\mathcal S^\kappa_\e(t):=\overline{E^\kappa_\e}-\frac{t}{\e}e_{_{N+1}}\quad t\geq 0,
\end{equation}
and the associated space-time tube
\[
\mathcal S^\kappa_\e := \left\{ (x,z,t)\in\R^{N+1}\times[0,+\infty): (x,z)\in\mathcal S^\kappa_\e(t) \right\}.
\]
By the Euler--Lagrange equation \eqref{euler:lagrange1} for $E^\kappa_\e$, \eqref{ijr} and translation invariance, the solitons satisfy
\begin{equation}\label{mob:EL-soliton}
H^{\Phi^\kappa}_{\partial^*\mathcal S^\kappa_\e(t)}
=\frac1\e\partial_z\Phi^\kappa(\nu)
=\frac{\nu\cdot e_{_{N+1}}}{\e\,\mathfrak M^\kappa(\nu)}
\quad\text{on $\partial^*\mathcal S^\kappa_\e(t)\cap\{z>-t/\e\}$}.
\end{equation}

We introduce the $(\mathfrak M^\kappa)^o$-distance
\begin{equation}\label{d:Mk}
d^\kappa_\e(x,z,t) := \operatorname{dist}^{(\mathfrak M^\kappa)^o} \big((x,z),\mathcal S^\kappa_\e(t)\big)
\end{equation}
and, with $\mathcal H_s:=\{(x,z)\in\R^{N+1}:z=s\}$, set
\[
\mathcal U^\kappa_\e := \left\{ (x,z,t)\in\R^{N+1}\times(0,+\infty):\, z>-\dfrac{t}{\e},\, 0<d^\kappa_\e(x,z,t)<\operatorname{dist}^{(\mathfrak M^\kappa)^o}\big((x,z),\mathcal H_{-t/\e}\big) \right\}.
\]
As usual for distance functions associated with a smooth strictly convex norm, $\mathfrak M^\kappa(\nabla d^\kappa_\e)=1$ a.e. in $\mathcal U^\kappa_\e$--see \eqref{Eikonal}.

\begin{proposition}[Translative solitons are super-flows]\label{mob:prop:superflow}
For every fixed $\kappa\in(0,\kappa_0)$ and $\e>0$, define $\mathcal z^\kappa_\e:=\nabla\Phi^\kappa(\nabla d^\kappa_\e)$ a.e. in $\mathcal U^\kappa_\e$. Then,
\[
\partial_t d^\kappa_\e \geq \Div\mathcal z^\kappa_\e = \Div\big(\nabla\Phi^\kappa(\nabla d^\kappa_\e)\big) \qquad\text{in }\mathcal U^\kappa_\e
\]
in the sense of distributions. Moreover, there exists a constant $C_\kappa>0$, depending only on $\Phi^\kappa$ and $\mathfrak M^\kappa$, such that for every $\sigma>0$,
\[
\Div\mathcal z^\kappa_\e \leq \frac{C_\kappa}{\sigma}\quad\text{in the sense of distributions in $\mathcal U^\kappa_\e\cap\{d^\kappa_\e>\sigma\}$.}
\]
\end{proposition}

\begin{proof}
The proof follows the argument of Proposition~\ref{prop:superflow}, with the substantial simplification that, at fixed $\kappa$, no approximation of the anisotropy and no
$(n,\delta,R)$-regularization are needed.

\medskip
\noindent\textit{Step 1. Viscosity inequality.} Fix $(x_0,z_0,t_0)\in\mathcal U^\kappa_\e$ and let $\varphi\in C^2$ be such that $d^\kappa_\e-\varphi$ has a local minimum at
$(x_0,z_0,t_0)$. Let $(\widehat x,\widehat z)\in\partial\mathcal S^\kappa_\e(t_0)$ be a point of $(\mathfrak M^\kappa)^o$-minimal distance from $(x_0,z_0)$. By the definition of
$\mathcal U^\kappa_\e$, we have  $\widehat z>-\frac{t_0}{\e}$. Arguing exactly as in Step~1 of the proof of Proposition~\ref{prop:superflow}, the translated test function
$\overline\varphi(x,z,t):=\varphi(x+x_0-\widehat x,z+z_0-\widehat z,t){-\varphi(x_0,z_0,t_0)}$ touches from below the signed $(\mathfrak M^\kappa)^o$-distance function
$\operatorname{sd}^{\kappa}_\e(x,z,t):=\operatorname{sd}^{\,(\mathfrak M^\kappa)^o}\big((x,z), \mathcal S^\kappa_\e(t)\big)$ from $\mathcal S^\kappa_\e(t)$ at $(\widehat x,\widehat
z,t_0)$.

The contact regularity and the smooth uniform ellipticity of $\Phi^\kappa$ imply that $\partial\mathcal S^\kappa_\e(t_0)$, and hence its signed $(\mathfrak M^\kappa)^o$-distance
$\operatorname{sd}^{\kappa}_\e(x,z,t)$, are smooth in a neighborhood of the touching point. Let $\nu$ denote the outer Euclidean unit normal to $\partial\mathcal S^\kappa_\e(t_0)$
at $(\widehat x,\widehat z)$. At the touching point $(\widehat x,\widehat z,t_0)$, $\nabla\operatorname{sd}^{\kappa}_\e=\nu/\mathfrak M^\kappa(\nu)$.

Since the sets $\mathcal S^\kappa_\e(t)$ translate with velocity $-e_{_{N+1}}/\e$, and using the Euler--Lagrange equation for the solitons \eqref{mob:EL-soliton}, we have at the
same point
\[
\partial_t\operatorname{sd}^{\kappa}_\e = \frac{\nu\cdot e_{_{N+1}}} {\e\,\mathfrak M^\kappa(\nu)} = \frac1\e\partial_z\Phi^\kappa(\nu) = H^{\Phi^\kappa}_{\partial\mathcal S^\kappa_\e(t_0)}= \Div\big(
\nabla\Phi^\kappa(\nu)\big)= \Div\left( \nabla\Phi^\kappa \big( \nabla\operatorname{sd}^{\kappa}_\e \big)\right),
\]
where, in the first identity we used the transport equation, and in last equality we used the fact that  $\nabla\Phi^\kappa$ is $0$-homogeneous. Thus,
\[
\begin{aligned}
\partial_t\varphi(x_0,z_0,t_0) &= \partial_t\overline\varphi(\widehat x,\widehat z,t_0) = \partial_t\operatorname{sd}^{\kappa}_\e (\widehat x,\widehat z,t_0)= \Div\left( \nabla\Phi^\kappa \big(
\nabla\operatorname{sd}^{\kappa}_\e \big)\right)(\widehat x,\widehat z,t_0)\\
&\geq \Div\big( \nabla\Phi^\kappa(\nabla\overline\varphi) \big) (\widehat x,\widehat z,t_0) = \Div\big( \nabla\Phi^\kappa(\nabla\varphi) \big) (x_0,z_0,t_0),
\end{aligned}
\]
proving that $\partial_t d^\kappa_\e\geq\Div(\nabla\Phi^\kappa(\nabla d^\kappa_\e))$ holds in the viscosity sense.

\medskip
\noindent\textit{Step 2. Distributional inequality.} As in Step~2 of Proposition~\ref{prop:superflow}, the translating character of $\mathcal S^\kappa_\e$ and the smoothness of
$(\mathfrak M^\kappa)^o$ imply that $d^\kappa_\e$ is locally semiconcave in space-time away from $\{d^\kappa_\e=0\}$. By Aleksandrov's theorem, $d^\kappa_\e$ is therefore twice
differentiable a.e., and by the viscosity inequality of Step~1 we have
\[
\partial_t d^\kappa_\e \geq \nabla^2\Phi^\kappa(\nabla d^\kappa_\e): \nabla^2 d^\kappa_\e \qquad\text{a.e. in }\mathcal U^\kappa_\e.
\]

The mollification argument of Step~3 of Proposition~\ref{prop:superflow} applies verbatim. Since the singular part of the distributional Hessian of a semiconcave function is
nonpositive and $\nabla^2\Phi^\kappa\geq0$, we obtain $\Div(\nabla\Phi^\kappa(\nabla d^\kappa_\e))\leq\nabla^2\Phi^\kappa(\nabla d^\kappa_\e):\nabla^2d^\kappa_\e$ in the sense of
distributions. Combining the two inequalities proves $\partial_t d^\kappa_\e\geq\Div(\nabla\Phi^\kappa(\nabla d^\kappa_\e))$ in the sense of distributions in $\mathcal
U^\kappa_\e$.

\medskip
\noindent\textit{Step 3. Upper bound on the divergence.} Fix $\sigma>0$. By the standard semiconcavity estimate for the $(\mathfrak M^\kappa)^o$-distance, proved exactly as in
Remark~\ref{phiC2}, we have
\[
\nabla^2d^\kappa_\e(\cdot,\cdot,t) \leq \frac{c_\kappa}{\sigma}\,\mathrm{Id}
\]
in the sense of distributions in $\{d^\kappa_\e(\cdot,\cdot,t)>\sigma\}$, where $c_\kappa$ depends only on $\mathfrak M^\kappa$.

Since $\mathfrak M^\kappa(\nabla d^\kappa_\e)=1$, the gradient $\nabla d^\kappa_\e$ takes values a.e. in the compact hypersurface $\{p\in\R^{N+1}:\mathfrak M^\kappa(p)=1\}$. Hence
\[
P_\kappa := \sup_{\mathfrak M^\kappa(p)=1} \operatorname{tr}\nabla^2\Phi^\kappa(p) <+\infty.
\]

Using the semiconcavity estimate and once again the fact that the singular part of $\nabla^2d^\kappa_\e$ is nonpositive, we obtain
\[
\Div\big( \nabla\Phi^\kappa(\nabla d^\kappa_\e) \big) \leq \frac{c_\kappa P_\kappa}{\sigma}
\]
distributionally in $\mathcal U^\kappa_\e\cap\{d^\kappa_\e>\sigma\}$.
\end{proof}

We next establish the analogue of Proposition~\ref{prop:paraboloidi}. Recall that, by construction, $\Phi^\kappa$ is symmetric with respect to the last variable. We use the
intrinsic paraboloids introduced in \eqref{abs:paraboloidi}-\eqref{abs:defg}, with $\Phi^\kappa$ in place of $\Psi$, and we set
\[
\phi^\kappa(\xi):=\Phi^\kappa(\xi,0),\qquad \alpha_\kappa:=\min_{p\in\S^N}\frac{\Phi^\kappa(p)}{\mathfrak M^\kappa(p)}>0, \qquad \eta_\kappa:=\frac{\alpha_\kappa}{4N}.
\]

\begin{proposition}[Intrinsic paraboloid inclusion]\label{mob:prop:paraboloidi}
Let $\kappa\in(0,\kappa_0)$ and $\e>0$. If $\overline{W^{\phi^\kappa}_r(x_0)}\subset\mathrm{int}(E_0)$, then
\[
\mathcal C^{\Phi^\kappa}_{r,\eta,\e}(x_0)\subset E^\kappa_\e
\]
for every $\eta\in(0,\eta_\kappa]$. Symmetrically, if $\overline{W^{\phi^\kappa}_r(x_0)}\subset\R^N\setminus\overline{E_0}$, then
\[
\mathcal C^{\Phi^\kappa}_{r,\eta,\e}(x_0) \subset\overline{\R^{N+1}\setminus E^\kappa_\e}
\]
for every $\eta\in(0,\eta_\kappa]$.

Setting
\begin{equation}\label{mob:def-tilde-paraboloid}
\widetilde{\mathcal C}^{\Phi^\kappa}_{r,\eta,\e}(x_0) := \mathcal C^{\Phi^\kappa}_{r,\eta,\e}(x_0) \cup \left(\overline{W^{\phi^\kappa}_r(x_0)}\times(-\infty,0]\right),
\end{equation}
for every fixed $\eta>0$, every sequence $\e_j\to0^+$ and every sequence $t_j\searrow0$ we have
\begin{equation}\label{mob:translated-paraboloid-limit}
\widetilde{\mathcal C}^{\Phi^\kappa}_{r,\eta,\e_j}(x_0) -\frac{t_j}{\e_j}e_{_{N+1}} \stackrel{\mathcal K}{\longrightarrow} \overline{W^{\phi^\kappa}_r(x_0)}\times\R \qquad\text{in $\R^{N+1}$.}
\end{equation}
\end{proposition}

\begin{proof}
By Remark~\ref{rm:complementare}, it is enough to consider the case $\overline{W^{\phi^\kappa}_r(x_0)}\subset\mathrm{int}(E_0)$. We apply Lemma~\ref{lm:paraboloidi} with
$\Psi=\Phi^\kappa$ and $g=0$. Notice first that the trace condition \eqref{abs:trace-hyp} is automatically satisfied.

It remains only to check \eqref{abs:vertical}. Let $\eta\in(0,\eta_\kappa]$ and, at a point where $g^{\Phi^\kappa}_{r,\eta,\e}>0$, set $p=(\nabla g^{\Phi^\kappa}_{r,\eta,\e},-1)$.
By the compatibility relation between $\Phi^\kappa$ and $\mathfrak M^\kappa$,
\[
\Phi^\kappa(p)\partial_z\Phi^\kappa(p) = -\frac{\Phi^\kappa(p)}{\mathfrak M^\kappa(p)} \leq-\alpha_\kappa = -4N\eta_\kappa < -2N\eta_\kappa.
\]
Consequently, $\e^{-1}\Phi^\kappa(p)\partial_z\Phi^\kappa(p)<-2N\eta_\kappa/\e$, which is precisely \eqref{abs:vertical}. Lemma~\ref{lm:paraboloidi} therefore yields $\mathcal
C^{\Phi^\kappa}_{r,\eta,\e}(x_0)\subset E^\kappa_\e$ for every $\eta\in(0,\eta_\kappa]$, and the complementary statement follows from Remark~\ref{rm:complementare}.

We now prove \eqref{mob:translated-paraboloid-limit}. Set $\lambda:=2\eta/\e$ and recall from \eqref{abs:defg} that the intrinsic paraboloid associated with $\Phi^\kappa$ is
defined in terms of $G:=\frac12\Phi^\kappa(\cdot,-1)^2$. Consider
\[
G_\lambda(x):=\frac1{\lambda^2}G(\lambda x) =\frac12\Phi^\kappa\left(x,-\frac1\lambda\right)^2.
\]
As $\lambda\to+\infty$, $G_\lambda\to\frac12(\phi^\kappa)^2$ locally uniformly in $\R^N$. Hence, by the standard convergence properties of Legendre transforms,
$\lambda^{-2}G^*(\lambda x)=G_\lambda^*(x)\to\frac12((\phi^\kappa)^o(x))^2$ locally uniformly in $\R^N$. Since $\{(\phi^\kappa)^o=r\}$ is compact, it follows from \eqref{abs:defa}
that $a^{\Phi^\kappa}_{\lambda,r}/\lambda^2\to r^2/2$. Recalling \eqref{abs:defg} and the identity $\e\lambda=2\eta$, we obtain
\[
\e\,g^{\Phi^\kappa}_{r,\eta,\e}(x) \longrightarrow \eta\big(r^2-((\phi^\kappa)^o(x))^2\big) \qquad\text{locally uniformly in $\R^N$ as $\e\to0^+$.}
\]

By \eqref{abs:footprint} and \eqref{mob:def-tilde-paraboloid}, the projection onto $\R^N$ of $\widetilde{\mathcal C}^{\Phi^\kappa}_{r,\eta,\e}(x_0)$ is contained in
$\overline{W^{\phi^\kappa}_r(x_0)}$ for every $\e>0$. Hence every Kuratowski upper limit in \eqref{mob:translated-paraboloid-limit} is contained in
$\overline{W^{\phi^\kappa}_r(x_0)}\times\R$.

Conversely, let $x\in W^{\phi^\kappa}_r(x_0)$. Then
\[
\e_j g^{\Phi^\kappa}_{r,\eta,\e_j}(x-x_0)-t_j \longrightarrow \eta\big(r^2-((\phi^\kappa)^o(x-x_0))^2\big)>0.
\]
It follows that $g^{\Phi^\kappa}_{r,\eta,\e_j}(x-x_0)-t_j/\e_j\to+\infty$. We deduce that for every $z\in\R$ the point $(x,z)$ belongs to $\widetilde{\mathcal
C}^{\Phi^\kappa}_{r,\eta,\e_j}(x_0)-t_j e_{_{N+1}}/\e_j$ for all sufficiently large $j$. Thus $W^{\phi^\kappa}_r(x_0)\times\R$ is contained in the Kuratowski lower limit.  This
proves \eqref{mob:translated-paraboloid-limit}.
\end{proof}

We next establish the analogue of Proposition~\ref{prop:wulff} at fixed $\kappa$. Recall that the regularity and uniform convexity of $\Phi^\kappa$ and $\mathfrak M^\kappa$ imply
that $\mathfrak M^\kappa$ is $\Phi^\kappa$-regular, see Section~\ref{subsec:cf}.

Fix $\kappa\in(0,\kappa_0)$ and set
\[
\Lambda_\kappa
:=
\max_{(\mathfrak M^\kappa)^o(q)=1}
\left[
\nabla^2\Phi^\kappa
\bigl(\nabla(\mathfrak M^\kappa)^o(q)\bigr)
:
\nabla^2(\mathfrak M^\kappa)^o(q)
\right],
\]
the maximum of the $\Phi^\kappa$-curvature of $\partial\mathcal W_1^{\mathfrak M^\kappa}$.

\begin{proposition}\label{mob:prop:wulff}
Let $\eta>0$ be fixed and set
\[
\e_0^\kappa(\eta)
:=
\frac{\eta^2}
{16\bigl(\mathfrak M^\kappa(e_{_{N+1}})\bigr)^2}.
\]
Assume
\[
\overline{\mathcal W_R^{\mathfrak M^\kappa}(x_0,z_0)}
\cap\mathcal S^\kappa_\e(t_0)=\varnothing,
\qquad
\overline{\mathcal W_R^{\mathfrak M^\kappa}(x_0,z_0)}
\subset
\left\{(x,z):z>-\frac{\eta}{2\e}\right\},
\]
for some $\e\in(0,\e_0^\kappa(\eta)\land\tfrac1{16R^2})$ and $t_0>\eta$. Setting
\[
R(s):=\sqrt{R^2-2\Lambda_\kappa s},
\]
we have
\[
\mathcal S^\kappa_\e(t_0+s)
\cap
\overline{\mathcal W_{R(s)}^{\mathfrak M^\kappa}
(x_0,z_0)}
=\varnothing
\qquad
\text{for all }
s\in\left(0,\frac{R^2}{2\Lambda_\kappa}\right).
\]
\end{proposition}

\begin{proof}
We first observe that the tube generated by $\overline{\mathcal W_{R(s)}^{\mathfrak M^\kappa} (x_0,z_0)}$ is a $(\Phi^\kappa,\mathfrak M^\kappa)$-superflow until its extinction
time. Indeed, outside the evolving Wulff shape its $(\mathfrak M^\kappa)^o$-distance is
\[
d(x,z,s)
=(\mathfrak M^\kappa)^o(x-x_0,z-z_0)-R(s).
\]
By the definition of $\Lambda_\kappa$ and homogeneity,
\[
\Div\bigl(\nabla\Phi^\kappa(\nabla d)\bigr)
\leq
\frac{\Lambda_\kappa}
{(\mathfrak M^\kappa)^o(x-x_0,z-z_0)}
\leq
\frac{\Lambda_\kappa}{R(s)}
=\partial_s d.
\]

Set
\[
\Delta(s)
:=
\dist^{(\mathfrak M^\kappa)^o}
\left(
\mathcal S^\kappa_\e(t_0+s),
\overline{\mathcal W_{R(s)}^{\mathfrak M^\kappa}
(x_0,z_0)}
\right),
\qquad
\Delta:=\Delta(0)\land\frac1{2\sqrt\e}.
\]
As in the proof of Proposition~\ref{prop:wulff}, $\Delta(\cdot)$ is continuous, and we claim that $\Delta(s)\geq\Delta$ throughout $[0,R^2/(2\Lambda_\kappa))$.

The localization required in that proof is available here as well. Indeed,
\[
\dist^{(\mathfrak M^\kappa)^o}
\bigl((x,z),\mathcal H_{-t/\e}\bigr)
=
\frac{z+t/\e}{\mathfrak M^\kappa(e_{_{N+1}})}
>
\frac{\eta}
{2\e\,\mathfrak M^\kappa(e_{_{N+1}})}
>
\frac1{\sqrt\e}
\]
whenever $z>-\eta/(2\e)$ and $t>\eta$.

Consequently, at a time $\bar\sigma$ such that $\Delta(\bar\sigma)=\Delta$, choosing $r>0$ as in the proof of Proposition~\ref{prop:wulff}, we obtain
\[
\overline{
\mathcal W_{R(\bar\sigma)+r}^{\mathfrak M^\kappa}
(x_0,z_0)}
\times\{t_0+\bar\sigma\}
\subset\mathcal U^\kappa_\e.
\]
In fact, on this enlarged Wulff shape,
\[
0<d^\kappa_\e(\cdot,t_0+\bar\sigma)
\leq 2R+\Delta+r<\frac1{\sqrt\e},
\]
and its points still satisfy $z>-\eta/(2\e)$, provided that $r$ is sufficiently small. Since $\mathcal U^\kappa_\e$ is open, the same inclusion holds on a short time interval
starting at $t_0+\bar\sigma$.

We may therefore repeat the localized comparison and continuation argument of Proposition~\ref{prop:wulff}, using Proposition~\ref{mob:prop:superflow} in place of
Proposition~\ref{prop:superflow}. This proves $\Delta(s)\geq\Delta$ and concludes the proof.
\end{proof}

Fix $\eta>0$ and let $(x_0,z_0)\notin\mathcal S^\kappa_\e(t_0)$, with $t_0>\eta$. Assume that
\begin{equation}\label{mob:eq*}
\mathcal W_{d^\kappa_\e(x_0,z_0,t_0)}
^{\mathfrak M^\kappa}(x_0,z_0)
\subset
\left\{(x,z):z>-\frac{\eta}{2\e}\right\},
\qquad
0<\e<
\e_0^\kappa(\eta)\land
\frac1{16\bigl(d^\kappa_\e(x_0,z_0,t_0)\bigr)^2}.
\end{equation}
For every $0<R<d^\kappa_\e(x_0,z_0,t_0)$, by Proposition~\ref{mob:prop:wulff} we have
\[
d^\kappa_\e(x_0,z_0,t_0+s)
\geq\sqrt{R^2-2\Lambda_\kappa s}
\qquad
\text{for }
0<s<\frac{R^2}{2\Lambda_\kappa}.
\]
Letting $R\nearrow d^\kappa_\e(x_0,z_0,t_0)$, we obtain
\begin{equation}\label{mob:eq:helly}
d^\kappa_\e(x_0,z_0,t_0+s)
\geq
\sqrt{
\bigl(d^\kappa_\e(x_0,z_0,t_0)\bigr)^2
-2\Lambda_\kappa s}
\end{equation}
for
\[
0<s<
\frac{\bigl(d^\kappa_\e(x_0,z_0,t_0)\bigr)^2}
{2\Lambda_\kappa},
\]
provided that \eqref{mob:eq*} is satisfied.

We are now ready to prove the analogue of Proposition~\ref{prop:helly} at fixed $\kappa$.

\begin{proposition}\label{mob:prop:helly}
Fix $\kappa\in(0,\kappa_0)$, and let $\e_j\to0^+$ be such that
\[
\mathcal S^\kappa_{\e_j} \stackrel{\mathcal K}{\longrightarrow} \mathcal S^\kappa \qquad\text{in }\R^{N+2}.
\]
Setting
\[
\mathcal S^\kappa(t) := \{(x,z)\in\R^{N+1}:(x,z,t)\in\mathcal S^\kappa\},
\]
the following properties hold:
\begin{itemize}
\item[(a)]
\[
\mathrm{int}(E_0)\times\R \subset \mathcal S^\kappa(0) \subset \overline{E_0}\times\R.
\]
\item[(b)] $\mathcal S^\kappa(\cdot)$ is Kuratowski left-continuous.
\item[(c)] If $\mathcal S^\kappa(t)=\varnothing$, then $\mathcal S^\kappa(s)=\varnothing$ for every $s>t$. In particular, the extinction time
\[
T_\kappa^* := \inf\{t>0:\mathcal S^\kappa(s)=\varnothing \text{ for every }s\geq t\}
\]
is well defined.
\item[(d)] There exists a countable set $\mathcal N_\kappa\subset(0,T_\kappa^*)$ such that, for every $t\in(0,T_\kappa^*)\setminus\mathcal N_\kappa$,
\[
d^\kappa_{\e_j}(\cdot,\cdot,t) \longrightarrow \dist^{(\mathfrak M^\kappa)^o} \big(\cdot,\mathcal S^\kappa(t)\big)
\]
locally uniformly in $\R^{N+1}$. Equivalently,
\[
\mathcal S^\kappa_{\e_j}(t) \stackrel{\mathcal K}{\longrightarrow} \mathcal S^\kappa(t).
\]
\item[(e)] For every $t>T_\kappa^*$,
\[
\mathcal S^\kappa_{\e_j}(t) \stackrel{\mathcal K}{\longrightarrow} \varnothing.
\]
\item[(f)] $\mathcal S^\kappa$ is a $(\Phi^\kappa,\mathfrak M^\kappa)$-superflow with initial datum $\overline{E_0}\times\R$.
\end{itemize}
\end{proposition}

\begin{proof}

By \eqref{mob:eq:helly}, the proof of properties \textup{(b)}--\textup{(e)} is exactly the argument of Proposition~\ref{prop:helly}, with $d^\e$, $\Phi^o$, and $N$ replaced
respectively by $d^\kappa_\e$, $(\mathfrak M^\kappa)^o$, and $\Lambda_\kappa$. Indeed, at fixed $\kappa$, the localization conditions in \eqref{mob:eq*} hold for all sufficiently
small $\e$ when the centers and radii remain bounded and the starting times stay bounded away from zero.

We now prove \textup{(a)}. Let $\overline{W^{\phi^\kappa}_r(x_0)}\subset\mathrm{int}(E_0)$, and fix $\eta\in(0,\eta_\kappa)$ (recall that $\phi^\kappa=\Phi^\kappa(\cdot,0)$). By
Proposition~\ref{mob:prop:paraboloidi}, together with the prescribed trace below $z=0$, \(\widetilde{\mathcal C}^{\Phi^\kappa}_{r,\eta,\e_j}(x_0) \subset E^\kappa_{\e_j}\) and
hence, for every sequence $t_j\searrow0$, \( \widetilde{\mathcal C}^{\Phi^\kappa}_{r,\eta,\e_j}(x_0) -\frac{t_j}{\e_j}e_{_{N+1}} \subset \mathcal S^\kappa_{\e_j}(t_j) \). The
translated-paraboloid convergence proved in Proposition~\ref{mob:prop:paraboloidi} yields
\[
\overline{W^{\phi^\kappa}_r(x_0)}\times\R \subset\mathcal S^\kappa(0).
\]
Since such Wulff shapes exhaust $\mathrm{int}(E_0)$, we obtain $\mathrm{int}(E_0)\times\R\subset\mathcal S^\kappa(0)$.

For the reverse inclusion, let $\overline{W^{\phi^\kappa}_r(x_0)} \subset\R^N\setminus\overline{E_0}$. Arguing as in the proof of Proposition~\ref{prop:helly}, we use the
complementary statement of Proposition~\ref{mob:prop:paraboloidi}, with $r/2$ in place of $r$ and any fixed $\eta\in(0,\eta_\kappa)$. The strict inclusion
$\overline{W^{\phi^\kappa}_{r/2}(x_0)} \Subset W^{\phi^\kappa}_r(x_0)$, together with the profile estimate used in the proof of that proposition, shows that the translated sets \(
\widetilde{\mathcal C}^{\Phi^\kappa}_{r/2,\eta,\e_j}(x_0) -\frac{t}{\e_j}e_{_{N+1}} \) remain at a positive distance from $\mathcal S^\kappa_{\e_j}(t)$, uniformly for $j$ large and
$t\geq0$. Taking $t=t_j\searrow0$ and using again the translated-paraboloid convergence, we obtain
\[
\big( \overline{W^{\phi^\kappa}_{r/2}(x_0)}\times\R \big) \cap\mathcal S^\kappa(0)=\varnothing.
\]
Hence $\mathcal S^\kappa(0)\subset\overline{E_0}\times\R$, proving \textup{(a)}.

It remains to prove \textup{(f)}. First, observe that if $K\Subset(\R^{N+1}\times(0,T_\kappa^*))\setminus\mathcal S^\kappa$, then there exists $\sigma>0$ such that
$d^\kappa\geq2\sigma$ on $K$ and, in turn,  for $j$ sufficiently large, $K\subset\mathcal U^\kappa_{\e_j}$ and $d^\kappa_{\e_j}\geq\sigma$ there. Moreover, recall that for almost
every $t$, $d^\kappa_{\e_j}(\cdot,\cdot,t)\to d^\kappa(\cdot,\cdot,t)$ locally uniformly. Therefore, since the spatial semiconcavity estimates are uniform in $j$ (see
Remark~\ref{phiC2}),
\[
\nabla d^\kappa_{\e_j} \longrightarrow \nabla d^\kappa \qquad\text{a.e. on }K.
\]
By dominated convergence we then obtain
\[
\nabla\Phi^\kappa(\nabla d^\kappa_{\e_j}) \longrightarrow \nabla\Phi^\kappa(\nabla d^\kappa) \qquad\text{in }L^1(K).
\]
Due to the latter convergence (and the arbitrariness of $K$), we can pass to the limit the distributional inequality and the divergence bound established in
Proposition~\ref{mob:prop:superflow} to establish  \textup{(f)}.
\end{proof}

We are now ready to prove the analogue of Theorem~\ref{thm:main1} for the elliptic regularization scheme with arbitrary mobility. We first present the result in terms of iterated
(Kuratowski) limits and then, in Corollary~\ref{mob:cor:diagonal} below, we show that the same limiting non-fattening evolution can be recovered along  a suitable diagonal
one-parameter family of approximating solitons.

\begin{theorem}\label{mob:thm:main1}
Let $\phi,\mathfrak m$, $\Phi^\kappa,\mathfrak M^\kappa$, $E^\kappa_\e$ and $\mathcal S^\kappa_\e$ be as above. Assume that $E_0$ is non-fattening for the $(\phi,\mathfrak
m)$-level-set flow in the sense of Definition~\ref{def:fat}, and let $(E(t))_{t>0}$ be the generalized $(\phi,\mathfrak m)$-motion starting from $E_0$ in the sense of
Definition~\ref{def:gen-motion}. Denote by
\[
E:=\big\{(x,t)\in\R^N\times(0,+\infty):\,x\in E(t)\big\}
\]
the associated space-time tube and by $\mathcal E$ the space-time tube built upon the cylindrical evolution $\mathcal E(t):=E(t)\times\R$, i.e.,
\[
\mathcal E:=\big\{(x,z,t)\in\R^{N+1}\times(0,+\infty):\,x\in E(t)\big\}.
\]

Let $\kappa_n\to0^+$ and $\e_j\to0^+$ be arbitrary sequences. Then, up to passing to a subsequence of $\{\e_j\}$, not relabeled, the following Kuratowski limits hold:
\begin{itemize}
\item[(a)]
\[
\lim_{n\to+\infty} \left( \lim_{j\to+\infty}\mathcal S^{\kappa_n}_{\e_j} \right) =\mathcal E \qquad\text{in }\R^{N+1}\times(0,+\infty).
\]

\item[(b)] There exists a countable set $\mathcal N\subset(0,+\infty)$, possibly depending on the above sequences and on the extracted subsequence, such that
\[
\lim_{n\to+\infty} \left( \lim_{j\to+\infty}\mathcal S^{\kappa_n}_{\e_j}(t) \right) =E(t)\times\R \qquad\text{in }\R^{N+1}
\]
for every $t\notin\mathcal N$.

\item[(c)]
\[
\lim_{n\to+\infty} \left( \lim_{j\to+\infty}\mathcal R_{\e_j}(E^{\kappa_n}_{\e_j}) \right) =E \qquad\text{in }\R^N\times(0,+\infty).
\]
\end{itemize}
\end{theorem}

\begin{remark}\label{mob:remark:complementary}
Under the same assumptions and with the same notation as in Theorem~\ref{mob:thm:main1}, set
\[
\mathcal S_{\e}^{\kappa,\complement}(t) := \overline{\R^{N+1}\setminus\mathcal S^\kappa_\e(t)} = \R^{N+1}\setminus\mathrm{int}\big(\mathcal S^\kappa_\e(t)\big),
\]
and denote by $\mathcal S_{\e}^{\kappa,\complement}$ the corresponding space-time tube. Choosing the diagonal subsequence in the proof of Theorem~\ref{mob:thm:main1} also for the
complementary solitons, and applying the theorem with $\R^N\setminus E_0$ in place of $E_0$, after possibly extracting a further (not relabeled) subsequence $\e_j$,  we obtain
\[
\begin{split}
&\lim_{n\to+\infty}\left(\lim_{j\to+\infty} \mathcal S_{\e_j}^{\kappa_n,\complement}\right) = \big(\R^{N+1}\times(0,+\infty)\big)\setminus\mathrm{int}(\mathcal E),\\
&\lim_{n\to+\infty}\left(\lim_{j\to+\infty} \mathcal S_{\e_j}^{\kappa_n,\complement}(t)\right) = \big(\R^N\setminus\mathrm{int}(E(t))\big)\times\R \qquad\text{for every }t\notin\mathcal N,\\
&\lim_{n\to+\infty}\left(\lim_{j\to+\infty} \left[ \big(\R^N\times(0,+\infty)\big) \setminus\mathcal R_{\e_j}\big(\mathrm{int}(E^{\kappa_n}_{\e_j})\big) \right]\right) =
\big(\R^N\times(0,+\infty)\big)\setminus\mathrm{int}(E),
\end{split}
\]
where, after possibly enlarging it, $\mathcal N$ is still countable and all the limits are understood in the Kuratowski sense.

Combining these convergences with Theorem~\ref{mob:thm:main1}, and setting
\[
\Gamma^\kappa_\e:=\mathcal R_\e(\pa E^\kappa_\e), \qquad \Gamma:=\big\{(x,t)\in\R^N\times(0,+\infty):\,u(x,t)=0\big\},
\]
where $u$ is as in Definition~\ref{def:fat}, we also deduce
\[
\begin{array}{rcll}
\displaystyle \lim_{n\to+\infty}\left(\lim_{j\to+\infty} \pa\mathcal S^{\kappa_n}_{\e_j}(t)\right) &=& \pa E(t)\times\R &\quad\text{in $\R^{N+1}$ for every $t\notin\mathcal N$},\\[3mm]
\displaystyle \lim_{n\to+\infty}\left(\lim_{j\to+\infty} \Gamma^{\kappa_n}_{\e_j}\right) &=& \Gamma &\quad\text{in $\R^N\times(0,+\infty)$}.
\end{array}
\]
\end{remark}

\begin{proof}[Proof of Theorem~\ref{mob:thm:main1}]
Let $\kappa_n\to0^+$ and $\e_j\to0^+$. By Proposition~\ref{mob:prop:helly} and a standard diagonal extraction, we may assume that, for every $n$, $\mathcal
S^{\kappa_n}_{\e_j}\stackrel{\mathcal K}{\longrightarrow}\mathcal S^{\kappa_n}$ in $\R^{N+2}$ as $j\to+\infty$. We perform the same diagonal extraction for the complementary
solitons.

Fix for the moment $n$ and write simply $\kappa=\kappa_n$. Arguing with the complementary solitons exactly as in the proof of Theorem~\ref{thm:main1}, we obtain a
$(\Phi^\kappa,\mathfrak M^\kappa)$-subflow $\mathcal A^\kappa$ such that
\begin{equation}\label{mob:main:AS}
\mathcal A^\kappa\subset\mathcal S^\kappa, \qquad \mathrm{int}(E_0)\times\R \subset\mathcal A^\kappa(0) \subset\mathcal S^\kappa(0) \subset\overline{E_0}\times\R.
\end{equation}

Let $u_0$ be the non-fattening level-set representation of $E_0$ given by the assumption, and let $u$ be the corresponding $(\phi,\mathfrak m)$-level-set flow. Set
$\phi^\kappa(\xi):=\Phi^\kappa(\xi,0)$ and $\mathfrak m^\kappa(\xi):=\mathfrak M^\kappa(\xi,0)$, and let $u_\kappa$ be the $(\phi^\kappa,\mathfrak m^\kappa)$-level-set flow
starting from $u_0$. By uniqueness, $\widehat u_\kappa(x,z,t):=u_\kappa(x,t)$ is the $(\Phi^\kappa,\mathfrak M^\kappa)$-level-set flow starting from $\widehat u_0(x,z):=u_0(x)$.

Since $\{u_0<0\}\subset E_0\subset\{u_0\leq0\}$, and the first and last sets are respectively open and closed, by \eqref{mob:main:AS} we have
\[
\{\widehat u_0<0\}\subset\mathrm{int}(E_0)\times\R \subset\mathcal A^\kappa(0)\subset\mathcal S^\kappa(0) \subset\overline{E_0}\times\R\subset\{\widehat u_0\leq0\}.
\]
Hence, by the minimality and maximality properties in Theorem~\ref{th:phiregularlevelset},
\begin{equation}\label{mob:main:sandwich}
\{\widehat u_\kappa<0\} \subset\mathcal A^\kappa \subset\mathcal S^\kappa \subset\{\widehat u_\kappa\leq0\}.
\end{equation}

By Proposition~\ref{prop:extended-norms}, $\phi^{\kappa_n}\to\phi$ and $\mathfrak m^{\kappa_n}\to\mathfrak m$ locally uniformly. Therefore, Theorem~\ref{th:levelset-stability}
yields
\begin{equation}\label{mob:main:uk-conv}
u_{\kappa_n}\longrightarrow u \qquad\text{locally uniformly in }\R^N\times[0,+\infty).
\end{equation}
If $\mathcal S$ is any subsequential Kuratowski limit of $\mathcal S^{\kappa_n}$, from \eqref{mob:main:sandwich} and \eqref{mob:main:uk-conv} we obtain
\[
\{\widehat u<0\}\subset\mathcal S\subset\{\widehat u\leq0\}, \qquad \widehat u(x,z,t):=u(x,t).
\]
By the non-fattening assumption we have $\overline{\{\widehat u<0\}}=\{\widehat u\leq0\}$ in $\R^{N+1}\times(0,+\infty)$. Since $\mathcal S$ is closed, it follows that $\mathcal
S=\{\widehat u\leq0\}=\mathcal E$ for positive times. Hence every subsequential limit is the same and $\mathcal S^{\kappa_n}\stackrel{\mathcal K}{\longrightarrow}\mathcal E$ as
$n\to+\infty$, proving \textup{(a)}.

We next prove \textup{(b)}. Set $d_\phi(x,t):=\dist^{\phi^o}(x,E(t))$. By comparison with the $(\phi,\mathfrak m)$-evolution of $\phi$-Wulff shapes, as in
\cite[Section~3.3]{CMNP19}, we obtain, whenever $d_\phi(x,t)>0$ and $s>t$ is sufficiently close to $t$,
\begin{equation}\label{mob:main:quasimon}
d_\phi(x,s)\geq d_\phi(x,t)-\frac{C}{d_\phi(x,t)}(s-t),
\end{equation}
where $C>0$ depends only on $\phi$ and $\mathfrak m$. Indeed, by the corresponding Wulff-shape estimate we first obtain $d_\phi(x,s)\geq\sqrt{d_\phi(x,t)^2-C(s-t)}$, and
\eqref{mob:main:quasimon} follows from the elementary inequality $\sqrt{R^2-a}\geq R-a/R$.

More precisely, for every $\delta>0$ there exists $\tau_\delta>0$ such that
\[
d_\phi(x,s)\geq d_\phi(x,t)-\frac{C}{\delta}(s-t)
\]
whenever $d_\phi(x,t)\geq\delta$ and $0<s-t<\tau_\delta$. Thus, as long as $d_\phi(x,\cdot)\geq\delta$, the map $t\mapsto d_\phi(x,t)+(C/\delta)t$ is nondecreasing. This is
precisely the local one-sided quasi-monotonicity used in the proof of \cite[Proposition~4.4]{CMP17}. Combined with the uniform $1$-Lipschitz continuity in space of
$d_\phi(\cdot,t)$ with respect to $\phi^o$, it implies that there exists a countable set $\mathcal N_0\subset(0,+\infty)$ such that
\begin{equation}\label{mob:main:distance-cont}
d_\phi(\cdot,t)\longrightarrow d_\phi(\cdot,t_0) \qquad\text{locally uniformly in $\R^N$ as $t\to t_0$}
\end{equation}
for every $t_0\notin\mathcal N_0$, which is clearly equivalent to the Kuratowski continuity of $E(\cdot)$ at $t_0$. At every such time we also have
\begin{equation}\label{mob:main:slice-nonfat}
\overline{\{u(\cdot,t)<0\}}=\{u(\cdot,t)\leq0\}=E(t).
\end{equation}
Indeed, otherwise a $\phi$-Wulff neighborhood of a point of $E(t)$ would be disjoint from $\{u(\cdot,t)<0\}$. Comparison for the complementary evolution would keep a smaller Wulff
neighborhood disjoint from $\{u<0\}$ for a nontrivial interval immediately after $t$; the space-time non-fattening property would then make it disjoint from $E(s)$ for $s>t$
sufficiently close to $t$, contradicting \eqref{mob:main:distance-cont}.

Fix $t\notin\mathcal N_0$. Taking time-$t$ sections in \eqref{mob:main:sandwich}, and using \eqref{mob:main:uk-conv} and \eqref{mob:main:slice-nonfat}, we get
\[
\{u_{\kappa_n}(\cdot,t)<0\}\times\R \subset\mathcal S^{\kappa_n}(t) \subset\{u_{\kappa_n}(\cdot,t)\leq0\}\times\R, \qquad \mathcal S^{\kappa_n}(t)\stackrel{\mathcal K}{\longrightarrow}E(t)\times\R.
\]
For the Kuratowski lower inclusion, it is enough to approximate every $x\in E(t)$ by $x_h\to x$ with $u(x_h,t)<0$, using \eqref{mob:main:slice-nonfat}, and then let
$h=h(n)\to+\infty$ sufficiently slowly in view of \eqref{mob:main:uk-conv}.

For every $n$, Proposition~\ref{mob:prop:helly} provides a countable set $\mathcal N_n$ outside which $\mathcal S^{\kappa_n}_{\e_j}(t)\stackrel{\mathcal K}{\longrightarrow}\mathcal
S^{\kappa_n}(t)$ as $j\to+\infty$. Hence $\mathcal N:=\mathcal N_0\cup\bigcup_{n\geq1}\mathcal N_n$ is countable and the assertion in \textup{(b)} follows.

Finally, \textup{(c)} follows by exactly the same argument as in Step~2 of the proof of Theorem~\ref{thm:main1}. Indeed, for fixed $n$ the convergence of $\mathcal
S^{\kappa_n}_{\e_j}$ implies the Kuratowski convergence of $\mathcal R_{\e_j}(E^{\kappa_n}_{\e_j})$ to the space-time trace of $\mathcal S^{\kappa_n}$ on $\{z=0\}$; i.e.,
\[
\mathcal R_{\e_j}(E^{\kappa_n}_{\e_j})\stackrel{\mathcal K}{\longrightarrow}F_n:=\big\{(x,t)\in\R^N\times(0,+\infty):\,(x,0,t)\in\mathcal S^{\kappa_n}\big\} \qquad\text{in }\R^N\times(0,+\infty).
\]
By \eqref{mob:main:sandwich},
\[
\{u_{\kappa_n}<0\}\subset F_n\subset\{u_{\kappa_n}\leq0\}.
\]
Hence, by \eqref{mob:main:uk-conv} and the space-time non-fattening property $\overline{\{u<0\}}=\{u\leq0\}=E$, we have $F_n\stackrel{\mathcal K}{\longrightarrow}E$ as
$n\to+\infty$.
\end{proof}

\begin{corollary}\label{mob:cor:diagonal}
Under the same assumptions on $\phi$, $\mathfrak m$, and $E_0$ as in Theorem~\ref{mob:thm:main1}, let $\Phi^\kappa$, $\mathfrak M^\kappa$, $E^\kappa_\e$, $\mathcal S^\kappa_\e$,
$E$, and $\mathcal E$ be as therein. Then there exists a function $\kappa:(0,1)\longrightarrow(0,\kappa_0)$ such that  $\kappa(\e)\longrightarrow0$ as $\e\to0^+$,
\begin{equation}\label{contr:ke}
    \e\,(L_{\kappa(\e)})^2\xrightarrow{\e\to 0^+}0,\qquad L_\kappa\coloneqq \big(1+\|\mathfrak{M}^\kappa\|_{L^\infty(\S^N)}+\operatorname{Lip}(\mathfrak{M}^\kappa)\big)
\end{equation}

and which may be chosen piecewise constant  and non-decreasing, with the following properties:
\begin{itemize}
\item[(a)]
\[
\mathcal S_{\e}^{\kappa(\e)} \stackrel{\mathcal K}{\longrightarrow} \mathcal E \qquad\text{in }\R^{N+1}\times(0,+\infty);
\]

\item[(b)] there exists a countable set $\mathcal N\subset(0,+\infty)$ such that
\[
\mathcal S_{\e}^{\kappa(\e)}(t) \stackrel{\mathcal K}{\longrightarrow} E(t)\times\R \qquad\text{in }\R^{N+1}
\]
for every $t\notin\mathcal N$;

\item[(c)]
\[
\mathcal R_{\e}\big(E_{\e}^{\kappa(\e)}\big) \stackrel{\mathcal K}{\longrightarrow} E \qquad\text{in }\R^N\times(0,+\infty).
\]
\end{itemize}
Moreover, $\kappa(\e)$ may be chosen so that the same conclusions hold simultaneously for the complementary construction associated with $\R^N\setminus E_0$.
\end{corollary}
Note that the function $\kappa(\e)$ in the above statement is not ``universal'', but it depends on $E_0$ (and possibly on the choice of $E^\kappa_\e$ in the minimization problem).
\begin{proof}[Proof of Corollary~\ref{mob:cor:diagonal}]
For fixed $\kappa\in(0,\kappa_0)$, let $\mathscr C_\kappa$ and $\mathscr D_\kappa$ be the families of Kuratowski cluster limits, as $\e\to0^+$, of $\mathcal S^\kappa_\e$ and
$\mathcal R_\e(E^\kappa_\e)$, respectively, and let $d_1,d_2$ be metrics inducing the corresponding local Kuratowski topologies. Arguing as in the proof of
Theorem~\ref{mob:thm:main1}, we obtain
\begin{equation}\label{mob:diag:clusters0}
\{\widehat u_\kappa<0\}\subset\mathcal S \qquad\text{for all }\mathcal S\in \mathscr C_\kappa
\end{equation}
and
\begin{equation}\label{mob:diag:clusters}
\sup_{\mathcal S\in\mathscr C_\kappa}d_1(\mathcal S,\mathcal E)\to0, \qquad \sup_{F\in\mathscr D_\kappa}d_2(F,E)\to0 \qquad\text{as }\kappa\to0^+.
\end{equation}
Moreover, compactness at fixed $\kappa$ implies
\begin{equation}\label{mob:diag:tail}
\dist_{d_1}(\mathcal S^\kappa_\e,\mathscr C_\kappa)\to0, \qquad \dist_{d_2}(\mathcal R_\e(E^\kappa_\e),\mathscr D_\kappa)\to0 \qquad\text{as }\e\to0^+,
\end{equation}
and
\begin{equation}\label{ke1:con}
   \e\,L_\kappa\to 0\quad \text{as $\e\to 0^+$,}
\end{equation}
with $L_\kappa$ as in \eqref{contr:ke}.

We also claim that, for fixed $\kappa$ and  every compact set $K\Subset \{(x,z,t)\in\R^{N+1}\times(0,+\infty):u_\kappa(x,t)<0\}$, we have
\begin{equation}\label{mob:diag:interior}
\lim_{\e\to0^+} \sup_{(x,z,t)\in K} \dist^{(\mathfrak M^\kappa)^o} \big((x,z),\mathcal S^\kappa_\e(t)\big)=0.
\end{equation}
Indeed, otherwise, for some $r>0$ we could find $\e_j\to0^+$ and $(x_j,z_j,t_j)\longrightarrow(x,z,t)\in K$ such that $\mathcal W_r^{\mathfrak M^\kappa}(x_j,z_j) \cap\mathcal
S^\kappa_{\e_j}(t_j)=\varnothing$. Since $\kappa$ is fixed, by Wulff-shape avoidance we have $\mathcal W_{\sqrt{r^2-2\Lambda_\kappa s}}^{\mathfrak M^\kappa}(x_j,z_j) \cap\mathcal
S^\kappa_{\e_j}(t_j+s)=\varnothing$ as long as $r^2-2\Lambda_\kappa s>0$. Choose $s>0$ and then $\eta>0$ so small that $u_\kappa(x,t+s)<0$ and
$\sqrt{r^2-2\Lambda_\kappa\tau}>\frac{2r}{3} \quad\text{for every }\tau\in(s-\eta,s+\eta)$. For $j$ sufficiently large, $\overline{\mathcal W_{r/2}^{\mathfrak M^\kappa}(x,z)}
\subset \mathcal W_{2r/3}^{\mathfrak M^\kappa}(x_j,z_j)$, and hence a fixed space-time neighborhood of $\overline{\mathcal W_{r/2}^{\mathfrak M^\kappa}(x,z)} \times\{t+s\}$ is
disjoint from $\mathcal S^\kappa_{\e_j}$. Passing, up to a further subsequence, to the limit $\mathcal S^\kappa_{\e_j} \stackrel{\mathcal K}{\longrightarrow}\mathcal S$, with
$\mathcal S\in\mathscr C_\kappa$, we obtain $\mathcal W_{r/2}^{\mathfrak M^\kappa}(x,z) \cap\mathcal S(t+s)=\varnothing$, contradicting \eqref{mob:diag:clusters0}, since
$u_\kappa(x,t+s)<0$.

Choose now $\kappa_n\downarrow0$  so that the two quantities in \eqref{mob:diag:clusters}, as well as their analogues for the complementary construction of
Remark~\ref{mob:remark:complementary}, are smaller than $2^{-n}$, and set
\[
K_n:= \left\{(x,z,t):|x|,|z|\leq n,\quad \frac1n\leq t\leq n,\quad u_{\kappa_n}(x,t)\leq-\frac1n\right\}.
\]
By \eqref{ke1:con}, \eqref{mob:diag:tail} and \eqref{mob:diag:interior}, imposing simultaneously the analogous conditions for the complementary construction, we may choose
$\e_n\downarrow0$ so that, whenever $0<\e\leq\e_n$,
\begin{equation}\label{cfe}
\e_n L_{\kappa_n}^2\leq 2^{-n},
\end{equation}
\[
\dist_{d_1}(\mathcal S^{\kappa_n}_\e,\mathscr C_{\kappa_n})<2^{-n}, \qquad \dist_{d_2}(\mathcal R_\e(E^{\kappa_n}_\e),\mathscr D_{\kappa_n})<2^{-n},
\]
and
\[
\sup_{(x,z,t)\in K_n}
\dist^{(\mathfrak M^{\kappa_n})^o}
\big((x,z),\mathcal S^{\kappa_n}_\e(t)\big)
<
\frac{1}{n\|\mathfrak M^{\kappa_n}\|_{L^\infty(\S^N)}}
\qquad\text{for all }0<\e\leq\e_n.
\]
Since $|(\xi,\zeta)| \leq \|\mathfrak M^{\kappa_n}\|_{L^\infty(\S^N)} (\mathfrak M^{\kappa_n})^o(\xi,\zeta)$, this implies
\begin{equation}\label{mob:diag:e_n}
    \B_{1/n}(x,z)\cap\mathcal S^{\kappa_n}_\e(t)
\neq\varnothing
\qquad\text{for every $(x,z,t)\in K_n$,}
\end{equation}
 where we recall that $\B_r(x_0,z_0)\subset \R^{N+1}$ is the Euclidean ball of radius $r$ and center $(x_0,z_0)$. Now define
\[
\kappa(\e):=\kappa_n \qquad\text{if }\e_{n+1}<\e\leq\e_n,
\]
and set $\kappa(\e):=\kappa_1$ for $\e_1<\e<1$. Then $\kappa(\e)\to0$ as $\e\to0^+$, and by \eqref{cfe},  \eqref{mob:diag:clusters} together with the preceding estimates we obtain
\eqref{contr:ke}, \textup{(a)} and \textup{(c)}, simultaneously for the original and complementary constructions.

It remains to prove \textup{(b)}. Let $\mathcal N$ be the union of the countable exceptional sets arising in the proof of Theorem~\ref{mob:thm:main1} for $E_0$ and its complement.
Thus, $\overline{\{u(\cdot,t)<0\}}=E(t)$ for every $t\notin\mathcal N$. The Kuratowski upper inclusion follows immediately from \textup{(a)}. Conversely, fix $t\notin\mathcal N$,
$(x,z)\in E(t)\times\R$, and $\eta>0$. Choose $y$ such that $|y-x|<\eta$ and  $u(y,t)<0$.  By \eqref{mob:main:uk-conv}, for all sufficiently large $n$,
$u_{\kappa_n}(y,t)<-\frac1n$, and hence $(y,z,t)\in K_n$. If $\e_{n+1}<\e\leq\e_n$, \eqref{mob:diag:e_n} provides $(y_\e,z_\e)\in\mathcal S^{\kappa(\e)}_\e(t)$ such that
$(y_\e,z_\e)\in\B_{1/n}(y,z)$, and thus \((y_\e,z_\e)\longrightarrow(y,z)\). Letting first $\e\to0^+$ and then $\eta\to0^+$, we obtain the Kuratowski lower inclusion. This proves
\textup{(b)}, and the same argument applies to the complementary construction.
\end{proof}

\begin{remark}\label{mob:remark:diagonal-complementary}
By the last assertion of Corollary~\ref{mob:cor:diagonal}, all the complementary convergences in Remark~\ref{mob:remark:complementary} hold along the same one-parameter regime
$\kappa=\kappa(\e)$. Consequently, after possibly enlarging the countable exceptional set $\mathcal N$,
\[
\pa\mathcal S_{\e}^{\kappa(\e)}(t) \stackrel{\mathcal K}{\longrightarrow} \pa E(t)\times\R \qquad\text{in }\R^{N+1}
\]
for every $t\notin\mathcal N$. Moreover, setting as in Remark~\ref{mob:remark:complementary}
\[
\Gamma^\kappa_\e:=\mathcal R_\e(\pa E^\kappa_\e), \qquad \Gamma:=\big\{(x,t)\in\R^N\times(0,+\infty):u(x,t)=0\big\},
\]
we have
\[
\Gamma^{\kappa(\e)}_\e \stackrel{\mathcal K}{\longrightarrow} \Gamma \qquad\text{in }\R^N\times(0,+\infty).
\]
\end{remark}

 \subsection{Energy estimates}
 In the following section, we revisit the arguments of Section \ref{sec:ellip:varifold} in the present setting of general mobility, highlighting the necessary modifications.

 Let $E_0\subset \R^N_x$ be a set of finite perimeter. For every $\e>0$ and $\kappa\in (0,\kappa_0)$, we let $E_{\e}^\kappa$ be a minimizer of \eqref{mob:minpb}. We also set
 \begin{equation}\label{mob:def:abbrevio}
 \mu_{\e,\kappa}:=\mu_{E_\e^\kappa}=\H^N\res \partial^* E_\e^\kappa.
 \end{equation}
 Owing to Proposition \ref{prop:Euler_lagrange}, for every $Y\in C^1_c(\{z>0\};\R^{N+1})$, we have
 \begin{equation}\label{mob:EFVEEe}
 \frac1\e \int e^{-z/\e}\, \Phi^\kappa(\nu)\,Y\cdot e_{_{N+1}}\,d\mu_{\e,\kappa}(x,z) =\int e^{-z/\e} B_{\Phi^\kappa}(\nu) : DY d\mu_{\e,\kappa}(x,z),
 \end{equation}
 and by the Euler-Lagrange \eqref{euler:lagrange1}, and recalling that (see \eqref{come:e:ovvia})
 \begin{equation}\label{mobprop:extnorm}
 \nabla \Phi^\kappa(v)\cdot e_{_{N+1}}=\frac{v\cdot e_{_{N+1}}}{\mathfrak{M}^\kappa(v)},\quad\text{for all $v\in \R^{N+1}\setminus\{0\}$,}
 \end{equation}
 we have that $E_\e^\kappa$ admits $\Phi^\kappa$-anisotropic mean curvature satisfying the soliton equation
 \begin{equation}\label{mob:eqsoliton}
 \mathbf{H}^{\Phi^\kappa}_{\partial^* E_\e^\kappa}+\frac{1}{\e}\frac{(\nu\cdot e_{_{N+1}})}{\mathfrak{M}^\kappa(\nu)}\nu=0\quad\text{on $\partial^* E_\e^\kappa\cap \{z>0\}$.}
 \end{equation}
 where $\nu=\nu_{E_\e^\kappa}$ denotes the outer unit normal to $E_\e^\kappa$.

 \noindent From now on, for the sake  of notation, in the following we  write $\mathbf{H}^{\Phi^\kappa}=\mathbf{H}^{\Phi^\kappa}_{\partial^* E_\e^\kappa}$.

 We start with the counterpart of Lemma  \ref{lem:stimEeHp}. Recall that $H^+=\{(x,z)\in \R^{N+1}:\,z>0\}$, and $\mathcal{A}_\e^{\Phi^\kappa}$ is given by
 \eqref{defAF}-\eqref{localized:Apsi} with $\Phi^\kappa$ in place of $\Psi$; we also define the norm $\phi^\kappa$ on $\R^N$ as
 \begin{equation}\label{phi:piccolo:kappa}
 \phi^\kappa(\xi)\coloneqq \Phi^\kappa(\xi,0),\quad\xi\in \R^N,
 \end{equation}
 hence, owing to \eqref{za:tesi},
 \begin{equation}\label{servira:forse}
 \phi^\kappa\xrightarrow{\kappa\to 0^+}\phi\quad\text{locally uniformly in $\R^N$.}
 \end{equation}

\begin{lemma}\label{mob:lem:stimEeHp}
    Let $E_{\e}^{\kappa}$ be a minimizer of \eqref{mob:minpb}. Then, for every $\e>0$ and $\kappa\in (0,\kappa_0)$, there holds
\begin{equation}\label{new:ti:rendi:conto}
    \int_{{\{z>-1\}}} e^{-z/\e}\,\Phi^\kappa(\nu)\,d\mu_{\e,\kappa}(x,z)=\mathcal{A}_\e^{\Phi^\kappa}(E_{\e}^{\kappa})\leq \mathcal{A}_\e^{\Phi^\kappa}(E_0\times \R),
\end{equation}
 where $\nu(x,z)=\nu_{E_{\e}^{\kappa}}(x,z)$ is the normal to $E_{\e}^{\kappa}$,   and
\begin{equation} \label{mob:en-uniform-bound}
    \frac{1}{\e}\mathcal{A}^{\Phi^\kappa}_\e(E_{\e}^{\kappa};H^+)\leq P^{\phi^\kappa}(E_0).
\end{equation}
\end{lemma}

\begin{proof}
   We argue analogously as in the proof of Lemma \ref{lem:stimEeHp}: using $E_0\times \R$  an admissible competitor for the minimization  problem \eqref{mob:minpb}, Equation
   \eqref{new:ti:rendi:conto} immediately follows; then observe that, by \eqref{phi:piccolo:kappa}, $\Phi^\kappa(\nu_{\partial^* E_0\times \R})=\phi^\kappa(\nu_{\partial^*E_0})$,
   and so we may compute
\begin{equation*}
\begin{split}
     \frac{1}{\e}\mathcal{A}^{\Phi^\kappa}_{\e}(E_0\times \R)&= \frac{1}{\e} \Big(\int_{-1}^{0}e^{-z/\e}dz\Big)\,\int_{\partial^*E_0} \phi^\kappa(\nu)d\mathcal{H}^{N-1}+\frac{1}{\e}\Big(\int_0^\infty e^{-z/\e}dz\Big)\int_{\partial^*E_0} \phi^\kappa(\nu)d\mathcal{H}^{N-1}
     \\
     &
     =\frac{1}{\e} \Big(\int_{-1}^{0}e^{-z/\e}dz\Big)\,P^{\phi^\kappa}(E_0)+P^{\phi^\kappa}(E_0),
     \end{split}
\end{equation*}
and again, as $E_\e^\kappa=E_0\times (-\infty,0)$ in $\{z<0\}$ by \eqref{mob:minpb}, we have
\begin{equation*}
    \frac{1}{\e} \mathcal{A}^{\Phi^\kappa}_{\e}(E_\e^\kappa)\geq \frac{1}{\e}\Big(\int_{-1}^0e^{-z/\e}dz\Big)\,P^{\phi^\kappa}(E_0)+\frac{1}{\e}\mathcal{A}^{\Phi^\kappa}_{\e}(E_\e^\kappa;H^+).
\end{equation*}
Thereby combining the two expressions above with \eqref{new:ti:rendi:conto}, we get \eqref{mob:en-uniform-bound}, that is our thesis.
\end{proof}

 Let us now prove the following
\begin{lemma}\label{mob:lem:possible:test}\label{moblem:1}\label{mob:lem:2}
Let $E_\e^\kappa$ be a minimizer of \eqref{mob:minpb}. For every Borel set $A\subset(0,\infty)$ of finite measure we have
\[
\int_{\R^N\times A}\left(\Phi^\kappa(\nu)-\frac{(\nu\cdot e_{_{N+1}})^2}{\mathfrak M^\kappa(\nu)}\right)\,d\mu_{\e,\kappa}\leq |A|P^{\phi^\kappa}(E_0),
\]
and
\begin{equation}\label{mob:stima:eneH2}
 \int_{H^+} \e |\mathbf{H}^{\Phi^\kappa}|^2  \, \mathfrak M^\kappa(\nu) d\mu_{\e,\kappa} \leq P^{\phi^\kappa}(E_0).
 \end{equation}
\end{lemma}
\begin{proof}
By \eqref{mobprop:extnorm}, \eqref{def:Bpsi}, and \eqref{meno:ovvia}, we have $0\leq B_{\Phi^\kappa}(\nu)e_{_{N+1}}\cdot e_{_{N+1}}=\Phi^\kappa(\nu)-\frac{(\nu\cdot e_{_{N+1}})^2}{\mathfrak
M^\kappa(\nu)}\leq\Phi^\kappa(\nu).$ Using the same vertical variations as in Lemma~\ref{lem:possible:test} and arguing similarly,  we obtain
\begin{equation}\label{mob:alEL}
\int_{H^+}e^{-z/\e}\left(\Phi^\kappa(\nu)-\frac{(\nu\cdot e_{_{N+1}})^2}{\mathfrak M^\kappa(\nu)}\right)\,\xi'(z)\,d\mu_{\e,\kappa}
=\frac1\e\int_{H^+}e^{-z/\e}\Phi^\kappa(\nu)\,\xi(z)\,d\mu_{\e,\kappa},
\end{equation}
for $\xi\in C^1([0,\infty))$ with $\xi(0)=0$ and both $\xi$ and $\xi'$ bounded. As in the proof of \eqref{altre:eulero:lagrange1}, we may take $\xi(z)=e^{z/\e}\eta(z)-\eta(0)$ for
$\eta\in C_c^1([0,\infty))$ nonnegative to get
\begin{equation}\label{mob:AEL1}
\frac1\e\int_{H^+}\eta(z)\frac{(\nu\cdot e_{_{N+1}})^2}{\mathfrak M^\kappa(\nu)}\,d\mu_{\e,\kappa}
-\int_{H^+}\eta'(z)\left(\Phi^\kappa(\nu)-\frac{(\nu\cdot e_{_{N+1}})^2}{\mathfrak M^\kappa(\nu)}\right)\,d\mu_{\e,\kappa}
=\frac{\eta(0)}{\e}\mathcal A_\e^{\Phi^\kappa}(E_\e^\kappa;H^+)\leq \eta(0) P^{\phi^\kappa}(E_0),
\end{equation}
where we used also \eqref{mob:en-uniform-bound}. We can now proceed exactly as in the proof of Lemma~\ref{lem:possible:test} and use also that $(\nu\cdot e_{_{N+1}})^2/(\e\mathfrak
M^\kappa(\nu))=\e|\mathbf H^{\Phi^\kappa}|^2\mathfrak M^\kappa(\nu)$ (by \eqref{mob:eqsoliton}) to get  the two asserted estimates.
\end{proof}

 In the remaining of this section, we will need a crucial control on the area of the non-vertical part $\Sigma^{(\e,\kappa),\gamma}$ of the minimizers $E_\e^\kappa$, defined as
 \begin{equation}\label{def:noncilindrico}
 \Sigma^{(\e,\kappa),\gamma}:=\big\{ (x,z)\in \partial^* E_\e^\kappa\cap H^+: (\nu\cdot e_{_{N+1}})^2\geq \gamma\big\}.
 \end{equation}
 for every $\gamma\in (0,1]$. This is proved in the following lemma, which is inspired by \cite[Lemma 4.1]{Tashiro}:

 \begin{lemma}\label{lem:quantif:cil}
 Under the notations and assumptions above,   for every $\e>0, \kappa \in (0,\kappa_0)$ and every $\gamma\in (0,1]$, we have
 \begin{equation}\label{bello:quantif}
 \H^{N}(\Sigma^{(\e,\kappa),\gamma})\leq \e\,\gamma^{-1}\,\|\mathfrak{M}^\kappa\|_{L^\infty(\S^N)}\,P^{\phi^\kappa}(E_0)\xrightarrow{\e\to 0^+}0.
 \end{equation}
 \end{lemma}

 \begin{proof}
 Owing to \eqref{mob:eqsoliton} and Chebyshev's inequality, we have
 \begin{equation*}
 \begin{split}
 \H^{N}(\Sigma^{(\e,\kappa),\gamma})&=\mathcal{H}^N \Big(\{(x,z)\in \partial^* E_\e^\kappa {\cap H^+}:|\mathbf{H}^{\Phi^\kappa}|^2\geq \gamma/(\e^2\,(\mathfrak M^\kappa(\nu))^2)\big\} \Big)\\
 &\leq \e\,\gamma^{-1}\,\int_{\partial^* E_\e^\kappa {\cap H^+}} \e\,|\mathbf{H}^{\Phi^\kappa}|^2\,(\mathfrak M^\kappa(\nu))^2\,d\H^{N} \\
 &\leq \e\,\gamma^{-1}\,\|\mathfrak{M}^\kappa\|_{L^\infty(\S^N)}\,\int_{H^+} \e\,|\mathbf{H}^{\Phi^\kappa}|^2\,\mathfrak M^\kappa(\nu)\,d\mu_{\e,\kappa}\stackrel{\eqref{mob:stima:eneH2}}{\leq} \e\,\gamma^{-1}\,\|\mathfrak{M}^\kappa\|_{L^\infty(\S^N)}\,P^{\phi^\kappa }(E_0).
 \end{split}
 \end{equation*}
 \end{proof}

 We then  prove the key energy estimates for the minimizer $E_\e^\kappa$

 \begin{proposition}\label{mob:prop:energyestimate}
 For every measurable set $A\subset (0,\infty)$ with $|A|<\infty$, the following perimeter estimates hold:
 \begin{equation} \label{mob:eq:energyestimate}
 P^{\Phi^\kappa}(E_\e^\kappa;\R^N\times A) \leq  (|A|+\e) \,P^{\phi^\kappa }(E_0),
 \end{equation}
 \begin{equation}\label{dimezzo:perimetro}
 P(E_\e^\kappa;\,\R^N\times A)\leq \Big[2\,\e\,\|\mathfrak{M}^\kappa\|_{L^\infty(\S^N)}+C_\phi\,(|A|+\e)\Big]\,P^{\phi^\kappa}(E_0),
 \end{equation}
 \begin{equation}\label{perimetro:mk}
 P^{\mathfrak M^\kappa}(E_\e^\kappa;\R^N\times A)\leq C_{\mathfrak{m},\phi}\,\Big[\e\,(L_\kappa)^{4}+  (|A|+\e)\,\Big]\,P^{\phi^\kappa }(E_0),
 \end{equation}
 where $L_\kappa$ is given by \eqref{contr:ke}, and    \begin{equation}\label{mob:eq:rescaled:energy} P\big(\mathcal{R}_\e(E_\e^\kappa);\R^N\times A\big)\leq
 C_{\mathfrak{m},\phi}\,\Big(\e\,L_\kappa^{2}+(\e^2+|A|)+(\e^2+|A|)^{1/2}\Big)\,P^{\phi^\kappa }(E_0).
 \end{equation}
 for all $\e\in(0,1)$ and all $\kappa<\tilde \kappa$, where $\tilde \kappa$  is a small parameter depending on $\mathfrak{m},\phi$, with $C_\phi$  a constant depending only on
 $N,\phi$, and $C_{\mathfrak{m},\phi}$ depending only  on $N,\phi,\mathfrak{m}$.
 \end{proposition}

 \begin{proof}

 Let $A\subset (0,\infty)$, and recalling \eqref{mob:def:abbrevio},  we split
 \begin{equation*}
 \begin{split}
 P^{\Phi^\kappa}(E_\e^\kappa;\R^N \times A) &=  \int_{\R^N\times A} \left(1-\frac{(\nu \cdot e_{_{N+1}})^2}{\Phi^\kappa(\nu)\,\mathfrak M^\kappa(\nu)} \right)\Phi^\kappa(\nu) d\mu_{\e,\kappa} + \int_{\R^N\times A} \frac{(\nu \cdot e_{_{N+1}})^2}{\Phi^\kappa(\nu)\,\mathfrak M^\kappa(\nu)} \Phi^\kappa(\nu) d\mu_{\e,\kappa}\\
 &\overset{\eqref{mob:eqsoliton}}{=}\int_{\R^N\times A} \left(1-\frac{(\nu \cdot e_{_{N+1}})^2}{\Phi^\kappa(\nu)\mathfrak M^\kappa(\nu)} \right)\Phi^\kappa(\nu) d\mu_{\e,\kappa} + \int_{\R^N\times A} \e^2|\mathbf{H}^{\Phi^\kappa}|^2 \mathfrak M^\kappa(\nu) d\mu_{\e,\kappa}.
 \end{split}
 \end{equation*}
Thus by the two bounds in Lemma~\ref{mob:lem:2} we obtain \eqref{mob:eq:energyestimate}.

 Now let $\gamma\in (0,1/2]$, and observe that, if $\nu=(\nu',\nu\cdot e_{_{N+1}})\in \mathbb{S}^N$ ($\nu'\in \R^{N}$) is such that $|\nu\cdot e_{_{N+1}}|^2\leq \gamma$, then
 $1=|\nu|^2\leq 2\,|\nu'|^2\leq 2$; thus by \eqref{za:tesi} and the equivalence of norms, there exists  $\kappa_1$ depending on $\phi$ small enough such that
 \[
 \Phi^\kappa (\nu)\geq \frac{1}{2}\phi(\nu')\geq  C_\phi\,|\nu'|\geq C'_\phi\,|\nu|,\quad \text{for all }\,\nu\in \S^N:\,|\nu\cdot e_{_{N+1}}|^2\leq \gamma,
 \]
 and for some constant $C_\phi$ depending on $\phi$, but independent of $\kappa$ and $\gamma$.  It follows that, if we split
 \[
 \big(\R^N\times A\big)=\bigg(\big(\R^N\times A\big)\cap \Sigma^{(\e,\kappa),\gamma}\bigg)\cup\bigg(\big(\R^N\times A\big)\setminus\Sigma^{(\e,\kappa),\gamma}\bigg)=: A^\gamma_1\cup A^\gamma_2,
 \]
 then by using the above information,  the definition of $\Sigma^{(\e,\kappa),\gamma}$ in \eqref{def:noncilindrico},  \eqref{bello:quantif} and \eqref{mob:eq:energyestimate} we
 find
 \begin{equation*}
 \begin{split}
 \mu_{\e,\kappa}(\R^N\times A)&=\mu_{\e,\kappa}(A^\gamma_1)+\mu_{\e,\kappa}(A^\gamma_2)\leq \mathcal{H}^N(\Sigma^{(\e,\kappa),\gamma})+ C_\phi\,P^{\Phi^\kappa}(E_\e^\kappa;\,A^\gamma_2) \\
 &\leq\Big[\e\,\gamma^{-1}\,\|\mathfrak{M}^\kappa\|_{L^\infty(\S^N)}+C_\phi\,(|A|+\e)\Big]\,P^{\phi^\kappa}(E_0).
 \end{split}
 \end{equation*}
 from which we deduce \eqref{dimezzo:perimetro} after choosing $\gamma=1/2$. Next we observe that, thanks to \eqref{za:tesi}, there exists $\tilde\kappa<\kappa_1$ depending on
 $\mathfrak{m},\phi$ such that $\mathfrak{M}^\kappa(\nu',0)\leq 2\,\mathfrak{m}(\nu')$ for all $\kappa<\tilde \kappa$, so we have
 \begin{equation*}
 \mathfrak{M}^\kappa(\nu',\nu\cdot e_{_{N+1}})\leq 2 \,\mathfrak{m}(\nu')+\mathrm{Lip}(\mathfrak{M}^\kappa)\,|\nu\cdot e_{_{N+1}}|\quad\text{for all $\kappa\in (0,\tilde \kappa)$.}
 \end{equation*}
 Therefore, we may estimate
 \begin{equation*}
 \begin{split}
 P^{\mathfrak M^\kappa}(E_\e^\kappa,\R^N\times A)&=\int_{\R^N\times A}\mathfrak M^\kappa(\nu) d\mu_{\e,\kappa} =\int_{A_1^\gamma}\mathfrak M^\kappa(\nu)d\mu_{\e,\kappa}+\int_{A_2^\gamma}\mathfrak M^\kappa(\nu) d\mu_{\e,\kappa} \\
 &\leq \H^N(\Sigma^{(\e,\kappa),\gamma})\,\|\mathfrak{M}^\kappa\|_{L^\infty(\S^N)}+2\,\int_{A_2^\gamma}\mathfrak{m}(\nu')\,d\mu_{\e,\kappa}+\mathrm{Lip}(\mathfrak{M}^\kappa) \int_{A_2^\gamma}|\nu\cdot e_{_{N+1}}|d\mu_{\e,\kappa} \\
 &\stackrel{\eqref{bello:quantif}}{\leq }\e\,\|\mathfrak{M}^\kappa\|_{L^\infty(\S^N)}^{2}\,\gamma^{-1}\,P^{\phi^\kappa}(E_0)+\left(2\,\|\mathfrak{m}\|_{L^\infty(\S^{N-1})}+\mathrm{Lip}(\mathfrak{M}^\kappa)\,\gamma^{1/2}\right)\,P(E_\e^\kappa,\R^N\times A) \\
 \stackrel{\eqref{dimezzo:perimetro}}{\leq}\Bigg\{\e\,\|\mathfrak{M}^\kappa\|_{L^\infty(\S^N)}^{2}&\,\gamma^{-1}+\Big(2\,\e\,\|\mathfrak{M}^\kappa\|_{L^\infty(\S^N)}+C_\phi\,(|A|+\e)\Big)\,\Big(2\,\|\mathfrak{m}\|_{L^\infty(\S^{N-1})}+\mathrm{Lip}(\mathfrak{M}^\kappa)\,\gamma^{1/2}\Big)\Bigg\}P^{\phi^\kappa}(E_0)
 \end{split}
 \end{equation*}
 where in the second inequality we used $|\nu\cdot e_{_{N+1}}|^2\leq\gamma $  on $A_2^\gamma$, and in the last inequality we used \eqref{dimezzo:perimetro}. By choosing
 $\gamma=\min\{\mathrm{Lip}(\mathfrak{M}^\kappa)^{-2},1/2\}$, recalling the definition $L_\kappa=1+\|\mathfrak{M}^\kappa\|_{L^\infty(\S^N)}+\operatorname{Lip}(\mathfrak M^\kappa)$
 and via simple algebraic manipulation we obtain \eqref{perimetro:mk}.

 Next, via a change of variables and \eqref{mob:eqsoliton}, we get
 \begin{equation*}
 \begin{split}
 P\big(\mathcal{R}_\e(E_\e^\kappa); \R^N\times A \big)&=\int_{\R^N\times \frac{A}{\e}}\sqrt{\e^2\,|\nu'|^2+|\nu\cdot e_{_{N+1}}|^2}\,d\mu_{\e,\kappa}=\e\,\int_{\R^N\times \frac{A}{\e}}\sqrt{|\nu'|^2+(\mathfrak M^\kappa(\nu))^2\,|\mathbf{H}^{\Phi^\kappa}|^2}\,d\mu_{\e,\kappa} \\
 &\leq \e\, P\Big(E_\e^\kappa;\, \R^N\times \frac{A}{\e} \Big)+\e\,\int_{\R^N\times \frac{A}{\e}}\mathfrak M^\kappa(\nu)\,|\mathbf{H}^{\Phi^\kappa}|\,d\mu_{\e,\kappa}\eqqcolon (I)+(II)
 \end{split}
 \end{equation*}
 By \eqref{dimezzo:perimetro} with $A/\e$  in place of $A$
 \begin{equation*}
 (I)\leq \Big[2\,\e^2\,\|\mathfrak{M}^\kappa\|_{L^\infty(\S^N)}+C_\phi(|A|+\e^2)\Big]\, P^{\phi^\kappa}(E_0),
 \end{equation*}
 while by  H\"older's inequality, \eqref{mob:stima:eneH2} and \eqref{perimetro:mk} with $A/\e$ in place of $A$, we get
 \begin{equation*}
 \begin{split}
 (II) &\leq  \Big(\e\,\int_{\R^N\times A/\e}\mathfrak M^\kappa(\nu) d\mu_{\e,\kappa} \Big)^{1/2}\Big(\e\,\int_{H^+}\mathfrak M^\kappa(\nu)\,|\mathbf{H}^{\Phi^\kappa}|^2\,d\mu_{\e,\kappa}\Big)^{1/2} \\
 &\leq \Big(\e\,P^{\mathfrak{M}^\kappa}({E_\e^\kappa;}\R^N\times A/\e)\Big)^{1/2}\,\big(P^{\phi^\kappa}(E_0)\big)^{1/2} \leq C_{\mathfrak{m},\phi}\,\Big[\e^2\,(L_\kappa)^{4}+  \,(|A|+\e^2)\,\Big]^{1/2}P^{\phi^\kappa}(E_0)
 \end{split}
 \end{equation*}
 Combining the three estimates above and performing some straightforward algebraic manipulations, we finally obtain \eqref{mob:eq:rescaled:energy}.

 \end{proof}

\subsection{Convergence of rescalings, translating solitons and space-time tracks}
\label{mob:completion:bvconvergence}

Throughout the remainder of this section, we assume that
\begin{equation*}
    \text{$E_0$ is a set of finite perimeter that is \textit{non-fattening in measure} for the pair $(\phi,\mathfrak{m})$}
\end{equation*}
in the sense of Definition \ref{def:fat}. We let $E^\kappa_\e$ be a minimizer of \eqref{mob:minpb}, and $\mathcal S^\kappa_\e$ be the corresponding solitons as in \eqref{Ssolliton}. We denote by $\{\widehat E(t)\}_{t\ge 0}$ the corresponding (unique) flat flow starting from $E_0$ as in Section \ref{subses:flat:flow}, and set
\[
 \widehat E:=\{(x,t)\in\R^N\times(0,\infty):x\in\widehat E(t)\}.
\]
Besides $\phi^\kappa$ in \eqref{phi:piccolo:kappa}, let us recall the definition $\mathfrak m^\kappa(p):=\mathfrak M^\kappa(p,0)$ and  $L_\kappa:=1+\|\mathfrak
M^\kappa\|_{L^\infty(\S^N)} +\Lip(\mathfrak M^\kappa).$ In particular, $\phi^\kappa\to\phi$ and $\mathfrak m^\kappa\to\mathfrak m$ locally uniformly. Let us now prove the
counterpart of Theorem  \ref{thm:altro:main}.

\begin{theorem}\label{mob:completion:diagonal}
Under the notation and assumptions above, let $\kappa(\e)$ be as in Corollary~\ref{mob:cor:diagonal}, chosen simultaneously for the original and complementary constructions. Then
the following limits hold:
\begin{itemize}
\item[(a)] there exists a set $\mathcal N\subset(0,+\infty)$
of null measure such that
\begin{equation}\label{mob:completion:solitonBV}
\mathcal S_\e^{\kappa(\e)}(t)
\xrightarrow{\e\to0^+}
\widehat E(t)\times\R
\quad\text{weakly-* in $BV_{\mathrm{loc}}(\R^{N+1})$
for every $t\in(0,+\infty)\setminus\mathcal N$;}
\end{equation}

\item[(b)] the rescaled sets satisfy
\begin{equation}\label{mob:completion:rescaledBV}
\mathcal R_\e\big(E_\e^{\kappa(\e)}\big)
\xrightarrow{\e\to0^+}
\widehat E
\quad\text{weakly-* in
$BV_{\mathrm{loc}}(\R^N\times(0,+\infty))$,}
\end{equation}
and we have the estimate
\begin{equation}\label{mob:completion:tubebound}
P\big(\widehat E;\R^N\times A\big)
\leq
C_{\mathfrak m,\phi}
\bigl(|A|+|A|^{1/2}\bigr)P^\phi(E_0)
\end{equation}
for every Borel set $A\subset(0,+\infty)$ with $|A|<+\infty$, where $C_{\mathfrak m,\phi}$ depends only on $N,\phi$ and $\mathfrak m$.
\end{itemize}
\end{theorem}

\begin{proof}
Along the diagonal of Corollary~\ref{mob:cor:diagonal} we have $\e L_{\kappa(\e)}^2\to0$ by \eqref{contr:ke}, while $P^{\phi^{\kappa(\e)}}(E_0)\to P^\phi(E_0)$ by
\eqref{servira:forse}. For fixed $t>0$ and a bounded interval $I\subset\R$, apply \eqref{dimezzo:perimetro} to $I+t/\e\subset(0,\infty)$ for $\e$ small. Since this interval has
length $|I|$, we obtain a uniform perimeter bound for the translated solitons on $\R^N\times I$. By \eqref{mob:eq:rescaled:energy} we obtain the corresponding local $BV$ bound for
the rescaled sets in positive time.

By Corollary~\ref{mob:cor:diagonal} and Remark~\ref{mob:remark:diagonal-complementary} we may realize the Kuratowski limits of the sets and their complements along the same
diagonal, outside a countable set of times for the soliton slices. Taking into account Remark~\ref{rem:Lmin} we can  apply Lemma~\ref{lemma:L1Kur} to both limits exactly as in the
proof of Theorem~\ref{thm:altro:main}. In turn, non-fattening in measure allows us to identify  every subsequential $L^1_{\mathrm{loc}}$ limit with $\widehat E(t)\times\R$ for a.e.
$t>0$, and with $\widehat E$ for the rescaled sets. This proves \eqref{mob:completion:solitonBV} and \eqref{mob:completion:rescaledBV} along the whole diagonal. Finally,
\eqref{mob:completion:tubebound} follows by lower semicontinuity in \eqref{mob:eq:rescaled:energy} for any open set $A$ and by outer approximation for any Borel $A$.
\end{proof}

\subsection{General transport and BV identities}
We are in the position to prove the counterpart of Theorem \ref{trasportobv} in the present setting with general mobility.

\begin{theorem}[Transport and BV identities]\label{mob:trasportobv}
Let $\phi,\mathfrak m$ be norms on $\R^N$ (possibly crystalline). Let $E_0\subset\R^N$ have finite perimeter and be non-fattening in measure with respect to $(\phi,\mathfrak m)$, and
let $\{\widehat E(t)\}_{t\geq0}$ be the corresponding flat flow.

Set
\[
\widehat E
:=\{(x,t)\in\R^N\times(0,\infty):x\in\widehat E(t)\},
\qquad
d\mu:=d\mu_{\widehat E(t)}\,dt,
\qquad
d\sigma:=d\mu^\phi_{\widehat E(t)}\,dt.
\]

Then there exist scalar fields $v\in L^2(\mu)$ and $h\in L^\infty(\sigma)$, a vector field $\mathbf b\in L^2(\sigma;\R^N)$, and a finite nonnegative Radon measure $\Lambda$ on
$\R^N\times(0,\infty)$, with
\begin{equation}\label{mob:h-bounds}
0<
\frac{\mathfrak m(\nu_{\widehat E(t)}(x))}
     {\phi(\nu_{\widehat E(t)}(x))}
\leq h(x,t)\leq\max_{\xi\in\mathbb S^{N-1}}
\frac{\mathfrak m(\xi)}{\phi(\xi)}
\qquad\text{for }\sigma\text{-a.e. }(x,t),
\end{equation}
such that the following properties hold.
\begin{itemize}
\item[i)]
The tube $\widehat E$ has locally finite perimeter, and
\begin{equation}\label{mob:tubofase}
D_{x,t}\chi_{\widehat E}
=(-\nu_{\widehat E(t)},v)\,\mu,
\qquad
\mathbf b\cdot\nu_{\widehat E(t)}=v
\quad\mu\text{-a.e.}
\end{equation}
In particular, for every $0\leq t_1<t_2$ and every $C^1$ test function $\eta$ with compact spatial support uniformly on $[t_1,t_2]$, one has
\begin{equation}\label{mob:trasportofase}
\int_{\widehat E(t_2)}\eta(x,t_2)\,dx
-\int_{\widehat E(t_1)}\eta(x,t_1)\,dx
=
\int_{t_1}^{t_2}\int_{\widehat E(t)}
\partial_t\eta\,dx\,dt+\int_{t_1}^{t_2}\int
\eta v\,d\mu_{\widehat E(t)}\,dt.
\end{equation}

\item[ii)]
The energy transport equation
\begin{equation}\label{mob:bilancioenergia}
\partial_t\sigma+\operatorname{div}_x(\mathbf b\sigma)
=-\Lambda
\end{equation}
holds in the sense of distributions on $\R^N\times(0,\infty)$, and the dissipation estimates hold:
\begin{equation}\label{mob:dissipazione}
\frac{\bigl(\phi^o(\mathbf b)\bigr)^2}{h}\,\sigma
\leq\Lambda,
\qquad
\Lambda(\R^N\times(0,\infty))
\leq P^\phi(E_0).
\end{equation}
More precisely, for every $0\leq t_1<t_2$ and every $C^1$ test function $\eta$ as above, one has
\begin{equation}\label{mob:bilanciointervalli}
\begin{aligned}
&\mu^\phi_{\widehat E(t_2)}(\eta(\cdot,t_2))
-\mu^\phi_{\widehat E(t_1)}(\eta(\cdot,t_1))\\
&\qquad=
\int_{t_1}^{t_2}\int
(\partial_t\eta+\mathbf b\cdot\nabla_x\eta)
\,d\mu^\phi_{\widehat E(t)}\,dt
-\int_{\R^N\times(t_1,t_2]}\eta\,d\Lambda.
\end{aligned}
\end{equation}
\end{itemize}
\end{theorem}

\begin{proof}
We follow the proof of Theorem~\ref{trasportobv} step by step, highlighting the necessary modifications.
    \\
    \noindent\textit{Step 1: Regularized norms and minimizers.}

Fix a starting time $t_0\geq0$, and set $G_0=\widehat E(t_0)$. In the present setting, the norms $\Phi^\kappa$ and $\mathfrak M^\kappa$ provided by
Proposition~\ref{prop:extended-norms} are already smooth away from the origin and uniformly convex. Thus, the additional approximation introduced in Step 1 of the proof of
Theorem~\ref{trasportobv} is not needed. For every $\e>0$ and $\kappa\in(0,\kappa_0)$, let $G_\e^\kappa$ minimize \eqref{mob:minpb} with initial datum $G_0$. Recall also that
$\phi^\kappa\to\phi$ locally uniformly, and hence
\[
P^{\phi^\kappa}(G_0)\longrightarrow P^\phi(G_0)
\qquad\text{as $\kappa\to0^+$.}
\]

\medskip
\noindent\textit{Step 2: Uniform estimates and convergence of the rescaled sets.}

In computations on $G_\e^\kappa$, we write $\nu=(\nu',\nu\cdot e_{_{N+1}})$ for its outer normal and abbreviate $H^{\Phi^\kappa}=H^{\Phi^\kappa}_{\partial^*G_\e^\kappa}$. By
\eqref{mob:eqsoliton},
\begin{equation}\label{mob:solitonbis}
H^{\Phi^\kappa}
=
\frac{\nu\cdot e_{_{N+1}}}
     {\e\,\mathfrak M^\kappa(\nu)}
\qquad\text{on $\partial^*G_\e^\kappa\cap H^+$.}
\end{equation}

For every Borel set $A\subset(0,+\infty)$ of finite measure, by \eqref{mob:stima:eneH2}, \eqref{mob:eq:energyestimate}, and \eqref{perimetro:mk} we obtain
\begin{equation}\label{mob:1}
\begin{split}
&\e\int_{H^+}|H^{\Phi^\kappa}|^2
       \mathfrak M^\kappa(\nu)\,d\mu_{G_\e^\kappa}
       \leq P^{\phi^\kappa}(G_0),\\
&P^{\Phi^\kappa}(G_\e^\kappa;\R^N\times A)
       \leq (|A|+\e)P^{\phi^\kappa}(G_0),\\
&P^{\mathfrak M^\kappa}(G_\e^\kappa;\R^N\times A)
       \leq C_{\mathfrak m,\phi}
       \bigl(\e L_\kappa^{4}+|A|+\e\bigr)
       P^{\phi^\kappa}(G_0),
\end{split}
\end{equation}
where $L_\kappa$ is defined in \eqref{contr:ke}. Moreover, \eqref{mob:eq:rescaled:energy} yields
\begin{equation}\label{mob:2}
P\bigl(\mathcal R_\e(G_\e^\kappa);\R^N\times A\bigr)
\leq C_{\mathfrak m,\phi}
\Bigl(\e L_\kappa^{2}+|A|+\e^2
      +(|A|+\e^2)^{1/2}\Bigr)
P^{\phi^\kappa}(G_0).
\end{equation}
The constant $C_{\mathfrak m,\phi}$ depends only on $N,\phi,\mathfrak m$, and the estimates hold for $\e\in(0,1)$ and all sufficiently small $\kappa$.

Now choose $\e_j\searrow0$ and set $\kappa_j:=\kappa(\e_j)$, where $\kappa(\e)$ is provided by Corollary~\ref{mob:cor:diagonal} applied with initial datum $G_0$. Then
Theorem~\ref{mob:completion:diagonal} yields
\begin{equation}\label{mob:3}
\mathcal R_{\e_j}(G_{\e_j}^{\kappa_j})+(0,t_0)
\longrightarrow\widehat E
\quad\text{weakly-* in }
BV_{\loc}\bigl(\R^N\times(t_0,+\infty)\bigr).
\end{equation}

Passing to the limit in \eqref{mob:2}, and recalling that $\e\,(L_{\kappa(\e)})^2\to 0$ by \eqref{contr:ke}, by lower semicontinuity of the perimeter we infer that the limiting
tube $\widehat{E}$ has locally finite perimeter.

\medskip
\noindent\textit{Step 3: Volume and energy identities before passage to the limit.}

In the following steps, all limits are taken along $(\e_j,\kappa_j)$, or along further subsequences. For notational simplicity, we also set
\[
G_j:=G_{\e_j}^{\kappa_j},\qquad
\nu_j:=\nu_{G_j},\qquad
H_j:=H^{\Phi^{\kappa_j}}_{\partial^*G_j}.
\]
Let $\eta\in C_c^1(\R^N\times(t_0,\infty))$. Changing variables and applying the Gauss--Green formula for sets of finite perimeter yields
\begin{equation}\label{mob:4}
\begin{split}
-\partial_t\chi_{\mathcal R_{\e_j}(G_j)+(0,t_0)}(\eta)
&=\int_{\mathcal R_{\e_j}(G_j)+(0,t_0)}
       \partial_t\eta\,dx\,dt=\int_{G_j\cap H^+}
       \partial_z[\eta(x,t_0+\e_jz)]\,dx\,dz\\
&=\int_{H^+}\eta(x,t_0+\e_jz)
       (\nu_j\cdot e_{_{N+1}})\,d\mu_{G_j}=\e_j\int_{H^+}\eta(x,t_0+\e_jz)H_j\,
       d\mu_{G_j}^{\mathfrak M^{\kappa_j}},
\end{split}
\end{equation}
where the last equality follows from \eqref{mob:solitonbis}. Consequently, by Cauchy--Schwarz, for every nonnegative $\eta\in C_c(\R^N\times(t_0,\infty))$ we have
\begin{equation}\label{mob:5}
\begin{split}
&\left|
\partial_t\chi_{\mathcal R_{\e_j}(G_j)+(0,t_0)}
\right|(\eta)\\
&\quad\leq
\left(
\e_j\int_{H^+}\eta(x,t_0+\e_jz)\,
d\mu_{G_j}^{\mathfrak M^{\kappa_j}}
\right)^{1/2}
\left(
\e_j\int_{H^+}\eta(x,t_0+\e_jz)|H_j|^2\,
d\mu_{G_j}^{\mathfrak M^{\kappa_j}}
\right)^{1/2}.
\end{split}
\end{equation}

Next, the first variation identity reads
\[
\int_{H^+}B_{\Phi^{\kappa_j}}(\nu_j):DY\,d\mu_{G_j}
=
\int_{H^+}H_j\,\nu_j\cdot Y\,d\mu_{G_j},
\qquad Y\in C_c^1(H^+;\R^{N+1}).
\]
Choosing $Y(x,z)=\eta(x,t_0+\e_jz)e_{_{N+1}}$, and recalling \eqref{def:Bpsi} and \eqref{mobprop:extnorm}, we compute
\[
\begin{split}
B_{\Phi^{\kappa_j}}(\nu_j):DY
&=\e_j\left(
\Phi^{\kappa_j}(\nu_j)
-\frac{(\nu_j\cdot e_{_{N+1}})^2}
       {\mathfrak M^{\kappa_j}(\nu_j)}
\right)\partial_t\eta
-(\nu_j\cdot e_{_{N+1}})
\nabla_\xi\Phi^{\kappa_j}(\nu_j)\cdot\nabla_x\eta.
\end{split}
\]
Here and in the following energy identity, $\eta$ and its derivatives are evaluated at $(x,t_0+\e_jz)$. Using again \eqref{mob:solitonbis}, and substituting into the first
variation identity, we obtain
\begin{equation}\label{mob:6}
\begin{split}
0=
\e_j\int_{H^+}\partial_t\eta\,
d\mu_{G_j}^{\Phi^{\kappa_j}}-\e_j\int_{H^+}
\Big(
H_j\nabla_\xi\Phi^{\kappa_j}(\nu_j)\cdot\nabla_x\eta
+|H_j|^2\eta
+\e_j^2|H_j|^2\partial_t\eta
\Big)\,
d\mu_{G_j}^{\mathfrak M^{\kappa_j}}.
\end{split}
\end{equation}

Finally, by convexity and one-homogeneity of $\Phi^{\kappa_j}$, for every $w\in\R^N$ we have
\[
\nabla_\xi\Phi^{\kappa_j}(\nu_j)\cdot w
\leq \Phi^{\kappa_j}(w,0)
=\phi^{\kappa_j}(w).
\]
By the definition of the dual norm, we obtain
\begin{equation}\label{mob:7}
\bigl(\phi^{\kappa_j}\bigr)^o
\bigl(\nabla_\xi\Phi^{\kappa_j}(\nu_j)\bigr)\leq1.
\end{equation}
In particular, since $\phi^{\kappa_j}\to\phi$ locally uniformly, we also obtain
\[
\bigl|\nabla_\xi\Phi^{\kappa_j}(\nu_j)\bigr|
\leq C_\phi
\]
for all sufficiently large $j$, with $C_\phi$ independent of $j$.

\medskip
\noindent\textit{Step 4: Limiting measures and transport identity before restarting.} Set
\[
\mathcal T_j(x,z):=(x,t_0+\e_jz),
\]
and define the following nonnegative Radon measures on $\R^N\times(t_0,\infty)$ by
\begin{equation}\label{mob:measures}
\begin{split}
\lambda_j
&:=\e_j(\mathcal T_j)_\#
\left(\mu_{G_j}^{\Phi^{\kappa_j}}\mres H^+\right),\\
\rho_j
&:=\e_j(\mathcal T_j)_\#
\left(\mu_{G_j}^{\mathfrak M^{\kappa_j}}\mres H^+\right),\\
\Lambda_j
&:=\e_j(\mathcal T_j)_\#
\left(|H_j|^2\mu_{G_j}^{\mathfrak M^{\kappa_j}}
\mres H^+\right).
\end{split}
\end{equation}
Rescaling the estimates in \eqref{mob:1}, we obtain
\begin{equation}\label{mob:questaserve}
\begin{split}
&\Lambda_j\bigl(\R^N\times(t_0,\infty)\bigr)
\leq P^{\phi^{\kappa_j}}(G_0),\\
&\lambda_j(\R^N\times I)
\leq (|I|+\e_j^2)P^{\phi^{\kappa_j}}(G_0),\\
&\rho_j(\R^N\times I)
\leq C_{\mathfrak m,\phi}
\left((\e_jL_{\kappa_j}^{2})^2+|I|+\e_j^2\right)
P^{\phi^{\kappa_j}}(G_0)
\end{split}
\end{equation}
for every bounded open interval $I\Subset(t_0,\infty)$. Since $\e_jL_{\kappa_j}^{2}\to0$ by \eqref{contr:ke} and $P^{\phi^{\kappa_j}}(G_0)\to P^\phi(G_0)$, we may extract a
non-relabeled subsequence such that
\[
\lambda_j\weakstar\lambda,\qquad
\rho_j\weakstar\rho,\qquad
\Lambda_j\weakstar\Lambda
\quad\text{locally in }\R^N\times(t_0,\infty),
\]
where $\lambda,\rho,\Lambda$ are nonnegative Radon measures. Moreover,
\begin{equation}\label{mob:forseserve}
\Lambda\bigl(\R^N\times(t_0,\infty)\bigr)
\leq P^\phi(G_0),
\qquad
\lambda(\R^N\times I)\leq |I|P^\phi(G_0).
\end{equation}
In particular, the time marginal of $\lambda$ is bounded above by $P^\phi(G_0)\mathcal L^1\mres(t_0,\infty)$. The disintegration theorem therefore yields a measurable family of
nonnegative Radon measures $\{\lambda_t\}_{t>t_0}$ such that
\begin{equation}\label{mob:9}
d\lambda=d\lambda_t(x)\,dt,
\qquad
\lambda_t(\R^N)\leq P^\phi(G_0)
\quad\text{for a.e. }t>t_0.
\end{equation}

We next compare $\rho$ with $\lambda$. Writing $\nu_j=(\nu_j',\nu_j\cdot e_{_{N+1}})$, we introduce the equatorial measures
\begin{equation}\label{mob:equatorialmeasures}
\begin{split}
\overline\lambda_j
&:=\e_j(\mathcal T_j)_\#
\left(\phi^{\kappa_j}(\nu_j')\mu_{G_j}\mres H^+\right),\\
\overline\rho_j
&:=\e_j(\mathcal T_j)_\#
\left(\mathfrak m^{\kappa_j}(\nu_j')
\mu_{G_j}\mres H^+\right),
\end{split}
\end{equation}
where we recall that  $\phi^\kappa(\xi)=\Phi^\kappa(\xi,0)$ and  $\mathfrak m^\kappa(\xi):=\mathfrak M^\kappa(\xi,0)$.

For $\nu_{j}=(\nu_{j}',\nu_{j}\cdot e_{_{N+1}})$, by convexity and \eqref{mobprop:extnorm} and \eqref{mob:solitonbis}
\[
0\leq \Phi^{\kappa_{j}}(\nu_j)-\phi^{\kappa_{j}}(\nu'_j)
\leq (\nu_j\cdot e_{_{N+1}})
     \nabla\Phi^{\kappa_{j}}(\nu_j)\cdot e_{_{N+1}}
=\frac{(\nu_j\cdot e_{_{N+1}})^2}{\mathfrak M^{\kappa_{j}}(\nu_j)}=\e_j^2\,|H_j|^2\,\mathfrak{M}^{\kappa_j}(\nu_j),
\]
whereas by the Lipschitz continuity of $\mathfrak{M}^{\kappa_{j}}$, the definition of $L_{\kappa_{j}}$ in \eqref{contr:ke}, and \eqref{mob:solitonbis},
\begin{equation}\label{Mm:Lh}
0
\leq
\mathfrak M^{\kappa_j}(\nu_j)
-\mathfrak m^{\kappa_j}(\nu_j')\leq L_{\kappa_j}|\nu_j\cdot e_{_{N+1}}|
=\e_jL_{\kappa_j}|H_j|
  \mathfrak M^{\kappa_j}(\nu_j).
\end{equation}
These inequalities yield
\begin{equation}\label{l:bl:eq}
\begin{split}
(\lambda_j-\overline\lambda_j)(\eta)
&\leq \e_j^3\int_{H^+}(\eta\circ\mathcal T_j)|H_j|^2
\,d\mu_{G_j}^{\mathfrak M^{\kappa_j}}
=\e_j^2\,\Lambda_j(\eta),
\end{split}
\end{equation}
and by using Cauchy-Schwarz inequality,
\[
\begin{split}
(\rho_j-\overline\rho_j)(\eta)
&\leq\e_j^2L_{\kappa_j}
\int_{H^+}(\eta\circ\mathcal T_j)|H_j|
\,d\mu_{G_j}^{\mathfrak M^{\kappa_j}}\leq \e_jL_{\kappa_j}
\bigl(\rho_j(\eta)\Lambda_j(\eta)\bigr)^{1/2}.
\end{split}
\]
for every nonnegative $\eta\in C_c(\R^N\times(t_0,\infty))$.

By \eqref{mob:questaserve} and \eqref{contr:ke}, the right-hand sides of the two estimates above tend to zero. Consequently,
\begin{equation}\label{co:conv}
\overline\lambda_j\weakstar\lambda,
\qquad
\overline\rho_j\weakstar\rho.
\end{equation}

Next, we set
\begin{equation}\label{mob:beta}
\begin{split}
\alpha_\kappa
&:=\min_{\xi\in\mathbb S^{N-1}}
\frac{\mathfrak m^\kappa(\xi)}{\phi^\kappa(\xi)},
\qquad
\beta_\kappa
:=\max_{\xi\in\mathbb S^{N-1}}
\frac{\mathfrak m^\kappa(\xi)}{\phi^\kappa(\xi)},\\
\alpha
&:=\min_{\xi\in\mathbb S^{N-1}}
\frac{\mathfrak m(\xi)}{\phi(\xi)},
\qquad
\beta
:=\max_{\xi\in\mathbb S^{N-1}}
\frac{\mathfrak m(\xi)}{\phi(\xi)}.
\end{split}
\end{equation}
Since $\phi^\kappa\to\phi$ and $\mathfrak m^\kappa\to\mathfrak m$ locally uniformly, we have $\alpha_\kappa\to\alpha$ and $\beta_\kappa\to\beta$, and by homogeneity
$0<\alpha\leq\beta<\infty$. We then have the trivial inequalities
\[
\alpha_{\kappa_j}\overline\lambda_j
\leq\overline\rho_j
\leq\beta_{\kappa_j}\overline\lambda_j,
\]
and therefore, passing to the limit, we obtain
\[
\alpha\lambda\leq\rho\leq\beta\lambda.
\]
Applying  the Radon--Nikodym theorem, we obtain
\begin{equation}\label{mob:rho-density}
h:=\frac{d\rho}{d\lambda}\in L^\infty(\lambda),
\qquad
\alpha\leq h(x,t)\leq\beta
\quad\text{for }\lambda\text{-a.e. }(x,t),\quad d\rho=h\,d\lambda=h(x,t)\,d\lambda_t(x)\,dt.
\end{equation}

 Next we introduce the vector-valued transport measures
\begin{equation}\label{transp:meas}
J_j:=-\e_j(\mathcal T_j)_\#
\left(
H_j\nabla_\xi\Phi^{\kappa_j}(\nu_j)
\mu_{G_j}^{\mathfrak M^{\kappa_j}}\mres H^+
\right).
\end{equation}
By \eqref{mob:7} and Cauchy--Schwarz, for every $X\in C_c(\R^N\times(t_0,\infty);\R^N)$,
\begin{equation}\label{mob:CS}
\begin{split}
\left|\int X\cdot dJ_j\right|
&\leq
\left(\int\phi^{\kappa_j}(X)^2\,d\rho_j\right)^{1/2}
\Lambda_j(\operatorname{spt}X)^{1/2}\\
&\leq C_\phi P^{\phi^{\kappa_j}}(G_0)^{1/2}
\left(\int|X|^2\,d\rho_j\right)^{1/2}.
\end{split}
\end{equation}
From the local mass bounds for $\rho_j$ we therefore obtain uniform bounds for $|J_j|$ on compact sets. Up to a further subsequence, we obtain an $\R^N$-valued Radon measure $J$
such that
\begin{equation}\label{JjJ:t0}
J_j\weakstar J
\quad\text{locally in }\R^N\times(t_0,\infty).
\end{equation}
Passing to the limit in \eqref{mob:CS} and recalling $\rho\leq\beta\lambda$, we find
\[
\left|\int X\cdot dJ\right|
\leq C_\phi P^\phi(G_0)^{1/2}
\left(\int|X|^2\,d\rho\right)^{1/2}
\leq C_\phi\bigl(\beta P^\phi(G_0)\bigr)^{1/2}
\left(\int|X|^2\,d\lambda\right)^{1/2}.
\]
By density of compactly supported continuous vector fields in $L^2(\lambda;\R^N)$ and the Riesz representation theorem, there exists $\mathbf b\in L^2(\lambda;\R^N)$ such that
\begin{equation}\label{mob:limitingJ}
J=\mathbf b\,\lambda=\mathbf b\,\lambda_t\,dt.
\end{equation}

 By the duality inequality $|v\cdot w|\leq (\phi^{\kappa})^o(v)\,\phi^{\kappa}(w)$, \eqref{mob:7} and Young's inequality,
\[
-2H_j\nabla_\xi\Phi^{\kappa_j}(\nu_j)\cdot w
-\phi^{\kappa_j}(w)^2
\leq |H_j|^2
\qquad\text{for every }w\in\R^N.
\]
Hence, for every $X\in C_c(\R^N\times(t_0,\infty);\R^N)$ and every nonnegative $\eta\in C_c(\R^N\times(t_0,\infty))$,
\[
2\int\eta X\cdot dJ_j
-\int\eta\,\phi^{\kappa_j}(X)^2\,d\rho_j
\leq\int\eta\,d\Lambda_j.
\]
Passing to the limit and using \eqref{mob:rho-density} and \eqref{mob:limitingJ}, we obtain
\begin{equation}\label{mob:mii}
\int\eta\bigl(2\mathbf b\cdot X-h\phi(X)^2\bigr)\,d\lambda
\leq\int\eta\,d\Lambda.
\end{equation}
Since $h>0$ $\lambda$-a.e., by duality we have
\[
\sup_{w\in\R^N}
\left\{2\mathbf b(x,t)\cdot w-h(x,t)\phi(w)^2\right\}
=
\frac{\bigl(\phi^o(\mathbf b(x,t))\bigr)^2}{h(x,t)}
\quad\text{for }\lambda\text{-a.e. }(x,t).
\]
Taking the supremum in \eqref{mob:mii}, using density of continuous vector fields in $L^2(\eta\lambda;\R^N)$, we conclude that
\begin{equation}\label{mob:10}
\frac{\bigl(\phi^o(\mathbf b)\bigr)^2}{h}\,
\lambda_t\,dt\leq\Lambda.
\end{equation}

We now rewrite \eqref{mob:6} in terms of the rescaled and transport measures \eqref{mob:measures}, \eqref{transp:meas}:
\[
\int\partial_t\eta\,d\lambda_j
+\int\nabla_x\eta\cdot dJ_j
=
\int\eta\,d\Lambda_j
+\e_j^2\int\partial_t\eta\,d\Lambda_j.
\]
The last term tends to zero, since its absolute value is bounded by $\e_j^2P^{\phi^{\kappa_j}}(G_0) \|\partial_t\eta\|_\infty$. Passing to the limit, we obtain
\begin{equation}\label{mob:11}
\int_{t_0}^\infty\int
\bigl(\partial_t\eta+\mathbf b\cdot\nabla_x\eta\bigr)
\,d\lambda_t\,dt
=
\int\eta\,d\Lambda
\end{equation}
for every $\eta\in C_c^1(\R^N\times(t_0,\infty))$.

Likewise, by \eqref{mob:3}, lower semicontinuity of the total variation of the time derivative, \eqref{mob:5} and \eqref{mob:rho-density} we obtain
\begin{equation}\label{mob:12}
\begin{split}
|\partial_t\chi_{\widehat E}|(\eta)
\leq
\left(\int\eta\,d\rho\right)^{1/2}
\left(\int\eta\,d\Lambda\right)^{1/2}=
\left(\int_{t_0}^\infty\int\eta h\,d\lambda_t\,dt\right)^{1/2}
\left(\int\eta\,d\Lambda\right)^{1/2}
\end{split}
\end{equation}
for every nonnegative $\eta\in C_c(\R^N\times(t_0,\infty))$.

We now show that the limiting energy dominates the anisotropic perimeter of the slices. Set $F_j:=\mathcal R_{\e_j}(G_j)+(0,t_0)$. For every nonnegative $\eta\in
C_c(\R^N\times(t_0,\infty))$, arguing as in \eqref{as:in}, \eqref{13}, we find
\begin{equation}\label{mob:as:in}
\begin{split}
\int_{t_0}^\infty\int
\eta\,d\mu^\phi_{\widehat E(t)}\,dt
&=
\int_{\partial^*\widehat E}
\eta\,\phi(\nu'_{\widehat E})\,d\Ha^N\leq
\liminf_{j\to\infty}
\int_{\partial^*F_j}
\eta\,\phi^{\kappa_j}(\nu'_{F_j})\,d\Ha^N\\
&=
\lim_{j\to\infty}\int\eta\,d\overline\lambda_j
=
\int_{t_0}^\infty\int\eta\,d\lambda_t\,dt.
\end{split}
\end{equation}
Specifically, we used the slicing formula, Reshetnyak's lower semicontinuity theorem together with \eqref{mob:3} and $\phi^{\kappa_j}\to\phi$, the change of variables $\mathcal
T_j$, \eqref{co:conv} and \eqref{mob:9}. The arbitrariness of $\eta$ yields
\begin{equation}\label{mob:13}
\lambda_t\geq\mu^\phi_{\widehat E(t)}
\qquad\text{for a.e. }t>t_0.
\end{equation}
In a similar fashion, since $\mathfrak m^\kappa\to \mathfrak m$ locally uniformly,
\[
\begin{aligned}
\int_{t_0}^\infty\int
\eta\,d\mu^{\mathfrak m}_{\widehat E(t)}\,dt
&=
\int_{\partial^*\widehat E}
\eta\,\mathfrak m(\nu'_{\widehat E})\,d\Ha^N\leq
\liminf_{j\to\infty}
\int_{\partial^*F_j}
\eta\,\mathfrak m^{\kappa_j}(\nu'_{F_j})\,d\Ha^N\\
&=
\liminf_{j\to\infty}
\e_j\int_{\partial^*G_j\cap H^+}
\eta(x,t_0+\e_jz)\,
\mathfrak m^{\kappa_j}(\nu_j')\,d\Ha^N\\
&=
\lim_{j\to\infty}\int\eta\,d\overline\rho_j
=\int\eta\,d\rho=
\int_{t_0}^\infty\int
\eta(x,t)h(x,t)\,d\lambda_t(x)\,dt.
\end{aligned}
\]
which, by the arbitrariness of $\eta$, shows that
\begin{equation}\label{la:usero}
h(\cdot,t)\lambda_t
\geq\mu^{\mathfrak m}_{\widehat E(t)}
\qquad\text{for a.e. }t>t_0.
\end{equation}

\medskip
\noindent\textit{Step 5: Choice of the precise representative, and monotonicity of the total mass.}

The argument of Step 5 of the proof of Theorem~\ref{trasportobv} applies with only a change in the constant in \eqref{stimetta}. Indeed, since $h\leq\beta$, by \eqref{mob:10} we
have $\bigl(\phi^o(\mathbf b)\bigr)^2\lambda\leq\beta\Lambda$ as measures.

Thus, by equivalence of norms, Cauchy--Schwarz, \eqref{mob:9}, and \eqref{mob:forseserve}, for every $T>t_0$ we obtain
\begin{equation}\label{mob:stimetta}
\begin{split}
\int_{t_0}^T\int|\mathbf b|\,d\lambda_t\,dt
&\leq C_\phi
\left(\int_{t_0}^T\lambda_t(\R^N)\,dt\right)^{1/2}
\left(\beta\Lambda(\R^N\times(t_0,T))\right)^{1/2}\\
&\leq C_\phi\sqrt{\beta}\,
P^\phi(G_0)\sqrt{T-t_0}.
\end{split}
\end{equation}
Using \eqref{mob:11} and \eqref{mob:13}, we may therefore repeat that argument to choose a representative $\{\lambda_t\}_{t\geq t_0}$ which is weakly-* left-continuous on
$(t_0,\infty)$ and satisfies
\[
\lambda_{t_0}=\mu^\phi_{G_0},
\qquad
\lambda_t\weakstar\mu^\phi_{G_0}
\quad\text{as }t\downarrow t_0.
\]
Moreover,
\begin{equation}\label{mob:14}
\lambda_t\geq\mu^\phi_{\widehat E(t)},\qquad\lambda_t(\R^N)
=P^\phi(G_0)-\Lambda(\R^N\times(t_0,t))
\qquad\text{for every }t\geq t_0.
\end{equation}
In particular, $t\mapsto\lambda_t(\R^N)$ is nonincreasing, left-continuous on $(t_0,\infty)$, and converges to $P^\phi(G_0)$ as $t\downarrow t_0$.

By the same argument we also obtain, for every $T>t_0$ and $\eta\in C_c^1(\R^N\times\R)$,
\begin{equation}\label{mob:thfund-general}
\begin{split}
\int_{t_0}^T\int
\bigl(\partial_t\eta+\mathbf b\cdot\nabla_x\eta\bigr)
\,d\lambda_t\,dt
={}&\int_{\R^N\times(t_0,T)}\eta\,d\Lambda+\lambda_T\bigl(\eta(\cdot,T)\bigr)
-\mu^\phi_{G_0}\bigl(\eta(\cdot,t_0)\bigr).
\end{split}
\end{equation}

Finally, applying the same reasoning with arbitrary starting time $t_0\geq0$, we conclude that $t\mapsto\mu^\phi_{\widehat E(t)}$ is weakly-* right-continuous on $[0,\infty)$, and
that $t\mapsto P^\phi(\widehat E(t))$ is nonincreasing and right-continuous.

\medskip
\noindent\textit{Step 6: Restarting and matching.}

We use the same restarting procedure as in Step 6 of the proof of Theorem~\ref{trasportobv}. Fix $\delta>0$, set $t_0^\delta=0$, and proceed inductively. Once $t_i^\delta$ has been
defined, apply Steps 1--5 with starting time $t_0=t_i^\delta$ and initial datum $G_0=\widehat E(t_i^\delta)$. Denote the resulting measures and fields by $\lambda_t^{(i)}$,
$\Lambda^{(i)}$, $\mathbf b^{(i)}$, and $h^{(i)}$. Define
\begin{equation}\label{mob:15}
t_{i+1}^\delta
:=
\inf\left\{
t>t_i^\delta:
\lambda_t^{(i)}(\R^N)
\geq P^\phi(\widehat E(t))+\delta
\right\},
\end{equation}
with value $+\infty$ if the set is empty. By \eqref{mob:14} and the right continuity of $t\mapsto P^\phi(\widehat E(t))$, the argument in Step 6 of the proof of
Theorem~\ref{trasportobv} applies unchanged: every finite restarting time is attained, $t_{i+1}^\delta>t_i^\delta$, and
\[
P^\phi(\widehat E(t_{i+1}^\delta))+\delta
\leq\lambda_{t_{i+1}^\delta}^{(i)}(\R^N)
\leq P^\phi(\widehat E(t_i^\delta)).
\]
Consequently, the number of restarts satisfies $N_\delta\leq\frac{P^\phi(\widehat E(0))}{\delta}$. Then we set $t_{N_\delta+1}^\delta=\infty$,
$I_i^\delta=[t_i^\delta,t_{i+1}^\delta)$, and
\[
\lambda_t^\delta:=\lambda_t^{(i)}\quad\text{for  $t\in I_i^\delta$, }
\]
and we also define
\[
\mathbf b^\delta:=\mathbf b^{(i)},
\qquad
h^\delta:=h^{(i)}
\]
$\lambda_t^\delta\,dt$-almost everywhere on the corresponding space-time strip. We remark that the values of these fields at restarting times are immaterial, since
$\lambda_t^\delta\,dt$ does not charge time slices. Define also the mobility measure
\begin{equation}\label{mob:rho-restarted}
d\rho^\delta(x,t)
:=h^\delta(x,t)\,d\lambda_t^\delta(x)\,dt.
\end{equation}
By \eqref{mob:rho-density} and \eqref{la:usero},
\begin{equation}\label{mob:rho-restarted-bounds}
\begin{split}
&\alpha\leq h^\delta\leq\beta
\qquad\text{for }\lambda_t^\delta\,dt\text{-a.e. }(x,t),\\
&\alpha\lambda_t^\delta\,dt
\leq\rho^\delta\leq\beta\lambda_t^\delta\,dt,
\qquad
\mu^{\mathfrak m}_{\widehat E(t)}\,dt\leq\rho^\delta.
\end{split}
\end{equation}
As in the natural-mobility case, define the total dissipation, including the losses at restarting times, by
\begin{equation}\label{mob:16}
\int\eta\,d\Lambda^\delta
\coloneqq
\sum_{i=0}^{N_\delta}
\int_{\R^N\times(t_i^\delta,t_{i+1}^\delta)}
\eta\,d\Lambda^{(i)}+
\sum_{i=1}^{N_\delta}
\left(
\lambda_{t_i^\delta}^{(i-1)}
-\mu^\phi_{\widehat E(t_i^\delta)}
\right)
\bigl(\eta(\cdot,t_i^\delta)\bigr).
\end{equation}
Each measure in the second sum is nonnegative by \eqref{mob:14}. Summing \eqref{mob:thfund-general} over the intervals therefore gives
\begin{equation}\label{mob:17}
\int_0^\infty\int
\bigl(\partial_t\eta+\mathbf b^\delta\cdot\nabla_x\eta\bigr)
\,d\lambda_t^\delta\,dt
=
\int\eta\,d\Lambda^\delta
\end{equation}
for every $\eta\in C_c^1(\R^N\times(0,\infty))$. At this point, we apply the same telescoping argument as  in \eqref{18}, and get
\begin{equation}\label{mob:18}
\Lambda^\delta(\R^N\times(0,\infty))
\leq P^\phi(\widehat E(0)),\quad\text{and}\quad\frac{\bigl(\phi^o(\mathbf b^\delta)\bigr)^2}{h^\delta}\,
\lambda_t^\delta\,dt
\leq\Lambda^\delta,
\end{equation}
the second inequality stemming from \eqref{mob:10}.

Next, by applying \eqref{mob:12} on each open interval
\begin{equation}\label{mob:19}
\begin{split}
|\partial_t\chi_{\widehat E}|(\eta)
\leq
\left(\int\eta\,d\rho^\delta\right)^{1/2}
\left(\int\eta\,d\Lambda^\delta\right)^{1/2}=
\left(
\int_0^\infty\int\eta h^\delta\,d\lambda_t^\delta\,dt
\right)^{1/2}
\left(\int\eta\,d\Lambda^\delta\right)^{1/2}
\end{split}
\end{equation}
for every nonnegative $\eta\in C_c(\R^N\times(0,\infty))$.

 Finally, the matching estimate is unchanged: indeed,  by construction and \eqref{mob:14},
\[
\mu^\phi_{\widehat E(t)}\leq\lambda_t^\delta,
\qquad
\lambda_t^\delta(\R^N)
\leq P^\phi(\widehat E(t))+\delta
\qquad\text{for every }t\geq0.
\]
Recalling that $d\sigma=d\mu^\phi_{\widehat E(t)}\,dt$, we conclude that
\begin{equation}\label{mob:20}
\left|\lambda_t^\delta\,dt-\sigma\right|
\bigl(\R^N\times(0,T)\bigr)
\leq T\delta
\qquad\text{for every }T<\infty.
\end{equation}
Moreover, \eqref{mob:rho-restarted-bounds} and \eqref{mob:14} give the uniform bound
\[
\rho^\delta(\R^N\times(0,T))
\leq\beta T P^\phi(\widehat E(0)),
\]
which provides local compactness of these mobility measures as $\delta\downarrow0$.

\medskip
\noindent\textit{Step 7: The matched limit.}

Let $\delta\downarrow0$. By \eqref{mob:18}, the bound $h^\delta\leq\beta$, equivalence of norms, and Cauchy--Schwarz, for every $X\in C_c(\R^N\times(0,\infty);\R^N)$ we have
\[
\left|
\int_0^\infty\int
X\cdot\mathbf b^\delta\,d\lambda_t^\delta\,dt
\right|
\leq
C_\phi\bigl(\beta P^\phi(\widehat E(0))\bigr)^{1/2}
\left(
\int_0^\infty\int|X|^2\,d\lambda_t^\delta\,dt
\right)^{1/2}.
\]
Combining this estimate with \eqref{mob:20}, we obtain locally uniform bounds on the total variations of the vector measures $\mathbf b^\delta\lambda_t^\delta\,dt$. The bounds in
Step 6 also provide local compactness of $\rho^\delta$ and $\Lambda^\delta$. Thus, after extraction, all three sequences admit weak-star limits.

By \eqref{mob:20}, the integral on the right-hand side converges to $\int|X|^2\,d\sigma$. The same density and Riesz representation argument as in Step 7 of the proof of
Theorem~\ref{trasportobv} therefore identifies the limiting vector measure with $\mathbf b\sigma$, for some $\mathbf b\in L^2(\sigma;\R^N)$. Hence, along the chosen subsequence,
\[
\mathbf b^\delta\lambda_t^\delta\,dt
\weakstar\mathbf b\sigma,
\qquad
\rho^\delta\weakstar\rho,
\qquad
\Lambda^\delta\weakstar\Lambda
\quad\text{locally in }\R^N\times(0,\infty).
\]
The limit measure $\Lambda$ is nonnegative, and by lower semicontinuity we have
\[
\Lambda(\R^N\times(0,\infty))
\leq P^\phi(\widehat E(0)).
\]

Passing to the limit in \eqref{mob:rho-restarted-bounds}, using \eqref{mob:20}, we obtain
\[
\alpha\sigma\leq\rho\leq\beta\sigma,
\qquad
\mu^{\mathfrak m}_{\widehat E(t)}\,dt\leq\rho.
\]
Thus, by the Radon--Nikodym theorem we obtain
\begin{equation}\label{mob:matched:rho}
h:=\frac{d\rho}{d\sigma}\in L^\infty(\sigma),
\qquad
d\rho=h\,d\sigma,
\qquad
\alpha\leq h\leq\beta
\quad\text{for }\sigma\text{-a.e. }(x,t).
\end{equation}
Moreover, since $d\mu^{\mathfrak m}_{\widehat E(t)}\,dt =\mathfrak m(\nu_{\widehat E(t)})\,d\mu$ and $ d\sigma=\phi(\nu_{\widehat E(t)})\,d\mu$, the second measure inequality above
yields the sharper lower bound
\begin{equation}\label{mob:matched:h}
h(x,t)\geq
\frac{\mathfrak m(\nu_{\widehat E(t)}(x))}
     {\phi(\nu_{\widehat E(t)}(x))}
\qquad\text{for }\sigma\text{-a.e. }(x,t).
\end{equation}
which proves \eqref{mob:h-bounds}. Now, in order to inherit the localized quadratic estimate, we repeat the Young inequality argument used in Step 4: for every nonnegative $\eta\in
C_c(\R^N\times(0,\infty))$ and every $X\in C_c(\R^N\times(0,\infty);\R^N)$, by \eqref{mob:18} we have
\[
2\int_0^\infty\int
\eta\,\mathbf b^\delta\cdot X\,d\lambda_t^\delta\,dt
-\int\eta\,\phi(X)^2\,d\rho^\delta
\leq\int\eta\,d\Lambda^\delta.
\]
Passing to the limit and using \eqref{mob:matched:rho}, we obtain
\[
\int\eta\bigl(2\mathbf b\cdot X-h\phi(X)^2\bigr)
\,d\sigma
\leq\int\eta\,d\Lambda.
\]
Taking the supremum over $X$, exactly as in Step 4, yields
\begin{equation}\label{mob:matched:dissipation}
\frac{\bigl(\phi^o(\mathbf b)\bigr)^2}{h}\,
\sigma\leq\Lambda.
\end{equation}

We can also pass to the limit in \eqref{mob:17}. For every $\eta\in C_c^1(\R^N\times(0,\infty))$, we obtain
\begin{equation}\label{mob:restart:l1}
\int
\bigl(\partial_t\eta+\mathbf b\cdot\nabla_x\eta\bigr)
\,d\sigma
=
\int\eta\,d\Lambda,
\end{equation}
which proves \eqref{mob:bilancioenergia}.

Finally, passing to the limit in \eqref{mob:19}, we obtain
\[
|\partial_t\chi_{\widehat E}|(\eta)
\leq
\left(\int\eta\,d\rho\right)^{1/2}
\left(\int\eta\,d\Lambda\right)^{1/2}
=
\left(\int\eta h\,d\sigma\right)^{1/2}
\left(\int\eta\,d\Lambda\right)^{1/2}
\]
for every nonnegative $\eta\in C_c(\R^N\times(0,\infty))$. By approximation, it follows that
\[
|\partial_t\chi_{\widehat E}|(A)
\leq\bigl(\rho(A)\Lambda(A)\bigr)^{1/2}
\]
for every bounded Borel set $A\subset\R^N\times(0,\infty)$. Since $\rho\leq\beta\sigma$, this implies $\partial_t\chi_{\widehat E}\ll\sigma$. Recalling that
$\sigma=\phi(\nu_{\widehat E(t)})\mu$ and that the norm $\phi$ is bounded above and below on the unit sphere, we obtain a function $v\in L^1_{\loc}(\mu)$ such that
\[
\partial_t\chi_{\widehat E}=v\mu.
\]
The slicing and Gauss--Green argument from Step 7 of the proof of Theorem~\ref{trasportobv} applies unchanged and we obtain
\[
D_x\chi_{\widehat E}
=-\nu_{\widehat E(t)}\mu.
\]
Consequently,
\[
D_{x,t}\chi_{\widehat E}
=(-\nu_{\widehat E(t)},v)\mu,
\qquad
|D_{x,t}\chi_{\widehat E}|
=\sqrt{1+v^2}\,\mu.
\]

\medskip
\noindent\textit{Step 8: Conclusion.}

The argument of Step 8 of the proof of Theorem~\ref{trasportobv} applies unchanged, using \eqref{mob:restart:l1} and the identities established at the end of Step 7. In particular,
it yields
\[
\mathbf b\cdot\nu_{\widehat E(t)}=v
\qquad\text{for }\mu\text{-a.e. }(x,t).
\]
Since $\mathbf b\in L^2(\sigma;\R^N)$ and $d\sigma=\phi(\nu_{\widehat E(t)})\,d\mu$, we also obtain $v\in L^2(\mu)$.

Finally, using the temporal continuity properties of the perimeter measures established in Step 5, together with \eqref{eq:flatflowholder}, we obtain by the same argument
\eqref{mob:bilanciointervalli} and \eqref{mob:trasportofase}, completing the proof.
\end{proof}

\begin{theorem}[Brakke and canonical BV realization]
\label{mob:brakkebv}
Under the notation and assumptions of Theorem~\ref{mob:trasportobv}, suppose additionally that $\phi\in C^1(\R^N\setminus\{0\})$ is strictly convex. Consider the unit-density
varifolds
\[
V_t=\var(\partial^*\widehat E(t)),
\qquad t\geq0.
\]
Then
\begin{equation}\label{mob:identifica:h}
h(x,t)=
\frac{\mathfrak m(\nu_{\widehat E(t)}(x))}
     {\phi(\nu_{\widehat E(t)}(x))}
\qquad\text{for }\sigma\text{-a.e. }(x,t).
\end{equation}
Moreover, for a.e. $t>0$, $V_t$ has generalized anisotropic mean curvature $\mathbf H^\phi_{\partial^*\widehat E(t)}$, and
\begin{equation}\label{mob:identifica}
\begin{split}
&\mathbf H^\phi_{\partial^*\widehat E(t)}
\perp\partial^*\widehat E(t),\quad\mathbf b
=h\bigl(\nu_{\widehat E(t)}\cdot
\mathbf H^\phi_{\partial^*\widehat E(t)}\bigr)
\nabla\phi(\nu_{\widehat E(t)}),\\
&v
=\mathfrak m(\nu_{\widehat E(t)})
\bigl(\nu_{\widehat E(t)}\cdot
\mathbf H^\phi_{\partial^*\widehat E(t)}\bigr)
\qquad
\mathcal H^{N-1}\text{-a.e. on }\partial^*\widehat E(t).
\end{split}
\end{equation}
In particular,
\begin{equation}\label{mob:curvatura}
h\bigl|\mathbf H^\phi_{\partial^*\widehat E(t)}\bigr|^2
\,\sigma
=
\mathfrak m(\nu_{\widehat E(t)})
\bigl|\mathbf H^\phi_{\partial^*\widehat E(t)}\bigr|^2
\,\mu
\leq\Lambda.
\end{equation}

As a consequence, the family $\{V_t\}_{t\geq0}$ is an integral (unit-density) $(\phi,\mathfrak m)$-anisotropic Brakke flow in the sense of Definition~\ref{def:brakke}; that is, for
every $0\leq t_1<t_2$ and every nonnegative $C^1$ test function $\eta$ with compact spatial support uniformly on $[t_1,t_2]$, one has
\begin{equation}\label{mob:brakke}
\begin{aligned}
&\mu^\phi_{\widehat E(t_2)}(\eta(\cdot,t_2))
-\mu^\phi_{\widehat E(t_1)}(\eta(\cdot,t_1))\\
&\quad\leq
\int_{t_1}^{t_2}\int
\Big[
\partial_t\eta
+h\bigl(\nu_{\widehat E(t)}\cdot
\mathbf H^\phi_{\partial^*\widehat E(t)}\bigr)
\nabla\phi(\nu_{\widehat E(t)})\cdot\nabla_x\eta
-\eta h\bigl|\mathbf H^\phi_{\partial^*\widehat E(t)}\bigr|^2
\Big]\,d\mu^\phi_{\widehat E(t)}\,dt.
\end{aligned}
\end{equation}
In particular, the optimal energy/dissipation inequality
\begin{equation}\label{mob:daje}
P^\phi(\widehat E(t_2))
+\int_{\R^N\times(t_1,t_2)}
\bigl|\mathbf H^\phi_{\partial^*\widehat E(t)}\bigr|^2
\,d\mu^{\mathfrak m}_{\widehat E(t)}dt
\leq P^\phi(\widehat E(t_1))
\end{equation}
holds for every $0\leq t_1<t_2$.

Finally, the flow is a canonical BV flow; that is,
\begin{equation}\label{mob:fase}
\begin{split}
&\int_{\widehat E(t_2)}\eta(x,t_2)\,dx
-\int_{\widehat E(t_1)}\eta(x,t_1)\,dx\\
&\quad=
\int_{t_1}^{t_2}\int_{\widehat E(t)}
\partial_t\eta\,dx\,dt
+\int_{t_1}^{t_2}\int
\eta\,\bigl(\nu_{\widehat E(t)}\cdot
\mathbf H^\phi_{\partial^*\widehat E(t)}\bigr)
\,d\mu^{\mathfrak m}_{\widehat E(t)}\,dt
\end{split}
\end{equation}
for every $0\leq t_1<t_2$ and every $C^1$ test function $\eta$ with compact spatial support uniformly on $[t_1,t_2]$.
\end{theorem}

\begin{proof}
As in the proof of Theorem \ref{brakkebv}, we start by observing that it is enough to identify the fields $h$ and $\mathbf b$ obtained in Theorem~\ref{mob:trasportobv}. More
precisely, we will show that \eqref{mob:identifica:h} holds and that, for almost every $t>0$, the varifold $\var(\partial^*\widehat E(t))$ admits generalized anisotropic mean
curvature, that this curvature is normal, and that the formula for $\mathbf b$ in \eqref{mob:identifica} holds.

Indeed, combining the identity $v=\mathbf b\cdot\nu_{\widehat E(t)}$ in \eqref{mob:tubofase} with Euler's identity for $\phi$ and \eqref{mob:identifica:h}, we obtain the formula
for $v$ in \eqref{mob:identifica}. Moreover, the identification of $\mathbf b$ and the Eikonal equation $\phi^o(\nabla\phi(\nu))=1$ yield
\[
\frac{\bigl(\phi^o(\mathbf b)\bigr)^2}{h}
=
h\bigl|\nu_{\widehat E(t)}\cdot
\mathbf H^\phi_{\partial^*\widehat E(t)}\bigr|^2
=
h\bigl|\mathbf H^\phi_{\partial^*\widehat E(t)}\bigr|^2,
\]
and thus \eqref{mob:dissipazione} implies \eqref{mob:curvatura}. In turn, using the latter and substituting the formula for $\mathbf b$ into \eqref{mob:bilanciointervalli}, we
obtain \eqref{mob:brakke}. The energy/dissipation inequality \eqref{mob:daje} then follows by using spatial cutoff functions converging to one. Finally, inserting the formula for
$v$ into \eqref{mob:trasportofase}, and using
\[
d\mu^{\mathfrak m}_{\widehat E(t)}
=\mathfrak m(\nu_{\widehat E(t)})\,d\mu_{\widehat E(t)}
=h\,d\mu^\phi_{\widehat E(t)},
\]
we obtain \eqref{mob:fase}.

To establish these identifications, we follow the proof of Theorem~\ref{brakkebv}, using the minimizers and approximating sequences selected in the proof of
Theorem~\ref{mob:trasportobv}. For each starting datum, we retain the notation $G_j:=G_{\e_j}^{\kappa_j}$, with $\kappa_j=\kappa(\e_j)$. All limits are taken along the
corresponding sequence $(\e_j,\kappa_j)$, or along further subsequences, which are not relabeled.

\medskip
\noindent\textit{Step 1: Varifolds associated with the slices of the rescaled minimizers.}

Fix $t_0\geq0$, set $G_0=\widehat E(t_0)$, and consider the sequence $G_j=G_{\e_j}^{\kappa_j}$ used in the corresponding construction of Theorem~\ref{mob:trasportobv}. For almost
every $t>t_0$, let
\[
V_t^j
:=\var\!\left(
\partial^*\{x\in\R^N:(x,(t-t_0)/\e_j)\in G_j\}
\right).
\]
By the same slicing and coarea argument as in Step~1 of the proof of Theorem~\ref{brakkebv}, and recalling \eqref{mob:equatorialmeasures}, we have
\begin{equation}\label{mob:sliceweights}
\|V_t^j\|^{\phi^{\kappa_j}}\,dt=\overline\lambda_j,
\qquad
\|V_t^j\|^{\mathfrak m^{\kappa_j}}\,dt
=\overline\rho_j.
\end{equation}

Since $\phi^{\kappa_j}\to\phi$ uniformly on $\mathbb S^{N-1}$, there exists $c_\phi>0$ such that $\phi^{\kappa_j}(\xi)\geq c_\phi|\xi|$ for every $\xi\in\R^N$ and all sufficiently
large $j$. Thus, using $\overline\lambda_j\leq\lambda_j$ and \eqref{mob:questaserve}, we obtain
\begin{equation}\label{step:1:massbound}
\int_I\|V_t^j\|(\R^N)\,dt
\leq C_\phi\lambda_j(\R^N\times I)
\leq C_\phi(|I|+\e_j^2)P^{\phi^{\kappa_j}}(G_0)
\end{equation}
for every bounded open interval $I\Subset(t_0,\infty)$. Arguing by compactness and disintegration as in Step 1 of the proof of Theorem~\ref{brakkebv}, we therefore obtain, along a
further subsequence,
\[
V_t^j\,dt\weakstar V_t\,dt
\quad\text{on }
\R^N\times(t_0,\infty)\times G(N,N-1).
\]
%Here $V_t$ denotes the limiting family, which will be identified below with $\var(\partial^*\widehat E(t))$.
Since $\phi^{\kappa_j}\to\phi$ and $\mathfrak m^{\kappa_j}\to\mathfrak m$ uniformly on $\mathbb S^{N-1}$, we may pass to the limit in \eqref{mob:sliceweights}. Recalling
\eqref{co:conv} and \eqref{mob:rho-density}, we conclude that
\begin{equation}\label{mob:pesoblocco}
\|V_t\|^\phi=\lambda_t,
\qquad
\|V_t\|^{\mathfrak m}=h(\cdot,t)\lambda_t
\qquad\text{for a.e. }t>t_0.
\end{equation}

\medskip
\noindent\textit{Step 2: The approximate slice curvature vector measure.}

Recall the notation
\[
\nu_j=\nu_{G_j},\qquad
H_j=H^{\Phi^{\kappa_j}}_{\partial^*G_j},
\qquad
\mathcal T_j(x,z)=(x,t_0+\e_jz).
\]
We introduce the $\R^N$-valued Radon measure $\mathcal J_j$ on $\R^N\times(t_0,\infty)\times G(N,N-1)$ defined by
\begin{equation}\label{mob:mathcalJj}
\int Y\cdot d\mathcal J_j
:=
-\e_j\int_{H^+}
H_j\frac{\nu_j'}{|\nu_j'|}
\cdot Y\!\left(\mathcal T_j(x,z),
                  \frac{\nu_j'}{|\nu_j'|}\right)
\mathfrak m^{\kappa_j}(\nu_j')\,d\mu_{G_j}
\end{equation}
for continuous compactly supported vector fields $Y$. As in the natural-mobility case, the integrand is set to zero when $\nu_j'=0$, and unit normals in the arguments of $Y$ are
identified with their projective classes.

Since $\mathfrak m^{\kappa_j}(\nu_j') \leq\mathfrak M^{\kappa_j}(\nu_j)$, by Cauchy--Schwarz  and \eqref{mob:1} we obtain
\begin{equation}
\begin{split}
\left|\int Y\cdot d\mathcal J_j\right|
&\leq
\e_j\int_{H^+}|H_j|
\left|Y\!\left(\mathcal T_j(x,z),
               \frac{\nu_j'}{|\nu_j'|}\right)\right|
\mathfrak m^{\kappa_j}(\nu_j')\,d\mu_{G_j}\\
&\leq
\left(
\e_j\int_{H^+}|H_j|^2
\mathfrak M^{\kappa_j}(\nu_j)\,d\mu_{G_j}
\right)^{1/2}
\left(
\e_j\int_{H^+}
\left|Y\!\left(\mathcal T_j(x,z),
               \frac{\nu_j'}{|\nu_j'|}\right)\right|^2
\mathfrak m^{\kappa_j}(\nu_j')\,d\mu_{G_j}
\right)^{1/2}\\
&=
\left(
\e_j\int_{H^+}|H_j|^2
\mathfrak M^{\kappa_j}(\nu_j)\,d\mu_{G_j}
\right)^{1/2}
\left(
\int_{t_0}^\infty\int
|Y|^2\mathfrak m^{\kappa_j}(\nu)\,dV_t^j\,dt
\right)^{1/2}\\
&\leq
P^{\phi^{\kappa_j}}(G_0)^{1/2}
\left(
\int_{t_0}^\infty\int
|Y|^2\mathfrak m^{\kappa_j}(\nu)\,dV_t^j\,dt
\right)^{1/2}.
\end{split}
\end{equation}
By the convergence established in Step 1 and the uniform convergence of $\mathfrak m^{\kappa_j}$ on $\mathbb S^{N-1}$, the compactness and Riesz representation argument of Step 2
of the proof of Theorem~\ref{brakkebv} applies unchanged. Thus, along a further subsequence,
\begin{equation}\label{mob:limitecurva}
\mathcal J_j\weakstar
\mathcal J=\mathbf h\,\mathfrak m(\nu)V_t\,dt,
\qquad
\mathbf h\in
L^2\bigl(\mathfrak m(\nu)V_t\,dt;\R^N\bigr).
\end{equation}
The same argument leading to \eqref{disscurva}, with $\mathfrak m(\nu)$ in place of $\phi(\nu)$, yields
\begin{equation}\label{mob:disscurva}
\int_{t_0}^\infty\int
\eta(x,t)|\mathbf h(x,t,\nu)|^2\,
\mathfrak m(\nu)\,dV_t(x,\nu)\,dt
\leq\int\eta\,d\Lambda
\end{equation}
for every nonnegative $\eta\in C_c(\R^N\times(t_0,\infty))$. Likewise, testing \eqref{mob:mathcalJj} with $(\mathrm{Id}_N-\nu\otimes\nu)Y$ and passing to the limit, we obtain
\begin{equation}\label{mob:normalita}
(\mathrm{Id}_N-\nu\otimes\nu)\mathbf h(x,t,\nu)=0
\qquad V_t\,dt\text{-a.e.}
\end{equation}

It remains to relate $\mathbf h$ to the transport field $\mathbf b$ constructed in Theorem~\ref{mob:trasportobv}. To this aim, define the equatorial flux $\overline J_j$ by
\begin{equation}\label{equatorial:transp}
\begin{split}
\int X\cdot d\overline J_j
&:=
\int
\bigl(\nabla\phi(\nu)\cdot X(x,t)\bigr)
\,\nu\cdot d\mathcal J_j(x,t,\nu)\\
&=
-\e_j\int_{H^+}
H_j\mathfrak m^{\kappa_j}(\nu_j')
\nabla\phi(\nu_j')\cdot X(\mathcal T_j(x,z))
\,d\mu_{G_j},
\end{split}
\end{equation}
where $X\in C_c(\R^N\times(t_0,\infty);\R^N)$. Here and below, occurrences of $\nabla\phi(\nu_j')$ are interpreted as zero when $\nu_j'=0$.

We claim that
\begin{equation}\label{mob:fluxerror}
J_j-\overline J_j\longrightarrow0
\qquad\text{locally in total variation.}
\end{equation}
where $J_j$ is the flux defined in \eqref{transp:meas}. Observe that, by \eqref{JjJ:t0}, \eqref{mob:limitingJ} this implies
\begin{equation}\label{bar:JjJ:t0}
    \overline J_j\weakstar J=\mathbf{b}\,\lambda_t\,dt
\quad\text{locally in }\R^N\times(t_0,\infty).
\end{equation}

To prove \eqref{mob:fluxerror} we split the difference of the flux coefficients as
\[
\begin{split}
\mathfrak M^{\kappa_j}(\nu_j)
 \nabla_\xi\Phi^{\kappa_j}(\nu_j)
&-\mathfrak m^{\kappa_j}(\nu_j')\nabla\phi(\nu_j')
\\
&\,=
\mathfrak M^{\kappa_j}(\nu_j)
\bigl(\nabla_\xi\Phi^{\kappa_j}(\nu_j)
      -\nabla\phi(\nu_j')\bigr)+
\bigl(\mathfrak M^{\kappa_j}(\nu_j)
      -\mathfrak m^{\kappa_j}(\nu_j')\bigr)
\nabla\phi(\nu_j').
\end{split}
\]
Applying Cauchy--Schwarz to the first term and using \eqref{Mm:Lh} and the boundedness of $\nabla\phi$ for the second, we obtain, for every compact
$K\subset\R^N\times(t_0,\infty)$,
\[
\begin{split}
|J_j-\overline J_j|(K)
&\leq \Lambda_j(K)^{1/2}
\left(
\e_j\int_{\mathcal T_j^{-1}(K)}
\left|
\nabla_\xi\Phi^{\kappa_j}(\nu_j)
-\nabla\phi(\nu_j')
\right|^2
\,d\mu_{G_j}^{\mathfrak M^{\kappa_j}}
\right)^{1/2}
+C_\phi\e_jL_{\kappa_j}\Lambda_j(K)
\\
&\leq \Lambda_j(K)^{1/2}
\left(
\e_j\int_{\mathcal T_j^{-1}(K)}
\left|
\nabla_\xi\Phi^{\kappa_j}(\nu_j)
-\nabla\phi(\nu_j')
\right|^2
\,d\mu_{G_j}^{\mathfrak M^{\kappa_j}}
\right)^{1/2}
+C_\phi\e_jL_{\kappa_j}P^{\phi^{\kappa_j}}(G_0),
\end{split}
\]
where the last inequality follows from \eqref{mob:1}. The last term tends to zero by \eqref{contr:ke} and the uniform boundedness of $P^{\phi^{\kappa_j}}(G_0)\to P^\phi(G_0)$.

Next choose a bounded open interval $I\Subset(t_0,\infty)$ such that $K\subset\R^N\times I$. Splitting the integral into the regions $|\nu_j'|\geq1/2$ and $|\nu_j'|<1/2$, and using
\eqref{mob:7}, we obtain
\begin{equation}\label{prove:step2}
\begin{split}
&\e_j\int_{\mathcal T_j^{-1}(K)}
\left|\nabla_\xi\Phi^{\kappa_j}(\nu_j)
-\nabla\phi(\nu_j')\right|^2
\,d\mu_{G_j}^{\mathfrak M^{\kappa_j}}
\\
&\quad\leq
\e_j P^{\mathfrak M^{\kappa_j}}
\left(G_j;\R^N\times\frac{I-t_0}{\e_j}\right)
\sup_{\substack{\nu\in\Sf^N\\|\nu'|\geq1/2}}
\left|\nabla_\xi\Phi^{\kappa_j}(\nu)
-\nabla\phi(\nu')\right|^2
+C_\phi\e_j
\int_{H^+\cap\{|\nu_j'|<1/2\}}
\,d\mu_{G_j}^{\mathfrak M^{\kappa_j}}
\\
&\quad\leq
C_{\mathfrak m,\phi}
\bigl((\e_jL_{\kappa_j}^{2})^2+|I|+\e_j^2\bigr)
P^{\phi^{\kappa_j}}(G_0)
\sup_{\substack{\nu\in\Sf^N\\|\nu'|\geq1/2}}
\left|\nabla_\xi\Phi^{\kappa_j}(\nu)
-\nabla\phi(\nu')\right|^2
\\
&\qquad\quad+
C_\phi\e_jL_{\kappa_j}^2
\int_{H^+}
\frac{(\nu_j\cdot e_{_{N+1}})^2}
     {\mathfrak M^{\kappa_j}(\nu_j)}
\,d\mu_{G_j}
\\
&\quad\leq
C_{\mathfrak m,\phi}
\bigl((\e_jL_{\kappa_j}^{2})^2+|I|+\e_j^2\bigr)
P^{\phi^{\kappa_j}}(G_0)
\sup_{\substack{\nu\in\Sf^N\\|\nu'|\geq1/2}}
\left|\nabla_\xi\Phi^{\kappa_j}(\nu)
-\nabla\phi(\nu')\right|^2
\\
&\qquad\quad+
C_\phi(\e_jL_{\kappa_j})^2
P^{\phi^{\kappa_j}}(G_0)
\longrightarrow0.
\end{split}
\end{equation}
Here we used \eqref{mob:1}, the bound $\mathfrak M^{\kappa_j}(\nu_j)\leq L_{\kappa_j}$, and the fact that $|\nu_j'|<1/2$ implies $(\nu_j\cdot e_{_{N+1}})^2>3/4$, followed by
\eqref{mob:solitonbis} and \eqref{mob:1}. The convergence follows from the fact that
\begin{equation*}
    \nabla_{\xi}\Phi^\kappa(\nu)\to \nabla \phi(\nu')\quad\text{locally uniformly in $ \S^N\cap \{|\nu'|>0\}$,}
\end{equation*}
due to \eqref{za:tesi2}, by \eqref{contr:ke}, and the uniform boundedness of $P^{\phi^{\kappa_j}}(G_0)$. Together with the uniform bound on $\Lambda_j(K)$, this proves
\eqref{mob:fluxerror}.

Now, for every $X\in C_c(\R^N\times(t_0,\infty);\R^N)$, by \eqref{equatorial:transp} and \eqref{mob:limitecurva}, we may write
\[
\begin{split}
&\int X\cdot dJ_j
-\int_{t_0}^\infty\int
(\mathbf h\cdot\nu)\nabla\phi(\nu)\cdot X\,
\mathfrak m(\nu)\,dV_t\,dt
\\
&\quad=
\int X\cdot d(J_j-\overline J_j)
+\int
\bigl(\nabla\phi(\nu)\cdot X(x,t)\bigr)\nu
\cdot d(\mathcal J_j-\mathcal J).
\end{split}
\]
The first term tends to zero by \eqref{mob:fluxerror}, since
\[
\left|\int X\cdot d(J_j-\overline J_j)\right|
\leq
\|X\|_\infty
|J_j-\overline J_j|(\mathrm{spt}\,X)
\longrightarrow0.
\]
The second term tends to zero by \eqref{mob:limitecurva}, because the field $(\nabla\phi(\nu)\cdot X(x,t))\nu$ is continuous, compactly supported, and well defined on the
Grassmannian. Finally, by \eqref{JjJ:t0} and \eqref{mob:limitingJ},
\[
\int X\cdot dJ_j
\longrightarrow
\int_{t_0}^\infty\int
\mathbf b\cdot X\,d\lambda_t\,dt.
\]
Combining these convergences, we finally obtain
\begin{equation}\label{mob:flussocurva}
\int_{t_0}^\infty\int
\mathbf b\cdot X\,d\lambda_t\,dt
=
\int_{t_0}^\infty\int
(\mathbf h\cdot\nu)\nabla\phi(\nu)\cdot X\,
\mathfrak m(\nu)\,dV_t\,dt.
\end{equation}

\medskip
\noindent\textit{Step 3: Horizontal first variation.}

Let $X\in C_c^1(\R^N\times(t_0,\infty);\R^N)$. We apply \eqref{euler:lagrange0}, with $g=0$ and $\Psi=\Phi^{\kappa_j}$, to the horizontal field
\[
Y(x,z)=\bigl(X(\mathcal T_j(x,z)),0\bigr).
\]
Recall that, by \eqref{come:e:ovvia} and \eqref{mob:solitonbis}, we have $\partial_z\Phi^{\kappa_j}(\nu_j) =(\nu_j\cdot e_{_{N+1}})/ \mathfrak M^{\kappa_j}(\nu_j) =\e_jH_j$. Thus,
arguing as in Step 3 of the proof of Theorem~\ref{brakkebv}, we obtain
\begin{equation}\label{mob:varorizz}
\begin{split}
&\e_j\int_{H^+}
\bigl[
\Phi^{\kappa_j}(\nu_j)\mathrm{Id}_N
-\nu_j'\otimes\nabla_\xi\Phi^{\kappa_j}(\nu_j)
\bigr]:D_xX\,d\mu_{G_j}
\\
&\qquad
-\e_j^3\int_{H^+}
H_j\nu_j'\cdot\partial_tX\,d\mu_{G_j}
=
\e_j\int_{H^+}
H_j\nu_j'\cdot X\,d\mu_{G_j},
\end{split}
\end{equation}
where $X$ and its derivatives are evaluated at $\mathcal T_j(x,z)$.

We first compare the matrix in the first integral with the one defining the first variation of the slice varifolds. By one-homogeneity and \eqref{def:Bpsi},
\[
\begin{split}
&\Phi^{\kappa_j}(\nu_j)\mathrm{Id}_N
-\nu_j'\otimes\nabla_\xi\Phi^{\kappa_j}(\nu_j)
-|\nu_j'|B_\phi\!\left(\frac{\nu_j'}{|\nu_j'|}\right)
\\
&\quad=
\bigl(\Phi^{\kappa_j}(\nu_j)-\phi(\nu_j')\bigr)
\mathrm{Id}_N
-\nu_j'\otimes
\bigl(\nabla_\xi\Phi^{\kappa_j}(\nu_j)
      -\nabla\phi(\nu_j')\bigr),
\end{split}
\]
where the terms containing $\nu_j'/|\nu_j'|$ or $\nabla\phi(\nu_j')$ are interpreted as zero when $\nu_j'=0$.

Set $K=\mathrm{spt}\,X$. Using the estimate for $\lambda_j-\overline\lambda_j$ given by \eqref{l:bl:eq}, together with $\phi^{\kappa_j}\to \phi$ uniformly on $\mathbb S^{N-1}$, we
obtain
\[
\begin{split}
\e_j&\,\int_{\mathcal T_j^{-1}(K)}
\bigl|\Phi^{\kappa_j}(\nu_j)-\phi(\nu_j')\bigr|
\,d\mu_{G_j}
\\
&\leq \e_j\int_{\mathcal T_j^{-1}(K)}
\bigl|\Phi^{\kappa_j}(\nu_j)-\phi^{\kappa_j}(\nu_j')\bigr|
\,d\mu_{G_j}+\e_j\int_{\mathcal T_j^{-1}(K)}
\bigl|\phi^{\kappa_j}(\nu'_j)-\phi(\nu_j')\bigr|
\,d\mu_{G_j}
\\
&\quad\leq
\e_j^2\Lambda_j(K)
+C_\phi
\|\phi^{\kappa_j}-\phi\|_{L^\infty(\mathbb S^{N-1})}
\lambda_j(K)
\longrightarrow0.
\end{split}
\]
To estimate the remaining term, recall that $\mathfrak m^{\kappa_j}\to\mathfrak m$ uniformly on $\mathbb S^{N-1}$, so for all sufficiently large $j$, we have $|\nu_j'| \leq
C_{\mathfrak m}\mathfrak m^{\kappa_j}(\nu_j') \leq C_{\mathfrak m}\mathfrak M^{\kappa_j}(\nu_j)$. Consequently, by Cauchy--Schwarz,
\[
\begin{split}
&\e_j\int_{\mathcal T_j^{-1}(K)}
|\nu_j'|
\bigl|\nabla_\xi\Phi^{\kappa_j}(\nu_j)
      -\nabla\phi(\nu_j')\bigr|\,d\mu_{G_j}
\\
&\quad\leq
C_{\mathfrak m}\rho_j(K)^{1/2}
\left(
\e_j\int_{\mathcal T_j^{-1}(K)}
\bigl|\nabla_\xi\Phi^{\kappa_j}(\nu_j)
      -\nabla\phi(\nu_j')\bigr|^2
\,d\mu_{G_j}^{\mathfrak M^{\kappa_j}}
\right)^{1/2}
\longrightarrow0.
\end{split}
\]
Here the convergence of the integral in parentheses was proved in Step 2 (see \eqref{prove:step2}), while the uniform bounds for $\lambda_j(K)$, $\rho_j(K)$, and $\Lambda_j(K)$
follow from \eqref{mob:questaserve} and  \eqref{contr:ke}.

By the coarea formula, these estimates show that replacing the first integral in \eqref{mob:varorizz} by
\begin{equation}\label{mob:slice1stvar}
\int_{t_0}^\infty\int
B_\phi(\nu):D_xX\,dV_t^j\,dt
\end{equation}
introduces an error tending to zero.

The term containing $\partial_tX$ is treated as in the natural-mobility case. Indeed, using the bound for $|\nu_j'|$ above and Cauchy--Schwarz,
\[
\left|
\e_j^3\int_{H^+}
H_j\nu_j'\cdot\partial_tX\,d\mu_{G_j}
\right|
\leq
C_{\mathfrak m}\e_j^2\|\partial_tX\|_\infty
\Lambda_j(K)^{1/2}\rho_j(K)^{1/2}
\longrightarrow0.
\]

Finally, by \eqref{mob:mathcalJj} and one-homogeneity of $\mathfrak m^{\kappa_j}$, the right-hand side of \eqref{mob:varorizz} can be written as
\[
\e_j\int_{H^+}
H_j\nu_j'\cdot X(\mathcal T_j(x,z))\,d\mu_{G_j}
=
-\int\frac{X(x,t)}{\mathfrak m^{\kappa_j}(\nu)}
\cdot d\mathcal J_j.
\]
The functions $1/\mathfrak m^{\kappa_j}$ converge uniformly to $1/\mathfrak m$ on the Grassmannian, so by using \eqref{mob:limitecurva}, we obtain
\[
-\int\frac{X(x,t)}{\mathfrak m^{\kappa_j}(\nu)}
\cdot d\mathcal J_j
\longrightarrow
-\int\frac{X(x,t)}{\mathfrak m(\nu)}
\cdot d\mathcal J
=
-\int_{t_0}^\infty\int
\mathbf h\cdot X\,dV_t\,dt.
\]
Since $B_\phi$ is continuous on the Grassmannian, we may also pass to the limit in \eqref{mob:slice1stvar}, using the convergence established in Step 1. We thus conclude that
\begin{equation}\label{mob:limvar}
\int_{t_0}^\infty\int
B_\phi(\nu):D_xX\,dV_t\,dt
=
-\int_{t_0}^\infty\int
\mathbf h\cdot X\,dV_t\,dt
\end{equation}
for every $X\in C_c^1(\R^N\times(t_0,\infty);\R^N)$.

\medskip
\noindent\textit{Step 4: Restarting and matching.}

We now carry the varifolds and fields constructed above through the restarting procedure of Theorem~\ref{mob:trasportobv}.

Fix $\delta>0$, and recall the times $0=t_0^\delta<\cdots<t_{N_\delta}^\delta <t_{N_\delta+1}^\delta=\infty$ and the intervals $I_i^\delta=[t_i^\delta,t_{i+1}^\delta)$ constructed
in Step 6 of that proof. For each $i$, apply Steps 1--3 of the current proof with $t_0=t_i^\delta$ and $G_0=\widehat E(t_i^\delta)$, along the sequence already selected for that
starting datum. Passing to further subsequences preserves the previously obtained measures and fields, and provides the additional limits $V_t^{(i)}$ and $\mathbf h^{(i)}$.

As in Step 4 of the proof of Theorem~\ref{brakkebv}, we join these families by setting
\[
V_t^\delta:=V_t^{(i)},
\qquad
\mathbf h^\delta(x,t,\nu):=\mathbf h^{(i)}(x,t,\nu)
\qquad\text{for }t_i^\delta<t<t_{i+1}^\delta.
\]
Again, values at the finitely many endpoints are immaterial. By \eqref{mob:pesoblocco} and \eqref{mob:rho-restarted}, we have
\begin{equation}\label{mob:pesobloccodelta}
\|V_t^\delta\|^\phi=\lambda_t^\delta
\quad\text{for a.e. }t>0,
\qquad
\|V_t^\delta\|^{\mathfrak m}\,dt
=\rho^\delta=h^\delta\lambda_t^\delta\,dt.
\end{equation}

The time cutoff argument of Step 4 of the proof of Theorem~\ref{brakkebv} applies unchanged to \eqref{mob:limvar} and \eqref{mob:flussocurva}. Indeed, for a cutoff $\zeta$
depending only on time, $D_x(\zeta X)=\zeta D_xX$, and the mass and quadratic bounds from Steps 1--2, together with \eqref{mob:10}, allow us to remove the cutoffs by dominated
convergence. Summing over the intervals, we obtain
\begin{equation}\label{mob:1stvardelta}
\int_0^\infty\int
B_\phi(\nu):D_xX\,dV_t^\delta\,dt
=
-\int_0^\infty\int
\mathbf h^\delta\cdot X\,dV_t^\delta\,dt
\end{equation}
and
\begin{equation}\label{mob:id-b-delta}
\int_0^\infty\int
\mathbf b^\delta\cdot X\,d\lambda_t^\delta\,dt
=
\int_0^\infty\int
(\mathbf h^\delta\cdot\nu)
\nabla\phi(\nu)\cdot X\,
\mathfrak m(\nu)\,dV_t^\delta\,dt
\end{equation}
for every $X\in C_c^1(\R^N\times(0,\infty);\R^N)$.

Likewise, applying \eqref{mob:disscurva} on each open interval and recalling \eqref{mob:16}, we obtain
\begin{equation}\label{mob:raccordivar}
\begin{split}
\int_0^\infty\int
\eta|\mathbf h^\delta|^2\mathfrak m(\nu)
\,dV_t^\delta\,dt
&\leq
\sum_{i=0}^{N_\delta}
\int_{\R^N\times(t_i^\delta,t_{i+1}^\delta)}
\eta\,d\Lambda^{(i)}\leq\int\eta\,d\Lambda^\delta
\end{split}
\end{equation}
for every nonnegative $\eta\in C_c(\R^N\times(0,\infty))$. Here the last inequality follows from the nonnegativity of the additional measures at the restarting times.

We now let $\delta\downarrow0$ along the sequence selected in Step 7 of the proof of Theorem~\ref{mob:trasportobv}. Along this sequence,
\begin{equation}\label{mob:limitozzi}
\begin{split}
&\lambda_t^\delta\,dt\longrightarrow\sigma
\quad\text{in total variation on }
\R^N\times(0,T)
\quad\text{for every }T<\infty,
\\
&\mathbf b^\delta\lambda_t^\delta\,dt
\weakstar\mathbf b\sigma,
\qquad
\rho^\delta\weakstar\rho=h\sigma,
\qquad
\Lambda^\delta\weakstar\Lambda,
\end{split}
\end{equation}
where the first convergence follows from \eqref{mob:20}, and the identification of $\rho$ is given by \eqref{mob:matched:rho}. All further extractions below are taken from this
sequence.

By \eqref{mob:pesobloccodelta} and \eqref{mob:14}, for every bounded time interval $I\subset(0,\infty)$,
\[
\int_I\|V_t^\delta\|(\R^N)\,dt
\leq C_\phi\int_I\lambda_t^\delta(\R^N)\,dt
\leq C_\phi|I|P^\phi(\widehat E(0)).
\]
Moreover, by \eqref{mob:raccordivar} and \eqref{mob:18} we obtain a uniform quadratic bound for $\mathbf h^\delta$ with respect to $\mathfrak m(\nu)V_t^\delta\,dt$. Thus, the same
compactness, disintegration, and Riesz representation arguments as in Step 4 of the proof of Theorem~\ref{brakkebv} yield, after a further extraction,
\begin{equation}\label{mob:varifoldlimits}
\begin{split}
&V_t^\delta\,dt\weakstar V_t\,dt,
\\
&\mathbf h^\delta\mathfrak m(\nu)V_t^\delta\,dt
\weakstar
\mathbf h\mathfrak m(\nu)V_t\,dt,
\end{split}
\end{equation}
on $\R^N\times(0,\infty)\times G(N,N-1)$, with
\[
\mathbf h\in
L^2\bigl(\mathfrak m(\nu)V_t\,dt;\R^N\bigr).
\]

We next identify the two weighted mass measures. Since $\phi$ and $\mathfrak m$ are continuous on the Grassmannian, passing to the limit in \eqref{mob:pesobloccodelta}, using
\eqref{mob:limitozzi} and \eqref{mob:varifoldlimits}, we obtain
\[
\|V_t\|^\phi\,dt=\sigma,
\qquad
\|V_t\|^{\mathfrak m}\,dt=\rho=h\sigma.
\]
Recalling that $d\sigma=d\mu^\phi_{\widehat E(t)}\,dt$, we conclude that
\begin{equation}\label{mob:massalimite}
\|V_t\|^\phi=\mu^\phi_{\widehat E(t)},
\qquad
\|V_t\|^{\mathfrak m}
=h(\cdot,t)\mu^\phi_{\widehat E(t)}
\qquad\text{for a.e. }t>0.
\end{equation}

By the same argument leading to \eqref{mob:disscurva}, using \eqref{mob:raccordivar}, \eqref{mob:varifoldlimits}, and $\Lambda^\delta\weakstar\Lambda$, we obtain
\begin{equation}\label{mob:pesolimite}
\int_0^\infty\int
\eta(x,t)|\mathbf h(x,t,\nu)|^2
\mathfrak m(\nu)\,dV_t(x,\nu)\,dt
\leq\int\eta\,d\Lambda
\end{equation}
for every nonnegative $\eta\in C_c(\R^N\times(0,\infty))$. Normality of the $\phi$-mean curvature is also preserved. Indeed, by \eqref{mob:normalita} we have
\[
\int_0^\infty\int
\mathbf h^\delta\cdot
(\mathrm{Id}_N-\nu\otimes\nu)Y\,
\mathfrak m(\nu)\,dV_t^\delta\,dt=0
\]
for every continuous compactly supported vector field $Y$ on the product space. Passing to the limit by the second convergence in \eqref{mob:varifoldlimits}, and using the
arbitrariness of $Y$, we obtain
\begin{equation}\label{mob:matchingnormalita}
(\mathrm{Id}_N-\nu\otimes\nu)\mathbf h(x,t,\nu)=0
\qquad\text{for }V_t\,dt\text{-a.e. }(x,t,\nu).
\end{equation}

Finally, we pass to the limit in \eqref{mob:1stvardelta} and \eqref{mob:id-b-delta}. For the curvature term in the first identity, we test the vector measure convergence in
\eqref{mob:varifoldlimits} with $X(x,t)/\mathfrak m(\nu)$. For the second identity, we use the continuous field $(\nabla\phi(\nu)\cdot X(x,t))\nu$, which is well defined on the
Grassmannian, together with the convergence of $\mathbf b^\delta\lambda_t^\delta\,dt$ in \eqref{mob:limitozzi}. We therefore obtain
\begin{equation}\label{mob:primavar}
\int_0^\infty\int
B_\phi(\nu):D_xX\,dV_t\,dt
=
-\int_0^\infty\int
\mathbf h\cdot X\,dV_t\,dt
\end{equation}
and
\begin{equation}\label{mob:id-b}
\int\mathbf b\cdot X\,d\sigma
=
\int_0^\infty\int
(\mathbf h\cdot\nu)\nabla\phi(\nu)\cdot X\,
\mathfrak m(\nu)\,dV_t\,dt
\end{equation}
for every $X\in C_c^1(\R^N\times(0,\infty);\R^N)$.

\medskip
\noindent\textit{Step 5: Rectifiability and identification of the curvature.}

By \eqref{mob:massalimite}, \eqref{mob:pesolimite}, \eqref{mob:matchingnormalita}, and \eqref{mob:primavar}, the argument of Step 5 of the proof of Theorem~\ref{brakkebv} applies
unchanged. Indeed, since $\mathfrak m$ is bounded below on $\mathbb S^{N-1}$, by \eqref{mob:pesolimite} we obtain the same square-integrability of $\mathbf h$ required there. Thus,
there exists a negligible set $\mathcal N\subset(0,\infty)$ such that, for every $t\notin\mathcal N$, $V_t$ is rectifiable and
\begin{equation}\label{mob:densitauno}
V_t
=\mu_{\widehat E(t)}
\otimes
\delta_{\operatorname{Tan}(\partial^*\widehat E(t),x)}.
\end{equation}
Moreover,
\[
\mathbf H^\phi_{\partial^*\widehat E(t)}(x)
:=
\mathbf h\bigl(x,t,[\nu_{\widehat E(t)}(x)]\bigr)
\]
is its generalized anisotropic mean curvature, and by \eqref{mob:matchingnormalita} we have
\[
\mathbf H^\phi_{\partial^*\widehat E(t)}
\perp\partial^*\widehat E(t)
\qquad
\mu_{\widehat E(t)}\text{-a.e.}
\]

We now identify $h$. By \eqref{mob:massalimite} and \eqref{mob:densitauno},
\[
h\,\sigma
=\|V_t\|^{\mathfrak m}\,dt
=\mathfrak m(\nu_{\widehat E(t)})
\,\mu_{\widehat E(t)}\,dt
=
\frac{\mathfrak m(\nu_{\widehat E(t)})}
     {\phi(\nu_{\widehat E(t)})}\,\sigma.
\]
This proves \eqref{mob:identifica:h}. Finally, substituting \eqref{mob:densitauno} and this identification into \eqref{mob:id-b}, we obtain
\[
\int\mathbf b\cdot X\,d\sigma
=
\int
h\bigl(
\mathbf H^\phi_{\partial^*\widehat E(t)}
\cdot\nu_{\widehat E(t)}
\bigr)
\nabla\phi(\nu_{\widehat E(t)})\cdot X
\,d\sigma
\]
for every $X\in C_c^1(\R^N\times(0,\infty);\R^N)$. By the arbitrariness of $X$, this proves the formula for $\mathbf b$ in \eqref{mob:identifica} $\sigma$-a.e. As explained at the
beginning of the proof, and redefining $V_t=\var(\partial^*\widehat E(t))$ also for $t\in\mathcal N$ and $t=0$, this concludes the proof.

\end{proof}

 \appendix

\section{Uniformly convex norms and approximation}\label{sec:normeunifconvesse}

In this section we report a few properties of uniformly convex norms and an approximation procedure. We recall that a norm $\psi$ in $\R^d$ is \textit{uniformly convex} if the set
\[
B^\psi\coloneqq\{x\in\R^d:\,\psi(x)\leq1\}
\quad\text{is of class $C^2$ and uniformly convex,}
\]
that is, its principal curvatures are bounded away from zero.

Observe also that $B^\psi$ coincides with the closed Wulff shape of the dual norm, namely $B^\psi=\overline{W_1^{\psi^o}}$. We report the following characterization of uniform
convexity.

\begin{lemma}\label{lem:varie:unifconv}
Let $\psi$ be a norm on $\R^d$. The following conditions are equivalent:
\begin{itemize}
    \item[(a)] $\psi$ is uniformly convex;
    \item[(b)] $\psi\in C^2(\R^d\setminus\{0\})$, and there exist $0<\lambda\leq\Lambda$ such that, for all $x\neq0$,
    \begin{equation}\label{b:unif:ell}
        \lambda\,\mathrm{Id}
        \leq\nabla^2(\psi^2)(x)
        \leq\Lambda\,\mathrm{Id}
        \quad\text{in the sense of matrices on $\R^d$;}
    \end{equation}
    \item[(c)] $\psi,\psi^o\in C^2(\R^d\setminus\{0\})$;
    \item[(d)] $\psi\in C^2(\R^d\setminus\{0\})$ and there exist $0<\lambda_0\leq\Lambda_0$ such that, for all $x\in\S^{d-1}$ and all $\eta\in x^\perp$,
    \[
        \lambda_0|\eta|^2
        \leq\nabla^2\psi(x)\eta\cdot\eta
        \leq\Lambda_0|\eta|^2;
    \]
    \item[(e)] $\psi\in C^2(\R^d\setminus\{0\})$, and for all $x\neq0$,
    \[
        \mathrm{Ker}\,\nabla^2\psi(x)
        =\mathrm{span}\{x\}.
    \]
\end{itemize}
If in addition $\psi\in C^k(\R^d\setminus\{0\})$, $k\geq2$, then also $\psi^o\in C^k(\R^d\setminus\{0\})$.
\end{lemma}

\begin{proof}
(a)$\iff$(b) is the content of \cite[Section~3]{CFV16}. (b)$\implies$(d) is a consequence of \cite[Lemma~3.1]{Antetal}.

To prove (d)$\implies$(b), by the $0$-homogeneity of $\nabla^2(\psi^2/2)$, it suffices to study the case $|x|=1$. Let $\xi\in\R^d\setminus\{0\}$ and decompose it as $\xi=\alpha
x+\eta$, with $\alpha\in\R$ and $\eta\in x^\perp$. Since
\[
\nabla^2(\psi^2/2)(x)
=\psi(x)\nabla^2\psi(x)
+\nabla\psi(x)\otimes\nabla\psi(x),
\]
and since $\nabla^2\psi(x)x=0$ and $\nabla\psi(x)\cdot x=\psi(x)$ by homogeneity, we have
\[
\begin{split}
\nabla^2(\psi^2/2)(x)\xi\cdot\xi
&=\psi(x)\nabla^2\psi(x)\xi\cdot\xi
  +|\nabla\psi(x)\cdot(\alpha x+\eta)|^2\\
&=\psi(x)\nabla^2\psi(x)\eta\cdot\eta
  +|\alpha\psi(x)+\nabla\psi(x)\cdot\eta|^2\\
&\geq\lambda_0\psi(x)|\eta|^2
  +|\alpha\psi(x)+\nabla\psi(x)\cdot\eta|^2.
\end{split}
\]
Since $\xi\neq0$, either $\eta\neq0$, or $\eta=0$ and $\alpha\neq0$. In both cases the right-hand side is positive. Thus $\nabla^2(\psi^2/2)(x)$ is positive definite for every
$x\in\S^{d-1}$. By continuity and compactness, its eigenvalues are uniformly bounded above and away from zero on the sphere. Statement (b) then follows by homogeneity.

To prove (d)$\iff$(e), by homogeneity it suffices to consider $x\in\S^{d-1}$. Since $\psi$ is convex, $\nabla^2\psi(x)$ is positive semidefinite, and by homogeneity we have
$\nabla^2\psi(x)x=0$. Hence (e) is equivalent to the positive definiteness of $\nabla^2\psi(x)$ on $x^\perp$ for every such $x$. The uniform bounds in (d) then follow by continuity
and compactness.

In order to prove (c)$\implies$(b), first observe that $\psi^o\in C^1(\R^d\setminus\{0\})$ implies that $B^\psi$ is strictly convex, see \cite[Corollary~1.7.3]{Schneider} or
\cite[eq.~(3.6)]{CiaSal09}. We may therefore apply \cite[Lemma~3.1]{CiaSal09} to conclude that the maps $\nabla(\psi^2/2)$ and $\nabla((\psi^o)^2/2)$ are inverse to each other. In
particular,
\begin{equation}\label{eq:ciasal}
\bigg(
\nabla\Big(\frac{\psi^2}{2}\Big)
\circ
\nabla\Big(\frac{(\psi^o)^2}{2}\Big)
\bigg)(x)=x
\end{equation}
for all $x\in\R^d$. Differentiating this identity and setting $y=y(x):=\nabla((\psi^o)^2/2)(x)$, we get
\begin{equation}\label{eq:ciasa2}
\nabla^2\Big(\frac{\psi^2}{2}\Big)(y)\,
\nabla^2\Big(\frac{(\psi^o)^2}{2}\Big)(x)
=\mathrm{Id}
\end{equation}
for all $x\neq0$. Both matrices are symmetric and positive semidefinite, and \eqref{eq:ciasa2} implies that they are positive definite. Moreover, by continuity and
zero-homogeneity, $\nabla^2((\psi^o)^2/2)$ is uniformly bounded outside the origin. Since $x\mapsto y(x)$ maps $\R^d\setminus\{0\}$ onto itself, by \eqref{eq:ciasa2} we obtain a
uniform positive lower bound for $\nabla^2(\psi^2/2)$. The corresponding upper bound follows again from continuity and zero-homogeneity. This proves (b).

Finally, we prove (b)$\implies$(c). By the equivalence (b)$\iff$(a), $B^\psi$ is uniformly convex and of class $C^2$. Hence, by \cite[Corollary~1.7.3]{Schneider}, see also
\cite[eq.~(3.6)]{CiaSal09}, we have $\psi^o\in C^1(\R^d\setminus\{0\})$. By \cite[Lemma~3.1]{CiaSal09}, the same inverse relations between the two gradient maps hold.

By \eqref{b:unif:ell}, the map $\nabla(\psi^2/2)$ has an invertible Jacobian at every point outside the origin. The inverse function theorem therefore implies that its inverse
$\nabla((\psi^o)^2/2)$ is of class $C^1$ outside the origin. Consequently, $(\psi^o)^2$, and hence $\psi^o$, belongs to $C^2(\R^d\setminus\{0\})$.

If in addition $\psi\in C^k(\R^d\setminus\{0\})$, by the same argument we obtain
\[
\nabla\Big(\frac{(\psi^o)^2}{2}\Big)
\in C^{k-1}(\R^d\setminus\{0\}),
\]
and therefore $\psi^o\in C^k(\R^d\setminus\{0\})$.
\end{proof}

Now let $\psi$ be a  (not necessarily smooth) norm on $\R^d$. Here we provide an approximation procedure for $\psi$ following the argument of \cite{Convolution-on-liegroups}. In
what follows we denote by  $O(d), SO(d)$ the orthogonal and the special orthogonal group of $\R^{d\times d}$, respectively.

Let also $\mathfrak{h}$ be the (normalized) Haar measure on $SO(d)$. We begin by constructing a smooth approximate identity on $SO(d)$, following \cite[Section
2]{Convolution-on-liegroups}, with a few additional properties suited to our purposes.

Fix $\eta\in C^\infty_c\big((-1,1);\R_+\big)$, with $\eta(0)>0$, and for $\epsilon>0$ define $ \rho_\epsilon:SO(d)\to \R_+$ by
\begin{equation}\label{rho:epslon}
    \rho_\epsilon(R)\coloneqq \frac{1}{\bar Z_\epsilon} \eta\left(\frac{|R-\mathrm{Id}|^2}{\epsilon^2} \right),\quad \bar Z_\epsilon=\int_{SO(d)}\eta\left(\frac{|R-\mathrm{Id}|^2}{\epsilon^2} \right)\,d\mathfrak{h}(R),
\end{equation}
so that
\begin{equation}\label{properties:ofrhoe}
  \rho_\epsilon\in C^\infty\big(SO(d);\R_+\big) ,\quad  \int_{SO(d)} \rho_\epsilon(R)\,d\mathfrak{h}(R)=1,\quad \mathrm{spt}\rho_\epsilon\subseteq\{R\in SO(d):\,|R-\mathrm{Id}|\leq\epsilon\}.
\end{equation}
Recall that $|A|=\sqrt{\mathrm{tr}(A^{^{\top}}A)}$ denotes the Hilbert-Schmidt norm of a matrix $A$. In particular, for every $Q\in O(d)$, we have
$|QRQ^{^{\top}}-\mathrm{Id}|=|Q(R-\mathrm{Id})Q^{^{\top}}|=|R-\mathrm{Id}|$, and thus
\begin{equation}\label{rho:e:propr2}
    \rho_\epsilon(QRQ^{^{\top}})=\rho_\epsilon(R)\quad\text{for all $R\in SO(d)$ and all $Q\in O(d)$.}
\end{equation}
We also claim that
\begin{equation}\label{R:int:SO}
    \int_{SO(d)}R\,\rho_\epsilon(R)\,d\mathfrak{h}(R)=\mathfrak{a}_\epsilon\,\mathrm{Id}\quad \text{with $\mathfrak{a}_\epsilon\xrightarrow{\epsilon\to 0^+}1$.}
\end{equation}

To prove it, set $A_\epsilon\coloneqq \int_{SO(d)}R\,\rho_\epsilon(R)\,d\mathfrak{h}(R)$, and observe that by \eqref{rho:e:propr2}, for all $Q\in O(d)$ we have
\begin{equation*}
    Q\,A_\epsilon Q^{^\top}=\int_{SO(d)}QRQ^{^\top}\rho_\epsilon(R)\,d\mathfrak{h}(R)=\int_{SO(d)}QRQ^{^\top}\rho_\epsilon\big(QRQ^{^\top}\big)\,d\mathfrak{h}(R)=A_\epsilon,
\end{equation*}
where in the last identity we used that, setting $g_Q:SO(d)\to SO(d)$, $g_Q(R)\coloneqq QRQ^{^{\top}}$, the pushforward
\begin{equation}\label{pushforward}
    (g_Q)_\#\mathfrak{h}=\mathfrak{h}\quad\text{in the sense of measure on $SO(d)$,}
\end{equation}
 by uniqueness of the Haar measure, as $(g_Q)_\#\mathfrak{h}$ is still left-invariant (and hence it is a normalized Haar measure). Therefore the matrix $A_\epsilon=QA_\epsilon
 Q^{^{\top}}$ for all $Q\in O(d)$, and this implies $A_\epsilon= \mathfrak{a}_\epsilon\,\mathrm{Id}$ for some $\mathfrak{a}_\epsilon \in \R$. Then, by \eqref{properties:ofrhoe}, we
 have
\begin{equation*}
    \begin{split}
        |(\mathfrak a_\epsilon\,-1)\,\mathrm{Id}|=|A_\epsilon-\mathrm{Id}|=\bigg|\int_{SO(d)}(R-\mathrm{Id}) \,\rho_\epsilon(R)\,d\mathfrak{h}(R) \bigg|\leq \int_{|R-\mathrm{Id}|<\epsilon} |R-\mathrm{Id}|\,\rho_\epsilon\,d\mathfrak{h}(R)<\epsilon,
    \end{split}
\end{equation*}
and this proves \eqref{R:int:SO}.

 For sufficiently small $\epsilon>0$, we have $\mathfrak a_\epsilon>0$. Given a norm $\psi$ on $\R^d$,  we now define the approximating norm via convolution w.r.t. the Haar
 measure, that is
\begin{equation}\label{altra:convoluzione}
     \mathcal{J}_\epsilon[\psi](\xi)\coloneqq \frac{1}{\mathfrak a_\epsilon}\int_{SO(d)} \psi(R\,\xi)\,\rho_\epsilon(R)\,d\mathfrak{h}(R),\quad \xi\in \R^d.
 \end{equation}
We have the following

\begin{lemma}\label{lemma:primaapprox}
    Let $\psi$ be a norm on $\R^d$. Then $\psi_\epsilon\coloneqq\mathcal{J}_\epsilon[\psi]$ is a norm of class $C^\infty(\R^d\setminus\{0\})$ such that
    \begin{equation}\label{toprove:unif}
        \psi_\epsilon\xrightarrow{\epsilon\to 0^+}\psi\quad\text{locally uniformly in $\R^d$.}
    \end{equation}
\end{lemma}

\begin{proof}
   It is straightforward to verify that $\psi_\epsilon$ is a norm. Concerning its regularity, by homogeneity it suffices to prove that $\psi_\epsilon\in C^\infty(\S^{d-1})$. So let
   $\xi_0\in \S^{d-1}$, fix $e\in \S^{d-1}$, and consider a smooth map
   \begin{equation*}
       \mathcal{Q}:\mathcal{U}(\xi_0)\to SO(d)\quad\text{ satisfying}\quad \mathcal{Q}(\xi)e=\xi\quad\text{for all $\xi\in \mathcal{U}(\xi_0)$,}
   \end{equation*}
   where $\mathcal{U}(\xi_0)$ is an open neighborhood of $\xi_0$ on $\S^{d-1}$. The map $\mathcal{Q}$  can be constructed by applying the local section theorem \cite[Theorem
   4.6]{Lee-smooth-manifolds} to the smooth submersion $\pi_e:SO(d)\to \S^{d-1}$, $\pi_e(R)\coloneqq Re$. Therefore, via the change of variables $R\mapsto R\,\mathcal{Q}(\xi)$, and
   using the right-invariance of the Haar measure, we get
\begin{equation*}
    \psi_\epsilon(\xi)=\frac{1}{\mathfrak{a}_\epsilon}\int_{SO(d)}\psi\big(R\,\mathcal{Q} (\xi)e\big)\,\rho_\epsilon(R)\,d\mathfrak{h}(R)=\frac{1}{\mathfrak{a}_\epsilon}\int_{SO(d)}\psi(Re)\,\rho_\epsilon\big(R\,\mathcal{Q}(\xi)^{-1}\big)\,d\mathfrak{h}(R),\quad \xi\in \mathcal{U}(\xi_0).
\end{equation*}
Owing to the smoothness of $\mathcal{Q}$ and  of $\rho_\epsilon$, it follows that $\psi_\epsilon\in C^\infty(\mathcal{U}(\xi_0))$ and thus $\psi_\epsilon\in C^\infty(\S^{d-1})$ by
the arbitrariness of $\xi_0\in \S^{d-1}$. Next, by homogeneity, it suffices to prove \eqref{toprove:unif} on $\S^{d-1}$; by \eqref{properties:ofrhoe},
\begin{equation*}
    \begin{split}
      \sup_{\xi\in \S^{d-1}} \left| \int_{SO(d)} \psi(R\xi)\,\rho_\epsilon(R)\,d\mathfrak{h}(R)-\psi(\xi) \right|\leq \sup_{\xi\in \S^{d-1}}\int_{|R-\mathrm{Id}|<\epsilon} \big|\psi(R\xi)-\psi(\mathrm{Id}\xi)\big|\,\rho_\epsilon(R)\,d\mathfrak{h}(R) \leq C(d)\,\epsilon\,\mathrm{Lip}(\psi),
    \end{split}
\end{equation*}
which goes to zero as $\epsilon\to 0^+$. This concludes the proof.
\end{proof}

 We also record that $\mathcal{J}$ preserves symmetry with respect to orthogonal transformations, as well as local linearity on the sphere. This is the content of

\begin{proposition}\label{prop:nuovaproprieta}
    Let $\psi$ be a norm on $\R^d$, let $\psi_\epsilon\coloneqq \mathcal{J}_\epsilon[\psi]$, and let $Q_0\in O(d)$ and $e\in \S^{d-1}$ be given.
    \begin{itemize}
        \item[(a)] if $\psi(Q_0\xi)=\psi(\xi)$ for all $\xi\in \R^d$, then also $\psi_\epsilon(Q_0\xi)=\psi_\epsilon(\xi)$ for all $\xi\in \R^d$.
        \item[(b)] If there exist $c_0>0$ and an open set $\mathcal{U}\subset \S^{d-1}$ such that $$\psi(\xi)=c_0\,\xi\cdot e\quad \text{for all $\xi\in \mathcal{U}$},$$
        then setting $\mathcal{U}_\epsilon=\{w\in \mathcal{U}:\,\mathrm{dist}(w,\partial \mathcal{U})>\epsilon\}$, we have
        \begin{equation}
            \psi_\epsilon(\xi)=c_0\,\xi\cdot e\quad\text{for all $\xi \in \mathcal{U}_\epsilon$.}
        \end{equation}
    \end{itemize}
\end{proposition}

Distances on $\mathbb S^{d-1}$ in this statement are Euclidean chord distances. Indeed $|R\xi-\xi|\leq\|R-I\|_{\mathrm{op}} \leq\|R-I\|_{\mathrm{HS}}$ for $|\xi|=1$, so no
additional dimensional constant is needed in the definition of $\mathcal U_\epsilon$.

\begin{proof}
    (a). Replacing $\xi$ with $Q_0^{^{\top}}\xi$, we deduce that $\psi(Q_0^{^{\top}}\xi)=\psi(\xi)$. Therefore
    \begin{equation*}
        \int_{SO(d)}\psi(RQ_0\,\xi)\,\rho_\epsilon(R)\,d\mathfrak{h}(R)=\int_{SO(d)}\psi(Q_0^{^{\top}}RQ_0\xi)\,\rho_\epsilon(R)\,d\mathfrak{h}(R)=\int_{SO(d)} \psi(S\xi)\,\rho_\epsilon(S)\,d\mathfrak{h}(S)
    \end{equation*}
    where in the last equality we used \eqref{pushforward} (with $Q=Q_0^{^{\top}}$), and \eqref{rho:e:propr2}. This proves (a). \\
   (b). We first observe that, if $\xi\in \mathcal{U}_\epsilon$ then, recalling \eqref{properties:ofrhoe}, we have $R\xi\in \mathcal{U}$ for all $R\in \mathrm{spt}\,\rho_\epsilon$.
   From \eqref{R:int:SO}, it then follows that
   \begin{equation*}
       \begin{split}
           \psi_\epsilon(\xi)=\frac{1}{\mathfrak{a}_\epsilon}\int_{\mathrm{spt}\,\rho_\epsilon}\psi(R\xi)\,\rho_\epsilon(R)\,d\mathfrak{h}(R)=\frac{c_0}{\mathfrak{a}_\epsilon}\bigg(\int_{\mathrm{spt}\,\rho_\epsilon}R\,\rho_\epsilon(R)\,d\mathfrak{h}(R)\bigg)\,\xi\cdot e=c_0\,\xi\cdot e,
       \end{split}
   \end{equation*}
for all $\xi\in \mathcal{U}_\epsilon$. This concludes the proof.
\end{proof}

 \section{Extended norms}\label{app:ext:norms}

Let $M$ be a norm on $\R^N_x\times\R_z$, symmetric with respect to $z$, namely
\begin{equation}\label{Mz:sym}
M(x,-z)=M(x,z)\qquad\text{for all $(x,z)\in\R^{N+1}$.}
\end{equation}
We begin with some preliminary remarks. First observe that, since $M$ is a norm,
\begin{equation}\label{Mz:sym1}
M(-x,z)=M(x,-z)=M(x,z).
\end{equation}
Thus $x\mapsto M(x,z)$ is convex and symmetric for every $z\in\R$, which implies $M(0,z)\leq M(x,z)$. Now set $c_M:=M(0,1)^{-1}$ and, for $0\leq q<c_M$, consider the sets
\begin{equation}\label{Kq:def}
K_q:=\{x\in\R^N:M(x,q)\leq1\}.
\end{equation}
They are compact, convex and symmetric, and contain the origin in their interior, since $M(0,q)=q/c_M<1$. Therefore their Minkowski gauges
\begin{equation}\label{mink:gauge}
\mu_q(x):=\inf\{\lambda>0:x\in\lambda K_q\}
\end{equation}
are norms on $\R^N$. For $M$ as above, we define
\begin{equation}\label{Tau:M}
\mathcal T[M](x,z):=c_M|z|+\int_{|z|}^{+\infty}\left(c_M-\frac{t}{M(x,t)}\right)\,dt.
\end{equation}
Observe that the integrand is nonnegative, since $M(x,t)\geq M(0,t)=t/c_M$ for $t>0$. Let us now show, that  $\mathcal T[M]$ is a norm.

\begin{lemma}\label{prop:TauMnorm}
Let $M$ be a norm on $\R^{N+1}$ satisfying \eqref{Mz:sym} and such that
\begin{equation}\label{ass:finite}
\int_0^{+\infty}\left(c_M-\frac{t}{M(x,t)}\right)\,dt<+\infty\qquad\text{for all $x\neq0$.}
\end{equation}
Then $\Psi:=\mathcal T[M]$ is a norm on $\R^{N+1}$, symmetric with respect to $z$, and satisfies
\begin{equation}\label{temp:repr}
\Psi(x,z)=\int_0^{c_M}\max\{|z|,q\mu_q(x)\}\,dq.
\end{equation}
Moreover,
\begin{equation}\label{its:clear}
\partial_z\Psi(x,z)=\frac{z}{M(x,z)}\qquad\text{for all $(x,z)\neq(0,0)$.}
\end{equation}
\end{lemma}

\begin{proof}
For $t>0$, set $h_x(t):=t/M(x,t)$. By convexity and symmetry, $s\mapsto M(sx,1)$ is nondecreasing for $s\geq0$. Hence $h_x(t)=1/M(x/t,1)$ is nondecreasing, with $0<h_x(t)\leq c_M$.

We now prove \eqref{temp:repr}. For every $0<q<c_M$ and $t>0$, by homogeneity we have
\[
h_x(t)\geq q
\iff M\left(\frac qt x,q\right)\leq1
\iff \frac qt x\in K_q
\iff t\geq q\mu_q(x).
\]
Thus,
\[
\int_{|z|}^{+\infty}(c_M-h_x(t))\,dt
=\int_0^{c_M}\int_{|z|}^{+\infty}\chi_{\{t<q\mu_q(x)\}}\,dt\,dq
=\int_0^{c_M}(q\mu_q(x)-|z|)_+\,dq.
\]
Substituting into \eqref{Tau:M}, and writing $c_M|z|=\int_0^{c_M}|z|\,dq$, we obtain \eqref{temp:repr}, since $|z|+(q\mu_q(x)-|z|)_+=\max\{|z|,q\mu_q(x)\}$.

For every $0<q<c_M$, the integrand in \eqref{temp:repr} is convex, even and positively one-homogeneous in $(x,z)$. Hence $\Psi$ has the same properties. Moreover, it is finite by
\eqref{ass:finite} and positive outside the origin. Thus $\Psi$ is a norm. Finally, \eqref{its:clear} follows from \eqref{Tau:M} and the fundamental theorem of calculus.
\end{proof}

We next study the regularity and uniform convexity of the norm obtained through $\mathcal T$.

\begin{lemma}\label{lemma:MC2}
Let $M\in C^2(\R^{N+1}\setminus\{0\})$ be a uniformly convex norm satisfying \eqref{Mz:sym}. Then $\Psi=\mathcal T[M]\in C^2(\R^{N+1}\setminus\{0\})$ is a uniformly convex norm. If
in addition $M\in C^k(\R^{N+1}\setminus\{0\})$, for some $k\geq2$, then also $\Psi\in C^k(\R^{N+1}\setminus\{0\})$.
\end{lemma}

\begin{proof}
By symmetry, it suffices to study the case $z\geq 0$. We split the proof into two steps.

\medskip
\noindent\textit{Step 1. Smoothness of $\Psi=\mathcal T[M]$.} Suppose that $M\in C^k(\R^{N+1}\setminus\{0\})$, with $k\geq2$, and set $F(y):=1/M(y,1)$. Then $F\in C^k(\R^N)$,
$F(0)=c_M$ and, by \eqref{Mz:sym1}, $DF(0)=0$. Using homogeneity and the change of variables $s=1/t$ in \eqref{Tau:M}, we find
\begin{equation}\label{regolaritasezione}
\Psi(x,1)=c_M+\int_0^1\frac{F(0)-F(sx)}{s^2}\,ds.
\end{equation}
By Taylor's formula, the integrand extends continuously to $s=0$, with value $-\frac12D^2F(0)[x,x]$. In particular, this shows that the integral is finite and, in turn, that
\eqref{ass:finite} holds. Thus $\Psi$ is a norm by Lemma~\ref{prop:TauMnorm}.

We now differentiate the integrand in \eqref{regolaritasezione} with respect to $x$. For $s>0$,
\[
D_x\left(\frac{F(0)-F(sx)}{s^2}\right)
=-\frac{DF(sx)}s
=-\int_0^1D^2F(\theta sx)x\,d\theta,
\qquad
D_x^j\left(\frac{F(0)-F(sx)}{s^2}\right)
=-s^{j-2}D^jF(sx),
\]
for  $2\leq j\leq k$. These derivatives also extend continuously to $s=0$. The extensions are jointly continuous in $(x,s)$, and hence bounded when $x$ ranges over a compact set
and $s\in[0,1]$. We may therefore differentiate under the integral sign up to order $k$, obtaining $\Psi(\cdot,1)\in C^k(\R^N)$. By \eqref{its:clear},
\[
\Psi(x,z)=\Psi(x,1)+\int_1^z\frac{t}{M(x,t)}\,dt.
\]
The second (integral) term is jointly of class $C^k$ in $(x,z)$ away from $(0,0)$. Thus $\Psi\in C^k(\R^{N+1}\setminus\{0\})$.

\medskip
\noindent\textit{Step 2. Uniform convexity.} For every $0<q<c_M$, the set $K_q$ is a section of the uniformly convex $C^2$ unit ball of $M$, and is therefore itself of class $C^2$
and uniformly convex. By Lemma~\ref{lem:varie:unifconv}, its associated norm satisfies
\begin{equation}\label{pos:D2xtx}
\mu_q\in C^2(\R^N\setminus\{0\}),\qquad D_x^2\mu_q(x)\geq0,\qquad \ker D_x^2\mu_q(x)=\operatorname{span}\{x\}.
\end{equation}
Moreover, the map $(x,q)\mapsto\mu_q(x)$ is $C^2$ for $x\ne0$ and $0<q<c_M$, by the implicit function theorem applied to $M(x/\lambda,q)=1$, with unknown $\lambda$. Indeed, at
$\lambda=\mu_q(x)$, by convexity we have $\partial_\lambda M(x/\lambda,q)=-\lambda^{-1}D_xM(x/\lambda,q)\cdot(x/\lambda)\le-\lambda^{-1}(1-M(0,q))<0$.

Fix $(x,z)$ with $x\neq0$. By symmetry and strict convexity  we have $M(0,z)<M(x,z)$. Consequently, $D_xM(x,z)\cdot x\geq M(x,z)-M(0,z)>0$. Differentiating \eqref{its:clear} and
using Euler's identity, we obtain
\begin{equation}\label{dzz:positive}
\partial_{zz}^2\Psi(x,z)=\frac{D_xM(x,z)\cdot x}{M(x,z)^2}>0.
\end{equation}
The same strict inequality $M(x,z)>M(0,z)=z/c_M$ shows that $z/M(x,z)<c_M$. We may therefore choose $z/M(x,z)<a<b<c_M$. By continuity, $|s|/M(y,s)<a$ in a neighborhood of $(x,z)$,
which we take small enough that $y\neq0$. By the equivalence established in the proof of Lemma~\ref{prop:TauMnorm} we then have $|s|<q\mu_q(y)$ throughout this neighborhood,
simultaneously for every $q\in[a,b]$. Hence \eqref{temp:repr} becomes
\[
\Psi(y,s)=\int_a^b q\mu_q(y)\,dq+\int_{(0,c_M)\setminus[a,b]}\max\{|s|,q\mu_q(y)\}\,dq.
\]
The second term is convex and, locally, is the difference of two $C^2$ functions. Thus
\begin{equation}\label{controllohessiana}
D^2\Psi(x,z)[(\xi,\zeta),(\xi,\zeta)]\geq\int_a^b q\,D_x^2\mu_q(x)[\xi,\xi]\,dq.
\end{equation}
If the left-hand side vanishes, \eqref{pos:D2xtx} implies $\xi=\alpha x$ for some $\alpha\in\R$. On the other hand, by homogeneity we have $D^2\Psi(x,z)(x,z)=0$, and therefore
\[
D^2\Psi(x,z)[(\alpha x,\zeta),(\alpha x,\zeta)]=(\zeta-\alpha z)^2\partial_{zz}^2\Psi(x,z).
\]
By \eqref{dzz:positive}, this expression vanishes if and only if $\zeta=\alpha z$ and thus we have  proved
\begin{equation}\label{ker:D2Psi1}
\ker D^2\Psi(x,z)=\operatorname{span}\{(x,z)\}\qquad\text{for all $x\neq0$.}
\end{equation}

We are left to study the poles. Since $D_xM(0,t)=0$ for $t>0$, differentiating \eqref{Tau:M} twice with respect to $x$, we obtain
\[
D_x^2\Psi(0,1)=\int_1^\infty\frac{t\,D_x^2M(0,t)}{M(0,t)^2}\,dt
=\frac{D_x^2M(0,1)}{M(0,1)^2}\int_1^\infty t^{-2}\,dt
=\frac{D_x^2M(0,1)}{M(0,1)^2},
\]
where we used $M(0,t)=tM(0,1)$ and $D_x^2M(0,t)=t^{-1}D_x^2M(0,1)$. Differentiation under the integral is justified by Step 1, equivalently using \eqref{regolaritasezione}.
Moreover, \eqref{its:clear} yields $D_x\partial_z\Psi(0,1)=0$, while $\Psi(0,z)=c_Mz$ for $z>0$ implies $\partial_{zz}^2\Psi(0,1)=0$. Hence
\[
D^2\Psi(0,1)=
\begin{pmatrix}
M(0,1)^{-2}D_x^2M(0,1)&0\\
0&0
\end{pmatrix}.
\]
The horizontal block is positive definite by the uniform convexity of $M$. Therefore \eqref{ker:D2Psi1} holds also at $(0,1)$ and, by symmetry and homogeneity, at every $(0,z)$
with $z\neq0$. The characterization in Lemma~\ref{lem:varie:unifconv} completes the proof.
\end{proof}

We now have all the ingredients to prove the main result of this section. Before stating it, we briefly describe the construction. We first define the extended mobility
$\widetilde{\mathfrak M}^{\gamma}$ in \eqref{def:firstMk}. For every sufficiently small $\gamma>0$, it agrees with $\mathfrak m(x)$ in a neighborhood of the equator
$\{z=0\}\cap\S^N$, and with $|z|/\gamma$ in neighborhoods of the two poles. Between these regions, we interpolate between $\phi(x)$ and $|z|/\gamma$ in such a way that the
associated norms $\widetilde\Phi^\gamma :=\mathcal T[\widetilde{\mathfrak M}^{\gamma}]$ converge locally uniformly to $\phi(x)$ as $\gamma\to0^+$. Finally, we regularize
$\widetilde{\mathfrak M}^{\gamma}$ using the approximation procedure of Appendix~\ref{sec:normeunifconvesse}.

\begin{center}
    \includegraphics[width=12cm]{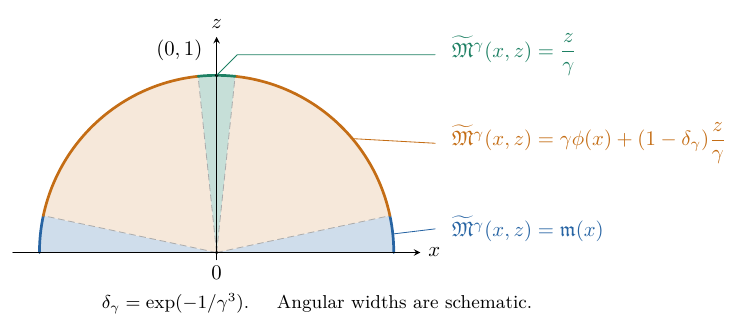}
\end{center}

\begin{proposition}\label{prop:extended-norms}
Let $\phi,\mathfrak m$ be norms on $\R^N$. Then there exist $\kappa_0>0$ and families of norms $\{\Phi^\kappa\}_{\kappa\in(0,\kappa_0)}$, $\{\mathfrak
M^\kappa\}_{\kappa\in(0,\kappa_0)}$ on $\R^{N+1}$, symmetric with respect to $z$, such that $\Phi^\kappa,\mathfrak M^\kappa\in C^\infty(\R^{N+1}\setminus\{0\})$ are uniformly
convex for every $\kappa$, and
\begin{equation}\label{come:e:ovvia}
\partial_z\Phi^\kappa(x,z)=\frac{z}{\mathfrak M^\kappa(x,z)}\qquad\text{for all $(x,z)\neq(0,0)$.}
\end{equation}
Moreover,
\begin{equation}\label{meno:ovvia}
\frac{z^2}{\Phi^\kappa(x,z)\,\mathfrak M^\kappa(x,z)}\leq1\qquad\text{for all $(x,z)\in\S^N$,}
\end{equation}
with equality if and only if $(x,z)=(0,\pm1)$, and
\begin{equation}\label{za:tesi}
\begin{aligned}
\mathfrak M^\kappa(x,0)&\longrightarrow\mathfrak m(x)&&\text{locally uniformly in $\R^N$,}\\
\Phi^\kappa(x,z)&\longrightarrow\phi(x)&&\text{locally uniformly in $\R^{N+1}$,}
\end{aligned}
\qquad\text{as $\kappa\to0^+$.}
\end{equation}
If in addition $\phi\in C^1(\R^N\setminus\{0\})$, then
\begin{equation}\label{za:tesi2}
\Phi^\kappa(x,z)\longrightarrow\phi(x)\qquad\text{in $C^1_{\loc}(\R^{N+1}\setminus\{x=0\})$.}
\end{equation}
\end{proposition}

\begin{proof}
\noindent\textit{Step 1. A first approximation.} Choose $\gamma_0>0$ so small that $\gamma\phi(x)<\mathfrak m(x)$ for every $x\neq0$ and $0<\gamma<\gamma_0$. Set
$\delta_\gamma:=\exp(-1/\gamma^3)$ and define
\begin{equation}\label{def:firstMk}
\widetilde{\mathfrak M}^{\gamma}(x,z):=\max\left\{\mathfrak m(x),\,\gamma\phi(x)+(1-\delta_\gamma)\frac{|z|}{\gamma},\,\frac{|z|}{\gamma}\right\}.
\end{equation}
This is a norm, symmetric with respect to $z$, and
\begin{equation}\label{prim:Mkequivm}
\widetilde{\mathfrak M}^{\gamma}(x,0)=\mathfrak m(x)\qquad\text{for all $x\in\R^N$.}
\end{equation}
For fixed $x$, it agrees with $|z|/\gamma$ when $|z|$ is sufficiently large. Hence \eqref{ass:finite} holds, with $c_{\widetilde{\mathfrak M}^{\gamma}}=\gamma$, and
$\widetilde\Phi^\gamma:=\mathcal T[\widetilde{\mathfrak M}^{\gamma}]$ is a norm by Lemma~\ref{prop:TauMnorm}.

We now study the limit of $\widetilde\Phi^\gamma(x,0)$. To this aim, consider the norm
\[
L^\gamma(x,z):=\max\left\{\gamma\phi(x)+(1-\delta_\gamma)\frac{|z|}{\gamma},\,\frac{|z|}{\gamma}\right\}.
\]
It also satisfies \eqref{ass:finite}, with $c_{L^\gamma}=\gamma$. Its two defining terms agree when $|z|=\gamma^2\phi(x)/\delta_\gamma$ and the maximum is given by the second term
for $|z|>\gamma^2\phi(x)/\delta_\gamma$. Therefore, for $x\neq0$, by direct computation we find
\begin{equation}\label{integralebase}
\begin{aligned}
\mathcal T[L^\gamma](x,0)
&=\int_0^{\gamma^2\phi(x)/\delta_\gamma}\left(\gamma-\frac{t}{\gamma\phi(x)+(1-\delta_\gamma)t/\gamma}\right)\,dt\\
&=\frac{\gamma^3\phi(x)}{(1-\delta_\gamma)^2}\left(\log\frac1{\delta_\gamma}-1+\delta_\gamma\right)
=\phi(x)\frac{1-\gamma^3(1-\delta_\gamma)}{(1-\delta_\gamma)^2}.
\end{aligned}
\end{equation}
The identity holds also at $x=0$, where both sides vanish.

Observe now that $\widetilde{\mathfrak M}^{\gamma}=\max\{\mathfrak m,L^\gamma\}$ and $L^\gamma(x,t)\geq t/\gamma$ for $t>0$. Thus the two norms agree whenever $t\geq\gamma\mathfrak
m(x)$, and
\begin{equation}\label{confrontointegrali}
0\leq\widetilde\Phi^\gamma(x,0)-\mathcal T[L^\gamma](x,0)
=\int_0^{\gamma\mathfrak m(x)}t\left(\frac1{L^\gamma(x,t)}-\frac1{\widetilde{\mathfrak M}^{\gamma}(x,t)}\right)\,dt
\leq\gamma^2\mathfrak m(x),
\end{equation}
where the last inequality follows from $t/L^\gamma(x,t)\leq\gamma$. Coupling \eqref{integralebase} and \eqref{confrontointegrali}, we deduce that
$\widetilde\Phi^\gamma(x,0)\to\phi(x)$ locally uniformly as $\gamma\to 0^+$. Moreover, by \eqref{its:clear},
\[
0\leq\widetilde\Phi^\gamma(x,z)-\widetilde\Phi^\gamma(x,0)
=\int_0^{|z|}\frac{t}{\widetilde{\mathfrak M}^{\gamma}(x,t)}\,dt\leq\gamma|z|.
\]
Consequently,
\begin{equation}\label{prim:phikconv}
\widetilde\Phi^\gamma(x,z)\longrightarrow\phi(x)\qquad\text{locally uniformly in $\R^{N+1}$ as $\gamma\to0^+$.}
\end{equation}

\medskip
\noindent\textit{Step 2. A further approximation.} Fix $\gamma\in(0,\gamma_0)$, and let $\mathcal J_\epsilon$ be as in \eqref{altra:convoluzione}. For $\epsilon>0$ sufficiently
small, define
\begin{equation}\label{second:Mk}
\mathfrak M^{\gamma,\epsilon}(x,z):=\sqrt{\bigl(\mathcal J_\epsilon[\widetilde{\mathfrak M}^{\gamma}](x,z)\bigr)^2+\epsilon|(x,z)|^2}.
\end{equation}
By Lemma~\ref{lemma:primaapprox} and Proposition~\ref{prop:nuovaproprieta}(a), these are norms of class $C^\infty(\R^{N+1}\setminus\{0\})$, symmetric with respect to $z$, and
\begin{equation}\label{dovro:richiamar}
\mathfrak M^{\gamma,\epsilon}\longrightarrow\widetilde{\mathfrak M}^{\gamma}\qquad\text{locally uniformly in $\R^{N+1}$ as $\epsilon\to0^+$.}
\end{equation}
Since  by  \eqref{second:Mk}, $D^2((\mathfrak M^{\gamma,\epsilon})^2)\geq2\epsilon\,\Id$ outside the origin, we have that  $\mathfrak M^{\gamma,\epsilon}$ is uniformly convex by
Lemma~\ref{lem:varie:unifconv}. Setting $\Phi^{\gamma,\epsilon}:=\mathcal T[\mathfrak M^{\gamma,\epsilon}]$, Lemma~\ref{lemma:MC2} shows that $\Phi^{\gamma,\epsilon}$ is also a
uniformly convex norm of class $C^\infty(\R^{N+1}\setminus\{0\})$.

We claim that, for fixed $\gamma$,
\begin{equation}\label{la:quasi:tesi}
\Phi^{\gamma,\epsilon}\longrightarrow\widetilde\Phi^\gamma\qquad\text{locally uniformly in $\R^{N+1}$ as $\epsilon\to0^+$.}
\end{equation}
By homogeneity, it suffices to prove uniform convergence on $\S^N$. Observe that $\widetilde{\mathfrak M}^{\gamma}(x,z)=|z|/\gamma$ in an open cone containing the vertical axis
outside the origin. Choose $R_\gamma>1$ so large that the directions of all $(x,t)$ with $|x|\leq1$ and $t\geq R_\gamma$ lie in a compact subset of this cone on the upper
hemisphere. By Proposition~\ref{prop:nuovaproprieta}(b), for all sufficiently small $\epsilon$,
\begin{equation}\label{tagliofisso}
\mathcal J_\epsilon[\widetilde{\mathfrak M}^{\gamma}](x,t)=\frac t\gamma,\qquad
\mathfrak M^{\gamma,\epsilon}(x,t)=\sqrt{(\gamma^{-2}+\epsilon)t^2+\epsilon|x|^2}
\qquad\text{if $|x|\leq1$ and $t\geq R_\gamma$.}
\end{equation}
Moreover, $c_{\mathfrak M^{\gamma,\epsilon}}=(\gamma^{-2}+\epsilon)^{-1/2}\to\gamma$.

We compare the expressions in \eqref{Tau:M} for the two norms $\mathfrak M^{\gamma,\epsilon}$ and $\widetilde{\mathfrak M}^{\gamma}$, splitting the integrals at $R_\gamma$. By
\eqref{dovro:richiamar} and homogeneity, it easily follows that $\mathfrak M^{\gamma,\epsilon}/\widetilde{\mathfrak M}^{\gamma}\to1$ uniformly outside the origin. In turn, since
$t/\widetilde{\mathfrak M}^{\gamma}(x,t)\leq\gamma$ for $t\geq0$, it follows that
\[
\sup_{\substack{(x,t)\neq(0,0)\\t\geq0}}\left|\frac{t}{\mathfrak M^{\gamma,\epsilon}(x,t)}-\frac{t}{\widetilde{\mathfrak M}^{\gamma}(x,t)}\right|\leq \gamma \sup_{\substack{(x,t)\neq(0,0)\\t\geq0}}\left|\frac{\widetilde{\mathfrak M}^{\gamma}(x,t)}{\mathfrak M^{\gamma,\epsilon}(x,t)}-1\right|\longrightarrow0.
\]
Hence the difference of the integrals over $[|z|,R_\gamma]$ corresponding to the two norms tends to zero uniformly for $(x,z)\in\S^N$, as does the  term $c_{\mathfrak
M^{\gamma,\epsilon}}|z|- c_{\widetilde{\mathfrak M}^{\gamma}}|z|=((\gamma^{-2}+\epsilon)^{-1/2}-\gamma)|z|$.

On $[R_\gamma,+\infty)$, the integrand corresponding to $\widetilde{\mathfrak M}^{\gamma}$ vanishes. On the other hand, \eqref{tagliofisso} implies
\[
0\leq c_{\mathfrak M^{\gamma,\epsilon}}-\frac{t}{\mathfrak M^{\gamma,\epsilon}(x,t)}
=\frac1{\sqrt{\gamma^{-2}+\epsilon}}-\frac1{\sqrt{\gamma^{-2}+\epsilon+\epsilon|x|^2/t^2}}
\leq\frac{\gamma^3\epsilon|x|^2}{2t^2}.
\]
For $|x|\leq1$, its integral over $[R_\gamma,+\infty)$ is bounded by $\gamma^3\epsilon/(2R_\gamma)$, which tends to zero.  Thus, we have established \eqref{la:quasi:tesi}.

\medskip
\noindent\textit{Step 3. Conclusion.} By \eqref{dovro:richiamar} and \eqref{la:quasi:tesi}, for every $0<\kappa<\gamma_0$ we may choose $0<\epsilon_\kappa<\kappa$ such that
\[
\|\mathfrak M^{\kappa,\epsilon_\kappa}-\widetilde{\mathfrak M}^{\kappa}\|_{L^\infty(\S^N)}
+\|\Phi^{\kappa,\epsilon_\kappa}-\widetilde\Phi^\kappa\|_{L^\infty(\S^N)}\leq\kappa
\]
and set $\mathfrak M^\kappa:=\mathfrak M^{\kappa,\epsilon_\kappa}$ and $\Phi^\kappa:=\Phi^{\kappa,\epsilon_\kappa}$. Then,  \eqref{prim:Mkequivm}, \eqref{prim:phikconv} imply
\eqref{za:tesi}, while the required regularity, symmetry and uniform convexity properties follow from Step 2.  Recall also that  identity \eqref{come:e:ovvia} follows from
Lemma~\ref{prop:TauMnorm}. Next, note that by Euler's identity and \eqref{come:e:ovvia}, $\Phi^\kappa(x,z)=D_x\Phi^\kappa(x,z)\cdot x+\frac{z^2}{\mathfrak M^\kappa(x,z)}$, while by
strict convexity and symmetry for $x\neq0$ we have $D_x\Phi^\kappa(x,z)\cdot x\geq\Phi^\kappa(x,z)-\Phi^\kappa(0,z)>0$. Combined with the previous identity, this proves
\eqref{meno:ovvia}, including its equality case. If $\phi\in C^1(\R^N\setminus\{0\})$,  \eqref{za:tesi2} follows from \eqref{za:tesi} and
\cite[Theorem~25.7]{RockafellarConvexAnalysis}, applied locally.
\end{proof}

 \bibliographystyle{siam}
\bibliography{Ilmanen-ext}
\end{document}